Benjamin M$^c$Kay

# Introduction to Exterior Differential Systems

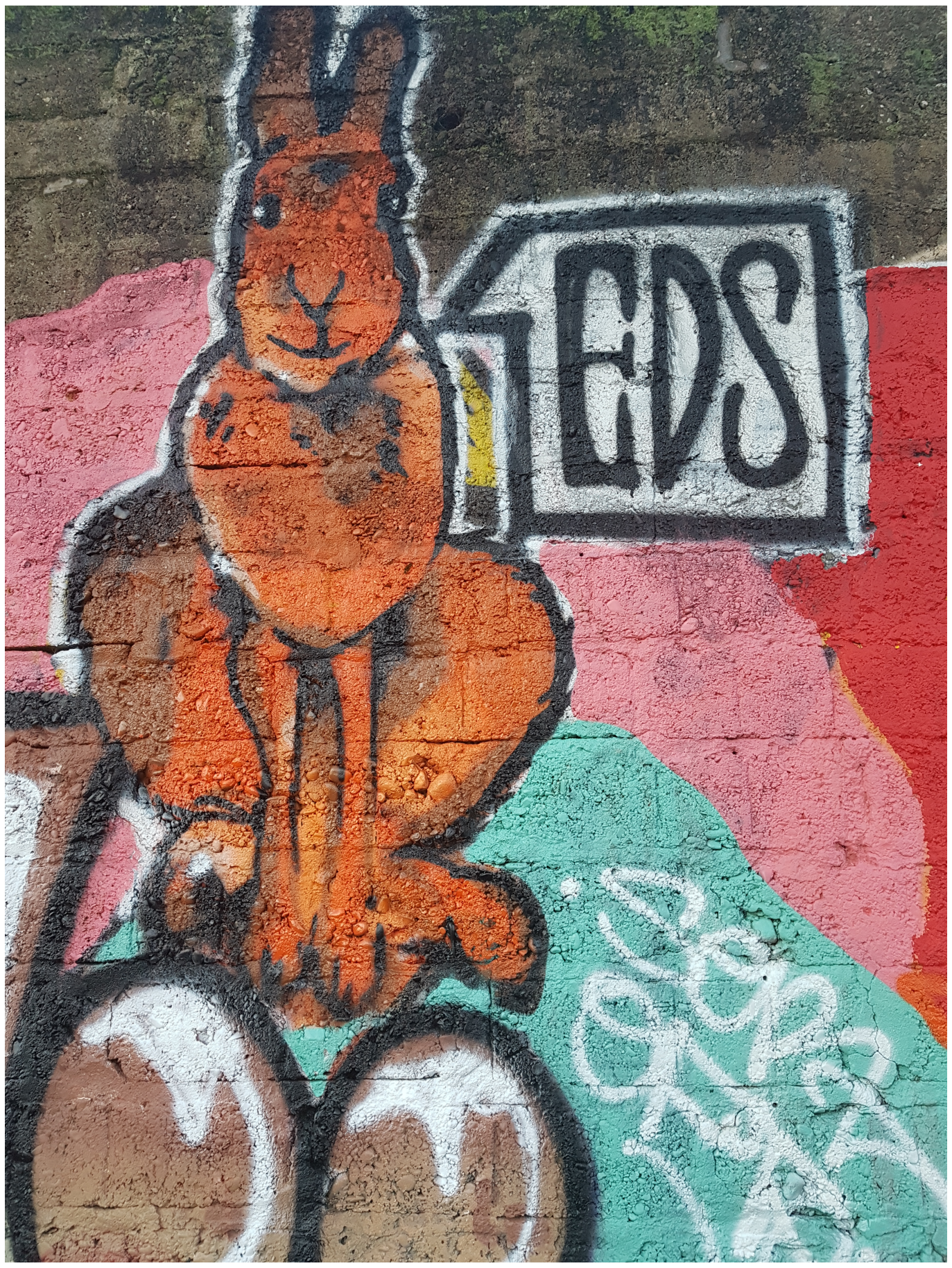
EDS

# Preface

To the reader who wants to dip a toe in the water: just read chapter 1. The reader who continues on to chapters 2 and 4 will pick up the rest of the tools. Subsequent chapters prove the theorems. We assume that the reader is familiar with elementary differential geometry on manifolds and with differential forms. These lectures explain how to apply the Cartan–Kähler theorem to problems in differential geometry. Given some differential equations, we want to decide if they are locally solvable. The Cartan–Kähler theorem gives a linear algebra test: if the equations pass the test, they are locally solvable. We give the necessary background on partial differential equations in appendices A, B, and the (not so necessary) background on moving frames in appendices D, F, G. The reader should be aware of [4], which we will follow closely, and also the canonical reference work [3] and the canonical textbook [19].

I wrote these notes for lectures at the Banach Center in Warsaw, at the University of Paris Sud, and at the University of Rome Tor Vergata. I thank Jan Gutt, Gianni Manno, Giovanni Moreno, Nefton Pali, and Filippo Bracci for invitations to give those lectures, and Francesco Russo for the opportunity to write these notes at the University of Catania.

# Contents

Chapter 1

# Exterior differential systems

We define exterior differential systems, and state the Cartan–Kähler theorem.

## Background material from differential geometry

By the expression *analytic*, we mean *real analytic*; see appendix A. We assume that all of our manifolds, maps vector fields and differential forms are analytic. (In a few exercises, we describe some such thing as *smooth*, i.e. $C^\infty$, not necessarily analytic.) A *submanifold* is an immersion of manifolds, defined up to diffeomorphism of the domain. Denote the dimension of a manifold $M$ by $\dim M$. The *codimension* of a submanifold $S$ of $M$ is $\dim M - \dim S$. A *hypersurface* in a manifold $M$ is a submanifold of codimension one. A *hyperplane* in a vector space $V$ is a linear subspace $W \subset V$ so that $\dim(V/W) = 1$.

## Differential equations encoded in differential forms

To express a differential equation $0 = f\left(x, u, \frac{\partial u}{\partial x}\right)$, add a variable $p$ to represent the derivative $\frac{\partial u}{\partial x}$. Let $\vartheta := du - p\,dx$, $\omega := dx$. Let $M$ be the manifold $M := \{\, (x, u, p) \mid 0 = f(x, u, p) \,\}$ (assuming it is a manifold). Any submanifold of $M$ on which $0 = \vartheta$ and $0 \neq \omega$ is locally the graph of a solution. It is easy to generalize this to any number of variables and any number of equations of any order.

## Exterior differential systems

An *integral manifold* of a collection of differential forms is a submanifold on which the forms vanish. An *exterior differential system* is an ideal $\mathcal{I} \subset \Omega^*$ of differential forms on a manifold $M$, closed under exterior derivative, which splits into a direct sum

$$\mathcal{I} = \mathcal{I}^1 \oplus \cdots \oplus \mathcal{I}^{\dim M}$$

of forms of various degrees: $\mathcal{I}^k := \mathcal{I} \cap \Omega^k$. Any collection of differential forms has the same integral manifolds as the exterior differential system it generates.

Some trivial examples: the exterior differential system generated by

*a.* $0$,

*b.* $\Omega^*$,

*c.* the pullbacks of all forms via a submersion,

*d.* $dx^1 \wedge dy^1 + dx^2 \wedge dy^2$ in $\mathbb{R}^4$,

*e.* $dy - z\,dx$ on $\mathbb{R}^3$.

**1.1** What are the integral manifolds of our trivial examples?

By definition $\mathcal{I}^0 = 0$, i.e. there are no nonzero functions in $\mathcal{I}$. If we wish some functions to vanish, we can replace our manifold $M$ by the zero locus of those functions (which might not be a manifold, a technical detail we will ignore).

## Integral elements

An *integral element* of an exterior differential system $\mathcal{I}$ is a linear subspace of a tangent space, on which all forms in $\mathcal{I}$ vanish. Every tangent space of an integral manifold is an integral element, but some integral elements of some exterior differential systems don't arise as tangent spaces of integral manifolds. A 1-dimensional integral element is an *integral line*. A 2-dimensional integral element is an *integral plane*.

**1.2** What are the integral elements of our trivial examples?

**1.3** Suppose that the 1-forms in an exterior differential system span a subspace of constant rank in each cotangent space. Prove that there is an integral curve tangent to each integral line.

**1.4** Give an example of an embedded integral manifold, whose every tangent space is a hyperplane in a unique integral element, but which is not a hypersurface in an integral manifold.

**1.5** Prove that a $k$-dimensional linear subspace of a tangent space is an integral element of an exterior differential system $\mathcal{I}$ just when all $k$-forms in $\mathcal{I}$ vanish on it.

The *polar equations* of an integral element $E$ are the linear functions

$$w \in T_m M \mapsto \vartheta(w, e_1, e_2, \dots, e_k)$$

where $\vartheta \in \mathcal{I}^{k+1}$ and $e_1, e_2, \dots, e_k \in E$, for any $k$. They vanish on a vector $w$ just when the span of $\{\, w \,\} \cup E$ is an integral element. The set of polar equations of $E$ is a linear subspace of $T_m^* M$. If an integral element $E$ lies in a larger one $F$, then all polar equations of $E$ occur among those of $F$: larger integral elements have more (or at least the same) polar equations.

**1.6** What are the polar equations of the integral elements of our trivial examples?

An *integral flag* of an exterior differential system is a sequence of nested integral elements

$$0 = E_0 \subset E_1 \subset E_2 \subset \cdots \subset E_p$$

with dimensions $0, 1, 2, \dots, p$. *Danger:* most authors require that a flag have subspaces of all dimensions; we *don't*: we only require that the subspaces have all dimensions $0, 1, 2, \dots, p$ up to some dimension $p$. Each subspace $E_0, E_1, \dots, E_p$ has successively larger space of polar equations. The *characters* $s_0, s_1, \dots, s_p$ of the flag are the increments in rank of polar equations: $E_k$ has polar equations of rank $s_0 + s_1 + \cdots + s_k$.

Consider all flags inside a given $E_p$. Polar equations remain linearly independent under small motions of a flag. So if a flag in $E_p$ has maximal dimensional polar equations, i.e. character sums $s_0, s_0 + s_1, \dots, s_0 + \cdots + s_p$, then so do all nearby flags in $E_p$. The *characters* of the integral element $E_p$ are those of any such flag.

**1.7** What are the characters of the integral flags of our trivial examples?

## Integral elements as points of the Grassmann bundle

The rank $p$ *Grassmann bundle* of a manifold $M$ is the set of all $p$-dimensional linear subspaces of tangent spaces of $M$.

**1.8** Recall how charts are defined on the Grassmann bundle. Prove that the Grassmann bundle is a fiber bundle.

The integral elements of an exterior differential system $\mathcal{I}$ form a subset of the Grassmann bundle. We would like to see if that subset is a submanifold, and predict its dimension.

## The Cartan–Kähler theorem

An integral element, say of dimension $p$ and with characters $s_0, \dots, s_p$, *predicts* the dimension $\dim M + s_1 + 2s_2 + \cdots + ps_p$. It predicts *correctly* if the nearby integral elements form a submanifold of the Grassmann bundle of the predicted dimension. (We will see that they always sit in a submanifold of the predicted dimension.) An integral element which correctly predicts dimension is *involutive*. An exterior differential system on a manifold $M$ is *involutive* if each component of $M$ contains an involutive maximal integral element. (We will see that involutivity of an exterior differential system implies involutivity of a dense open set of maximal integral elements.)

**Theorem 1.1** (Cartan–Kähler). *Every involutive integral element of any analytic exterior differential system is tangent to an analytic integral manifold.*

**1.9** Give an example of an exterior differential system whose characters are not constant.

## Example: Lagrangian submanifolds

We employ the Cartan–Kähler theorem to prove the existence of Lagrangian submanifolds of complex Euclidean space. Let

$$\vartheta := dx^1 \wedge dy^1 + dx^2 \wedge dy^2 + \cdots + dx^n \wedge dy^n.$$

Let $\mathcal{I}$ be the exterior differential system generated by $\vartheta$ on $M := \mathbb{R}^{2n}$. Writing spans of vectors in angle brackets,

| Flag | Polar equations | Characters |
|---|---|---|
| $E_0 = \{0\}$ | $\{0\}$ | $s_0 = 0$ |
| $E_1 = \langle \partial_{x^1} \rangle$ | $\langle dy^1 \rangle$ | $s_1 = 1$ |
| $\vdots$ | $\vdots$ | $\vdots$ |
| $E_n = \langle \partial_{x^1}, \partial_{x^2}, \dots, \partial_{x^n} \rangle$ | $\langle dy^1, dy^2, \dots, dy^n \rangle$ | $s_n = 1$ |

The flag predicts

$$\dim M + s_1 + 2\,s_2 + \cdots + n\,s_n = 2n + 1 + 2 + \cdots + n.$$

The nearby integral elements at a given point of $M$ are parameterized by $dy = a\,dx$, which we plug in to $\vartheta = 0$ to see that $a$ can be any symmetric matrix. So the space of integral elements has dimension

$$\dim M + \frac{n(n+1)}{2} = 2n + \frac{n(n+1)}{2},$$

correctly predicted. The Cartan–Kähler theorem proves the existence of Lagrangian submanifolds of complex Euclidean space, one (at least) through each point, tangent to each subspace $dy = a\,dx$, at least for any symmetric matrix $a$ close to 0.

**1.10** On a complex manifold $M$, take a Kähler form $\vartheta$ and a holomorphic volume form $\Psi$, i.e. closed forms expressed in local complex coordinates as

$$\vartheta = \frac{\sqrt{-1}}{2} g_{\mu\bar{\nu}} dz^{\mu} \wedge dz^{\bar{\nu}},$$
$$\Psi = f(z)\, dz^1 \wedge dz^2 \wedge \cdots \wedge dz^n,$$

with $f(z)$ a holomorphic function and $g_{\mu\bar{\nu}}$ a positive definite self-adjoint complex matrix of functions. Prove the existence of *special Lagrangian submanifolds*, i.e. integral manifolds of both $\vartheta$ and the imaginary part of $\Psi$.

## The last character

We always look for integral manifolds of a particular dimension $p$. For simplicity, we can assume that our exterior differential system contains all differential forms of degree $p+1$ and higher. In particular, the $p$-dimensional integral elements are maximal. The polar equations of any maximal integral element $E_p$ cut out precisely $E_p$, i.e. have rank $\dim M - p$ on $E_p$. We encounter $s_0, s_1, \dots, s_p$ polar equations at each increment, so rank $s_0 + s_1 + \dots + s_p$. We can find $s_p$ from the other characters:

$$s_0 + s_1 + \dots + s_{p-1} + s_p = \dim M - p.$$

For even greater simplicity, we take this as a definition for the final character $s_p$; we only *pretend* that our exterior differential system contains all differential forms of degree $p+1$ and higher.

## Example: harmonic functions

We will prove the existence of harmonic functions on the plane with given value and first derivatives at a given point. On $M = \mathbb{R}^5_{x,y,u,u_x,u_y}$, let $\mathcal{I}$ be the exterior differential system generated by

$$\begin{aligned} \vartheta &:= du - u_x\,dx - u_y\,dy, \\ \Theta &:= du_x \wedge dy - du_y \wedge dx. \end{aligned}$$

Note that

$$d\vartheta = dx \wedge du_x + dy \wedge du_y$$

also belongs to $\mathcal{I}$ because any exterior differential system is closed under exterior derivative. Any integral surface $X$ on which $0 \neq dx \wedge dy$ is locally the graph of a harmonic function $u = u(x,y)$ and its derivatives $u_x = \frac{\partial u}{\partial x}$, $u_x = \frac{\partial u}{\partial x}$.

Each integral plane $E_2$ (i.e. integral element of dimension 2) on which $dx \wedge dy \neq 0$ is given by equations

$$\begin{aligned} du_x &= u_{xx}dx + u_{xy}dy, \\ du_y &= u_{yx}dx + u_{yy}dy, \end{aligned}$$

for a unique choice of 4 constants $u_{xx}, u_{xy}, u_{yx}, u_{yy}$ subject to the 2 equations $u_{xy} = u_{yx}$ and $0 = u_{xx} + u_{yy}$. Hence integral planes at each point have dimension 2. The space of integral planes has dimension $\dim M + 2 = 5 + 2 = 7$.

Every integral line is the span $E_1 = \langle v \rangle$ of a vector

$$v = \left(\dot{x}, \dot{y}, u_x\dot{x} + u_y\dot{y}, u_{xx}\dot{x} + u_{xy}\dot{y}, u_{yx}\dot{x} + u_{yy}\dot{y}\right).$$

Compute

$$v \lrcorner \begin{pmatrix} d\vartheta \\ \Theta \end{pmatrix} = \begin{pmatrix} \dot{x}\,du_x + \dot{y}\,du_y - (u_{xx}\dot{x} + u_{xy}\dot{y})\;dx - (u_{yx}\dot{x} + u_{yy}\dot{y})\;dy \\ \dot{y}\,du_x - \dot{x}\,du_y + (u_{xx}\dot{x} + u_{xy}\dot{y})\;dy - (u_{yx}\dot{x} + u_{yy}\dot{y})\;dx \end{pmatrix}.$$

and

| Flag | Polar equations | Characters |
|---|---|---|
| $E_0=\{0\}$ | $\langle\vartheta\rangle$ | $s_0=1$ |
| $E_1=\langle v\rangle$ | $\langle\vartheta, v \lrcorner d\vartheta, v \lrcorner \Theta\rangle$ | $s_1=2$ |

We are only interested in finding integral surfaces; we compute the final character from

$$s_0+s_1+s_2=\dim M-2.$$

The characters are $s_0, s_1, s_2 = 1, 2, 0$ with predicted dimension $\dim M + s_1 + 2s_2 = 5+2+2\cdot 0 = 7$: involution. So harmonic functions exist near any point of the plane, with prescribed value and first derivatives at that point.

**1.11** What are the integral elements and integral manifolds of the exterior differential system generated by $du \wedge dx, dv \wedge dx, u\, du \wedge dv$ on $\mathbb{R}^4_{x,y,u,v}$?

## Background material from differential geometry

**1.12** Prove: constant rank linear equations, depending analytically on parameters, have a local solution depending analytically on those parameters.

**Lemma 1.2** (Cartan's lemma). *Suppose that $\xi_1, \xi_2, \dots, \xi_k$ are linearly independent analytic 1-forms on some manifold $M$. (They need not form a basis.) Suppose that $\alpha_1, \alpha_2, \dots, \alpha_k$ are analytic 1-forms satisfying $0 = \alpha_1 \wedge \xi_1 + \alpha_2 \wedge \xi_2 + \dots + \alpha_k \wedge \xi_k$. Then there is a unique symmetric matrix $A = (a_{ij})$ of analytic functions so that $\alpha_i = \sum a_{ij}\xi_j$.*

**1.13** Prove Cartan's lemma using the solution of problem 1.12.

## The Frobenius theorem

The *pullback* of an exterior differential system by a map is the exterior differential system generated by the pulled back forms. A exterior differential system is *locally generated* by forms with some property if every point lies in an open set so that the pullback to that open set is generated by forms with that property.

**Theorem 1.3** (Frobenius). *Suppose that $\mathcal{I}$ is an exterior differential system on a manifold $M$. The following are equivalent; if any, hence all, hold we say that $\mathcal{I}$ is* Frobenius.

a. *Every integral element of $\mathcal{I}$ lies in a unique $p$-dimensional integral element.*

b. *There is a constant so that $\mathcal{I}$ has a $p$-dimensional integral element at each point of $M$, with character $s_0$ equal to that constant, and all other characters $s_1, s_2, \dots, s_p$ vanishing.*

c. *$\mathcal{I}$ has an involutive $p$-dimensional integral element at each point of $M$. $\mathcal{I}$ is locally generated by $\dim M - p$ linearly independent 1-forms.*

d. *$\mathcal{I}$ is locally generated as an ideal of differential forms by* $\dim M - p$ *linearly independent* 1*-forms.*

e. *$\mathcal{I}$ is locally generated by linearly independent* 1*-forms* $\theta^1, \theta^2, \dots, \theta^q$ *with* $d\theta^i = -\sum_j \xi^i_j \wedge \theta^j$ *for some* 1*-forms* $\xi^i_j$.

f. *$\mathcal{I}$ is locally generated by linearly independent* 1*-forms. The vector fields on which all* 1*-forms in $\mathcal{I}$ vanish are closed under Lie bracket. They span a $p$-dimensional linear subspace in each tangent space of $M$.*

g. *Each point of $M$ has coordinates* $x^1, \dots, x^p, y^1, \dots, y^q$ *so that $\mathcal{I}$ is locally generated by* $dy^1, dy^2, \dots, dy^q$.

h. *The $p$-dimensional $\mathcal{I}$-integral manifolds form the leaves of a foliation $F$ of $M$. $\mathcal{I}$ is locally generated by the* 1*-forms vanishing on the leaves.*

**1.14** Use a statement of the Frobenius theorem from a differential geometry textbook, or lemma 7.5 on page 55, to prove the Frobenius theorem.

To solve both $u_x = f(x, y, u)$ and $u_y = g(x, y, u)$, note

$$u_{xy} = f_y + f_u u_y = u_{yx} = g_x + g_u u_x,$$

the compatibility condition between these equations. The exterior differential system generated by $\theta = du - f\,dx - g\,dy$ contains

$$d\theta = (f_y + f_u g - g_x - g_u f) dx \wedge dy + (f_u dx + g_u dy) \wedge \theta,$$

so takes the compatibility condition into account. By the Frobenius theorem, $f_u + f_u g = g_x + g_u f$ everywhere just when there are solutions with arbitrary initial condition $u(x_0, y_0) = u_0$.

## Example: surfaces with constant shape operator

In this section, we assume familiarity with appendix D. Use the Frobenius theorem:

**1.15** Prove that, for point $x \in \mathbb{E}^3$ and plane $P$ in $\mathbb{E}^3$ through $x$, any connected and embedded surface with zero shape operator passing through $x$ and tangent to $P$ at $x$ is an open subset of $P$.

**1.16** An oriented surface is $c_0$*-round*, for a constant $c_0 > 0$, if it has shape operator $I\!I(u, v) = c_0 u \cdot v$ at every point. Prove that, for each constant $c_0 \in \mathbb{R}$, point $x \in \mathbb{E}^3$ and oriented plane $P$ through $x$, there is a $c_0$-round oriented surface passing through $x$, with oriented tangent space at $x$ equal to $P$. If $c_0 \neq 0$, prove that any such surface, if connected and embedded, is an open subset of a sphere.

**1.17** Prove that any connected and embedded surface in $\mathbb{E}^3$ with constant principal curvatures is an open subset of a plane, sphere or cylinder.

**1.18** Prove that any connected embedded curve in $\mathbb{E}^3$ of constant curvature and torsion is an open subset of a line, circle or helix.

**1.19** An *umbilic point* is a point of a surface in $\mathbb{E}^3$ where the shape operator is a multiple of the dot product: $I\!I(u,v) = c\, u \cdot v$, for some number $c$. Find all connected and embedded surfaces in $\mathbb{E}^3$ which are everywhere umbilic, with perhaps varying value of $c$.

**1.20** Find every connected and embedded surface in $\mathbb{E}^3$ whose every geodesic is planar.

## Background material from differential geometry

A collection of 1-forms $\omega^1, \dots, \omega^n$ on a manifold $M$ *coframes* $M$, or is a *coframe* or *coframing* of $M$, if, at every point $m \in M$, $\omega^1_m, \dots, \omega^n_m$ are a basis of $T^*_m M$.

## Example: triply orthogonal webs

A *triply orthogonal web* is a triple of foliations by surfaces whose leaves are pairwise perpendicular.

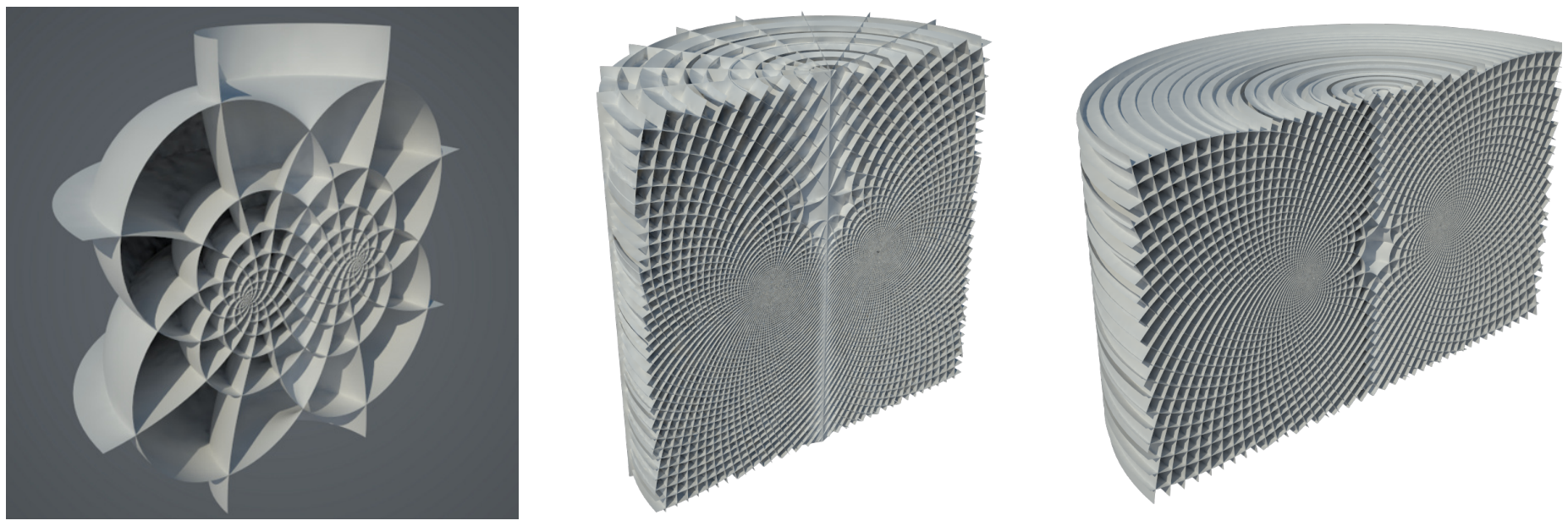

Images: Daniel Piker, 2015

**Theorem 1.4.** *There are infinitely many triply orthogonal webs, depending on 3 functions of 2 variables, defined near any given point of 3-dimensional Euclidean space.*

*Proof.* Picture a triply orthogonal web. Each leaf is perpendicular to a unique unit length 1-form $\eta_i$, up to $\pm$, which satisfies $0 = \eta_i \wedge d\eta_i$, by the Frobenius theorem. Denote 3-dimensional Euclidean space by $X$. Let $M$ be the set of all orthonormal bases of the tangent spaces of $X$, with obvious bundle map $x\colon M \to X$, so that each point of $M$ has the form $m = (x, e_1, e_2, e_3)$ for some $x \in X$ and orthonormal basis $e_1, e_2, e_3$ of $T_x X$. The *soldering* 1*-forms* $\omega_1, \omega_2, \omega_3$ on $M$ are defined by

$$(v \lrcorner \omega_i) e_i = x_* v.$$

Note: they are 1-forms on $M$, not on $X$. Let

$$\omega := \begin{pmatrix} \omega_1 \\ \omega_2 \\ \omega_3 \end{pmatrix}.$$

As explained in appendix D, there is a unique matrix-valued 1-form $\pi$, valued in antisymmetric $3 \times 3$ matrices and known as the *Levi-Civita connection*, so that $d\omega = -\pi \wedge \omega$, i.e.

$$d \begin{pmatrix} \omega_1 \\ \omega_2 \\ \omega_3 \end{pmatrix} = - \begin{pmatrix} 0 & -\pi_3 & \pi_2 \\ \pi_3 & 0 & -\pi_1 \\ -\pi_2 & \pi_1 & 0 \end{pmatrix} \wedge \begin{pmatrix} \omega_1 \\ \omega_2 \\ \omega_3 \end{pmatrix}.$$

Our triply orthogonal web is precisely a section $X \to M$ of the bundle map $M \to X$ on which $0 = \omega_i \wedge d\omega_i$ for all $i$, hence an integral 3-manifold of the exterior differential system $\mathcal{I}$ on $M$ generated by the closed 3-forms

$$\omega_1 \wedge d\omega_1, \quad \omega_2 \wedge d\omega_2, \quad \omega_3 \wedge d\omega_3.$$

Using the equations above, $\mathcal{I}$ is also generated by

$$\pi_3 \wedge \omega_{12}, \quad \pi_1 \wedge \omega_{23}, \quad \pi_2 \wedge \omega_{31}.$$

The integral manifolds coframed by $\omega_1, \omega_2, \omega_3$ are locally precisely the triply orthogonal webs. The reader can find the characters: $s_1, s_2, s_3 = 0, 3, 0$. The integral elements coframed by $\omega_1, \omega_2, \omega_3$ form a manifold of dimension 12 (parameterized by choice of a point of $M$ and 6 coefficients to determine values of $\pi_1, \pi_2, \pi_3$ at each point of $M$): involution. □

The reader familiar with Riemannian geometry will note that this proof works just as well for any 3-dimensional Riemannian manifold. For more on orthogonal webs in Euclidean space, see [7, 9, 35, 39].

**1.21*** A question of Alin Albu-Schäffer: suppose we are given an analytic foliation of an open subset of 3-dimensional Euclidean space by surfaces. Is there, at least locally, a triply orthogonal web which has this as one of its three foliations?

## Generality of integral manifolds

This section can be omitted without loss of continuity.

If I describe some family of submanifolds as integral manifolds of an exterior differential system, you might find a different description of the same submanifolds as integral manifolds of a different exterior differential system, perhaps on a different manifold, with different characters.

A function $f(x)$ of one variable is equivalent information to having a constant $f(0)$ and a function $df/dx$ of one variable.

We will solve $\nabla \times u = f - u$ for an unknown vector field $u$ in 3-dimensional Euclidean space two ways: an elementary approach on page 83, and following Cartan's strategy in problem 2.5 on page 21. The first approach uses 1 function of 1 variable and 2 functions of 2 variables. The second approach uses 3 constants, 3 functions of 1 variable and 2 functions of 2 variables.

Lagrangian submanifolds of $\mathbb{C}^n$ depend on 1 function of $n$ variables: those which are graphs $y = y(x)$ are of the form

$$y = \frac{\partial S}{\partial x}$$

for some potential function $S(x)$, unique up to adding a constant. On the other hand, the proof of the Cartan–Kähler theorem builds up each Lagrangian manifold from a choice of 1 function of 1 variable, 1 function of 2 variables, and so on.

We will see that Cartan's proof of the Cartan–Kähler theorem constructs the general solution of any involutive exterior differential system, using $s_0$ constants, $s_1$ functions of 1 variable, $s_2$ functions of 2 variables, and so on. We "trust" the last nonzero $s_p$ of an involutive exterior differential system to give the generality of integral manifolds: they depend on $s_p$ functions of $p$ variables, but we don't "trust" $s_0, s_1, \dots, s_{p-1}$.

Immersed plane curves are the integral curves of $\mathcal{I} = 0$ on $M = \mathbb{R}^2$. Any integral flag has $s_0, s_1 = 0, 1$. Immersed plane curves are also the integral curves of the ideal $\mathcal{J}$ generated by

$$\vartheta := \sin\phi \, dx - \cos\phi \, dy$$

on $M := \mathbb{R}^2_{x,y} \times S^1_\phi$. Here $s_0, s_1 = 1, 1$.

## What could go wrong?

On $M := \mathbb{R}^3_{x,y,z}$, take the ideals $\mathcal{I}$ generated by $dx \wedge dz, dy \wedge dz$ and $\mathcal{J}$ generated by $dx \wedge dz, dy \wedge (dz - y\,dx)$. At the origin, these differential forms are identical, so the integral elements and characters are the same at that point. But $\mathcal{I}$ has integral surfaces $z =$ constant, while $\mathcal{J}$ does not have integral surfaces.

**1.22** Compute characters and dimensions of spaces of integral elements for these. Show that the Cartan–Kähler theorem does not apply. Find their integral surfaces.

## Generalizations

The Cartan–Kähler theorem also holds for holomorphic exterior differential systems on complex manifolds, formal power series exterior differential systems, and Denjoy–Carleman exterior differential systems, with the same proof. In particular, any involutive smooth exterior differential system is also a formal power series system about each point, so has formal power series solutions, perhaps divergent. The theorem also holds for certain smooth systems [21, 22, 38].

# Chapter 2

# Tableaux

A *tableau* is a matrix encoding an exterior differential system, a computational tool to organize the linear algebra needed to uncover the characters.

Recall triply orthogonal webs on page 8: an exterior differential system generated by 1-forms $\theta^a$ and by 3-forms of the form

$$\omega^{12} \wedge \pi^3, \omega^{31} \wedge \pi^2, \omega^{23} \wedge \pi^1,$$

with coframe $\theta^a, \omega^i, \pi^\alpha$. Write these 3-forms as rows of a matrix wedge product:

$$\begin{pmatrix} \pi^3 & 0 & 0 \\ 0 & \pi^2 & 0 \\ 0 & 0 & \pi^1 \end{pmatrix} \wedge \begin{pmatrix} \omega^{12} \\ \omega^{13} \\ \omega^{23} \end{pmatrix}$$

## Definition

Take an exterior differential system $\mathcal{I}$ on a manifold $M$. We are looking for $p$-dimensional integral manifolds. Take a point $m \in M$. Let

$$\mathcal{I}_m := \{ \vartheta_m \in \Lambda^* T^*_m M \mid \vartheta \in \mathcal{I} \}.$$

We carry out all of our work below modulo the ideal $(\mathcal{I}^1_m) \subseteq \mathcal{I}_m$.

*a.* Take a basis $\omega^i, \pi^\alpha$ of $T^*_m M / \mathcal{I}^1_m$.

*b.* Write a set of forms generating $\mathcal{I}_m / (\mathcal{I}^1_m)$ in a column

$$\vartheta = \begin{pmatrix} \vartheta^1 \\ \vartheta^2 \\ \vdots \\ \vartheta^r \end{pmatrix}.$$

*c.* Let $\omega^{ij} := \omega^i \wedge \omega^j$, and so on. Let

$$\omega = \begin{pmatrix} \omega^1 \\ \hline \omega^2 \\ \omega^{12} \\ \hline \omega^3 \\ \omega^{13} \\ \omega^{23} \\ \omega^{123} \\ \hline \vdots \\ \hline \omega^p \\ \vdots \\ \omega^{1\cdots p} \end{pmatrix},$$

a column vector of all wedge products of the $\omega^i$. Arrange by grades: grade $j$ consists of all forms $\omega^{\cdots j}$ which are wedge products of 1-forms from among $\omega^1, \dots, \omega^j$, and must contain $\omega^j$. We sometimes mark grades with horizontal lines. Drop any entry of $\omega$ if it doesn't appear in the forms $\vartheta^a$ in our generating set, expanded in our basis.

*d.* Write out $\vartheta = \varpi \wedge \omega + \dots$ for a matrix $\varpi$, the *tableau*, so that each matrix entry is a linear combination of $\pi^\alpha$, and the $\dots$ consists of

1. terms with only $\omega$ in them, no $\pi^\alpha$ 1-form wedged into them, the *torsion* and
2. terms with two or more $\pi^\alpha$ 1-forms wedged into them, the *nonlinearity.*

We sometimes draw vertical lines in $\varpi$, marking out grades at widths matching the grade heights in $\omega$. The nonlinearity we assign grade $p$.

*e.* A *polar* is an entry of $\varpi$ linearly independent of all entries found in all earlier grades and above or to the left in the same grade. Highlight all polars. In practice, each polar is often one of the $\pi^\alpha$ of our basis; we can always change basis to arrange this.

*f.* Take any basis $\theta^a$ of $\mathcal{I}^1_m$; declare the basis elements to be polars of grade zero.

The *character* $s_j$ of the tableau is the number of polars in grade $j$.

Suppose $\mathcal{I}$ is an exterior differential system spanned by

1-forms $\theta^1, \theta^2, \theta^3$,
a 2-form $\theta^1 \wedge \omega^3 + \omega^1 \wedge \pi^1 + \omega^2 \wedge \pi^2 + \omega^3 \wedge \pi^3$ and
a 3-form $\pi^{123} - \omega^{12} \wedge \pi^3 + \omega^{13} \wedge \pi^2 - \omega^{23} \wedge \pi^1$,

for a coframing

$$\theta^1, \theta^2, \theta^3, \omega^1, \omega^2, \omega^3, \pi^1, \pi^2, \pi^3.$$

Dropping multiples of the $\theta^a$:

no 1-forms
a 2-form $\omega^1 \wedge \pi^1 + \omega^2 \wedge \pi^2 + \omega^3 \wedge \pi^3$ and
a 3-form $-\omega^{12} \wedge \pi^3 + \omega^{13} \wedge \pi^2 - \omega^{23} \wedge \pi^1 + \dots$.

So modulo $\theta^1, \theta^2, \theta^3$:

$$\vartheta = \begin{pmatrix} \omega^1 \wedge \pi^1 + \omega^2 \wedge \pi^2 + \omega^3 \wedge \pi^3 \\ -\omega^{12} \wedge \pi^3 + \omega^{13} \wedge \pi^2 - \omega^{23} \wedge \pi^1 \end{pmatrix} + \begin{pmatrix} 0 \\ \pi^{123} \end{pmatrix},$$

$$= -\begin{pmatrix} \pi^1 & \pi^2 & 0 & \pi^3 & 0 & 0 \\ 0 & 0 & \pi^3 & 0 & -\pi^2 & \pi^1 \end{pmatrix} \wedge \begin{pmatrix} \omega^1 \\ \omega^2 \\ \omega^{12} \\ \omega^3 \\ \omega^{13} \\ \omega^{23} \end{pmatrix} + \begin{pmatrix} 0 \\ \pi^{123} \end{pmatrix}.$$

$$\begin{matrix} s_1 & s_2 & s_3 \\ 1 & 2 & 0 \end{matrix}$$

Recall that Lagrangian submanifolds are integral manifolds of

$$\vartheta := dx^1 \wedge dy^1 + dx^2 \wedge dy^2 + \dots + dx^n \wedge dy^n,$$

$$= -\begin{pmatrix} \pi^1 & \pi^2 & \dots & \pi^n \end{pmatrix} \wedge \begin{pmatrix} \omega^1 \\ \omega^2 \\ \vdots \\ \omega^n \end{pmatrix}$$

$$\begin{matrix} s_1 & s_2 & \dots & s_n \\ 1 & 1 & \dots & 1 \end{matrix}$$

with $\omega^i := dx^i$, $\pi^i := dy^i$.

On $M = \mathbb{R}^5_{x,y,u,u_x,u_y}$, let $\mathcal{I}$ be generated by

$$\theta := du - u_x\,dx - u_y\,dy,$$
$$\vartheta := du_x \wedge dy - du_y \wedge dx,$$

and note

$$d\theta = -du_x \wedge dx - du_y \wedge dy$$

also belongs to $\mathcal{I}$. An integral surface $X$ on which $0 \neq dx \wedge dy$ is locally the graph of a harmonic function $u = u(x, y)$ and its derivatives $u_x = \frac{\partial u}{\partial x}$, $u_y = \frac{\partial u}{\partial y}$. Taking

$$\omega^1 := dx, \omega^2 := dy, \pi^1 := du_x, \pi^2 := -du_y,$$

gives tableau

$$\begin{pmatrix} d\theta \\ \vartheta \end{pmatrix} = -\begin{pmatrix} \pi^1 & -\pi^2 \\ \pi^2 & \pi^1 \end{pmatrix} \wedge \begin{pmatrix} \omega^1 \\ \omega^2 \end{pmatrix}.$$

**2.1** As in our example of triply orthogonal webs, construct an exterior differential system whose integral 4-manifolds are foliations of open subsets of 3-dimensional Euclidean space. Write out the tableau and find the characters.

## Torsion

The tableau is *adapted* to the flag

$$E_j := (0 = \theta^a = \pi^\alpha = \omega^{j+1} = \cdots = \omega^p).$$

Any flag has many tableau adapted to it.

**2.2** What are the polar equations of each $E_j$? Prove that the characters as defined above are the characters as defined in chapter 1, and the polars are a basis of the polar equations.

**2.3** Prove that torsion vanishes just when the flag is integral.

Take an exterior differential system generated by 1-forms $\theta^1, \theta^2$ with

$$d\begin{pmatrix} \theta^1 \\ \theta^2 \end{pmatrix} = -\begin{pmatrix} \pi^1 & 0 \\ \pi^2 & \pi^3 \end{pmatrix} \wedge \begin{pmatrix} \omega^1 \\ \omega^2 \end{pmatrix} + \begin{pmatrix} \omega^{12} \\ 0 \end{pmatrix} \mod \theta^1, \theta^2.$$

Note that $d\theta^1, d\theta^2$ generate $\mathcal{I}/(\mathcal{I}^1)$. The term $\omega^{12}$ is the torsion, as it has no $\pi$ in it.

Let $\pi'^1 := \pi^1 + \omega^2$:

$$d\begin{pmatrix}\theta^1\\ \theta^2\end{pmatrix} = -\begin{pmatrix}\pi'^1 & 0\\ \pi^2 & \pi^3\end{pmatrix} \wedge \begin{pmatrix}\omega^1\\ \omega^2\end{pmatrix} \quad \operatorname{mod} \theta^1, \theta^2;$$

we *absorb* the torsion.

Changing bases to arrange that torsion vanishes is *absorbing the torsion* A tableau can only examine integral elements coframed by the $\omega^i$. Torsion absorbs just when there is such an integral element.

The exterior differential system $\mathcal{I}$ is the true fundamental geometric object; the choice of tableau is like a choice of coordinate system: a magnifying glass with which to examine $\mathcal{I}$.

## Borrowing

Recall triply orthogonal webs on page 8 had exterior differential system generated by 1-forms $\theta^a$ and by 3-forms of the form

$$\omega^{12} \wedge \pi^3, \omega^{31} \wedge \pi^2, \omega^{23} \wedge \pi^1,$$

with coframe $\theta^a, \omega^i, \pi^\alpha$:

$$\begin{array}{c}\left(\begin{array}{|c|cc}\pi^3 & 0 & 0\\ 0 & \pi^2 & 0\\ 0 & 0 & \pi^1\end{array}\right) \wedge \begin{pmatrix}\overline{\omega^{12}}\\ \omega^{13}\\ \omega^{23}\end{pmatrix}\\ \begin{array}{cc} s_2 & s_3\\ 1 & 2\end{array}\end{array}$$

*Warning:* these $s_1, s_2, s_3$ are *not* the characters we computed in chapter 1.

Reconsider the same example, with new choices of 1-forms $\omega^i$. Let $\omega'^1 := \omega^1$, $\omega'^2 := \omega^2$, $\omega'^3 := \omega^1 - \omega^2 + \omega^3$. Write this as

$$\begin{pmatrix}\omega^1\\ \omega^2\\ \omega^3\end{pmatrix} = \begin{pmatrix}\omega'^1\\ \omega'^2\\ -\omega'^1 + \omega'^2 + \omega'^3\end{pmatrix}.$$

In these 1-forms, the tableau is:

$$
\begin{pmatrix}
\pi^3 & 0 & 0 \\
\pi^2 & 0 & 0 \\
\pi^1 & 0 & 0
\end{pmatrix}
\wedge
\begin{pmatrix}
\omega'^{12} \\
\hline
\omega'^{13} \\
\omega'^{23}
\end{pmatrix}
$$

$$
\begin{array}{cc}
s_2 & s_3 \\
3 & 0
\end{array}
$$

These *are* the characters we computed in chapter 1. We have *borrowed* polars from later grades into earlier. We change the choice of integral flag from

$$
\begin{aligned}
E_1 &= (0 = \omega^2 = \omega^3 = \theta^a = \pi^\alpha), \\
E_2 &= (0 = \omega^3 = \theta^a = \pi^\alpha), \\
E_3 &= (0 = \theta^a = \pi^\alpha),
\end{aligned}
$$

to

$$
\begin{aligned}
E'_1 &= (0 = \omega'^2 = \omega'^3 = \theta^a = \pi^\alpha) \\
E'_2 &= (0 = \omega'^3 = \theta^a = \pi^\alpha), \\
E_3 &= (0 = \theta^a = \pi^\alpha) \text{ unchanged.}
\end{aligned}
$$

**2.4** Prove that the characters depend only on the flag.

Borrowing as many polars as we can, and perhaps permuting the coframe, the integral flag has largest $s_1$ among all integral flags at that point. Subject to that $s_1$, it has largest $s_2$, and so on.

Take a row which represents a $k$-form. Permuting the $\omega^i$, we can get any polar in that row to appear in grade $k-1$: $\pi^\alpha \wedge \omega^{1\dots k-1}$. If there is another polar in that row, say in grade $\ell$, add a multiple of $\omega^k$ to $\omega^\ell$ to borrow it to grade $k$. Continue in this way: for one particular row, representing a $k$-form in $\mathcal{I}$, we arrange polars in successive grades, starting at grade $k-1$, all followed by any nonpolar entries in that row.

Write wedge products $\omega^{i_1\dots i_q}$ with $i_1 < \dots < i_q$. Order any two wedge products by last entry $i_q$, then by next to last, and so on. Borrow to have polars arising in sorted order before any other entries.

Since this occurs for some linear transformation of $\omega^i$, it also occurs for all linear transformations of $\omega^i$ except for those with certain minors vanishing. We can thus borrow simultaneously for all rows, by generic linear transformation.

$$\left(\begin{array}{c|c} \pi^1 & 0 \\ 0 & \pi^1 \end{array}\right) \wedge \left(\frac{\omega^{12}}{\omega^3}\right)$$

has a polar appearing in grade 2. Permuting indices 1 and 3:

$$\left(\begin{array}{c|c} 0 & \pi^1 \\ \pi^1 & 0 \end{array}\right) \wedge \left(\frac{\omega^{1}}{\omega^{23}}\right)$$

puts it in grade 1.

The torsion is absorbable just when there is an integral element $E = (\pi = p\omega)$. Absorb the torsion by subtracting $p\omega$ from $\pi$. Thus there is a torsion free tableau just when there is a generic torsion free tableau.

Take a tableau

$$\begin{pmatrix} \pi^1 & 0 \\ \pi^2 & \pi^3 \end{pmatrix} \wedge \begin{pmatrix} \omega^1 \\ \omega^2 \end{pmatrix} + \begin{pmatrix} \pi^4 \wedge \pi^2 \\ 0 \end{pmatrix}$$

with a polar in the nonlinearity. Add $\omega^2$ to $\pi^2$ to get the polar to appear in the tableau:

$$\begin{pmatrix} \pi^1 & \pi^4 \\ \pi'^2 & \pi^3 \end{pmatrix} \wedge \begin{pmatrix} \omega^1 \\ \omega^2 \end{pmatrix} - \begin{pmatrix} 0 \\ \omega^{12} \end{pmatrix} + \begin{pmatrix} \pi^4 \wedge \pi'^2 \\ 0 \end{pmatrix}.$$

This produces torsion, but we can absorb it.

$$\begin{pmatrix} \pi^1 & \pi^2 \end{pmatrix} \wedge \begin{pmatrix} \omega^1 \\ \omega^2 \end{pmatrix} + \pi^3 \wedge \pi^2$$

has a polar in the nonlinearity, but we absorb it by $\omega'^2 := \omega^2 - \pi^3$.

Some tableaux have nonabsorbable polars in the nonlinearity:

$$\begin{pmatrix} \pi^1 & \pi^2 \end{pmatrix} \wedge \begin{pmatrix} \omega^1 \\ \omega^2 \end{pmatrix} + \pi^3 \wedge \pi^4.$$

Any polar in the nonlinearity can be demoted to a new $\omega^i$ 1-form, but we won't need to do so.

## Integral elements

Picture a tableau

$$
\begin{pmatrix}
\pi^1 & \pi^4 & \pi^6 \\
\pi^2 & \pi^5 & \pi^1 - \pi^5 \\
\pi^3 & \pi^2 & \pi^1 + \pi^2 \\
\pi^1 - \pi^2 & 0 & \pi^3
\end{pmatrix}
\wedge
\begin{pmatrix}
\omega^1 \\ \omega^2 \\ \omega^3
\end{pmatrix}.
$$

Take any 3-dimensional linear subspace $E$ of a tangent space of $M$ coframed by $\omega^1, \omega^2, \omega^3$ and on which $\theta^a = 0$. Then $E$ has coefficients $\pi^1 = p^1_1 \omega^1 + p^1_2 \omega^2 + p^1_3 \omega^3$ and so on. Plug in to the tableau to find equations for integral elements. Since $\omega^{ij} = -\omega^{ji}$, each tableau entry in column $i$ has coefficient of $\omega^j$ exactly equal to the tableau entry in column $j$ coefficient of $\omega^i$:

$$
\begin{aligned}
p^1_2 &= p^4_1, & p^1_3 &= p^6_1, & p^4_3 &= p^6_2, \\
p^2_2 &= p^5_1, & p^2_3 &= p^1_1 - p^5_1, & p^5_3 &= p^1_2 - p^5_2, \\
p^3_2 &= p^2_1, & p^3_3 &= p^1_1 + p^2_1, & p^2_3 &= p^1_2 + p^2_2, \\
p^1_2 - p^2_2 &= 0, & p^1_3 - p^2_3 &= p^3_1, & 0 &= p^3_2.
\end{aligned}
$$

Again imagine an exterior differential system generated by 1-forms $\theta^a$ and by 3-forms of the form

$$
\omega^{12} \wedge \pi^1, \omega^{31} \wedge \pi^2, \omega^{23} \wedge \pi^3,
$$

with coframe $\theta^a, \omega^i, \pi^\alpha$. (E.g. triply orthogonal webs; see page 8.) The tableau:

$$
\begin{pmatrix}
0 & 0 & \pi^1 \\
0 & \pi^2 & 0 \\
\pi^3 & 0 & 0
\end{pmatrix}
\wedge
\begin{pmatrix}
\omega^{12} \\ \omega^{13} \\ \omega^{23}
\end{pmatrix}
$$

| $s_2$ | $s_3$ |
|---|---|
| 1 | 2 |

3-dimensional integral elements:

$$
\begin{pmatrix}
\pi^1 \\ \pi^2 \\ \pi^3
\end{pmatrix}
=
\begin{pmatrix}
0 & p^1_2 & p^1_3 \\
p^2_1 & 0 & p^2_3 \\
p^3_1 & p^3_2 & 0
\end{pmatrix}
\begin{pmatrix}
\omega^1 \\ \omega^2 \\ \omega^3
\end{pmatrix}
$$

a 6-dimensional space of integral elements at each point.

$$
s_1 + 2s_2 + 3s_3 = 0 + 2(1) + 3(2) = 8 > 6,
$$

involution fails.

The same example, but borrow:

$$\begin{pmatrix}\omega^1\\ \omega^2\\ \omega^3\end{pmatrix} = \begin{pmatrix}\omega'^1\\ \omega'^2\\ \omega'^1+\omega'^2+\omega'^3\end{pmatrix},$$

yielding tableau:

$$\begin{pmatrix} -\pi^1 & 0 & \pi^3 \\ -\pi^2 & \pi^2 & 0 \\ \pi^3 & 0 & 0 \end{pmatrix} \wedge \begin{pmatrix}\omega'^{12}\\ \omega'^{13}\\ \omega'^{23}\end{pmatrix}$$

$$\begin{matrix} s_2 & & s_3 \\ 3 & & 0 \end{matrix}$$

$$s_1 + 2s_2 + 3s_3 = 0 + 2(3) + 3(0) = 6,$$

involution: there are 3-dimensional integral manifolds.

**2.5** Write the equation $\nabla \times u = f - u$, which we unravel in detail on page 83, as an exterior differential system. Find a tableau. Can you absorb torsion? What submanifold contains the integral manifolds? Is the exterior differential system in involution on that submanifold?

## Example: Lie's third theorem

This section can be omitted without loss of continuity.

In this section we assume familiarity with Lie groups [33]. Lie's third theorem: every Lie algebra $\mathfrak{g}$, say of dimension $p$, is isomorphic to a Lie algebra of vector fields spanning every tangent space on a $p$-dimensional manifold. This theorem is a first step in constructing a Lie group with a given Lie algebra. Since we employ differential forms, it is easier for us to prove the dual statement about the dual 1-forms to those vector fields. A *Maurer–Cartan form* is a 1-form $\xi$ valued in a Lie algebra $\mathfrak{g}$, defined on an open subset $U \subset \mathbb{R}^p$, so that, at every point $x \in U$, $\xi_x \colon \mathbb{R}^p \to \mathfrak{g}$ is a linear isomorphism and $d\xi + \frac{1}{2}[\xi\xi] = 0$.

**2.6** Explain how a Maurer–Cartan form determines such vector fields and vice versa.

**Theorem 2.1** (Lie III). *Any Lie algebra has a Maurer–Cartan form.*

*Proof.* Choose a basis $e_1, \dots, e_p \in \mathfrak{g}$ and write out the Lie bracket in the basis: $[e_i e_j] = c^k_{ij} e_k$. Any such $\xi$ will then be $\xi = \xi^i e_i$, for a coframing $\xi^i = g^i_j(x) dx^j$, with $g^i_j(x)$ an invertible matrix for each $x \in U$. We want these $\xi^i$ to satisfy $0 = d\xi^i + \frac{1}{2} c^{jk}_i \xi^j \wedge \xi^k$. Let $M := \mathbb{R}^p \times \mathrm{GL}(p, \mathbb{R})$ with coordinates $x^i, g^i_j$. Take the exterior differential system $\mathcal{I}$ generated by the components $\vartheta^i$ of the 2-form

$$\vartheta = d(g\,dx) + \frac{1}{2}\left[g\,dx, g\,dx\right].$$

**2.7** Use the Jacobi identity, either in components or directly, to see that $0 = d\vartheta$.

Let $\omega^i := dx^i$.

**2.8** If we want $\vartheta^i = \pi^i_j \wedge \omega^j$, then $\pi^i_j = dg^i_j + q^i_{jk}(c) dx^k$ where $q = q(c)$ is some quadratic polynomial expression in the coefficients $c^k_{ij}$. Prove this by computing $q(c)$.

$$\begin{pmatrix} \vartheta^1 \\ \vartheta^2 \\ \vdots \\ \vartheta^p \end{pmatrix} = \begin{pmatrix} \pi^1_1 & \pi^1_2 & \dots & \pi^1_p \\ \pi^2_1 & \pi^2_2 & \dots & \pi^2_p \\ \vdots & \vdots & \dots & \vdots \\ \pi^p_1 & \pi^p_2 & \dots & \pi^p_p \end{pmatrix} \wedge \begin{pmatrix} \omega^1 \\ \omega^2 \\ \vdots \\ \omega^p \end{pmatrix}.$$

$$\begin{matrix} s_1 & s_2 & \dots & s_p \\ p & p & \dots & p \end{matrix}$$

Integral elements of dimension $p$: $\pi^i_j = p^i_{jk} \omega^k$ with $p^i_{jk} = p^i_{kj}$, so the space of integral elements has dimension $\dim M + p(p+1)/2$: involution. □

**2.9** Uniqueness: prove that any two Maurer–Cartan forms for the same Lie algebra are locally identified by a diffeomorphism.

## Example: surface invariants

This section can be omitted without loss of continuity.

We suppose that the reader has read appendix D. Does every quadratic form on a plane through a point arise as the shape operator on the tangent plane of some surface? If we had such a surface, its adapted frame bundle would satisfy

$$\begin{aligned} \omega_3 &= 0, \\ \gamma_{13} &= a_{11}\omega_1 + a_{12}\omega_2, \\ \gamma_{23} &= a_{21}\omega_1 + a_{22}\omega_2, \end{aligned}$$

with $a_{ij} = a_{ji}$. Let $V$ be the set of all symmetric $2 \times 2$ matrices, with typical element written as

$$a = \begin{pmatrix} a_{11} & a_{12} \\ a_{21} & a_{22} \end{pmatrix}.$$

On the manifold $V \times \ulcorner\mathbb{E}^3$, take the exterior differential system $\mathcal{I}$ generated by

$$\omega_3, \gamma_{13} - a_{11}\omega_1 - a_{12}\omega_2, \gamma_{23} - a_{21}\omega_1 - a_{22}\omega_2.$$

Every surface in $\mathbb{E}^3$ has frame bundle an integral manifold.

**2.10** Prove that any $\mathcal{I}$-integral manifold coframed by $\omega_1, \omega_2, \gamma_{12}$ is locally the frame bundle of a surface.

Summing over $i, j, k, \ell = 1, 2$, let

$$Da_{ij} := da_{ij} - a_{ik}\gamma_{kj} - a_{jk}\gamma_{ki},$$

$$d\begin{pmatrix} \omega_3 \\ \gamma_{13} - a_{11}\omega_1 - a_{12}\omega_2 \\ \gamma_{23} - a_{21}\omega_1 - a_{22}\omega_2 \end{pmatrix} = \begin{pmatrix} 0 & 0 & 0 \\ Da_{11} & Da_{12} & 0 \\ Da_{12} & Da_{22} & 0 \end{pmatrix} \wedge \begin{pmatrix} \omega_1 \\ \omega_2 \\ \gamma_{12} \end{pmatrix},$$

$$\begin{array}{ccc} s_1 & s_2 & s_3 \\ 2 & 1 & 0 \end{array}$$

Integral elements coframed by $\omega_1, \omega_2, \gamma_{12}$ are

$$Da_{ij} = a_{ijk}\omega_k + A_{ij}\gamma_{12}.$$

Plug into the tableau: integral elements have $a_{ijk}$ symmetric in all 3 indices, and $A_{ij} = 0$. Each integral element is identified with $I\!I\!I := e_3 a_{ijk}\omega_i\omega_j\omega_k$, the third fundamental form on any integral manifold arising from a surface. The space of integral elements is 4-dimensional at each point of our manifold, involution: there is an integral manifold coframed by $\omega_1, \omega_2, \gamma_{12}$ through every point of $V \times \ulcorner\mathbb{E}^3$.

**Theorem 2.2.** *Shape operators are arbitrary. To be precise, take a point of $\mathbb{E}^3$, a plane through that point, a symmetric bilinear form and a symmetric trilinear form on that plane, valued in the perpendicular line. There is a surface in $\mathbb{E}^3$ through that point, tangent to that plane, and with that bilinear form as shape operator and that trilinear form as third fundamental form.*

**2.11** Prove that there are surfaces in 3-dimensional Euclidean space preserved by no rigid motions except the identity.

**2.12** Are the shape operators of the leaves of a triply orthogonal web arbitrary?

Chapter 3

# Example: almost complex structures

This chapter can be omitted without loss of continuity. As an application of the Cartan–Kähler theorem, we demonstrate the integrability of torsion free almost complex structures.

## Example: the Cauchy–Riemann equations

This section can be omitted without loss of continuity.

Take several complex variables $z^1, \dots, z^n$, with real and imaginary parts $z^\mu = x^\mu + iy^\mu$. A *holomorphic* function is a complex valued function $f(z^1, \dots, z^n)$ satisfying the *Cauchy–Riemann equations* [8, 15, 18, 32, 34]

$$\frac{\partial f}{\partial x^\mu} = i\frac{\partial f}{\partial y^\mu}, \quad \mu = 1, 2, \dots, n.$$

In real and imaginary parts $f = u + iv$, this becomes

$$\begin{aligned} \frac{\partial u}{\partial x^\mu} &= \frac{\partial v}{\partial y^\mu}, \\ \frac{\partial v}{\partial x^\mu} &= -\frac{\partial u}{\partial y^\mu}. \end{aligned}$$

Let $M := \mathbb{R}^{4n+2}_{x^\mu, y^\mu, u, v, p_\mu, q_\mu}$ on which we take the exterior differential system $\mathcal{I}$ generated by

$$\begin{aligned} \theta^1 &:= du - p_\mu dx^\mu - q_\mu dy^\mu, \\ \theta^2 &:= dv + q_\mu dx^\mu - p_\mu dy^\mu. \end{aligned}$$

$$d\begin{pmatrix}\theta^1 \\ \theta^2\end{pmatrix} = - - \begin{pmatrix} dp_1 & dp_2 & \dots & dp_n & dq_1 & \dots & dq_n \\ -dq_1 & -dq_2 & \dots & -dq_n & dp_1 & \dots & dp_n \end{pmatrix} \wedge \begin{pmatrix} dx^1 \\ dx^2 \\ \vdots \\ dx^n \\ dy^1 \\ \vdots \\ dy^n \end{pmatrix}$$

$$\begin{matrix} s_1 & s_2 & \dots & s_n & s_{n+1} & \dots & s_{2n} \\ 2 & 2 & \dots & 2 & 0 & \dots & 0 \end{matrix}$$

So

$$\dim M + s_1 + 2s_2 + \cdots + 2ns_{2n} = 4n + 2 + n(n+1).$$

Integral elements at each point $(x, y, u, v, p, q)$ are

$$\begin{pmatrix} dp_\mu \\ dq_\mu \end{pmatrix} = \begin{pmatrix} r_{\mu\nu} & -s_{\mu\nu} \\ s_{\mu\nu} & r_{\mu\nu} \end{pmatrix} \begin{pmatrix} dx^\nu \\ dy^\nu \end{pmatrix}$$

with $r, s$ symmetric in lower indices. So the space of integral elements has dimension $4n + 2 + n(n+1)$, involution with general solution depending on 2 functions of $n$ variables.

For a complex notation, let $\omega^\mu := dx^\mu + i\,dy^\mu$, $\theta := \theta^1 + i\theta^2$, $\pi_\mu := dp_\mu + i\,dq_\mu$ and $\omega^{\bar\mu} := \overline{\omega^\mu} = dx^\mu - i\,dy^\mu$. A tableau in complex differential forms will have terms expressed in wedge products with $\omega^\mu, \omega^{\bar\mu}$. But for the Cauchy–Riemann equations we find

$$d\theta = -\begin{pmatrix} \pi^1 & \pi^2 & \dots & \pi^n \end{pmatrix} \wedge \begin{pmatrix} \omega^1 \\ \omega^2 \\ \vdots \\ \omega^n \end{pmatrix}$$

with no $\omega^{\bar\mu}$ terms. The characters count complex polars, so double them to count real and imaginary parts: $s_0, s_1, s_2, \dots, s_n = 2, 2, \dots, 2$. Integral elements are $\pi_\mu = p_{\mu\nu}\omega^\nu$, with complex numbers $p_{\mu\nu} = p_{\nu\mu}$, so $n(n+1)/2$ complex numbers, hence $n(n+1)$ real and imaginary parts.

## Example: almost complex structures

This section can be omitted without loss of continuity.

We can turn a real vector space into a complex vector space, just by picking a real linear map $J$ so that $J^2 = -I$. An *almost complex structure* on an even dimensional manifold $M$ is a choice of complex vector space structure on each tangent space, analytically varying, in that is described by an analytic map $J\colon TM \to TM$, acting linearly on each tangent space, with $J^2 = -I$.

Complex Euclidean space has its usual almost complex structure $J(x, y) = (-y, x)$, preserved by biholomorphisms (i.e. holomorphic diffeomorphisms) between open sets, as their derivatives are complex linear maps. Any complex manifold $M$ has an almost complex structure, defined by using holomorphic coordinates to identify locally with complex Euclidean space.

On an even dimensional manifold, a *complex coframing* is a collection of complex valued 1-forms $\omega^\mu$ so that their real and imaginary parts coframe. It is *complex linear* for an almost complex structure if each $\omega^\mu$ is a complex linear map on each tangent space, i.e. $\omega^\mu \circ J = \sqrt{-1}\omega^\mu$ for every $\mu$.

**3.1** Prove that every almost complex structure has, near each point, some complex linear complex coframing.

**3.2** Prove that every complex coframing is complex linear for a unique almost complex structure.

Complex coframings are a useful way to exhibit almost complex structures.

The complex coframing

$$\begin{aligned}\omega^1 &:= dz,\\ \omega^2 &:= dw - \bar{w}\, d\bar{z}\end{aligned}$$

yields a unique almost complex structure on the space parameterized by two complex variables $z, w$.

**3.3** Prove that any two complex coframings $\omega^\mu, \omega'^\mu$ yield the same almost complex structure just when $\omega'^\mu = g^\mu_\nu \omega^\nu$ for a unique matrix $g = (g^\mu_\nu)$ of functions.

Any complex differential form is expressed in any complex coframing $\omega^\mu$ as a sum of wedge products of $\omega^\mu$ and $\bar{\omega}^\mu$. Following convention, write $\bar{\omega}^\mu$ as $\omega^{\bar{\mu}}$. A $(p, q)$-*form* is a differential form expressed with $p$ factors of $\omega^\mu$ and $q$ of $\omega^{\bar{\mu}}$. For example, $\omega^\mu$ is $(1, 0)$ while $\omega^{\bar{\mu}}$ is $(0, 1)$. In particular,

$$d\omega^\mu = t^\mu_{\mu\sigma}\omega^\nu \wedge \omega^\sigma + t^\mu_{\nu\bar{\sigma}}\omega^\nu \wedge \omega^{\bar{\sigma}} + t^\mu_{\bar{\nu}\sigma}\omega^{\bar{\nu}} \wedge \omega^\sigma + t^\mu_{\bar{\nu}\bar{\sigma}}\omega^{\bar{\nu}} \wedge \omega^{\bar{\sigma}},$$

for unique complex valued functions $t$, antisymmetric in lower indices. The exterior derivative of any $(p, q)$-form is uniquely expressed as a sum of forms with $(p, q)$ raised by $(2, -1), (1, 0), (0, 1)$ or $(-1, 2)$. We thus split up

$$d = \tau + \partial + \bar{\partial} + \bar{\tau}:$$

$$\begin{aligned}\tau\omega^\mu &= 0, & \partial\omega^\mu &= t^\mu_{\nu\sigma}\omega^\nu \wedge \omega^\sigma,\\ \bar{\tau}\omega^\mu &= t^\mu_{\bar{\nu}\bar{\sigma}}\omega^{\bar{\nu}} \wedge \omega^{\bar{\sigma}}, & \bar{\partial}\omega^\mu &= t^\mu_{\nu\bar{\sigma}}\omega^\nu \wedge \omega^{\bar{\sigma}} + t^\mu_{\bar{\sigma}\nu}\omega^{\bar{\sigma}} \wedge \omega^\nu.\end{aligned}$$

Let $\omega := (\omega^\mu)$. Expanding out in the coframing, we find that $0 = \tau = \bar{\tau}$ on all differential forms if and only if $\bar{\tau}\omega = 0$. Change coframing: replace $\omega$ by some $g\omega$.

$$d(g\omega) = dg \wedge \omega + g d\omega,$$

expands out to

$$\begin{aligned}\tau(g\omega) &= 0, & \partial(g\omega) &= (\partial g) \wedge \omega + g\,\partial\omega,\\ \bar{\tau}(g\omega) &= g\,\bar{\tau}\omega, & \bar{\partial}(g\omega) &= (\bar{\partial} g) \wedge \omega + g\,\bar{\partial}\omega.\end{aligned}$$

Hence $\bar{\tau}\omega = 0$ just when $\bar{\tau}(g\omega) = 0$. So vanishing of $\bar{\tau}$ is a property of the almost complex structure.

**3.4** Prove that $\bar{\tau} = 0$ just when $\tau = 0$.

**3.5** Construct a tensor whose vanishing is equivalent to $\bar{\tau} = 0$.

The complex coframing

$$\begin{aligned}\omega^1 &:= dz,\\ \omega^2 &:= dw - |w|^2\, d\bar{z}\end{aligned}$$

has

$$\begin{aligned}\bar{\tau}\omega^1 &= 0,\\ \bar{\tau}\omega^2 &= -w\, d\bar{w} \wedge d\bar{z},\\ &= w\,\omega^{\bar{1}} \wedge \omega^{\bar{2}},\\ &\neq 0.\end{aligned}$$

Therefore the associated almost complex structure is not a complex structure.

A complex valued function $f\colon M \to \mathbb{C}$ on an almost complex manifold is *holomorphic* if $df$ is complex linear.

**Lemma 3.1.** *An almost complex manifold which admits holomorphic functions with arbitrary complex linear differentials at each point has $\bar{\tau} = 0$.*

*Proof.* Denote by $2n$ the dimension of $M$ as a real manifold. Since this problem is local, we can assume that $M$ has a global complex linear coframing $\omega^\mu$. Take the manifold $M' := M \times \mathbb{C}_z \times \mathbb{C}^n_Z$, and the exterior differential system generated by $dz - Z_\mu\omega^\mu$. Any integral manifold on which $\omega^\mu$ are complex-linearly independent is locally the graph of a holomorphic function, and vice versa. The tableau

$$\begin{aligned}d(dz - Z\omega) &= -dZ \wedge \omega - Z\, d\omega,\\ &= -DZ \wedge \omega - Z\,\bar{\tau}\omega,\end{aligned}$$

where

$$DZ_\mu := dZ_\mu + Z_\sigma(t^\sigma_{\nu\mu}\omega^\nu + 2t^\sigma_{\bar{\nu}\mu}\omega^{\bar{\nu}}).$$

The torsion, where $Z \neq 0$, consists of the expression $Z\,\bar{\tau}\omega = Z_\mu\bar{\tau}\omega^\mu$. ☐

**3.6** Find all holomorphic functions for the almost complex structure of the complex coframing

$$\begin{aligned}\omega^1 &:= dz,\\ \omega^2 &:= dw - \bar{w}\, d\bar{z}\end{aligned}$$

Any wedge product $\pi \wedge \omega$ of complex valued 1-forms can be rewritten as a real wedge product in real and imaginary parts

$$\begin{pmatrix}\pi^1 & -\pi^2\\ \pi^2 & \pi^1\end{pmatrix} \wedge \begin{pmatrix}\omega^1 & -\omega^2\\ \omega^2 & \omega^1\end{pmatrix}.$$

If $\pi \wedge \omega$ occurs in a tableau, at some grade, it contributes 2 linearly independent 1-forms:

$$\begin{pmatrix} \pi^1 & -\pi^2 \\ \pi^2 & \pi^1 \end{pmatrix} \wedge \begin{pmatrix} \omega^1 & -\omega^2 \\ \omega^2 & \omega^1 \end{pmatrix}.$$

Count with a complex tableau as if it were real linear, but double the characters.

**Theorem 3.2.** *An analytic almost complex structure has $\bar{\tau} = 0$ just when it arises from a complex structure.*

This theorem remains true with milder assumptions than analyticity [17].

*Proof.* On a complex manifold with holomorphic coordinates $z^\mu$, the 1-forms $\omega^\mu = dz^\mu$ are complex linear for the standard almost complex structure, and have $d\omega^\mu = 0$, so no torsion.

Take an almost complex structure $J$ which has $\bar{\tau} = 0$. Our problem is to construct local holomorphic coordinate functions locally identifying $J$ with the standard complex structure on complex Euclidean space. Again take the exterior differential system generated by $\theta := dz - Z\omega$:

$$d\theta = -\begin{pmatrix} DZ_1 & DZ_2 & \dots & DZ_n \end{pmatrix} \wedge \begin{pmatrix} \omega^1 \\ \omega^2 \\ \vdots \\ \omega^n \end{pmatrix}.$$

$$\begin{matrix} s_1 & s_2 & \dots & s_n \\ 2 & 2 & \dots & 2 \end{matrix}$$

Integral elements are $DZ_\mu = p_{\mu\nu}\omega^\nu$, $p_{\mu\nu} = p_{\nu\mu}$, $n(n+1)/2$ complex constants, so $n(n+1)$ real constants, involution. So there are holomorphic functions with arbitrary differentials at a point, i.e. local holomorphic coordinates. □

The complex coframing

$$\omega^1 = dz,$$
$$\omega^2 = dw - w\, d\bar{z},$$

determines a complex structure.

The expression $\bar{\tau}\omega^\mu = t^\mu_{\bar{\nu}\bar{\sigma}}\omega^{\bar{\nu}} \wedge \omega^{\bar{\sigma}}$ is antisymmetric in $\bar{\nu}, \bar{\sigma}$. It vanishes if $M$ has complex dimension 1, i.e. real dimension 2: every almost complex manifold of real dimension 2 is a Riemann surface.

In this example, we assume familiarity with matrix groups [33]. The manifold $SU_3$ is the collection of all $3 \times 3$ complex unitary matrices $z = (z_{\mu\bar\nu})$ of determinant 1. Write $z_{\bar\mu\nu}$ to mean $\bar z_{\mu\bar\nu}$, so that $z^*$ has entries $z^*_{\mu\bar\nu} = z_{\bar\nu\mu}$. Unitarity is $z_{\mu\bar\sigma} z_{\bar\nu\sigma} = \delta_{\mu\bar\nu}$. Note that $SU_3$ is a real submanifold, *not* a complex submanifold, of the $3 \times 3$ matrices, as this unitarity equation is not complex analytic.

Since $SU_3$ is a submanifold of matrices, each tangent vector $v$ to $SU_3$ is expressed as a matrix. The tangent space $T_I \, SU_3$ at the identity matrix, denoted $\mathfrak{su}_3$, is the set of all traceless skew adjoint $3 \times 3$ complex matrices.

It is traditional to write the identity function on any group as $g$, so $g(z) = z$. The *Maurer–Cartan* 1-form $\omega := g^{-1} dg$ is a 1-form on $SU_3$, valued in $\mathfrak{su}_3$, i.e. to each tangent vector $v \in T_z \, SU_3$, which we identify with a matrix $A$, $v \lrcorner \omega = z^{-1} A$. Write out $\omega$ as a matrix of complex valued 1-forms

$$\omega = \begin{pmatrix} \omega_{1\bar1} & \omega_{1\bar2} & \omega_{1\bar3} \\ \omega_{2\bar1} & \omega_{2\bar2} & \omega_{2\bar3} \\ \omega_{3\bar1} & \omega_{3\bar2} & \omega_{3\bar3} \end{pmatrix}$$

with $\omega_{\mu\bar\nu} = -\omega_{\bar\nu\mu}$.

**3.7** Prove that $\omega$ is invariant under left translation.

**3.8** Calculate that $d\omega = -\omega \wedge \omega$.

In matrix entries, $d\omega_{\mu\bar\nu} = -\omega_{\mu\bar\sigma} \wedge \omega_{\sigma\bar\nu}$. Consider the coframing

$$\omega_{1\bar1} + i\omega_{2\bar2}, \omega_{1\bar2}, \omega_{1\bar3}, \omega_{2\bar3},$$

**3.9** Take exterior derivatives and find torsion vanishing: $SU_3$ has a left invariant complex structure.

**3.10** On the other hand, if we conjugate one of the last three 1-forms in the coframing, prove that we arrive at a left invariant almost complex structure which is not complex.

**3.11*** Prove theorem 3.2 on the preceding page by complexifying variables locally, and applying the Frobenius theorem.

## Almost complex submanifolds

This section can be omitted without loss of continuity.

An *almost complex submanifold* of an almost complex manifold $M$ is a submanifold whose tangent planes are $J$-invariant. Suppose that $M$ has real dimension $2(p+q)$. Let's look for almost complex submanifolds of dimension $2p$. Take a complex linear coframing $\omega^1, \dots, \omega^p, \pi^1, \dots, \pi^q$. Take an almost complex submanifold of real dimension $2p$ on which $\omega^1, \dots, \omega^p$ have linearly independent real and imaginary parts. Then on that submanifold, $\pi = p\omega$ for a complex matrix $p$. Our almost complex submanifold is an integral manifold of the exterior differential system $\pi = p\omega$ on $M \times \mathbb{C}^{pq}$. We let $\theta := \pi - p\omega$. Then on any integral manifold $0 = \bar{\tau}\theta = \bar{\tau}\pi - p\,\bar{\tau}\omega$.

**3.12** Prove that, for any $p > 1$, if there is an almost complex submanifold of real dimension $2p$ tangent to any complex linear subspace of complex dimension $p$ in any tangent space, then $M$ is a complex manifold.

A *holomorphic curve*, often called a *pseudoholomorphic curve*, in an almost complex manifold $M$ is a map $C \to M$ from a Riemann surface, with complex linear differential.

**3.13** Prove the existence of embedded holomorphic disks, i.e. embedded holomorphic curves $C \to M$ where $C$ is the unit disk in the complex plane, in every analytic almost complex manifold, tangent to every complex line in every tangent space.

Chapter 4

# Prolongation

In our next example, we will see what to do when the Cartan–Kähler theorem does not apply to an exterior differential system.

## Notation

In this chapter, for convenience of notation, we drop our convention of writing $\omega^{12}$ to mean $\omega^1 \wedge \omega^2$ etc.

## Example: isometric immersion, involution fails

Take a surface $S$ with a Riemannian metric. Naturally we are curious if there is an isometric immersion $\phi\colon S \to \mathbb{E}^3$, i.e. a map preserving the lengths of all curves on $S$.

This surface (viewed from various angles)

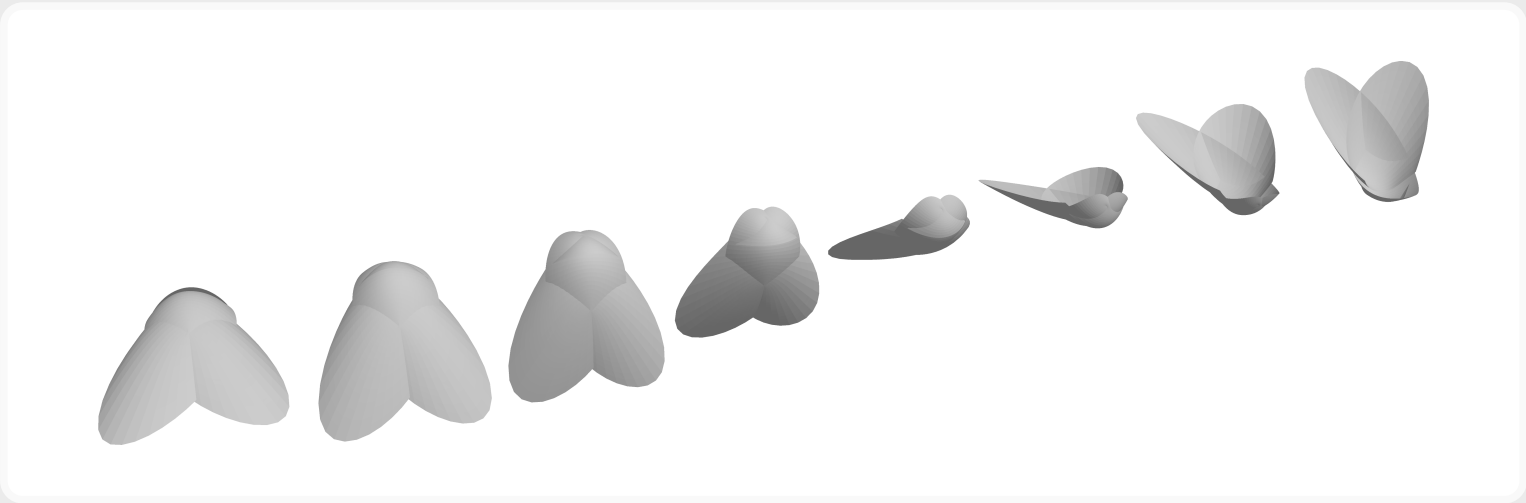

is the image of an isometric immersion of a piece of this paraboloid

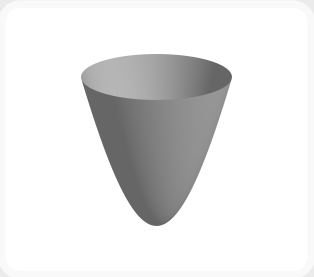

On the frame bundle $\ulcorner S$ of oriented orthonormal frames, denote the soldering forms as $\omega = \omega_1 + i\omega_2$ and the connection by 1-form) $\alpha$ so that $d\omega = i\alpha \wedge \omega$ and $d\alpha = (i/2)K\omega \wedge \bar{\omega}$. On ${}^{\ulcorner}\mathbb{E}^3$ there is a soldering 1-form $\omega'_i$ and a connection 1-form $\gamma'_{ij}$ so that $d\omega'_i = -\gamma'_{ij} \wedge \omega'_j$ and $d\gamma'_{ij} = -\gamma'_{ik} \wedge \gamma'_{kj}$. Let $\alpha' := \gamma'_{12}$.

Suppose that there is an isometric immersion $\phi\colon S \to \mathbb{E}^3$. Its *adapted frame bundle* $X := X_\phi \subset M := {}^{\ulcorner}S \times {}^{\ulcorner}\mathbb{E}^3$ is the set of all tuples

$$(x, e_1, e_2, x', e'_1, e'_2, e'_3)$$

where $x \in S$ with orthonormal frame $e_1, e_2 \in T_xS$ and $x' \in \mathbb{E}^3$ with orthonormal frames $e'_1, e'_2, e'_3 \in T_{x'}\mathbb{E}^3$, so that $\phi_* e_1 = e'_1$ and $\phi_* e_2 = e'_2$. Let $\mathcal{I}$ be the exterior differential system on $M$ generated by the 1-forms

$$\begin{pmatrix} \theta_1 \\ \theta_2 \\ \theta_3 \end{pmatrix} := \begin{pmatrix} \omega'_1 - \omega_1 \\ \omega'_2 - \omega_2 \\ \omega'_3 \end{pmatrix}.$$

Along $X$, all of these 1-forms vanish, while the 1-forms $\omega_1, \omega_2, \alpha$ coframe. Conversely, we will eventually prove that all integral manifolds coframed by $\omega_1, \omega_2, \alpha$ are locally frame bundles of isometric immersions. For the moment, we concentrate on asking whether we can apply the Cartan–Kähler theorem.

$$d\begin{pmatrix} \theta_1 \\ \theta_2 \\ \theta_3 \end{pmatrix} = -\begin{pmatrix} 0 & \alpha' - \alpha & 0 \\ -(\alpha' - \alpha) & 0 & 0 \\ -\gamma'_{13} & -\gamma'_{23} & 0 \end{pmatrix} \wedge \begin{pmatrix} \omega_1 \\ \omega_2 \\ \alpha \end{pmatrix} \mod \theta_1, \theta_2, \theta_3$$

$$\begin{matrix} s_1 & s_2 & s_3 \\ 2 & 1 & 0 \end{matrix}$$

Each 3-dimensional integral element has $\omega' = \omega$, so is determined by the linear equations giving $\gamma'_1, \gamma'_2, \gamma'_3$ in terms of $\omega_1, \omega_2, \alpha$ on which $d\theta = 0$:

$$\begin{pmatrix} \gamma'_{13} \\ \gamma'_{23} \\ \alpha' - \alpha \end{pmatrix} = \begin{pmatrix} a & b & 0 \\ b & c & 0 \\ 0 & 0 & 0 \end{pmatrix} \begin{pmatrix} \omega_1 \\ \omega_2 \\ \alpha \end{pmatrix}.$$

Therefore there is a 3-dimensional space of integral elements at each point. But $s_1 + 2s_2 = 4 > 3$: not involutive, so we can't apply the Cartan–Kähler theorem.

**4.1** Is the existence of a $p$-dimensional involutive integral element precisely the condition that $p + s_0 + s_1 + \dots + s_p \le \dim M$ on every integral element of dimension $p$?

There is another way to look at the failure of involution. In chapter 7, we will see that the Cartan–Kähler theorem can only apply if the generic integral line sits in an integral plane. The equations on integral lines are $\omega'_1 = \omega_1$, $\omega'_2 = \omega_2$, $\omega'_3 = 0$. On any integral plane $\alpha' = \alpha$. The generic integral line does not sit in an integral plane, because it doesn't have to satisfy $\alpha' = \alpha$.

## What to do?

On every integral element, we said that

$$\begin{pmatrix} \gamma'_{13} \\ \gamma'_{23} \\ \alpha' - \alpha \end{pmatrix} = \begin{pmatrix} a & b & 0 \\ b & c & 0 \\ 0 & 0 & 0 \end{pmatrix} \begin{pmatrix} \omega_1 \\ \omega_2 \\ \alpha \end{pmatrix}.$$

Make a new manifold $M' := M \times \mathbb{R}^3_{a,b,c}$, and on $M'$ let $\mathcal{I}'$ be the exterior differential system generated by

$$\begin{pmatrix} \theta_4 \\ \theta_5 \\ \theta_6 \end{pmatrix} := \begin{pmatrix} \gamma'_{13} \\ \gamma'_{23} \\ \alpha' - \alpha \end{pmatrix} - \begin{pmatrix} a & b & 0 \\ b & c & 0 \\ 0 & 0 & 0 \end{pmatrix} \begin{pmatrix} \omega_1 \\ \omega_2 \\ \alpha \end{pmatrix}.$$

## Prolongation

What should we do if there are no involutive integral elements? Take an exterior differential system $\mathcal{I}$ on a manifold $M$. Denote by $\mathrm{Gr}_p M$ the Grassmann bundle of $p$-dimensional linear subspaces of tangent spaces of $M$. Take the subset $M' \subset \mathrm{Gr}_p M$ of all $p$-dimensional integral elements. Cut out any point of $M'$ near which $M'$ is not a submanifold of $\mathrm{Gr}_p M$. Suppose we write out a tableau in 1-forms $\theta, \omega, \pi$. Write each linear subspace coframed by $\omega$ as the solutions of the linear equations $\pi = p\omega, \theta = q\omega$ for some constants $p, q$. So on an open subset of $\mathrm{Gr}_p M$, $p, q$ are functions valued in some vector spaces. Pull back the 1-forms $\theta, \omega, \pi$ to that open subset via the map $(m, E) \in \mathrm{Gr}_p M \mapsto m \in M$. The subset $M'$ is cut out by the equation $q = 0$. Let $\theta' := \pi - p\omega$. The exterior differential system $\mathcal{I}'$ on $M'$ generated by $\vartheta'$ is the *prolongation* of $\mathcal{I}$.

More abstractly, without choosing a local basis of 1-forms: a 1-form $\theta$ on $\mathrm{Gr}_p M$ is *contact* if, at each point $E \in \mathrm{Gr}_p M$, $\theta_E$ is the pullback of a cotangent vector vanishing on $E$. The ideal $\mathcal{I}'$ is generated by the restrictions to $M'$ of the contact 1-forms. Any integral manifold $X \subset M$ determines an integral manifold $x \in X \mapsto T_x X \in M'$. Recall that the *fibers* of a map are the preimages of points. All integral manifolds on $M'$ nowhere tangent to the fibers of $M' \to M$ arise in this way.

Roughly speaking, prolongation is differentiation:

**4.2** Write an exterior differential system associated to the partial differential equation $u_x = u_{yy}$. Consider the system of differential equations obtained by differentiating both sides of the equation once in each of the variable:

$$u_x = u_{yy}$$

$$(\dots)_x \swarrow \qquad \searrow (\dots)_y$$

$$u_{xx} = u_{yyx} \qquad\qquad u_{xy} = u_{yyy}$$

Explain why the prolongation is associated to that system.

Some theorems prove that, under some hypotheses, exterior differential systems become involutive after sufficiently many prolongations [3] p. 255 theorem 3.2, [26] p. 68 theorem 4.11.

## Example: isometric immersion, prolonging

Returning to our example of isometric immersion of surfaces, prolong:

$$\begin{pmatrix}\theta_4\\ \theta_5\\ \theta_6\end{pmatrix} := \begin{pmatrix}\gamma'_{13}\\ \gamma'_{23}\\ \alpha'-\alpha\end{pmatrix} - \begin{pmatrix}a & b & 0\\ b & c & 0\\ 0 & 0 & 0\end{pmatrix}\begin{pmatrix}\omega_1\\ \omega_2\\ \alpha\end{pmatrix}.$$

Note that $0 = d\theta_1, d\theta_2, d\theta_3$ modulo $\theta_4, \theta_5, \theta_6$, so we can forget about them.

Calculate the exterior derivatives:

$$d\begin{pmatrix}\theta_4\\ \theta_5\\ \theta_6\end{pmatrix} = -\begin{pmatrix}Da & Db & 0\\ Db & Dc & 0\\ 0 & 0 & 0\end{pmatrix}\wedge\begin{pmatrix}\omega_1\\ \omega_2\\ \alpha\end{pmatrix} + \begin{pmatrix}0\\ 0\\ t\omega_1\wedge\omega_2\end{pmatrix} \mod \theta_1,\dots,\theta_6$$

where

$$\begin{pmatrix}Da\\ Db\\ Dc\end{pmatrix} := \begin{pmatrix}da + 2b\alpha + a_1\omega_1 + a_2\omega_2,\\ db + (a-c)\alpha + b_1\omega_1 + b_2\omega_2,\\ dc + 2b\alpha + c_1\omega_1 + c_2\omega_2,\end{pmatrix}$$

where $a_1, a_2, b_1, b_2, c_1, c_2$ can be chosen as we like, as long as

$$a_2 = b_1,$$
$$b_2 = c_1.$$

(It is convenient to leave the freedom to choose these later.) The torsion is

$$t := ac - b^2 - K.$$

This torsion clearly has to vanish on any 3-dimensional $\mathcal{I}'$-integral element, i.e. every 3-dimensional $\mathcal{I}'$-integral element lives over the subset of $M'$ on which $K = ac - b^2$. To ensure that this subset is a submanifold, we let $M'_0 \subset M'$ be the set of points where this equation is satisfied and at least one of $a, b, c$ is not zero. Clearly $M'_0 \subset M'$ is a submanifold, on which we find $Da, Db, Dc$ linearly dependent. On $M'_0$:

$$d\begin{pmatrix}\theta_4\\ \theta_5\\ \theta_6\end{pmatrix} = -\begin{pmatrix}Da & Db & 0\\ Db & Dc & 0\\ 0 & 0 & 0\end{pmatrix}\wedge\begin{pmatrix}\omega_1\\ \omega_2\\ \alpha\end{pmatrix} \mod \theta_1,\dots,\theta_6.$$

| $s_1$ | $s_2$ | $s_3$ |
|---|---|---|
| 2 | 0 | 0 |

There are 2 dimensions of integral elements at each point, involution: there is an integral manifold through each point of $M'_0$, and in particular above every point of the surface. The prolongation exposes the hidden necessary condition for existence of a solution: the relation $K = ac - b^2$ between the curvature of the surface and the shape operator.

**4.3** Take any 3-dimensional integral manifold $X$ of the exterior differential system constructed above. Suppose that on $X$, $0 \neq \omega_1 \wedge \omega_2 \wedge \alpha$. Prove that every point of $X$ lies in some open subset $X_0 \subset X$ so that $X_0$ is an open subset of the adapted frame bundle of an isometric immersion $U \to \mathbb{E}^3$ of an open subset $U \subset S$.

We can easily generalize these computations:

**Theorem 4.1.** *Take any surface $S$ with Riemannian metric and a point $x_0 \in S$. Denote the Gauss curvature by $K$. Take any 3-dimensional manifold $X'$ with Riemannian metric, a point $x'_0$, and a linear isometric injection $F \colon T_{x_0}S \to T_{x'_0}X'$. Let $R'$ be the sectional curvature tensor of $X'$ on the image of $F$. Pick a nonzero quadratic form $q$ on the tangent plane $T_{x_0}S$ so that $\det q = K - R'$. Then there is an isometric immersion $f$ of some neighborhood of $x_0$ to $X'$, so that $f'(x_0) = F$ and so that $f$ induces shape operator $q$ at $x_0$.*

Write $dK = K_1\omega_1 + K_2\omega_2$. It is convenient to pick our $a_1, a_2, b_1, b_2, c_1, c_2$ to satisfy not just $a_2 = b_1, b_2 = c_1$, as required above, but also to satisfy

$$
\begin{aligned}
K_1 &= 2bb_1 - ac_1 - ca_1, \\
K_2 &= 2bb_2 - ac_2 - ca_2,
\end{aligned}
$$

so that $0 = c\,Da + a\,Dc - 2b\,Db$, a simple relation among the differential forms in the exterior differential system. We can always do this near any point of $M'_0$, by choice of $a_1, a_2, b_1, b_2, c_1, c_2$, using the fact that one of $a, b, c$ is not zero.

**4.4** Prove that every sufficiently small spherical cap on the unit sphere admits an isometric embedding to $\mathbb{E}^3$ *not* contained in a surface of revolution.

For more on isometric immersions, see [16].

Chapter 5

# Cartan's test

We explain how involutivity can be seen as having the generic integral line lie in an integral plane and so on.

## Generic points and generic submanifolds

When we say that the *generic* point of some manifold has some property, we mean that those points which fail to have that property lie in a closed, nowhere dense set. The *generic* submanifold has some property if the property holds for all submanifolds whose tangent spaces avoid some closed, nowhere dense subset of the Grassmann bundle. The definition of *generic* in analysis is usually more sophisticated; we could use a more sophisticated definition without changing any of our proofs or results, but we won't need to.

## Cartan's strategy

Draw a point, an integral curve through that point, an integral surface through that curve, and so on up to some required dimension. We fail unless the integral curve is *extendable*, i.e. lies inside an integral surface, and so on. Integral curves are not always extendable. Don't try to draw all integral curves, surfaces, and so on, but only the generic ones: ask whether the generic integral curve is extendable, and the generic integral surface, and so on.

## Cartan's test

An integral line is *extendable* if it lies in an integral plane, and so on. We can imagine that extendability of generic integral manifolds is related to extendability of generic integral elements.

**5.1** Prove that any integral element is extendable just when its codimension, in the space of integral elements, exceeds the rank of its polar equations.

A 0-dimensional integral element is *regular* if its polar equations have locally maximal rank. Locally maximal rank is the same as locally constant rank, as linearly independent polar equations remain linearly independent nearby.

An integral element is *ordinary* if it contains a regular hyperplane, and *regular* if in addition its polar equations have locally maximal rank. Regularity and ordinarity are difficult to test directly, so we will prove:

**Theorem 5.1** (Cartan's test)**.** *An integral element is involutive if and only if it is ordinary.*

**Corollary 5.2.** *In each component of the space of integral elements, either none are involutive or the generic one is involutive.*

Ordinarity occurs precisely when the generic integral line sits in an integral plane, and so on:

**5.2** Prove that an integral plane $E_2$ is ordinary just when, for any integral line $E_1 \subset E_2$, every integral line close enough to $E_1$ sits in an integral plane.

## Differentiating on the Grassmann bundle

Each differential form $\vartheta$ yields an equation $0 = \vartheta$, satisfied by the integral elements of the exterior differential system it generates. At any given point $m_0 \in M$, integral elements are points of the Grassmannian $\operatorname{Gr}_p T_{m_0} M$.

**5.3** Identify the polar equations at an integral element with the differentials of those equations on the Grassmannian.

**Corollary 5.3.** *Take an involutive integral element $E$ of any exterior differential system $\mathcal{I}$. Take a subsystem $\mathcal{I}' \subseteq \mathcal{I}$ and a flag in $E$. If the $\mathcal{I}'$-characters of that flag are the $\mathcal{I}$-characters of $E$, then $\mathcal{I}$ and $\mathcal{I}'$ have the same integral elements near $E$, and so the same integral manifolds.*

*Proof.* The $\mathcal{I}'$-integral elements lie in a submanifold of the Grassmann bundle, with codimension equal to the number of linearly independent polar equations, by problem 5.3. This submanifold has the same dimension as the set of $\mathcal{I}$-integral elements, which it contains. □

## Extending integral elements

**Lemma 5.4.** *Every ordinary integral element is involutive.*

*Proof.* Impose $p-1$ generic linear constraints; they cut down our given ordinary integral element $E$ to an integral line. They also cut down nearby integral elements to nearby integral lines. On the other hand, a line satisfying those generic constraints is integral just when it satisfies $s_0$ linearly independent polar equations. Pick another $p-1$ generic linear constraints. A line satisfying these will sum to our integral line to span an integral plane, just when it satisfies the $s_0 + s_1$ linearly independent polar equations of our integral line. By induction, the general $p$-dimensional integral element near $E$ is cut out by solving $ps_0 + (p-1)s_1 + \cdots + s_{p-1}$ equations. Solutions form an analytic manifold by problem 1.12 on page 6, of dimension $\dim M + s_1 + 2s_2 + \cdots + ps_p$. □

**5.4** Prove theorem 5.1.

**5.5** Prove that the generic linear subspace $E_k \subset E_p$ of dimension $k$ of an involutive integral element is involutive with the same characters $s_1, s_2, \dots, s_k$ as $E_p$, up to dimension $k$.

## View from a tableau

**5.6** Prove that, for any involutive integral element $E \subset T_{m_0}M$ of an exterior differential system $\mathcal{I}$, there is a tableau defined near $m_0$, with characters equal to those of $\mathcal{I}$, generating an exterior differential system with the same characters and integral elements as $\mathcal{I}$ near $E$.

**5.7** Prove in addition that locally (not just at a point) some change of coframing absorbs torsion.

Put all polars $\pi^\alpha$ of grade $i$ into the entries of a column vector $\pi^i$, with $s_i$ rows. Write the equations for integral elements by plugging $\pi^i = p^i_j \omega^j$ into the tableau. The *grade* of $p^i_j$ is $i - j$.

For the moment, suppose there is no nonlinearity. Plug into the tableau to get linear equations solving for all coefficients $p^i_j$ for $i < j$ in terms of various coefficients $p^\beta_i$, with a smaller subscript. Inductively, we solve for all negative grade coefficients in terms of semipositive grade coefficients (where *semipositive* means *not negative*). (There may be other equations as well, from entries that are not polars; ignore these.) We solve for $p^1_i$, for $i = 2, 3, \dots, p$, hence $(p-1)s_1$ equations. Similarly for the other grades, a total of

$$(p-1)s_1 + (p-2)s_2 + \dots + s_{p-1}$$

equations at least, on $p(s_1 + s_2 + \dots + s_p)$ variables. Each integral element at this point is determined by the values of its semipositive coefficients. The number of semipositive coefficients is the difference between number of variables and number of equations:

$$s_1 + 2s_2 + \dots + ps_p.$$

We assumed no nonlinearity. Each nonlinearity term puts a quadratic or higher order expression into each of these equations. By the implicit function theorem, near an integral element, say arranged to be at the origin, higher order terms do not alter the possibility to solve the equations locally analytically.

This tableau is only computed at one point. Extend the 1-forms in which it is expressed to be defined nearby as well. Careful: there might be differential forms in $\mathcal{I}$ which vanish at the point where we computed the tableau; we are ignoring them. So when we extend the forms of the tableau to nearby points, we obtain a tableau for some subsystem of $\mathcal{I}$, but maybe not all of $\mathcal{I}$.

The space of integral elements thus locally lies in a submanifold of the Grassmann bundle of dimension

$$\dim M + s_1 + 2s_2 + \dots + ps_p,$$

parameterized by choices of point of $M$ and semipositive grade coefficients. Near a given integral element, the space of integral elements is a submanifold of that dimension just when the tableau is involutive at that integral element. Involutivity holds just all other equations arising on integral elements are consequences of those above which solved for negative grade coefficient. The semipositive grade coefficients can then vary in some open set.

## Differential equations

All polars are linearly independent at our starting point, so we can assume that at some given point they are differentials of coordinate functions $\pi^\alpha = du^\alpha$, and that $\omega^i = dx^i$. We can extend the $\omega^i$ to anything we like at nearby points, so can assume $\omega^i = dx^i$ nearby. We can only arrange that the polars at nearby points are various multiples of $dx^i, du^\alpha$. Careful: at nearby points, some other tableau entries might become linearly independent, i.e. more polars, and so more equations. Careful: there need not be any actual integral elements near this one; the equations we generate are only necessary conditions for an integral element. Nonetheless, by the implicit function theorem again, each equation of negative grade coefficients of each polar $\pi^\alpha$ can be written as an equation of negative grade coefficients of some $du^\alpha$, and vice versa, hence as functions of the semipositive ones and the coordinates $x, u$.

Let $u^i$ be the column vector of functions $u^\alpha$ associated to polars $\pi^\alpha$ of grade $i$. The equations become differential equations

$$\frac{\partial u^i}{\partial x^{>i}} = \text{some function}\left(x, u, \frac{\partial u^j}{\partial x^{\le j}}\right),$$

for negative grade derivatives in terms of semipositive grade. We see one differential equation for each polar equation on each integral element in our flag.

So far, we have not assumed involutivity: every integral manifold of every exterior differential system satisfies such equations. Involutivity is the condition that the tableau yields no more differential equations on integral manifolds.

Choose any initial values

$$u^i(x^1, \dots, x^i, 0, \dots, 0),$$

so $s_i$ functions of $i$ variables, $i = 0, 1, \dots, p$. The differential equations become determined, so have local solutions near the origin by the Cauchy–Kovalevksaya theorem (theorem A.2 on page 80).

The first differential equation solves for $u^0$ along the $x^1$ axis, from initial values at the origin.

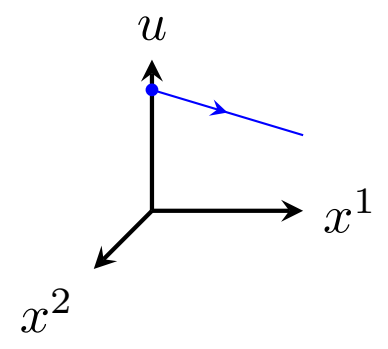

The second solves for $u^0$ and $u^1$, in the $x^1, x^2$ plane, from initial values along the $x^1$ axis.

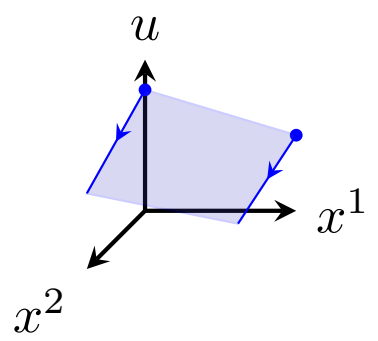

But is the first differential equation still satisfied by $u^0$, at any constant value of $x^2$?

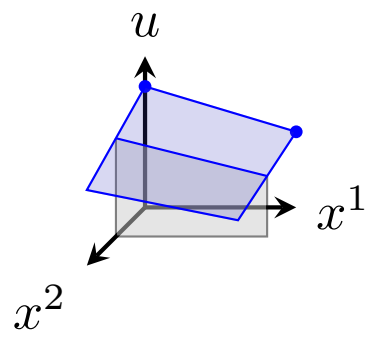

In other words, we need to see that these equations are compatible with one another.

**5.8** Write out these differential equations for the exterior differential system of the Cauchy–Riemann equations in two complex variables. Can you solve these differential equations, by specifying initial data as above?

**5.9** Prove that any involutive integral element of dimension $p$, whose characters all vanish except perhaps $s_{p-1}$ and $s_p$, is tangent to an integral manifold.

**5.10** Use the Cartan–Kähler theorem to prove that every involutive integral element of any exterior differential system lies tangent to a leaf of a foliation by integral manifolds, foliating some open set. How many functions of how many variables do such foliations depend on? Prove that any involutive integral manifold is covered in open sets, each of which is a leaf in a foliation by integral manifolds, foliating some open set.

## Compatibility and involutivity

Geometrically, we sweep a point into an integral curve, sweep that into a surface, and so on. Compatibility asks that the surface is an integral surface, and so on. Suppose that the differential equations are compatible, i.e. we can pick any initial values in the domain of our coordinates, and obtain an integral manifold. Then we can also pick initial values for the same differential equations nearby in those coordinates. The space of integral elements, parameterized by the semipositive derivatives

$$\frac{\partial u^j}{\partial x^{\le j}},$$

has the predicted dimension at all nearby points. So compatibility implies involutivity.

An incompatibility, arising as an additional differential equation, might vanish, perhaps to high order, at the particular point where we are working, but obstruct at nearby points. So we check the dimension of the space of integral elements nearby.

Incompatibilities can arise from commuting partial derivatives in our differential equations. Any exterior differential system is closed under exterior derivative, expressing the commutativity of first partial derivatives. So we expect that all incompatibilities are already present in the differential equations. So we expect that incompatibilities force the dimension of nearby integral elements below the predicted dimension. So we expect that compatibility is involutivity; chapter 7 gives a proof.

## Generality of integral manifolds

This section can be omitted without loss of continuity.

As above, Cartan's strategy constructs integral manifolds of an involutive system by solving a sequence of determined equations, with initial data $s_i$ functions of $i$ variables, $i = 0, 1, 2, \dots, p$. Count Taylor coefficients of initial data: the Taylor series of order $k$ of integral manifolds, at a chosen point in our coordinates, form a manifold of dimension

$$s_0 + \binom{k}{0} s_1 + \binom{k+1}{1} s_2 + \cdots + \binom{k+p-1}{p-1} s_p.$$

As we vary $k$, these dimensions of determine all of the characters.

Any other choice of determined systems of differential equations, giving rise to the same integral manifolds (or at least to the same Taylor series of integral manifolds at each order), injectively on each order of Taylor series, also has general solution depending on $s_0$ constants, $s_1$ functions of one variable, $s_2$ functions of two variables, and so on.

## Deformation

This section can be omitted without loss of continuity.

A *local deformation* of an integral manifold $X$ is an analytic map $\phi$ defined on an open subset of $X \times \mathbb{R}$ containing $X \times \{ 0 \}$ so that, for each constant $t$, $x \mapsto \phi(x, t)$ is an integral manifold, where defined, with $\phi(x, 0) = x$. The *velocity* of the deformation is the projection of $\left.\frac{\partial \phi}{\partial t}\right|_{t=0}$ to the normal bundle of $X$. By deforming initial data, any integral manifold of any determined system is covered in open sets admitting local deformation with velocity any solution of the linearization. We have not yet justified Cartan's strategy to prove the Cartan–Kähler theorem, but we can already see that, since the strategy is a sequence of determined problems, we can apply the same reasoning; any involutive analytic integral manifold of any exterior differential system is covered in open sets admitting local deformation with velocity any solution of the linearization. Local analytic velocity fields exist in the same generality as integral manifolds, solving linear determined problems.

The Holmgren uniqueness theorem [20] p. 80 proves that, on any analytic involutive integral manifold, any continuously differentiable velocity field is uniquely determined by its continuously differentiable initial data at each step in the sequence of linear determined problems. However, existence of continuously differentiable velocity fields is not guaranteed.

## Prolongation

This section can be omitted without loss of continuity.

When we prolong, we form new equations $\pi = p_j\omega^j$ for each polar. Only the semipositive coefficients $p^i_{\le i}$ can vary freely. (For noninvolutive systems, not all of them vary freely.)

$$d\begin{pmatrix}\pi^1 - p^1_i\omega^1\\ \pi^2 - p^2_i\omega^1\\ \vdots\\ \pi^p - p^p_i\omega^i\end{pmatrix} = -\begin{pmatrix} dp^1_1 & * & * & \dots & *\\ dp^2_1 & dp^2_2 & * & \dots & *\\ \vdots & \vdots & \vdots & \vdots & \vdots\\ dp^p_1 & dp^p_2 & \dots & \dots & dp^p_p\end{pmatrix}\wedge\begin{pmatrix}\omega^1\\ \omega^2\\ \vdots\\ \omega^p\end{pmatrix}$$

modulo torsion, since we haven't included $d\pi^i$ terms and $d\omega^i$ terms. The $*$ terms are $dp^i_j$ of negative grade, each solved for in terms of semipositive grade, so generate no polar. There is no nonlinearity. Some of these $dp$ might be linearly dependent, due to noninvolutivity, so we don't know where to highlight polars, but the polars lie in some spots in among these semipositive $dp$, on or below the diagonal.

Consider the characters $s'_j$ of the prolongation. There are $s_1 + s_2 + \dots + s_p$ rows to this tableau, one for each polar, each representing a 1-form in $\mathcal{I}'$. Grade zero in $\mathcal{I}'$ also includes $\mathcal{I}^1$, so

$$s'_0 = s_0 + s_1 + \dots + s_p.$$

In the first column, there is at most one $dp$ polar in each row:

$$s'_1 \le s_1 + s_2 + \dots + s_p.$$

In the second column,

$$s'_2 \le s_2 + s_3 + \dots + s_p,$$

and so on, with finally $s'_p \le s_p$: the last nonzero character cannot increase. These inequalities are equalities just for involutive exterior differential systems.

# Chapter 6

# The characteristic variety

We give a geometric description of the characteristics of the associated partial differential equations of exterior differential systems.

## Linearization

The reader unfamiliar with linearization of partial differential equations, or characteristics, might look at appendix B. Take an exterior differential system $\mathcal{I}$ on a manifold $M$, and an integral manifold $X$. Suppose that the flow of a vector field $v$ on $M$ moves $X$ through a family of integral manifolds. The tangent spaces of $X$ are carried by the flow of $v$ through integral elements of $\mathcal{I}$. Equivalently, the flow pulls back each form in $\mathcal{I}$ to vanish on $X$. So $0 = \mathcal{L}_v \vartheta|_X$ for any $\vartheta \in \mathcal{I}$.

**6.1** Prove that all vector fields $v$ tangent to $X$ satisfy this equation.

More generally, suppose that $E \subset T_m M$ is an integral element of $\mathcal{I}$. If a vector field $v$ on $M$ carries $E$ through a family of integral elements, then $0 = \mathcal{L}_v \vartheta|_E$ for each $\vartheta \in \mathcal{I}$.

**6.2** Compute $\mathcal{L}_v \vartheta|_E$ in coordinates.

Take any submanifold $X$. Suppose that a differential form $\vartheta$ vanishes on the linear subspace $E := T_m X$. Writing $X$ as the graph of some functions, the expression $\vartheta|_X$, as a nonlinear first order differential operator on those functions, has linearization $\vartheta \mapsto \mathcal{L}_v \vartheta|_E$. That linearization is applied to sections $v$ of the normal bundle $TM|_X / TX$.

**6.3** In coordinates, prove that the linearized operator at the origin of our coordinates depends only on the integral element $E = T_m X$, not on the choice of submanifold $X$.

Linearize any exterior differential system about any integral element by linearizing its differential forms, i.e. by linearizing the differential equations given by asking that those forms vanish on submanifolds. The linearization at a $p$-dimensional integral element depends only on $\mathcal{I}^p_m$: the set of all values $\vartheta_m$ of forms in $\mathcal{I}^p$.

**6.4** If $P := \langle E^{\perp} \rangle \subset \Lambda^* T^*_m M$, identify the linearized exterior differential system with $\left(\mathcal{I}_m + P^2\right)/P^2 \subset P/P^2 = E^{\perp} \otimes \Lambda^* E^*$.

In a tableau, $P$ is generated by the polars, so the linearization is precisely the same tableau, modulo the nonlinearity. So the linearization about an involutive integral element is involutive. On the other hand, the linearization at a noninvolutive integral element may or may not be involutive.

**6.5** Compute the linearization of $u_{xx} = u_{yy} + u_{zz} + u_x^2$ around $u = 0$ by setting up this equation as an exterior differential system.

**6.6** For a section $v$ of the normal bundle of $X$ which vanishes at a point $m$, compute the linearization $\mathcal{L}_v \vartheta|_E$.

## The characteristic variety

An integral element is *noncharacteristic* if it is a hyperplane in precisely one integral element. An integral manifold is *noncharacteristic* (or *characteristic*) if all of its tangent spaces are. A characteristic hypersurface in an integral manifold could perhaps lie in more than one integral manifold, a potential nontangential intersection of integral manifolds. Glue along such an intersection, to create a "crease" along an integral manifold. So we imagine that integral manifolds are more flexible along characteristic hypersurfaces, although there is no theorem to prove that.

**6.7** Find the characteristics of the wave equation as an exterior differential system. Show creasing of solutions along characteristics, and not along noncharacteristic curves.

**Lemma 6.1.** *Take an integral manifold $X$ of an exterior differential system $\mathcal{I}$, a point $x \in X$, and let $E := T_x X$. The characteristic variety $\Xi_x \subset \mathbb{P}E^*$ of the linearization is the set of characteristic hyperplanes in $E$.*

*Proof.* We compute the symbol by replacing derivatives $\frac{\partial v^a}{\partial x^i}$ by expressions $v^a \xi_i$:

$$\sigma(\xi)\, v = \xi \wedge (v \lrcorner \vartheta)|_E \, .$$

The linearization of $\mathcal{I}$ is just the sum of the linearizations of any forms $\vartheta \in \mathcal{I}$ spanning $\mathcal{I}_m$:

$$\sigma \colon \xi \in E^*, v \in T_m M/E \mapsto \mathcal{I}_m^{p*} \otimes \Lambda^p E^* .$$

So $\Xi_x$ is the set of $\xi \in E^*$ for which there is some section $v$ of the normal bundle of $X$ with $v(x) \neq 0$ and

$$\xi \wedge (v \lrcorner \vartheta)|_E = 0$$

for every $\vartheta \in \mathcal{I}_m^p$. This says precisely that the vector $v$ can be added to the hyperplane $E \cap (0 = \xi) \subset E$ to make an integral element enlarging the hyperplane. ∎

This lemma allows us to define the characteristic variety of any integral element, even if not tangent to any integral manifold: the *characteristic variety* of an integral element $E$ is the set of characteristic hyperplanes in $E$, denoted $\Xi_E \subset \mathbb{P}E^*$.

**6.8** For an involutive integral element, prove that the following are equivalent:

*a.* the integral element contains a noncharacteristic hyperplane,

*b.* every regular hyperplane is noncharacteristic,

*c.* the final character is zero.

**6.9** More equations, fewer characteristics: if $\mathcal{J} \subset \mathcal{I}$, prove that $\Xi_{\mathcal{I}} \subset \Xi_{\mathcal{J}}$.

**6.10** What is the characteristic variety of a Frobenius exterior differential system?

## Determined systems

A $p$-dimensional integral element $E \subset T_{m_0}M$ of an exterior differential system $\mathcal{I}$ is *determined* if $E$ contains a noncharacteristic hyperplane and $\mathcal{I}^p_m$ has constant dimension $s_0 + s_1 + \cdots + s_{p-1}$ for $m$ near $m_0$. Note that then $s_p = 0$.

**Theorem 6.2.** *Take an analytic exterior differential system. Every noncharacteristic analytic integral manifold $X$, whose every tangent space is a hyperplane in a determined integral element, is a hypersurface in an analytic integral manifold. Any two analytic integral manifolds containing $X$ as a hypersurface share a set, open in both, containing $X$.*

*Proof.* As above, the symbol matrix $\sigma(\xi)$ at each $p$-dimensional integral element $E$ is a linear map in

$$T_mM/E \to \mathcal{I}^{p*}_m \otimes \Lambda^p E^*,$$

square at each determined integral element $E$. The dimension of $\mathcal{I}^p_m$ cannot drop below $s_0 + s_1 + \cdots + s_p$, because it generates that many polar equations at $E$, and hence near $E$. Pick out that number of linearly independent $p$-forms from $\mathcal{I}$: they span $\mathcal{I}^p$ at every point nearby. The determined exterior differential system they generate has the same $p$-dimensional integral manifolds, polar equations, and symbol. □

**6.11** Prove that every real immersed curve in an almost complex manifold lies in an immersed holomorphic disk.

**6.12** Prove that if $E \subset E_+$ is characteristic, and $E_+$ has dimension $p$ then $s_p > 0$.

To compute the characteristic variety from a tableau:

*a.* Drop the nonlinearity.

*b.* Turn each polar $\pi^\alpha$ into a formal expression $v^\alpha \xi$ where $\xi = \xi_i \omega^i$.

*c.* Expand out the tableau, and collect up terms in each $\omega^I$.

*d.* Write out these terms as linear expressions in the variables $v^\alpha$, and reorganize those linear expressions in a matrix: the symbol matrix $\sigma(\xi)\, v = 0$.

*e.* The determinants of the minors of the symbol matrix are the equations of $\Xi$.

This works because the equations of the characteristic variety include $0 = \theta^a(v)\xi$, which forces $0 = \theta^a(v)$ for any nonzero $\xi$, so that the nonzero components of $v$ are $v^\alpha = \pi^\alpha(v)$. Compute $\xi \wedge (v \lrcorner \vartheta)$ by plugging in $v$ to each polar $\pi^\alpha$, yielding $v^\alpha$, and wedge in a factor of $\xi$.

For the tableau of triply orthogonal webs

$$
\begin{pmatrix} \pi^3 & 0 & 0 \\ 0 & \pi^2 & 0 \\ 0 & 0 & \pi^1 \end{pmatrix} \wedge \begin{pmatrix} \overline{\omega^{12}} \\ \omega^{13} \\ \omega^{23} \end{pmatrix}
$$

we find

$$
\begin{aligned}
\begin{pmatrix} v^3 & 0 & 0 \\ 0 & v^2 & 0 \\ 0 & 0 & v^1 \end{pmatrix} \xi \wedge \begin{pmatrix} \omega^{12} \\ \omega^{13} \\ \omega^{23} \end{pmatrix} &= \begin{pmatrix} v^3 & 0 & 0 \\ 0 & v^2 & 0 \\ 0 & 0 & v^1 \end{pmatrix} \begin{pmatrix} \xi_3 \\ -\xi_2 \\ \xi_1 \end{pmatrix} \omega^{123}, \\
&= \begin{pmatrix} 0 & 0 & \xi_3 \\ 0 & -\xi_2 & 0 \\ \xi_1 & 0 & 0 \end{pmatrix} \begin{pmatrix} v^1 \\ v^2 \\ v^3 \end{pmatrix} \omega^{123}, \\
&= \sigma(\xi)\, v \omega^{123}.
\end{aligned}
$$

So $\det \sigma(\xi) = \xi_1 \xi_2 \xi_3$: the characteristic variety $\Xi$ is the triple of lines $(\xi_1 = 0)$, $(\xi_2 = 0)$ and $(\xi_3 = 0)$. Any 3-dimensional integral element coframed by $\omega^1, \omega^2, \omega^3$ has characteristic hyperplanes $0 = \xi_1 \omega^1 + \xi_2 \omega^2 + \xi_3 \omega^3$ for any one of the three coefficients vanishing, i.e. containing any one of the three axes.

**Theorem 6.3.** *Take an embedded analytic surface $S \subset \mathbb{E}^3$. Pick orthonormal analytic vector fields $e_1, e_2, e_3$ along $S$, none tangent to $S$. Then there is a triply orthogonal web near $S$ in $\mathbb{E}^3$ with leaves perpendicular to $e_1, e_2, e_3$ at each point of $S$. Any two such agree near $S$.*

The reader familiar with Riemannian geometry may recognize that this proof works identically replacing $\mathbb{E}^3$ by any analytic Riemannian 3-manifold. Any symmetry of $S$ and $e_1, e_2, e_3$ on $S$ is shared by the triply orthogonal web, by uniqueness.

In proving Lie's third theorem, we had tableau

$$\begin{pmatrix} \pi^1_1 & \pi^1_2 & \dots & \pi^1_p \\ \pi^2_1 & \pi^2_2 & \dots & \pi^2_p \\ \vdots & \vdots & \dots & \vdots \\ \pi^p_1 & \pi^p_2 & \dots & \pi^p_p \end{pmatrix} \wedge \begin{pmatrix} \omega^1 \\ \omega^2 \\ \vdots \\ \omega^p \end{pmatrix},$$

i.e. $0 = \pi^i_j \wedge \omega^j$. We plug in $\pi^i_j = v^\alpha \xi$, to get

$$0 = \frac{1}{2}(v^i_j \xi_k - v^i_k \xi_j)\omega^{kj}.$$

Our linear expressions are $v^i_j \xi_k - v^i_k \xi_j$. One particular solution $v$ to these equations is $v^i_j := \xi_j$ for all $i$. So every hyperplane in every $p$-dimensional integral element is characteristic. This is not surprising: any diffeomorphism takes every Maurer–Cartan form to a Maurer–Cartan form; we can pick a diffeomorphism which is very close to the identity except along a hypersurface "bending" every Maurer–Cartan form "along" that hypersurface.

**6.13** Find the characteristic variety of $\nabla \times u = u - f$.

**6.14** Find the characteristic variety for the exterior differential system of harmonic functions in the plane. Show that it doesn't satisfy the conditions of the Cauchy–Kovalevskaya theorem as described in theorem 6.2 on page 49. (Harder: can you "fix it"?)

Chapter 7

# Proof of the Cartan–Kähler theorem

## Cartan's strategy II

This section can be omitted without loss of continuity.

Pick a hypersurface, a hypersurface in the hypersurface, and so on, a flag of submanifolds, stopping at codimension $p$. Imagine that the generic $p$-dimensional integral manifold locally intersects each submanifold of our flag transversally: intersecting the smallest in a point, the next smallest in an integral curve, the next smallest in an integral surface, and so on.

Construct an integral manifold by Cartan's strategy: draw a point inside the smallest submanifold of the flag, an integral curve through that point inside the next smallest, an integral surface through that curve inside the next smallest, and so on. Trouble: there might be many integral curves, passing through that point, lying on that flag submanifold. Inside the flag submanifold, pick a smaller submanifold, a *restraining manifold*, cutting down dimensions so that our point lies on a unique integral curve in the restraining manifold. We choose nested restraining manifolds, one in each flag submanifold, starting from the largest and going down in dimension. Start from the smallest and go up in dimension: pick a point in the smallest restraining manifold, sweep out an integral curve through it in the next restraining manifold, and so on. We will see that each restraining manifold is locally the choice of $s_k$ functions of $k$ variables.

Fix a flag. Construct integral manifolds by varying the restraining manifolds inside the submanifolds of that fixed flag. We will see that different selections of restraining manifolds give rise to a different integral manifold in the final stage. In this sense, integral manifolds depend on $s_0$ constants, $s_1$ functions of one variable, and so on.

We wrote out differential equations on page 42. The flag:

$$M_i = (0 = x^{i+1} = \cdots = x^p)$$

the variables over which to solve differential equations. The restraining manifolds are the initial data:

$$R_p = (u^p = u^p(x^1, \dots, x^p))$$

and, for $i = p-1, p-2, \dots, 2, 1, 0$,

$$R_i = R_{i+1} \cap M_i \cap (u^i = u^i(x^1, \dots, x^i, 0, \dots, 0)).$$

**7.1** For the equation $\nabla \times u = f - u$ from problem 2.5 on page 21, write flag and restraining manifolds in coordinates, and the associated differential equations.

## Chapter summary

**Theorem 7.1** (Cartan–Kähler I). *Take an analytic exterior differential system and a noncharacteristic analytic integral manifold $X$ with locally maximal rank polar equations in every tangent space. Then $X$ is a hypersurface in an analytic integral manifold,* locally unique *in that any two analytic integral manifolds in which $X$ is a hypersurface both contain a subset, open in both and containing $X$.*

A submanifold $R$ of $M$ *restrains* an integral manifold $X$ of an exterior differential system if the exterior differential system pulls back to $R$ to make $X$ noncharacteristic, with each tangent space of $X$ having locally maximal rank polar equations.

**Corollary 7.2.** *If an analytic manifold restrains an analytic integral manifold of an analytic exterior differential system, then the integral manifold is a hypersurface in a locally unique analytic integral submanifold of the restraining manifold.*

**Theorem 7.3** (Cartan–Kähler II). *Take an analytic exterior differential system $\mathcal{I}$ on a manifold $M$ and an analytic integral manifold $X$ with locally maximal rank polar equations. Take an analytic submanifold $R$ containing $X$ so that, for each $x \in X$,*

a. *all nonzero polar equations of $T_xX$ pull back to nonzero linear functions on $T_xR$ and*

b. *the vectors in $T_xR$ on which those polar equations vanish form a subspace containing $T_xX$ as a hyperplane.*

*Then $R$ restrains $X$, corollary 7.2 applies.*

*Proof.* When we pull back to $R$, the polar equations of integral elements near $T_xX$ also pull back without losing rank. So polar equations on $T_xX$ are locally maximal among all nearby integral elements tangent to $R$. At each $x \in X$, $T_xX$ is a hyperplane in a unique integral element in $T_xR$: the vanishing locus of the polar equations. So $X$ is noncharacteristic in $R$; apply theorem 7.1. □

**Corollary 7.4.** *If the tangent spaces of an analytic integral manifold of an analytic exterior differential system have locally maximal rank polar equations, then the following are equivalent:*

a. *that rank is less than the codimension of the integral manifold,*

b. *the integral manifold is covered in open sets, each of which is a hypersurface in an analytic integral manifold.*

*Proof.* In one tangent space $T_xM$, pick covectors, linearly independent modulo the polar equations of $T_xX$, so that the subspace on which they and the polar equations vanish is a hyperplane containing $T_xX$. That hyperplane can be any one that contains $T_xX$ and satisfies the polar equations of $T_xX$, so in particular we can make $T_xX$ noncharacteristic in it. Pick any submanifold $R$ normal to those covectors. □

## Extending an integral manifold

Proof of theorem 7.1 on the preceding page.

*Proof.* Let $p := 1 + \dim X$; denote $X$ as $X_{p-1}$. At each $x \in X_{p-1}$, there is a unique $p$-dimensional integral element $E_p \subset T_xM$ containing $E_{p-1} := T_xX_{p-1}$. By the solution of problem 6.8 on page 49, $s_p = 0$ on $E_p$.

By problems 5.6 and 5.7 on page 41, we can take a torsion-free tableau for $\mathcal{I}$, defined in some open set, adapted to a generic flag which includes $E_{p-1}$ and $E_p$. Each polar

$$\vartheta = \dots + \pi \wedge \omega^I + \dots$$

has $I \subseteq \{1, 2, \dots, p-1\}$, since $s_p = 0$. Let $J := \{1, 2, \dots, p-1\} - I$, and $\vartheta' := \vartheta \wedge \omega^J$. So $\vartheta'$ is $\vartheta$ "raised up" to become a $p$-form of grade $p-1$, shifting its polar to that grade. Similarly, if $\vartheta$ is one of the coframing 1-forms $\theta^a$, let

$$\vartheta' = \theta^a \wedge \omega^{1\dots p-1}.$$

Let $\mathcal{I}'$ be the exterior differential system generated by these $p$-forms $\vartheta'$.

**7.2** Prove that $X_{p-1}$ is a hypersurface in some $\mathcal{I}'$-integral manifold $X_p$.

We need to see that $\mathcal{I} = 0$ on $X_p$. The 1-forms $\omega^1, \dots, \omega^p$ coframe $X_p$ near $x$. Every $p$-form in $\mathcal{I}^p$ on $X_p$ is wedged up a form $\vartheta \in \mathcal{I}$ by some $\omega^J$. If $p \notin J$, then $\vartheta \wedge \omega^J \in \mathcal{I}'$ vanishes on $X_p$. So every $p$-form in $\mathcal{I}^p$ on $X_p$ is a multiple of $\omega^p$. We will see in lemma 7.6 on the following page that $\mathcal{I} = 0$ on $X_p$. □

## Background material from differential geometry

**Lemma 7.5.** *For any $(p+1)$-form $\vartheta$ and vector fields $v_0, \dots, v_p$,*

$$\begin{aligned} d\vartheta(v_0, \dots, v_p) = {} & (-1)^i v_i(\vartheta(v_0, \dots, \hat{v}_i, \dots, v_p)) \\ & + \sum_{i<j} (-1)^{i+j} \vartheta([v_i, v_j], v_0, \dots, \hat{v}_i, \dots, \hat{v}_j, \dots, v_p). \end{aligned}$$

**7.3** Prove lemma 7.5.

## Compatibility

Cartan's strategy picks out a determined subset of differential equations, sweeping an integral curve into a surface, and so on. Why is the surface integral? In other words, why are the differential equations sweeping out the surface compatible with those which swept out the curve?

**Lemma 7.6.** *On a connected p-dimensional manifold $X_+$, take a locally finitely generated analytic exterior differential system $\mathcal{I}$. Suppose that $\mathcal{I}$ has an analytic integral hypersurface. Suppose that $\mathcal{I}^p = \mathcal{I}^{p-1} \wedge \omega$ for some 1-form $\omega$. Then $\mathcal{I} = 0$.*

*Proof.*

**7.4** Why does it suffice to prove that $\mathcal{I}^{p-1} = 0$?

**7.5** Why does it suffice to assume that $\omega \neq 0$ everywhere on $X_+$?

Denote $\omega$ as $\omega^p$. Perhaps replacing $X_+$ by an open subset, choose analytic 1-forms $\omega^1, \dots, \omega^{p-1}$ on $X_+$ which coframe $X$ so that $\omega^1, \dots, \omega^p$ coframe $X_+$. Pick differential forms $\phi^a \in \mathcal{I}^{p-1}$ generating $\mathcal{I}$. In particular, $d\phi^a, \phi^a \wedge \omega^i \in \mathcal{I}^p = \mathcal{I}^{p-1} \wedge \omega^p$. Let $\phi$ be the column vector with entries $\phi^a$:

$$
\begin{aligned}
d\phi &= h\phi \wedge \omega^p, \\
\phi \wedge \omega^i &= h^i \phi \wedge \omega^p, i = 2, \dots, p-1
\end{aligned}
$$

for matrices $h, h^i$ whose entries are analytic functions.

Denote the vector fields dual to $\omega^1, \dots, \omega^p$ by $e_1, \dots, e_p$. Because $\omega^1, \dots, \omega^{p-1}$ coframe $X$, $e_p$ is not tangent to $X$. Write $e_{\hat{\imath}}$ to denote

$$e_1, \dots, \hat{e}_i, \dots, e_p.$$

**7.6** Use lemma 7.5 on the previous page to find a determined linear system satisfied by $f := \phi(e_{\hat{p}})$.

By the Cauchy–Kovalevskaya theorem (theorem A.2 on page 80), there is a unique solution $f$ near $X$ with given values on the noncharacteristic hypersurface $X$; these values are $f = 0$, so the solution is $f = 0$. So $\phi \wedge \omega^p = 0$, so $\phi \wedge \omega^i = 0$ for all $i$, so $\mathcal{I}^{p-1} = 0$. □

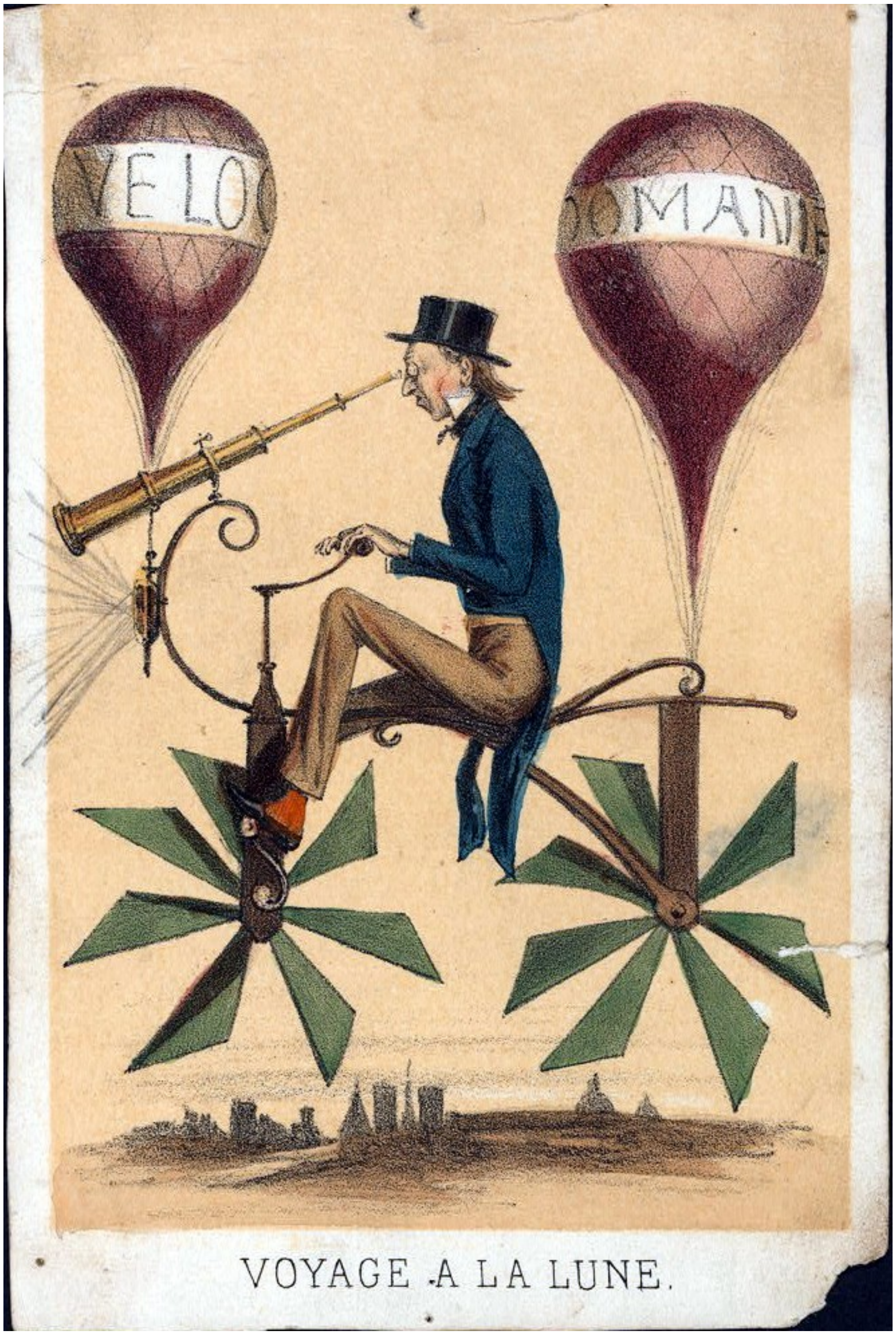
VELO
MANIE
VOYAGE A LA LUNE.

Chapter 8

# Cauchy characteristics

Often an exterior differential system can be written in a smaller number of variables than we would at first expect.

## Redefinition of exterior differential system

A *symmetry vector field* of an exterior differential system is a vector field whose flow preserves the system. It might seem natural to define a symmetry vector field as one whose flow permutes integral manifolds, but we have no test for this. The flow of a vector field might not be defined globally, while all of the forms in an exterior differential system are defined globally. Moreover, we want to test whether a vector field is a symmetry by local computation. So we consider a different concept of exterior differential system, in which the forms need only be defined locally.

An *exterior differential system* $\mathcal{I}$ on a manifold $M$ is an ideal $\mathcal{I}_U \subset \Omega^*_U$ of differential forms on each open set $U \subseteq M$ so that

*a.* $d$-closed: the exterior derivative takes $\mathcal{I}_U \to \mathcal{I}_U$ and

*b.* restricts: if $U \subseteq V \subseteq M$ are open sets, then restricting forms takes $\mathcal{I}_V \to \mathcal{I}_U$, and

*c.* glues: if $U = \bigcup_a U_a$, then a differential form belongs to $\mathcal{I}_U$ just when its restriction to each $U_a$ belongs to $\mathcal{I}_{U_a}$, and

*d.* graded: $\mathcal{I}_U = \mathcal{I}^1_U \oplus \mathcal{I}^2_U \oplus \dots \oplus \mathcal{I}^{\dim M}_U$, $\mathcal{I}^k_U := \mathcal{I}_U \cap \Omega^k_U$.

All of our theorems so far hold, with the same proofs, for this definition of exterior differential system.

**8.1** Give an example of an exterior differential system $\mathcal{J}$ in the sense of the old definition which is *not* the ideal $\mathcal{I}_M$ of an exterior differential system $\mathcal{I}$ in the sense of the new definition.

Take any collection of differential forms defined on various open subsets of a manifold. Without changing the submanifolds on which they vanish, we can add forms to our collection until we obtain an exterior differential system.

## Convergence

Differential forms *converge* when their component functions do in local coordinates, as analytic functions p. 97.

**Theorem 8.1.** *Every exterior differential system is closed under convergence.*

*Proof.* Take a convergent sequence $\vartheta_i \to \vartheta$, with $\vartheta_i \in \mathcal{I}$. Take coordinates with origin at some chosen point. A differential form $\vartheta = f_I dx^I$ is a vector valued map, valued in $\Lambda^* \mathbb{R}^n$. Germs of forms from $\mathcal{I}$ constitute a submodule of the germs of forms. By theorem C.4 on page 98, the germ of $\vartheta$ is among those germs. So $\vartheta \in \mathcal{I}_U$ for some open set $U$ around the chosen point. □

**8.2** Give an example of an exterior differential system, according to our old definition, not closed under convergence.

## Symmetries

**8.3** Prove that, for any vector field $v$ and differential form $\vartheta$ defined near some point, the point lies in an open set in which

$$e^{tv*}\vartheta = \sum \frac{t^k}{k!} \mathcal{L}_v^k \vartheta.$$

**8.4** For any exterior differential system $\mathcal{I}$ and vector field $v$, prove that the following are equivalent:

*a.* $v$ is a symmetry of $\mathcal{I}$,

*b.* $\mathcal{L}_v \mathcal{I}_U \subseteq \mathcal{I}_U$ for all open sets $U$,

*c.* $\mathcal{L}_v \mathcal{I}_U \subseteq \mathcal{I}_U$ for some open sets $U$ forming a basis for the topology of $M$.

**8.5** What are the symmetry vector fields of a Frobenius system?

**8.6** Prove that the symmetry vector fields of any exterior differential system form a Lie algebra.

**8.7** Give an example of a smooth exterior differential system $\mathcal{I}$ and a complete analytic vector field $v$, so that $\mathcal{L}_v \mathcal{I} \subseteq \mathcal{I}$ but the flow of $v$ does *not* preserve $\mathcal{I}$, *nor* permute integral manifolds.

**8.8** Give an example of an analytic exterior differential system $\mathcal{I}$ with no nonzero smooth symmetry vector field, so that the smooth exterior differential system it generates has a nonzero smooth symmetry vector field.

**8.9** Give an example of a complete vector field whose flow preserves the integral elements of an exterior differential system, but does not preserve the exterior differential system.

**8.10** For a smooth exterior differential system $\mathcal{I}$, prove *a.* if $\mathcal{I}$ is locally finitely generated, $v$ is a nowhere vanishing vector field, and $\mathcal{L}_v \mathcal{I} \subseteq \mathcal{I}$, then $v$ is a symmetry.

*b.* if $\mathcal{I}$ is closed under uniform convergence on compact sets with all derivatives, and $v$ is a symmetry, then $\mathcal{L}_v \mathcal{I} \subseteq \mathcal{I}$.

## Finite type

An exterior differential system has *finite type* if every point of $M$ lies in an open set $U \subseteq M$ on which there are finitely many forms $\vartheta_j \in \mathcal{I}_U$ so that, for any point $m \in U$, and form $\vartheta$ from $\mathcal{I}$ defined near $m$, there are forms $\phi_j$ defined near $m$, so that $\vartheta = \sum \phi_j \wedge \vartheta_j$ near $m$.

> Any exterior differential system generated by finitely many globally defined differential forms has finite type; this includes all of our examples.

**8.11** Give an example of an infinite type exterior differential system. Is your example in involution?

**8.12** Prove that, in the definition of finite type, we can always assume that these $\phi_j$ are functions, i.e. 0-forms.

**8.13** Prove compatibility (lemma 7.6 on page 56) for finite type exterior differential systems. Give an infinite type counterexample.

**8.14** Give an example of an infinite type involutive exterior differential system.

**8.15** Prove that, for any finite type exterior differential system, a vector field $v$ is a symmetry just when every point lies in an open set $U$ on which a finite set of forms $\vartheta_i$ generate $\mathcal{I}_U$ with $\mathcal{L}_v \vartheta_i \in \mathcal{I}_U$.

## Pointwise linear independence

This section can be omitted without loss of continuity.

An exterior differential system $\mathcal{I}$ is *bundled* if it is generated by forms in certain degrees, and in these degrees is locally spanned by pointwise linearly independent forms.

> Any exterior differential system generated by finitely many globally defined differential forms, all of the same degree, everywhere linearly independent, is bundled; this includes all of our examples.

**8.16** Prove that every bundled exterior differential system is of finite type.

**8.17** Give an example of an unbundled finite type exterior differential system $\mathcal{I}$ on a manifold $M$.

**8.18** Solve problem 8.4 on page 60 for bundled smooth systems.

**8.19** Give an example of an unbundled finite type involutive exterior differential system.

**8.20** We assume familiarity with vector bundles [6]. Prove that an exterior differential system $\mathcal{I}$ is bundled just when there are vector subbundles $I^p \subseteq \Lambda^p T^* M$ for various values of $p = p_1, p_2, \dots, p_\ell$ and $M$ is covered by open sets $U_a \subseteq M$ on which

*a.* $\mathcal{I}^p_{U_a}$ is the collection of differential forms from $I^p$ on $U_a$, $p = p_1, p_2, \dots, p_k$ and

*b.* $\mathcal{I}_{U_a}$ is generated by these sections.

## Cauchy characteristics

The exterior differential system generated by $dy - p\,dx$ on $M := \mathbb{R}^4_{x,y,p,q}$ doesn't make use of the variable $q$. We can build a quotient manifold $\bar{M} := \mathbb{R}^3_{x,y,p}$, map $(x, y, p, q) \in M \mapsto (x, y, p) \in \bar{M}$. Maximal integral manifolds in $M$ are locally the preimages of maximal integral manifolds of $dy - p\,dx$ on $\bar{M}$. Our aim in this section is to find "unused variables", and quotient them out.

A vector field $v$ is *Cauchy characteristic* for an exterior differential system $\mathcal{I}$ if $v \lrcorner \mathcal{I} \subseteq \mathcal{I}$.

**8.21** Prove that the Cauchy characteristic vector fields of an exterior differential system form a Lie subalgebra of the vector fields.

Denote by $\mathcal{I}_m$ the set of values $\vartheta_m \in \Lambda^* T_m M^*$ of forms $\vartheta \in \mathcal{I}_U$ for some open set $U$ with $m \in U$. A *Cauchy characteristic vector* of an exterior differential system $\mathcal{I}$ is a vector $v \in T_m M$ so that $v \lrcorner \mathcal{I}_m \subseteq \mathcal{I}_m$. The *rank* of Cauchy characteristic vectors at each point is their dimension as a vector space.

**8.22** We assume familiarity with vector bundles [6]. Prove that Cauchy characteristic vectors have constant rank just when they form a vector subbundle of the tangent bundle $TM$.

**8.23** Prove: if the Cauchy characteristic vectors have constant rank then a vector field is a Cauchy characteristic vector field if and only if its value at each point is a Cauchy characteristic vector.

**8.24** Give an example of an involutive exterior differential system, on a connnected manifold, whose Cauchy characteristic vectors do not have constant rank.

The *retracting space* of an exterior differential system $\mathcal{I}$ is the collection of 1-forms vanishing on its Cauchy characteristic vectors. A *Cauchy characteristic* is an integral submanifold of the retracting space. The *pullback* $\pi^*\bar{\mathcal{I}}$ by a map $\pi\colon M \to \bar{M}$ of an exterior differential system $\bar{\mathcal{I}}$ is the exterior differential system generated by pullbacks $\pi^*\bar{\vartheta}$ of differential forms $\bar{\vartheta}$ from $\bar{\mathcal{I}}$.

**8.25** If $\mathcal{I} = \pi^*\bar{\mathcal{I}}$ is the pullback by a submersion $\pi$, prove that a submanifold of $M$ is the $\pi$-preimage of an $\bar{\mathcal{I}}$-integral manifold just when it is an $\mathcal{I}$-integral manifold and contains any of fiber of $\pi$ it touches.

**8.26** How do characters behave when we pull back?

Clearly the pullback of a finite type system is finite type. The *pushforward* $\pi_*\mathcal{I}$ by a map $\pi\colon M \to \bar{M}$ of an exterior differential system $\mathcal{I}$ is the exterior differential system consisting of the forms $\bar{\vartheta}$ on $\bar{M}$ whose pullback $\pi^*\bar{\vartheta}$ lies in $\mathcal{I}$.

**8.27** Prove that $\bar{\mathcal{I}} \subseteq \pi_*\pi^*\bar{\mathcal{I}}$.

**8.28** Prove that $\pi^*\pi_*\mathcal{I} \subseteq \mathcal{I}$.

**8.29** Prove that the vectors on which $\pi_* = 0$ are Cauchy characteristic vectors for $\mathcal{I} := \pi^*\bar{\mathcal{I}}$, and that $\pi^*\Omega^1_{\bar{M}}$ lies in the retracting space of $\mathcal{I}$.

**8.30** Take an exterior differential system with characters $s_i$ at an integral element transverse to the fibers of a map. Suppose that its push forward by that map has characters $s_i'$. Prove that the restriction to each preimage of each integral manifold has characters $s_i - s_i'$.

**8.31** Suppose that $\mathcal{I}$ is an exterior differential system on a manifold $M$, and that the retracting space of $\mathcal{I}$ has constant rank. Prove:

*a.* The retracting space generates a Frobenius exterior differential system lying inside $\mathcal{I}$ and

*b.* every point of $M$ lies in an open set $U$ so that the retracting space of $\mathcal{I}_U$ is the pullback $\pi^*\Omega^1_{\bar{U}}$ of a surjective submersion $\pi\colon U \to \bar{U}$ to some manifold $\bar{U}$.

**8.32** Give an example of a finite type exterior differential system on a manifold $M$, and a submersion $\pi\colon M \to \bar{M}$ with Cauchy characteristic fibers, so that $\mathcal{I} \neq \pi^*\pi_*\mathcal{I}$.

**Theorem 8.2.** *Suppose that $\mathcal{I}$ is a finite type exterior differential system on a manifold $M$. Take a submersion $\pi\colon M \to \bar{M}$ with Cauchy characteristic fibers. Suppose that for any two components of any fiber of $\pi$, there is a diffeomorphism of $M$ preserving $\mathcal{I}$ and $\pi$ and interchanging these components. Then $\mathcal{I} = \pi^*\pi_*\mathcal{I}$ and $\pi_*\mathcal{I}$ has finite type.*

*Proof.* Let $\bar{\mathcal{I}} := \pi_* \mathcal{I}$. We know that $\pi^* \pi_* \mathcal{I} \subseteq \mathcal{I}$, i.e. that $\pi^* \bar{\mathcal{I}} \subseteq \mathcal{I}$.

If we can cover $M$ in open sets so that every element of $\mathcal{I}$ defined on one of those sets lies in $\pi^* \bar{\mathcal{I}}$, we glue. So we need only prove that, for any arbitrary point, elements of $\mathcal{I}$ defined near that point lie in $\pi^* \bar{\mathcal{I}}$, i.e. are multiples of pullbacks from $\bar{\mathcal{I}}$.

Take an open set $U \subseteq M$. Take a Cauchy characteristic vector field $v$ defined and nonzero near some point of $U$. Shrink $U$ if needed to arrange that $v$ is defined in $U$ and, by finite type, that the forms in $\mathcal{I}$ near any point of $U$ are generated by finitely many forms $\vartheta^i$. But $\mathcal{I}_U$ is $\mathcal{L}_v$-closed, so, after perhaps shrinking $U$ again,

$$\mathcal{L}_v \vartheta^i = f^i_j \vartheta^i,$$

for some functions $f^i_j$ on $U$; denote this equation $\mathcal{L}_v \vartheta = f \vartheta$. Pick a hypersurface $H \subset M$ through $m_0$ transverse to $v$. After perhaps shrinking $U$, we can define functions $g = (g^i_j)$ by $g = I$ along $H$ and

$$\mathcal{L}_v g = -g f.$$

If $\bar{\vartheta} := g\vartheta \in \mathcal{I}_U$, $\mathcal{L}_v \bar{\vartheta} = 0$, and $\vartheta = g^{-1} \bar{\vartheta}$ so $\mathcal{I}_U$ is generated by $v$-invariant forms.

Similarly, for any finite set of commuting and nonvanishing vector fields $v_1, \dots, v_k$, after perhaps shrinking $U$, we can generate $\mathcal{I}_U$ by forms invariant under all of these. In particular, after perhaps shrinking $U$, we can pick these vector fields to give a basis of local sections of the kernel of $\pi'$, i.e. of the vertical vector fields. Generators become invariant generators.

Take coordinates $x^i, y^a$ on $M$ so that these $x$ coordinates are pulled back from $\bar{M}$. Taking as $v$ the various $\partial_{y^a}$, we have seen that $\mathcal{I}_U$ is generated by invariant differential forms, i.e. forms which, expanded out as $\vartheta = f_{IA} dx^I \wedge dy^A$, have $f_{IA} = f_{IA}(x)$. We need to arrange that there are no $dy^a$ terms in our generators. Since $\vartheta \in \mathcal{I}_U$, we know that $v \lrcorner \vartheta \in \mathcal{I}_U$ for each $v = \partial_{y^a}$, and similarly if we wedge several $\partial_{y^a}$ vector fields in, so $f_{IA}(x) dx^I \in \mathcal{I}_U$, for each $dy^A$, so we can replace $\vartheta$ by the various $f_{IA}(x) dx^I$.

So $\mathcal{I}$ is locally generated by pulled back differential forms. Since Cauchy characteristic vector fields are symmetries, these pullback forms continue to generate on open sets invariant under Cauchy characteristics. □

**8.33** Finish the proof.

**8.34** Does theorem 8.2 on the preceding page generalize to infinite type?

**8.35** Prove that, for any surjective submersion $\pi \colon M \to \bar{M}$ and exterior differential system $\bar{\mathcal{I}}$ on $\bar{M}$, $\bar{\mathcal{I}} = \pi_* \pi^* \bar{\mathcal{I}}$.

**8.36** Under the hypotheses of theorem 8.2 on the previous page, prove that $\mathcal{I}$ is bundled just when $\pi_* \mathcal{I}$ is.

**8.37** Under the hypotheses of problem 8.35, prove that $\bar{\mathcal{I}}$ is bundled just when $\pi^* \bar{\mathcal{I}}$ is.

**8.38** Give an example of an exterior differential system which is pulled back via a map with connected Cauchy characteristic fibers, but not locally spanned by pointwise linearly independent forms.

**8.39** A exterior differential system is *local* if it is closed under locally finite sums. Prove a variant of theorem 8.2 on page 63, for smooth local exterior differential systems, using either the old or the new definition of exterior differential system. Give an analytic counterexample, for the old definition.

To find Cauchy characteristics: in some local coframing, write out a set of differential forms which span an exterior differential system. Any 1-forms of the coframing which do not appear in the spanning set of the system are dual to Cauchy characteristic vector fields.

## Example: surface invariants

Return to the study of surface invariants on page on page 22. Note that $\mathcal{I}$ has a Cauchy characteristic: the vector fields $v$ on which $\gamma_{12} \neq 0$ but

$$0 = \omega_1 = \omega_2 = \omega_3 = \gamma_{3i} - a_{ij}\omega_j = Da_{ij}.$$

We can quotient locally by this Cauchy characteristic, so that each integral manifold $X$ coframed by $\omega_1, \omega_2, \gamma_{12}$ projects to a surface $S$ in the quotient space of $M = \ulcorner\mathbb{E}^3 \times V$ by the orthogonal group of the plane. We can in addition quotient by a reflection $e_3 \mapsto -e_3, A \mapsto -A$. The quotient space $\bar{M}$ is the set of all choices of point $x \in \mathbb{E}^3$, plane $P$ through $x$, and symmetric quadratic form $A$ valued in the normal line to $P$ at $x$.

## Example: isometric immersion, Cauchy characteristics

Recall the isometric immersion notation:

- $S$ is a surface with Riemannian metric, with Gauss curvature $K$,
- $\ulcorner S$, $\ulcorner\mathbb{E}^3$ are the orthonormal frame bundles,
- Identify $\mathbb{R}^3_{a,b,c}$ with the set of $2 \times 2$ symmetric matrices

$$A = \begin{pmatrix} a & b \\ b & c \end{pmatrix}.$$

- $M'_0 \subset \ulcorner S \times \ulcorner\mathbb{E}^3 \times \mathbb{R}^3_{a,b,c}$ is the subset on which $A \neq 0$ and $K = \det A$.
- $\omega_1, \omega_2$ are the soldering forms on $S$, and $\alpha$ is the connection form,
- $\omega'_1, \omega'_2, \omega'_3$ are the soldering forms on $\ulcorner\mathbb{E}^3$, and $\gamma'_{12}, \gamma'_{23}, \gamma'_{31}$ are the connection forms; we let $\alpha' := \gamma'_{12}$.

We let

$$\begin{pmatrix} Da \\ Db \\ Dc \end{pmatrix} := \begin{pmatrix} da + 2b\alpha + a_1\omega_1 + a_2\omega_2, \\ db + (a-c)\alpha + b_1\omega_1 + b_2\omega_2, \\ dc + 2b\alpha + c_1\omega_1 + c_2\omega_2, \end{pmatrix}$$

with $a_1, a_2, b_1, b_2, c_1, c_2$ any functions chosen so that

$$\begin{aligned} a_2 &= b_1, \\ b_2 &= c_1, \\ 0 &= 2bb_1 - ac_1 - ca_1 - K_1, \\ 0 &= 2bb_2 - ac_2 - ca_2 - K_2. \end{aligned}$$

This ensures that $c\,Da + a\,Dc = 2b\,Db$. The ideal for isometric immersions on $M_0'$ is generated by the 1-forms

$$\begin{pmatrix} \gamma_{13}' \\ \gamma_{23}' \\ \alpha' - \alpha \end{pmatrix} - \begin{pmatrix} a & b \\ b & c \\ 0 & 0 \end{pmatrix} \begin{pmatrix} \omega_1 \\ \omega_2 \end{pmatrix},$$

which have exterior derivatives

$$-\begin{pmatrix} Da & Db & 0 \\ Db & Dc & 0 \\ 0 & 0 & 0 \end{pmatrix} \wedge \begin{pmatrix} \omega_1 \\ \omega_2 \\ \alpha \end{pmatrix} \quad \text{mod } \theta_1, \dots, \theta_6.$$

How do we spot Cauchy characteristics? On $M_0'$, the 1-forms

$$\omega_1, \omega_2, \omega_1' - \omega_1, \omega_2' - \omega_1, \omega_3', \alpha' - \alpha, \gamma_{13}' - (a\omega_1 + b\omega_2), \gamma_{23}' - (b\omega_1 + c\omega_2), Da, Db, Dc, \alpha + \alpha'$$

form a basis, except for the one relation $c\,Da + a\,Dc = 2b\,Db$. When we write out the tableau in our basis, we don't use the last basis element: $\alpha + \alpha'$. So in this basis, $\alpha + \alpha'$ is dual to a Cauchy characteristic vector field $v$, i.e. $v$ hooks to zero in every 1-form appearing in the tableau, and so hooks the exterior differential system into itself.

We can do a little better. Let $G$ be the group of all orthogonal $3 \times 3$ matrices $g$ preserving the vertical axis:

$$g = \begin{pmatrix} h & 0 \\ 0 & (-1)^k \end{pmatrix},$$

with $k = 0$ or 1, $h$ an orthogonal $2 \times 2$ matrix. Recall that $r_g^* \omega' = g^{-1}\omega'$ and $r_g^* \gamma' = g^{-1}\gamma' g$. Extend the $G$-action to $M_0'$:

$$r_g(x, e, x', e', A) = (x, eh, x', e'g, (-1)^k h^{-1} A h).$$

Problem D.9 on page 108 shows that this action preserves $M_0'$ and the exterior differential system. Check that

$$v = \left.\frac{d}{dt}\right|_{t=0} r_{g(t)},$$

where

$$g(t) = \begin{pmatrix} \cos t & -\sin t \\ \sin t & \cos t \end{pmatrix}.$$

The quotient space $\bar{M}_0'$ of $M_0'$ by the $G$-action is the space of choices of linear isometry $F\colon T_xS \to T_{x'}\mathbb{E}^3$ together with a quadratic form $q \neq 0$ on $T_xS$ so that $K = \det q$. By problem 8.33 on page 64, $\mathcal{I}'$ on $M_0'$ is pulled back from a unique exterior differential system $\bar{\mathcal{I}}'$ on $\bar{M}_0'$. But, unlike $M_0'$, the manifold $\bar{M}_0'$ does not have a canonical choice of coframing. We would struggle to write out the quotient exterior differential system $\bar{\mathcal{I}}'$. Except for $\alpha$, the other 1-forms in our tableau vanish on the Cauchy characteristics, so define a tableau for the quotient exterior differential system:

$$\begin{pmatrix} \theta_4 \\ \theta_5 \\ \theta_6 \end{pmatrix} = \begin{pmatrix} \gamma_1' \\ \gamma_2' \\ \alpha' - \alpha \end{pmatrix} - \begin{pmatrix} a & b \\ b & c \\ 0 & 0 \end{pmatrix} \begin{pmatrix} \omega_1 \\ \omega_2 \end{pmatrix},$$

$$d \begin{pmatrix} \theta_4 \\ \theta_5 \\ \theta_6 \end{pmatrix} = - \begin{pmatrix} Da & Db \\ Db & Dc \end{pmatrix} \wedge \begin{pmatrix} \omega_1 \\ \omega_2 \end{pmatrix} \quad \text{mod } \theta_1, \dots, \theta_6.$$

$$\begin{matrix} s_1 & s_2 \\ 2 & 0 \end{matrix}$$

## Example: isometric immersion, noncharacteristic data

The *asymptotic curves* of a surface are those on which the shape operator vanishes.

**8.40** Show that the characteristic curves of the isometric immersion problem on any integral surface are the asymptotic curves.

Since we have assumed that $a, b, c$ do not all simultaneously vanish, the shape operator is nowhere zero, so not every curve is asymptotic.

A *ribbon* is a choice of curve, its *spine*, and a *ruling line* at each point of the spine, perpendicular to the tangent line to the spine, and analytically varying along the spine. A ribbon is *nondegenerate* if the ruling line is nowhere perpendicular to the curvature vector of the spine. A nondegenerate ribbon has a unit tangent vector along its spine, and a curvature vector, and their cross product, so imposes a basis, orienting the perpendicular space to every ruling line.

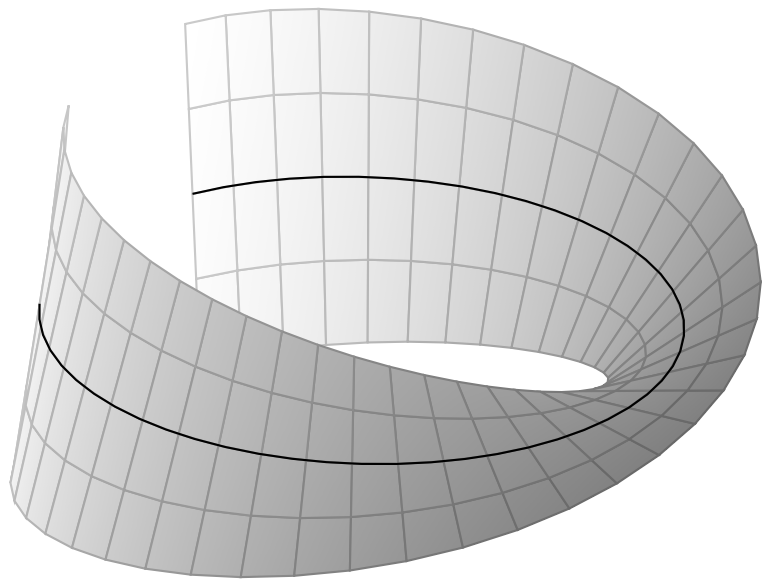

Take a surface $S$ with analytic Riemannian metric and an embedded connected analytic curve $C$ in $S$ with an orientation of $S$ defined along $C$. Take a nondegenerate ribbon in $\mathbb{E}^3$, with spine $C'$, and an analytic isometry $\iota\colon C \to C'$. Use the isometry to identify an orientation of $C$ with one of $C'$. Use the orientation defined along $C$ to extend $\iota_*$ to a linear isometry of each tangent plane to the perpendicular to each ruling line.

**Theorem 8.3.** *The isometry of curve to spine extends to a locally unique analytic isometric immersion normal to the ruling lines of its ribbon.*

Draw an infinitely long curve on the peel of an orange, accumulating only toward two points, so embedded in the sphere with those points deleted. Slice the peel close that curve, and lay out the peel tangent to an infinitely long ribbon: an isometric immersion.

Take a Möbius strip $S$ and a closed curve $C$ so that $S$ is not orientable in any neighborhood of $C$. Take a nondegenerate ribbon whose spine has the same length as $C$. An isometry of curve to spine exists, but cannot extend to an isometric immersion. It can't even extend to an isometry of tangent spaces of the surface to perpendicular planes to the ruling line, as the ruling lines are oriented.

*Proof.* Take $e'_1$ to be the unit tangent to $C'$ and $e'_3$ to be the unit tangent vector to the ruling lines. Locally extend $e'_1, e'_3$ into an orthonormal basis $e'_1, e'_2, e'_3$, giving an immersed curve in $\ulcorner\mathbb{E}^3$, on which $\gamma'_{31} = a\omega'_1$ and $\gamma'_{32} = b\omega'_1$ for some functions $a, b$. The ruling lines are nowhere perpendicular to the curvature vector, i.e. $e'_3$ is nowhere perpendicular, i.e. $a \neq 0$. Clearly we must then let

$$c := \frac{b^2 - K}{a}$$

at every point. Since $a \neq 0$, this curve in $M'$ is not characteristic:

$$a\omega_1^2 + 2b\omega_1\omega_2 + c\omega_2^2 = a\omega_1^2 \neq 0.$$

□

An isometric immersion arises as in this theorem precisely when it has nonzero shape operator everywhere.

Chapter 9

# Example: conformal maps of surfaces

We use the Cartan–Kähler theorem to prove the existence of local conformal maps between surfaces. We assume familiarity with appendix D.

## Conformal maps

A *conformal map* is a local diffeomorphism $\phi\colon S \to S'$ between surfaces in $\mathbb{E}^3$ preserving angles between curves.

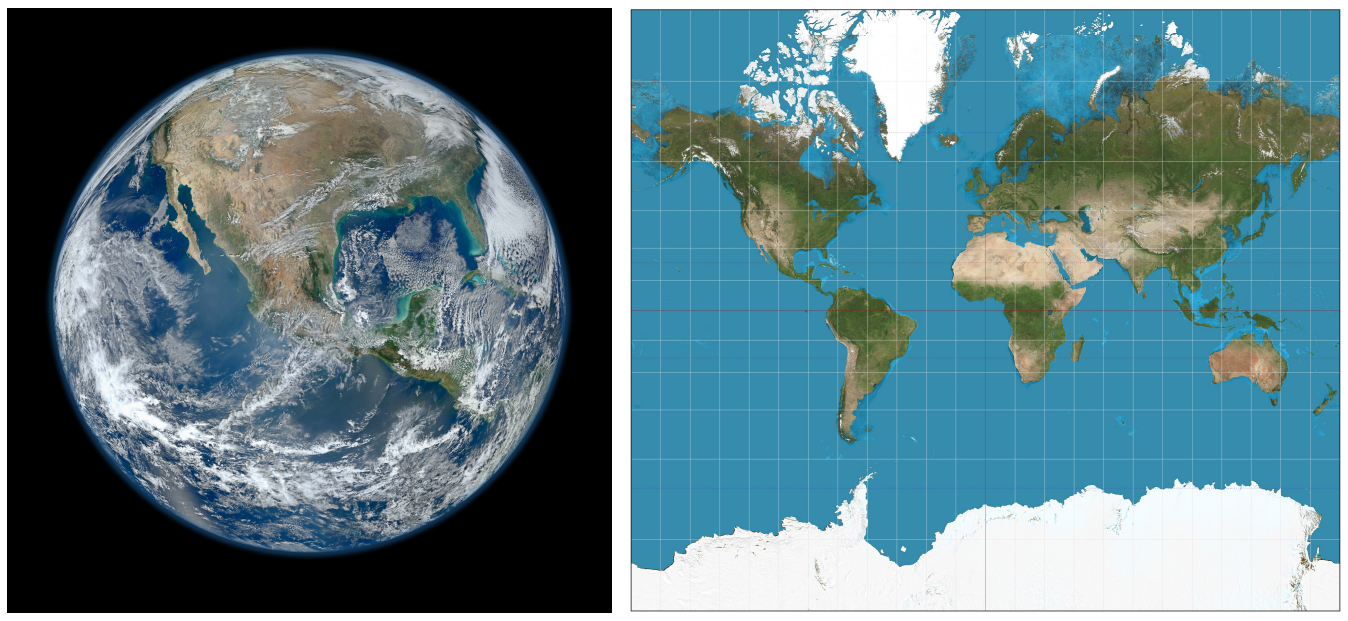

Mercator projection, a conformal map of an open subset of the sphere.
By Strebe - Own work, CC BY-SA 3.0,
https://commons.wikimedia.org/w/index.php?curid=16115307!

**Theorem 9.1.** *Given two surfaces $S, S'$ in $\mathbb{E}^3$, and points $x_0 \in S, x_0' \in S'$, there is a conformal map $\phi\colon U \to U'$ from an open set $U \subset S$ containing $x_0$ to an open set $U' \subset S'$ containing $x_0'$ so that $\phi(x_0) = x_0'$.*

See [10] for an elementary proof for smooth surfaces.

**9.1** Prove that a linear isomorphism of the plane preserves angles just it is uniquely expressed as a product of a rotation, a rescaling, and perhaps a reflection:

$$r \begin{pmatrix} \cos\theta & -\sin\theta \\ \sin\theta & \cos\theta \end{pmatrix}$$

or

$$r \begin{pmatrix} \cos\theta & -\sin\theta \\ \sin\theta & \cos\theta \end{pmatrix} \begin{pmatrix} 1 & 0 \\ 0 & -1 \end{pmatrix}.$$

**9.2** Given two surfaces of revolution, show that we can explicitly compute a rotationally invariant conformal map from one to the other by solving a first order ordinary differential equation for one function of one variable. If one of the surfaces is a cylinder or a sphere, explain how to reduce the construction of the conformal map to solving an integral, rather than an ordinary differential equation.

## The graph of a conformal map as an integral manifold

It is convenient to write the conformal scaling factor not as $r$ but instead as $e^{-u}$. Suppose that $\phi\colon S \to S'$ is a conformal map. Inside the 7-dimensional manifold $M := \ulcorner S \times \ulcorner S' \times \mathbb{R}_u$, consider the 3-dimensional submanifold $X$ consisting of points $(x, e, x', e', u)$ with $(x, e) \in \ulcorner S$ and $(x', e') \in \ulcorner S'$ so that

$$
\begin{aligned}
x' &= \phi(x), \\
e'_1 &= e^{-u}\phi'(x)e_1, \\
e'_2 &= e^{-u}\phi'(x)e_2.
\end{aligned}
$$

Recall the structure equations $d\omega = i\alpha \wedge \omega$ on $S$ and $d\omega' = i\alpha' \wedge \omega'$ on $S'$.

**9.3** On $X$, show that $\omega' = e^{-u}\omega$ and $du + i(\alpha' - \alpha) = u'\omega$ for a unique function $u'$ on $X$.

On the 7-dimensional manifold $M$, take the exterior differential system $\mathcal{I}$ generated by $\omega' - e^{-u}\omega$. Any conformal map $\phi\colon S \to S'$ has associated submanifold $X$ of $\ulcorner S \times \ulcorner S'$ an integral manifold. The tableau is

$$
d(\omega' - e^{-u}\omega) = e^{-u}(du + i(\alpha' - \alpha)) \wedge \omega.
$$

In real and imaginary parts, the tableau $du + i(\alpha' - \alpha)$ is

$$
\begin{pmatrix} du & -(\alpha' - \alpha) \\ \alpha' - \alpha & du \end{pmatrix}
$$

$$
\begin{array}{cc} s_1 & s_2 \\ 2 & 0 \end{array}
$$

The integral elements are the complex numbers $u'$ so that $du + i(\alpha' - \alpha) = u'\omega$, i.e. the real and imaginary parts of these numbers, $s_1 + 2s_2 = 2 + 2(0) = 2 = s$, involution: an integral manifold $X \subset M$ exists through any point of $M$.

Finally, we need to prove that our integral manifold is actually constructed from some conformal map $\phi$. Since $u$ is real, $du + i(\alpha' - \alpha)$ has real part $du$ and imaginary part $\alpha' - \alpha$. Note that $\mathcal{I}$ has a Cauchy characteristic: the vector fields $v$ on which $\alpha \neq 0$ but

$$
0 = \omega_1 = \omega_2 = \omega'_1 = \omega'_2 = du = \alpha' - \alpha.
$$

We can in addition quotient by a reflection $e_2 \mapsto -e_2, e'_2 \mapsto -e'_2$. The quotient space $\bar{M}$ is the set of all conformal linear maps from tangent spaces of $S$ to those of $S'$, as we are quotienting out by conformally changing the frame $e_1, e_2$ and correspondingly the frame $e'_1, e'_2$. So $\bar{M}$ is a 3-manifold, and $X$ projects to a surface $\bar{X}$ in $\bar{M}$. As $\omega_1, \omega_2, \alpha$ are linearly independent on $X$, $\bar{X}$ has $\omega_1, \omega_2$ still linearly independent, since both $\omega_1$ and $\omega_2$ vanish on the Cauchy characteristic vectors. Therefore $\bar{X}$ projects to $S$ by a local diffeomorphism. Similarly, $\bar{X}$ projects to $S'$ by a local diffeomorphism. So $\bar{X}$ injects into $S \times S'$, as the graph of a local diffeomorphism, say $\phi\colon S \to S'$.

We need to prove that $\phi$ is conformal. Take some tangent vector $v \in T_x S$. Pick a point

$$(x, e, x', e') \in X.$$

Write the vector $v$ as $v = v_1 e_1 + v_2 e_2$. We make a vector $\hat{v}$ on $X$ which projects to $v$, by asking that

$$\hat{v} \lrcorner \begin{pmatrix} \omega_1 \\ \omega_2 \\ \alpha \end{pmatrix} = \begin{pmatrix} v_1 \\ v_2 \\ 0 \end{pmatrix}.$$

Then $\hat{v}$ projects to $S'$ to a vector $v'$ with $v' = v'_1 e'_1 + v'_2 e'_2$ given by

$$\hat{v} \lrcorner \begin{pmatrix} \omega'_1 \\ \omega'_2 \end{pmatrix} = \begin{pmatrix} v'_1 \\ v'_2 \end{pmatrix}.$$

But on $X$, $\omega'_1 = e^{-u}\omega_1$, $\omega'_2 = e^{-u}\omega_2$ so

$$\begin{pmatrix} v'_1 \\ v'_2 \end{pmatrix} = e^{-u} \begin{pmatrix} v_1 \\ v_2 \end{pmatrix}.$$

In other words, $v' = e^{-u} v_1 e'_1 + e^{-u} v_2 e'_2$, so that $\phi$ is a conformal map.

## Characteristics

**Theorem 9.2.** *Take any two surfaces $S, S'$ in $\mathbb{E}^3$, and embedded curves $C \subset S$ and $C' \subset S'$. Every local diffeomorphism $\phi\colon C \to C'$ extends to a conformal map $\phi\colon U \to U'$ from an open set $U \subset S$ containing $C$ to an open set $U' \subset S'$ containing $C'$. Any two such maps agree on any connected open neighborhood of $C$ on which both are defined.*

Note that the curves $C, C'$ may be closed curves here; the result is global as regards the curves, but local in that we only construct a conformal map near the curves.

> Near each point of any surface, there is a conformal map to the plane, i.e. there are coordinates $x, y$, called *isothermal*, identifying angle measurements with those of the plane

*Proof.* The characteristic variety of our tableau:

$$\begin{pmatrix} \pi^1 & -\pi^2 \\ \pi^2 & \pi^1 \end{pmatrix} = \begin{pmatrix} du & -(\alpha' - \alpha) \\ \alpha' - \alpha & du \end{pmatrix}$$

Replace each polar $\pi^\alpha$ in the tableau with $v^\alpha \xi$:

$$0 = \begin{pmatrix} v^1\xi & v^2\xi \\ -v^2\xi & v^1\xi \end{pmatrix} \wedge \begin{pmatrix} \omega^1 \\ \omega^2 \end{pmatrix} = \begin{pmatrix} v^2\xi_1 - v^1\xi_2 \\ v^1\xi_1 + v^2\xi_2 \end{pmatrix} \omega^{12} = \begin{pmatrix} -\xi_2 & \xi_1 \\ \xi_1 & \xi_2 \end{pmatrix} \begin{pmatrix} v^1 \\ v^2 \end{pmatrix} \omega^{12}.$$

The characteristic variety equation is the determinant of the symbol matrix: $\xi_1^2 + \xi_2^2 = 0$, i.e. there are no real characteristics, an elliptic determined system after we mod out the Cauchy characteristics. So any integral curve in $\bar{M}$ lies in a unique integral surface. A curve in the 6-dimensional manifold $\bar{M}$ is a choice of curve $C$ in $S$, a curve $C'$ in $S'$, a diffeomorphism between the curves, and a conformal factor $e^{-u}$ along the curve. Such a curve is an integral curve of the system just when $\omega_i' = e^{-u}\omega_i$, i.e. $u$ is determined by the ratio of the lengths of the tangent vector to $C$ and that to $C'$. $\square$

Chapter 10

# Example: Weingarten surfaces

## Weingarten surfaces

A *Weingarten surface* is an oriented surface $S$ in $\mathbb{E}^3$ whose Gauss and mean curvature satisfy some relation, i.e.

$$(H, K): S \to W,$$

for some curve $W$ in the plane. For example, we could ask that $K = 1$ or $H = 0$. The "generic" surface has no such relation, as we saw on page 22. Every surface of revolution has such a relation: the Gauss and mean curvature are invariant under the revolution, so have values determined along any one meridian.

## Weingarten surfaces as integral manifolds

**10.1** Prove that, on any surface in $\mathbb{E}^3$, $H^2 \geq K$ with equality just at umbilic points.

If a surface consists entirely of umbilic points, it is a plane or sphere. So suppose that $W$ is a curve in the plane, lying in the open set of points $(x, y) \in \mathbb{R}^2$ so that $x^2 > y$, and $S$ is a Weingarten surface associated to $W$. On its frame bundle $\ulcorner S$ in $\ulcorner \mathbb{E}^3$, we have

$$
\begin{aligned}
\omega_3 &= 0, \\
\begin{pmatrix} \gamma_{13} \\ \gamma_{23} \end{pmatrix} &= \begin{pmatrix} a_{11} & a_{12} \\ a_{12} & a_{22} \end{pmatrix} \begin{pmatrix} \omega_1 \\ \omega_2 \end{pmatrix}, \\
0 &= a_{12} - a_{21}, \\
K &= a_{11}a_{22} - a_{12}^2, \\
H &= \frac{a_{11} + a_{22}}{2}, \\
(H, K) &\in W.
\end{aligned}
$$

Let $\hat{W}$ be the set of symmetric matrices

$$a = \begin{pmatrix} a_{11} & a_{12} \\ a_{12} & a_{22} \end{pmatrix}$$

so that

$$\left( \frac{\operatorname{tr} a}{2}, \det a \right) \in W.$$

Then $X := {}^{\ulcorner}S$ is a 3-dimensional integral manifold in $M := {}^{\ulcorner}\mathbb{E}^3 \times \hat{W}$ of

$$\begin{pmatrix} \theta_0 \\ \theta_1 \\ \theta_2 \end{pmatrix} = \begin{pmatrix} \omega_3 \\ \gamma_{13} - (a_{11}\omega_1 + a_{12}\omega_2) \\ \gamma_{23} - (a_{21}\omega_1 + a_{22}\omega_2) \end{pmatrix}$$

on which $\omega_1, \omega_2, \omega_{12}$ are linearly independent. Calculate the tableau:

$$d \begin{pmatrix} \theta_0 \\ \theta_1 \\ \theta_2 \end{pmatrix} = - \begin{pmatrix} 0 & 0 & 0 \\ \pi_1 & \pi_2 & 0 \\ \pi_2 & \pi_3 & 0 \end{pmatrix} \wedge \begin{pmatrix} \omega_1 \\ \omega_2 \\ \gamma_{12} \end{pmatrix}$$

where

$$\begin{pmatrix} \pi_1 \\ \pi_2 \\ \pi_3 \end{pmatrix} = d \begin{pmatrix} a_{11} \\ a_{12} \\ a_{22} \end{pmatrix} + \begin{pmatrix} 0 & 2 & 0 \\ -1 & 0 & 1 \\ 0 & -2 & 0 \end{pmatrix} \begin{pmatrix} a_{11} \\ a_{12} \\ a_{22} \end{pmatrix} \omega_{12}.$$

Locally, we can write $W$ as the set of solutions of an equation $f(x, y) = 0$ in the plane with $df \neq 0$. So on $S$,

$$0 = f\left(\frac{\operatorname{tr} a}{2}, \det a\right).$$

Let

$$f_H := \frac{\partial f}{\partial H}, f_K := \frac{\partial f}{\partial K}.$$

Compute out that this gives

$$0 = \left(\frac{da_{11} + da_{22}}{2}\right) f_H + (a_{22}da_{11} + a_{11}da_{22} - 2a_{12}da_{12}) f_K.$$

In terms of $\pi_1, \pi_2, \pi_3$ this relation is

$$0 = \left(\frac{\pi_1 + \pi_3}{2}\right) f_H + (a_{22}\pi_1 + a_{11}\pi_3 - 2a_{12}\pi_2) f_K$$

This equation has coefficients of $\pi_1, \pi_2, \pi_3$ given by

$$a_{11}f_K + f_H, a_{22}f_K + f_H, a_{12}f_K.$$

**10.2** Prove that all of these vanish, i.e. there is no addition linear relation among $\pi_1, \pi_2, \pi_3$, precisely when $a_{11} = a_{22}$ and $a_{12} = 0$, i.e. an umbilic point.

Since our curve $W$ lies inside $H^2 > K$ ("away from umbilic points"), our exterior differential system has characters $s_1 = 2, s_2 = 0, s_3 = 0$ so involution.

## Cauchy characteristics

The Cauchy characteristics are the rotations of frame tangent to the surface. The 6-dimensional quotient manifold $\bar{M}$ is the set of choices of point in $\mathbb{E}^3$, plane through that point, unit normal vector to that plane, and symmetric bilinear form on that plane, with half trace and determinant lying in $W$. On $\bar{M}$, the exterior differential system is determined. Each Weingarten surface $S$ gives an integral surface of that exterior differential system, mapping each point of $S$ to its tangent plane and shape operator on that tangent plane.

## Characteristic variety

The symbol matrix is

$$\begin{pmatrix} -\xi_2 & \xi_1 \\ \frac{f_H}{2}\xi_1 + f_K\left(a_{22}\xi_1 - a_{12}\xi_2\right) & \frac{f_H}{2}\xi_2 + f_K\left(a_{11}\xi_2 - a_{12}\xi_1\right) \end{pmatrix}$$

which has determinant

$$-f_K\left(a_{11}\xi_2^2 - 2a_{12}\xi_1\xi_2 + a_{22}\xi_1^2\right) - \frac{f_H}{2}\left(\xi_1^2 + \xi_2^2\right).$$

Recall that the characteristic variety consists of the hyperplanes

$$0 = \sum_i \xi_i\omega_i = 0,$$

satisfying these equations. A vector $v = v_1e_1 + v_2e_2$ lies in such a characteristic hyperplane just when $0 = \xi_1v_1 + \xi_2v_2$, so then, up to scaling

$$(\xi_1, \xi_2) = (v_2, -v_1).$$

Plug this in to see that the characteristics are

$$0 = f_K\left(a_{11}v_1^2 + 2a_{12}v_1v_2 + a_{22}v_2^2\right) + \frac{f_H}{2}\left(v_1^2 + v_2^2\right).$$

i.e. in classical notation,

$$0 = f_K \mathit{II} + \frac{f_H}{2} I.$$

So the characteristics are the curves on $S$ with velocity $v$ satisfying this quadratic equation. Since we have assumed that our surface contains no umbilic points and that $df \neq 0$, not every curve is characteristic.

## Initial data

Take a ribbon along a curve $C$ in $\mathbb{E}^3$. At each point of that curve, draw the perpendicular plane to the ruling line of the ribbon. On that plane, take a symmetric bilinear form $I\!I$ with two distinct eigenvalues, analytically varying along $C$. Let $H := \operatorname{tr} I\!I/2$, $K := \det I\!I$. The form $I\!I$ is *nondegenerate* if, on every tangent line to $C$, the homogeneous cubic form

$$dH\, I\!I + \frac{dK}{2} I$$

is not zero, where $I$ is the Euclidean inner product.

Consider what a noncharacteristic curve looks like in the 6-dimensional manifold $M = {}^{\ulcorner}\mathbb{E}^3 \times \hat{W}$. Such a curve consists of a ribbon $x(t), e_3(t)$, together with the additional data of $e_1, e_2$ and the values $a_{ij}$. The curve $x, e, a$ is an integral curve of the exterior differential system, i.e. $\omega_3 = 0$ and $\gamma_{i3} = a_{ij}\omega_j$. The equation $\omega_3 = 0$ is just the requirement that $e_3 \perp \dot{x}$, i.e. a ribbon. Since we can rotate the frame $e_1, e_2$, we can ask that $e_1$ be tangent to the curve $x(t)$. So then $\gamma_{i3} = a_{ij}\omega_j$ just when

$$\begin{aligned} I\!I(e_1, e_1) &= -e_3 \cdot \frac{de_1}{dt}, \\ I\!I(e_1, e_2) &= -e_3 \cdot \frac{de_2}{dt}. \end{aligned}$$

So the shape operator is partly determined by the ribbon. Noncharacteristicity is precisely that

$$e_3 \cdot \frac{de_1}{dt} f_K \neq \frac{f_H}{2}.$$

**10.3** Prove that the Weingarten equation $f = 0$ locally recovers the coefficient $a_{22}$, hence the entire shape operator, from the data of the ribbon.

Parameterize a curve $C$ by arc length as $x(s)$, and let $e_1 := \dot{x}$. Take a ribbon on that curve and write the direction of the ruling line as $e_3$. Let $e_2$ be the vector so that $e_1, e_2, e_3$ is a positively oriented orthonormal frame along $C$. A form $I\!I$ is *compatible* with the ribbon if $I\!I(e_1, e_1) = -e_3 \cdot \dot{e}_1$ and $I\!I(e_1, e_2) = -e_3 \cdot \dot{e}_2$. Compatibility is independent of the choice of arc length parameterization and of rotation of $e_1, e_2$. Define $H$ and $K$, the trace and determinant of $I\!I$. Define an immersed curve $W$: the image of the $s \mapsto (H(s), K(s))$.

**Theorem 10.1.** *Take an analytic ribbon along a connected curve $C$ in $\mathbb{E}^3$, and a symmetric bilinear form $I\!I$, defined in the perpendicular plane of the ruling line of the ribbon at one point of $C$, nondegenerate and compatible with the ribbon. Then $I\!I$ extends uniquely locally to be defined along an open subset of $C$, analytically varying, nondegenerate and compatible with the ribbon. If $I\!I$ extends to all of $C$, then there is an analytic Weingarten surface $S$ in $\mathbb{E}^3$ containing $C$ whose shape operator is $I\!I$ at each point of $C$ and whose Gauss and mean curvature lie in the image of $C$ in the plane under the map $(H, K)$. Any two such surfaces agree near $C$.*

**10.4** Which analytic curves in $\mathbb{E}^3$ are geodesics on minimal surfaces?

Appendix A

# The Cauchy–Kovalevskaya theorem

We prove the Cauchy–Kovalevskaya theorem: analytic determined systems of partial differential equations have local solutions, with arbitrary initial conditions.

## Formal Taylor series

For $x_1, x_2, \dots, x_n$ real variables and $a_1, a_2, \dots, a_n$ nonnegative integers, let

$$
\begin{aligned}
x &:= (x_1, x_2, \dots, x_n),\\
a &:= (a_1, a_2, \dots, a_n),\\
x^a &:= x_1^{a_1} \dots x_n^{a_n},\\
a! &:= a_1! \dots a_n! \text{ and}\\
\partial^a &:= \frac{\partial^{a_1}}{\partial x_1^{a_1}} \dots \frac{\partial^{a_n}}{\partial x_n^{a_n}}.
\end{aligned}
$$

A *formal Taylor series* is an expression $\sum c_a (x - x_0)^a$ with real constants $c_a$, not required to converge.

**A.1** Prove that the formal Taylor series $\sum_n (2n)! t^n$ diverges for $t \ne 0$.

We add, subtract, multiply, differentiate and compose in the obvious way: finitely many terms at a time. Crucially, each output term depends only on input terms of lower or equal order (or, for the derivative, on just one order higher), so on only *finitely* many input terms. When we add, multiply, differentiate or compose, each step is only adding or multiplying coefficients. In particular, the sum, product, derivative (in any variable) and composition of formal Taylor series with positive terms has positive terms.

A formal Taylor series $\sum b_a x^a$ *majorizes* another $\sum c_a x^a$ if $b_a \ge |c_a|$ for all $a$. If a convergent series majorizes another, the other is absolutely convergent. If $f$ majorizes $g$ and $u$ majorizes $v$ then $f \circ u$, $f + u$, $fu$, $\partial^a f$ majorizes $g \circ v$, $g + v$, $gv$, $\partial^a g$ respectively.

The geometric series

$$
\frac{1}{1-t} = 1 + t + t^2 + \dots
$$

has a Taylor series with positive coefficients. A formal Taylor series $c_0 + c_1 t + \dots$ with all $|c_j| < 1$ converges absolutely for $|t| < 1$, because the geometric series converges and majorizes it. Rescaling $t$ and our series, we see that if a formal Taylor series is majorized by a geometric series, then it converges near the origin, i.e. if $|c_j|$ are bounded by a geometric series, in other words $|c_j| \le Mr^j$ for some $M, r > 0$, then $c_0 + c_1 t + \dots$ converges absolutely near the origin. An *analytic function* is one which is locally the sum of a convergent Taylor series.

**Lemma A.1.** *A formal Taylor series converges to an analytic function just when it is majorized by a product of geometric series in one variable each.*

*Proof.* We give the proof for one variable around the point $x = 0$, and let the reader generalize. Take any analytic function $f(x)$ with convergent Taylor series $f(x) = \sum c_n x^n$. Since it converges absolutely for $x$ near 0, $\sum |c_n| r^n$ converges for $r$ near 0, so the terms of this series are bounded. Rescale to get a bound of 1, i.e. $|c_n| r^n < 1$ for all $n$. Therefore

$$|c_n| \le \frac{1}{r^n}$$

i.e. $f(x)$ is majorized by $1/(1-x)$. □

**A.2** Prove that any two formal Taylor series which converge near the origin to the same analytic function agree.

**A.3** Prove that an analytic function given by a convergent Taylor series around some point is also given by a convergent Taylor series around every nearby point.

**A.4** Prove that any two analytic functions defined on a connected open set which agree near some point agree everywhere.

The derivative and composition of analytic functions is analytic, and the analytic inverse and implicit function theorems hold with the usual proofs [24, 31].

## Complexification

Every Taylor series converging in some domain of real variables continues to converge for complex values of those variables, with small imaginary parts, by the same majorization. Hence every analytic function of real variables extends to an open set in a domain of complex variables. By elementary complex analysis [1], sums, differences, products and derivatives of analytic functions on an open set are analytic, as are quotients, wherever the denominator doesn't vanish, and the analytic inverse and implicit function theorems hold.

## Solving for Taylor series using a differential equation

Let's solve a simple partial differential equation $u_t = u_x$. Think about the vector field $X = \partial_t - \partial_x$ in the $tx$-plane.

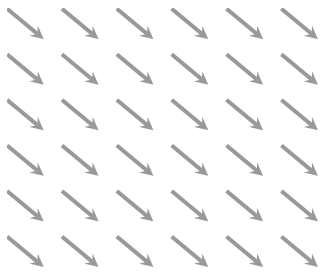

Clearly $Xu = 0$ just when $u$ satisfies our equation; but this happens just when $u$ is constant along the flow lines of $X$. The flow lines of $X$ are the diagonal lines, i.e. $t + x$ constant, so the solutions of our partial differential equation are $u(t, x) = U(t + x)$ for any function $U$. Solutions of partial differential equations can involve arbitrary functions, as in this example. Each solution $u(t, x)$ is determined by its value at "time" $t = 0$. If we pick our initial data function $U$ not smooth at a point, then $u(t, x)$ is not smooth along the flow line of that point. So solutions can be worse than the initial data from which we construct them.

**A.5** Suppose that $X$ is an analytic vector field on a manifold $M$ and that $H \subset M$ is an embedded analytic hypersurface and $X$ is not tangent to $H$ at any point of $H$. Pick analytic functions $v \colon H \to \mathbb{R}$ and $f \colon M \times \mathbb{R} \to \mathbb{R}$. Prove that there is an analytic function $u \colon \text{open} \subset M \to \mathbb{R}$ so that $u(m, 0) = v(m)$ for $m \in H$ and $Xu(m) = f(m, u(m))$ wherever $u$ is defined. Prove that any two such functions $u = u_1$ and $u = u_2$ agree near $H$.

Take the differential equation $u_t = u^2 + t$. Differentiate both sides:

$$
\begin{aligned}
u_{tt} &= 2u\, u_t + 1, \\
&= 2u \left(u^2 + t\right) + 1.
\end{aligned}
$$

Our differential equation allows us to replace differentiation in $t$ with some expression without any derivatives.

The differential equation $u_t = u^2 u_{xx} + \sin(u)$ allows us, using the chain rule, to replace any number of $t$ derivatives by complicated expressions in $x$ derivatives. We can calculate any Taylor coefficient in any power of $t$ as a function of finitely many Taylor coefficients involving only powers of $x$, i.e. Taylor coefficients of $u(x, 0)$. Danger: after we compute all of the Taylor coefficients, we have *no guarantee* that they sum up to a convergent Taylor series.

A *formal solution* of an analytic differential equation is a formal Taylor series so that expanding out the compositions term by term, all terms in the equation turn to zeroes. As above, for any formal Tayor series $u(0, x)$, there is an associated formal Taylor series solution $u(t, x)$ of each partial differential equation $u_t = f(t, x, u, u_x, u_{xx}, \dots, u_{x\dots x})$ for any smooth function or formal Taylor series $f$.

The heat equation $u_t = u_{xx}$ has a unique formal Taylor series solution with intial condition

$$u\big|_{t=0} = \frac{1}{1+x^2},$$

given by

$$u(t,x) = \sum_{j,k=0}^{\infty} (-1)^{j+k} (2(j+k))! \frac{t^j}{j!} \frac{x^{2k}}{(2k)!}$$

which doesn't converge anywhere except at $t = 0$.

## The method of majorants

**Theorem A.2** (Cauchy–Kovalevskaya)**.** *An analytic system of partial differential equations of the form $u_t = f(t,x,u,u_x)$, defined in an open set of $t \in \mathbb{R}, x \in \mathbb{R}^n, u \in \mathbb{R}^p$, with any analytic initial conditions $u = U(x)$ at $t = t_0$, has an analytic solution near $t = t_0$. Any two solutions agree near $t = t_0$.*

*Proof.* Instead of looking at Taylor series of solutions, look at Taylor series of equations. Start with a simpler problem: take a differential equation $u_t = f(t,u)$ with initial condition $u(0) = u_0$. Go back to our expressions for Taylor coefficients:

$$\begin{aligned}
u_t(0) &= f(0,u_0)\,, \\
u_{tt}(0) &= f_t(0,u_0) + f_u(0,u_0)\, u_t, \\
&= f_t(0,u_0) + f_u(0,u_0)\, f(0,u_0)\,, \\
u_{ttt}(0) &= \dots.
\end{aligned}$$

If $u_0 \ge 0$ and all of the Taylor coefficients of $f$ are $\ge 0$ then inductively all of the Taylor coefficients of $u$ are $\ge 0$. In other words, if $u_0 \ge 0$ and $f$ majorizes 0 then $u$ majorizes 0.

By the same reasoning, if we have two differential equations $u_t = f(t,u)$ and $v_t = g(t,v)$ with initial conditions $u(0) \ge v(0)$, and if $f$ majorizes $g$, then by the same induction, $u$ majorizes $v$. In particular, if the Taylor series of $u$ converges, then so does that of $v$.

**A.6** Check that the function

$$u(t) = 1 - \sqrt{1 + 2\log(1-t)}$$

satisfies the *toy* equation

$$u_t = \frac{1}{(1-t)(1-u)}$$

with initial condition $u = 0$ at $t = 0$.

Don't expand the unpleasant function $u(t)$ but expand the toy equation

$$
\begin{aligned}
u_t &= \frac{1}{(1-t)(1-u)}, \\
&= \left(1 + t + t^2 + \dots\right)\left(1 + u + u^2 + \dots\right), \\
&= \sum_{jk} t^j u^k :
\end{aligned}
$$

a geometric series. Suitable rescaling of $t$ and $u$ produces any convergent geometric series we like. By lemma A.1 on page 78, any analytic function $f(t,u)$ is majorized by a geometric series in $t,u$. So the equation $u_t = f(t,u)$ is majorized by the toy equation, after some rescaling of variables. So the solution $u(t)$ to $u_t = f(t,u)$ with $u(0) = 0$ is majorized by the toy solution, so has convergent Taylor series.

We will generalize this toy example several times, to get a larger class of examples of majorizing equations. The function

$$
u(t,x) = 1 - x - \sqrt{(1-x)^2 - 2t}
$$

satisfies

$$
u_t = \frac{1 + u_x}{1 - x}
$$

with $u = 0$ at $t = 0$. With suitable rescalings of $x, t, u$, this equation majorizes any equation of the form $u_t = f(x,u) + g(x,u)\,u_x$ for any analytic $f, g$, with $x, u \in \mathbb{R}$ in a suitable open set. Therefore all such equations have local analytic solutions so that $u = 0$ when $t = 0$.

To allow more $x$ variables: if $u(t,s)$ is the function above

$$
u(t,s) = 1 - s - \sqrt{(1-s)^2 - 2t}
$$

and we let $v(t,x) = u(t,s)$ where we let $s = \sum_i x_i$ where $x \in \mathbb{R}^n$, then

$$
v_t = \frac{1 + \frac{1}{n}\sum v_{x_i}}{1 - \sum x_i}
$$

and $v = 0$ at $t = 0$. Again, with a little rescaling, this equation majorizes any equation of the form $u_t = f(x,u) + g(x,u)\,u_x$ for any analytic $f, g$, and $t, u \in \mathbb{R}$, $x \in \mathbb{R}^n$ in some open set. Of course the same trick works even if we allow many $u$ functions, i.e. $u \in \mathbb{R}^p$.

To allow the $t$ variable to enter into the right hand side, we add a new equation $v_t = 1$ with initial condition $v = 0$ at $t = 0$, so we are just forcing $v = t$ everywhere. Then any system like $u_t = f(t,x,u) + g(t,x,u)u_x$ with $u = 0$ when $t = 0$ is equivalent to

$$
\begin{aligned}
u_t &= f(v,x,u) + g(v,x,u)u_x, \\
v_t &= 1,
\end{aligned}
$$

with $u = v = 0$ when $t = 0$. We have already seen that such systems have solutions near the origin,

If we want to solve a system of the form $u_t = f(t, x, u, u_x)$, invent a new variable $v$ and solve instead

$$
\begin{aligned}
u_t &= f(t, x, u, v)\,, \\
v_t &= f_x(t, x, u, v) + f_u(t, x, u, v)\, v + f_v(t, x, u, v)\, v_x,
\end{aligned}
$$

with $u = v = 0$ when $t = 0$, to force $v = u_x$.

Suppose that we want to solve a system $u_t = f(t, x, u, u_x)$ with $u = U(x)$ at $t = 0$ instead of $u = 0$ at $t = 0$. Replace the system by the system $v_t = f(t, x, v+U, v_x+U_x)$ with $v = 0$ at $t = 0$ and then let $u(t, x) = v(t, x)+U(x)$. □

Given a solution $u$ to the wave equation $u_{tt} = u_{xx}$ with initial conditions $u = U$ and $u_t = W$ at $t = 0$, let

$$
v(t, x) = \int W(x)\, dx + \int_0^t u_x(s, x)\, ds
$$

so that $v_t = u_x$. Differentiate under the integral sign and apply the fundamental theorem of calculus to find that $v_x = u_t$. Therefore $u, v$ satisfy

$$
\begin{aligned}
u_t &= v_x, \\
v_t &= u_x
\end{aligned}
$$

with initial conditions $u = U$ and $v = W$ at $t = 0$. Conversely, any solution to this system of equations recovers our solution $u$ to the wave equation. An analytic solution of the wave equation is uniquely determined by the initial height and velocity of the waves.

To solve a second order system

$$
u_{tt} = f(t, x, u, u_x, u_t, u_{xx}, u_{xt})\,,
$$

with initial conditions $u = U(x)$ and $u_x = W(x)$ at $t = 0$, let $p = u_x$ and $q = u_t$. Then

$$
\begin{aligned}
u_t &= q, \\
p_t &= q_x, \\
q_t &= f(t, x, u, p, q, p_x, q_x)\,,
\end{aligned}
$$

with initial conditions $u = U(x)$, $p = U'(x)$, $q = W(x)$ at $t = 0$. Any solution or formal solution of this system arises from a solution or formal solution to the original system and vice versa. This same trick reduces any system of any order to a first order system in more variables.

Is there a vector field $u$ so that $\nabla \times u = f$ for a given vector field $f$ in $\mathbb{E}^3$? Write $u = \left(u^1, u^2, u^3\right)$ and $f = \left(f^1, f^2, f^3\right)$ and let

$$u^1_{x_2} := \frac{\partial u^1}{\partial x_2}, \text{etc.}$$

Our equation $\nabla \times u = f$ expands out to

$$\begin{aligned} u^3_{x_2} - u^2_{x_3} &= f^1, \\ u^1_{x_3} - u^3_{x_1} &= f^2, \\ u^2_{x_1} - u^1_{x_2} &= f^3 : \end{aligned}$$

3 first order equations for 3 unknowns. By analogy with the Cauchy–Kovalevskaya theorem, we expect to find a unique solution with initial values of $u^1, u^2, u^3$ given on some surface in $\mathbb{E}^3_{x_1,x_2,x_3}$. This expectation is *not* correct. Taking $\nabla\cdot$ on both sides of our equations reveals $0 = \nabla \cdot f$. If $0 \neq \nabla \cdot f$ then there is no solution $u$. If $0 = \nabla \cdot f$ then the solutions are $u = u_0 + \nabla\phi$ for any one fixed solution $u_0$ and for any twice continuously differentiable function $\phi$. So the local solutions depend on one function $\phi$ of 3 variables, not on initial data from a surface in $\mathbb{E}^3$.

We might spot trouble on the horizon when we notice that we can't solve our equations for the derivatives with respect to $x_1$, since there are only two of them, and similarly for $x_2$ and for $x_3$. Nonetheless, we can easily make more complicated examples of systems of partial differential equations to which the Cauchy–Kovalevskaya theorem does not apply, with complicated "compatibility conditions" needed to ensure existence of a solution.

Is there a vector field $u$ so that $\nabla \times u = f - u$ for a given vector field $f$ in $\mathbb{E}^3$? Our equation $\nabla \times u = f - u$ expands out to

$$\begin{aligned} u^3_{x_2} - u^2_{x_3} &= f^1 - u^1, \\ u^1_{x_3} - u^3_{x_1} &= f^2 - u^2, \\ u^2_{x_1} - u^1_{x_2} &= f^3 - u^3 : \end{aligned}$$

3 first order equations for 3 unknowns. Again we might expect to find a unique solution with initial values of $u^1, u^2, u^3$ given on some surface in $\mathbb{E}^3_{x_1,x_2,x_3}$. Again this intuition is *not* correct. Taking $\nabla\cdot$ on both sides of our equations reveals $0 = \nabla \cdot f - \nabla \cdot u$, i.e.

$$u^1_{x_1} + u^2_{x_2} + u^3_{x_3} = f^1_{x_1} + f^2_{x_2} + f^3_{x_3},$$

a fourth first-order equation which is linearly independent of the other three.

Split the four equations as:

$$u^3_{x_2} - u^2_{x_3} = f^1 - u^1,$$

and

$$\begin{aligned} u^1_{x_3} - u^3_{x_1} &= f^2 - u^2, \\ u^2_{x_1} - u^1_{x_2} &= f^3 - u^3, \\ u^1_{x_1} + u^2_{x_2} + u^3_{x_3} &= f^1_{x_1} + f^2_{x_2} + f^3_{x_3}. \end{aligned}$$

Solve the first of these equations:

$$u^3_{x_2} - u^2_{x_3} = f^1 - u^1,$$

but just along the surface $x_1 = 0$, using the Cauchy–Kovalevskaya theorem in the variables $x_2, x_3$, starting with any choice of functions $u^1, u^2$ on the surface $x_1 = 0$ and any function $u^3$ on the curve $x_1 = x_2 = 0$. Apply the Cauchy–Kovalevskaya theorem to the last 3 equations:

$$\begin{aligned} u^1_{x_3} - u^3_{x_1} &= f^2 - u^2, \\ u^2_{x_1} - u^1_{x_2} &= f^3 - u^3, \\ u^1_{x_1} + u^2_{x_2} + u^3_{x_3} &= f^1_{x_1} + f^2_{x_2} + f^3_{x_3}. \end{aligned}$$

in all 3 variables $x_1, x_2, x_3$, by treating $x_1$ as the variable we differentiate in:

$$\begin{pmatrix} u^1 \\ u^2 \\ u^3 \end{pmatrix}_{x_1} = \begin{pmatrix} -u^2_{x_2} - u^3_{x_3} + \nabla \cdot f \\ u^1_{x_2} + f^3 - u^3, \\ u^1_{x_3} - f^2 + u^2 \end{pmatrix}.$$

The tricky question is whether the resulting functions actually solve the original problem. In fact, they do: we have solved the first equation

$$u^3_{x_2} - u^2_{x_3} = f^1 - u^1,$$

at $x_1 = 0$. But if we differentiate that expression in $x_1$, check that

$$\left(u^3_{x_2} - u^2_{x_3} - f^1 + u^1\right)_{x_1} = 0$$

modulo the 3 other equations, so vanishes by construction. The general solution to $\nabla \times u = f - u$ is given by picking one arbitrary function of one variable, and 2 arbitrary functions of two variables, again not what we would guess from familiarity with the Cauchy–Kovalevskaya theorem.

Again, we might spot trouble on the horizon when we notice that we can't solve our equations for the derivatives with respect to $x_1$, since there are only two of them, and similarly for $x_2$ and for $x_3$.

## Smooth counterexamples

For functions which are not analytic, there are numerous theorems to prove the existence of solutions, under a wide variety of hypotheses, but there are some surprising counterexamples:

*a.* The Cauchy–Riemann equations admit only analytic solutions: no solution has smooth nonanalytic initial data.

*b.* Some linear determined equations $Lu = f$ with analytic $L$ have no solutions for any smooth $f$ unless $f$ is analytic [25].

Appendix B

# Symbol and characteristic variety

We define the notion of symbol of a system of partial differential equations.

## The problem

The Cauchy–Kovalevskaya theorem says that we can solve a system of differential equations, if we can write it as solving for the highest derivative in a single variable $t$ in terms of all other derivatives, in all variables, and the values of independent and dependent variables. If we can't do this in the variables we are given, we might be able to do it after a change of variables. So we need a coordinate free description of "how many derivatives" are taken in a given direction.

> The Laplace operator in the plane: $u \mapsto u_{xx} + u_{yy}$. At first glance, it appears to take two derivatives along each coordinate axis, and no derivatives in any other directions. But the Laplace operator is invariant under rotation, so it actually "feels" the second derivatives in all directions.

## The symbol of a linear operator

If you differentiate a high frequency wave in a direction in which it oscillates, it gets bigger. If you differentiate in a direction in which the wave is constant, it gets killed. We detect derivatives by plugging in waves.

Take a large number $\lambda$ and any function $f(x)$ vanishing at the origin. Near the origin, the function $e^{i\lambda f}$ looks like a high frequency wave: expand $f$ in a Taylor series $f(x) = \sum \xi_i x_i + \dots$,

$$e^{i\lambda f} = e^{i\lambda \sum \xi_i x_i + \dots}.$$

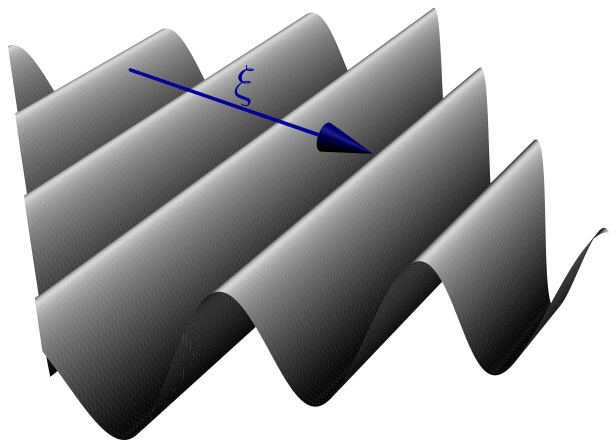

Any linear differential operator,

$$P = \sum_{|a| \le m} c_a(x) \partial^a$$

is a polynomial $P = P(x, \partial)$ in the operators $\partial_{x_i}$; the coefficients are functions of $x$. Write this polynomial

$$P(x, \xi) = \sum_{|a| \le m} c_a(x) \xi^a.$$

Let $P_m$ be the highest order terms, also called the *symbol*:

$$P_m(x, \xi) = \sum_{|a| = m} c_a(x) \xi^a,$$

also denoted $\sigma_P(\xi)$, for $\xi = \xi_i dx^i \in T^*M$ a 1-form. To see "how many derivatives" $P$ takes in some direction, expand:

$$e^{-i\lambda f} P\left[e^{i\lambda f} u\right] = \sum_{|a| \le m} e^{-i\lambda f} c_a(x) \partial^a \left(e^{i\lambda f} u\right),$$

Every time a derivative hits $e^{i\lambda f} u$, it either hits the exponential factor, pulling down a power of $\lambda$, or hits the $u$, so no power of $\lambda$. So the leading order term in $\lambda$ is

$$\begin{aligned} e^{-i\lambda f} P\left[e^{i\lambda f} u\right] &= i^m \lambda^m \sum_{|a| = m} c_a(x) \xi^a u + O(\lambda)^{m-1}, \\ &= i^m \lambda^m P_m(x, \xi) u + O(\lambda)^{m-1}. \end{aligned}$$

where $\xi = df(x)$. Without taking coordinates,

$$\sigma_P(df)\, u = \lim_{\lambda \to \infty} \frac{e^{-i\lambda f} P\left[e^{i\lambda f} u\right]}{i^m \lambda^m}$$

for any continuously differentiable function $f$ vanishing at $x$.

The heat operator in the plane is

$$Pu = u_t - u_{xx} - u_{yy}$$

for a function $u = u(t, x, y)$, and the symbol is $\sigma_P(c\,dt + a\,dx + b\,dy) = -a^2 - b^2$.

All along we could have allowed $u$ to be valued in, say, $\mathbb{R}^q$, and allow the coefficients $c_a(x)$ to be $p \times q$ matrices. Then $\sigma_P(\xi)$ is a $p \times q$ matrix. The matrix entries are functions of $x$ and polynomials of degree $m$ in $\xi$.

**B.1** Find the symbol of the Cauchy–Riemann operator

$$P \begin{bmatrix} u \\ v \end{bmatrix} = \begin{pmatrix} u_x - v_y \\ u_y + v_x \end{pmatrix}.$$

**B.2** If we are interested in a linear system of *equations* $0 = Pu$, rather than an operator $P$, the equations have the same solutions if we rescale both sides by some invertible matrix $g(x)$. Similarly, we could rescale the $u$ variable by a matrix $h(x)$. Check that

$$\sigma_{gPh}(\xi) = g\sigma_P(\xi)\,h.$$

If we change variables by a diffeomorphism $y = h(x)$, then $\partial_{x_i} = \sum_j \frac{\partial y_j}{\partial x_i} \partial_{y_j}$ and for any 1-form $\xi = \sum_j \xi_j dy_j = \sum_{ij} \xi_j \frac{\partial y_j}{\partial x_i} dx_i$, so when we substitute $\xi_i$ for $\partial_{y_i}$,

$$\sigma_{h_* P}(\xi) = \sigma_P(h^* \xi).$$

We knew this already: the symbol $\sigma_P(\xi)$ is defined independently of coordinates.

Taking exterior derivative $d$ as our differential operator,

$$e^{-i\lambda f} d(e^{i\lambda f}) = i\lambda df \wedge \omega + d\omega,$$

so the symbol is

$$\sigma_d(\xi) = \xi\wedge,$$

i.e. the symbol is the linear transformation $\omega \mapsto \xi \wedge \omega$, where $\xi$ is a 1-form.

## Background material from projective geometry

The *projective space* $\mathbb{P}V$ of a vector space $V$, over any field, is the set of lines through the origin of $V$ [29]. If $V$ has dimension $n+1$, we may write $\mathbb{P}V$ as $\mathbb{P}^n$. For any $v \in V$ with $v \neq 0$, $[v] \in \mathbb{P}V$ is the line through $v$. Any nonzero linear function $f$ on $V$ vanishes on a hyperplane in $V$. This hyperplane determines $f$ up to rescaling, and vice versa. So each point $[f]$ of $\mathbb{P}V^*$ is naturally identified with a hyperplane in $V$.

## The characteristic variety

The *characteristic variety* of a linear differential operator $P$ is

$$\Xi_{P,m} = \{\, [\xi] \in \mathbb{P}T^*_m M \mid \ker \sigma_P(\xi) \neq 0 \,\}\,.$$

Let $T^{\mathbb{C}*}_m M := T^*_m M \otimes \mathbb{C}$. The *complex characteristic variety* consists of complex lines in a complex vector space:

$$\Xi^{\mathbb{C}}_{P,m} = \left\{\, [\xi] \in \mathbb{P}T^{\mathbb{C}*}_m M \;\middle|\; \ker \sigma_P(\xi) \neq 0 \,\right\}.$$

A 1-form $\xi$ belongs to the characteristic variety just when the operator $P$ takes "exceptionally few" derivatives in the direction of $\xi$.

For the heat operator, the characteristic variety is the set of $[c, a, b] \in \mathbb{P}^2$ so that $a^2 + b^2 = 0$, i.e. the single point $[c, a, b] = [1, 0, 0] \in \mathbb{P}^2$, while complex characteristic variety consists of the complex number solutions to the same equations: $[c, a, \pm a]$, i.e. $b = ia$ and $b = -ia$, a union of two complex lines.

Take exterior derivative $d$ on $k$-forms as our differential operator. The symbol is $\sigma_d(\xi) = \xi\wedge$, so the characteristic variety is the set of all lines $[\xi]$ spanned by 1-forms $\xi \neq 0$ so that $\xi \wedge \omega = 0$ for some $\omega \neq 0$. If we work with only 0-forms $\omega$, then this forces $\xi = 0$: the characteristic variety is empty. If $k > 0$, then any $\xi$ has $\xi \wedge \omega = 0$ for some $\omega \neq 0$: take any $\omega$ of the form $\xi \wedge \eta$. So the characteristic variety of the exterior derivative on $k$-forms is

$$\Xi_{d,m} = \begin{cases} \varnothing, & \text{if } k = 0, \\ \mathbb{P}T^*_m M, & \text{if } k > 0. \end{cases}$$

The same calculation computes the complex characteristic variety:

$$\Xi^{\mathbb{C}}_{d,m} = \begin{cases} \varnothing, & \text{if } k = 0, \\ \mathbb{P}T^{\mathbb{C}*}_m M, & \text{if } k > 0. \end{cases}$$

For the heat operator, the characteristic variety at each point consists of the single hyperplane $dt = 0$. This hyperplane is the tangent plane to each surface representing space at constant time.

For a vector field $X$, thought of as a differential operator $Pu = Xu$, the characteristic variety is the set of $[\xi]$ so that $\xi(X) = 0$, i.e. the set of all hyperplanes containing the vector $X$.

For the wave operator, $Pu = u_{tt} - u_{xx} - u_{yy}$, the symbol is $\sigma_P(t, x, y, c, a, b) = c^2 - a^2 - b^2$. Up to rescaling, we can arrange $c = 1$, so the characteristic variety is a circle $[1, a, b] \in \mathbb{P}^2$ so that $a^2 + b^2 = 1$. As a family of hyperplanes, this circle is the set of hyperplanes tangent to the "light cone" $dt^2 = dx^2 + dy^2$.

Points of the complex characteristic variety are complex hyperplanes in the spaces of complexified tangent vectors.

**B.3** Find the characteristic variety of Maxwell's equations in the vacuum.

## Linearization

Take a nonlinear differential equation, say $u_t = u_{xx}u + \sin(u_x)$, and take a solution $u(t, x)$. Take a function $v(t, x)$ and plug in $u + \varepsilon v$ instead of $u$ to both sides, supposing that $u + \varepsilon v$ is also a solution, at least up to an error smaller than $\varepsilon$:

$$
\begin{aligned}
u_t + \varepsilon v_t &= (u_{xx} + \varepsilon v_{xx})(u + \varepsilon v) + \sin(u_x + \varepsilon v_x), \\
&= u u_{xx} + \sin(u_x) + \varepsilon\left(u_{xx}v + v_{xx}u + \cos(u_x)\, v_x\right) + O(\varepsilon)^2
\end{aligned}
$$

using the Taylor series of sin. Since $u$ satisfies the differential equation, we can cancel off $u_t = u u_{xx} + \sin(u_x)$ from both sides, and divide by $\varepsilon$ and send $\varepsilon \to 0$:

$$
v_t = u_{xx}v + v_{xx}u + \cos(u_x)\, v_x,
$$

a linear differential equation in $v$, but with coefficients which depend on the solution $u$ that we perturbed about.

More generally, for any expression $P[u] = F(x, u, u_x)$ (a differential operator, perhaps nonlinear), and any function $u$ for which $P[u] = 0$, the *linearization* of $P$ about $u$ is the linear differential operator

$$
P'[u]v = F_u(x, u, u_x)\, v + F_{u_x}(x, u, u_x)\, v_x,
$$

and similarly for higher order operators.

**B.4** Let $\Delta = \partial_{xx} + \partial_{yy}$. Linearize the equation $\Delta u = |du|^2$ around the solution $u = 1$.

**B.5** Let $\Delta = \partial_{xx} + \partial_{yy}$. Linearize the Liouville equation $\Delta \log u = -u^2$ around the solution $u = \frac{1}{2}\log\left(\frac{4}{(1-x^2-y^2)^2}\right)$.

When we take about a *differential equation* $0 = F(x, u, u_x)$, defined for $x \in \mathbb{R}^n, u \in \mathbb{R}^k, u_x \in \mathbb{R}^{nk}$ lying in some open set $(x, u, u_x) \in U \subset \mathbb{R}^{n+k+nk}$, we implicitly assume that $F$ is analytic and that the set $M$ of all $(x, u, p) \in U \subset \mathbb{R}^{n+k+nk}$ where $F(x, u, p) = 0$ is a manifold and that the map $(x, u, p) \in M \mapsto (x, u)$ is a submersion (so the differential equation doesn't, locally, constraint the value of the independent or dependent variables).

If we are given the values of $u$ and $u_x$ at just one point $x$, then we can calculate $F(x, u, u_x)$ at that one point, and if $0 = F(x, u, u_x)$ then we can also calculate the linearization $P'[u]$ at $x$, or in other words we can calculate the coefficients of the linearization at $x$:

$$F_u(x, u, u_x) \text{ and } F_{u_x}(x, u, u_x) .$$

Therefore we can calculate the *characteristic variety* $\Xi_{x,u,u_x}$ of our differential equation, by which we mean that of the linearization at $(x, u, u_x)$.

The linearization of $u_t = uu_{xx} + u_{yy} + u^2$ is $v_t = vu_{xx} + uv_{xx} + v_{yy} + 2uv$, so the characteristic variety is $0 = 0 + u\xi_x^2 + \xi_y^2 + 0$, if we write the components of a 1-form $\xi$ as $\xi = \xi_t \, dt + \xi_x \, dx + \xi_y \, dy$.

## Initial value problems

A differential equation is *determined* just when the symbol matrix $\sigma(\xi)$ is square, and invertible for some $\xi$.

**Lemma B.1.** *A differential equation $0 = F(x, u, u_x)$ can locally be written in the form of the Cauchy–Kovalevskaya theorem in some coordinates just when it is determined.*

*Proof.* If we linearize the equation around some point $(x, u, p)$, we get the linear equation

$$0 = F_u(x, u, u_x)v + F_p(x, u, u_x)v_x,$$

with symbol $\sigma(\xi) \, v = F_p(x, u, u_x)\xi v$. So if we write coordinates as $x^i, u^a, p_i^a$, then the symbol matrix is

$$\sigma(\xi) = F_{p_i}\xi_i.$$

We want to see when we can somehow smoothly solve, at least locally, the equation $0 = F(x, u, u_x)$ for some $u_{x_i}$ as a function of $x, u$ and the other $u_{x_j}$. By the implicit function theorem, we can do this just when $F_{p_i}$ is an invertible square matrix. Note that

$$\sigma(dx_i) = F_{p_i}.$$

On the other hand, if $\sigma(\xi)$ is an invertible square matrix for some $\xi$, we can linearly change variables to arrange that $\xi = dx_i$ and reverse our steps. □

The problem of specifying initial data for nonlinear equations is much more complicated than for linear equations. Take a hypersurface $H \subset \mathbb{R}^n$ of the $x$-variables and maps $u(x)$ and $p(x)$ defined along $H$ so that $p(x)v = u'(x)v$ for any tangent vector $v$ to $H$. The maps $u(x), p(x)$ are *noncharacteristic initial data* if $(x, u(x), p(x)) \in M$ for each $x \in H$ and $T_x H$ is noncharacteristic for the the linearization of our equation at each point $(x, u(x), p(x))$.

**Lemma B.2.** *A determined system of equations has some analytic noncharacteristic initial data through any point $(x, u, p)$.*

*Proof.* Since the symbol matrix is somewhere invertible, the characteristic variety $\Xi_m$ is not all of $\mathbb{P}T^*_m M$, where $m = (x, u, p)$. □

**Theorem B.3.** *Take a determined system of differential equations and noncharacteristic initial data along an embedded hypersurface, all analytic. There is an analytic solution to the differential equations agreeing with the initial data along the hypersurface. Any two such agree near the hypersurface and agree in any connected open set where both are defined.*

*Proof.* If we can prove uniqueness of local solutions, then we can glue them together to get a unique global solution. So it suffices to prove the result locally. By an analytic change of variables, arrange that our hypersurface is $t = 0$ in some $(t, x)$ variables. By the implicit function theorem, if we have a smooth map $F(t, x, u, p, q)$, we can locally smoothly solve $F(t, x, u, p, q) = 0$ for a variable $p$ as a function of the other variables just when $F_p$ is a matrix of full rank. We can rewrite the equation $F(t, x, u, u_t, u_x) = 0$ in the Cauchy–Kovaleskaya form, i.e. solve for $u_t$ as a function of the other variables, just when $F_{u_t}$ is a matrix of full rank, i.e. when $t = 0$ is not characteristic for the linearization of the differential equation. □

**Corollary B.4.** *A determined analytic system of differential equations has a local analytic solution near any point.*

Take the equation $u_t^2 + u_x^2 = 1$ with initial data $u = 0$ at $t = 0$. Differentiating the requirement that $u = 0$, we require $u_x = 0$ too at $t = 0$, as part of our initial data. The differential equation then says that $u_t = \pm 1$ at $t = 0$. The linearization at any point is

$$2u_t v_t + 2u_x v_x = 0.$$

Along our initial data, this is

$$2v_t = 0.$$

The linearization is determined, so the nonlinear equation is, and our theorem guarantees a solution near $t = 0$. Easier: we could have just solved

$$u_t = \pm\sqrt{1 - u_x^2}$$

without using the theorem above, and then any initial data $u(0, x)$ with $-1 < u_x < 1$ all along $t = 0$ would give rise to a unique solution $u(t, x)$ near $t = 0$.

## Higher order differential equations

Consider a linear differential operator $Pu = u_{xx}$; its symbol is $\sigma_P(\xi) = -\xi^2$. Take a first order operator with the same kernel:

$$Q\begin{pmatrix} u \\ v \end{pmatrix} = \begin{pmatrix} v_x \\ u_x \end{pmatrix} = \begin{pmatrix} 0 & 1 \\ 1 & 0 \end{pmatrix} \partial_x \begin{pmatrix} u \\ v \end{pmatrix}.$$

The symbol of $Q$ is

$$\sigma_Q(\xi) = \begin{pmatrix} 0 & \xi \\ \xi & 0 \end{pmatrix}.$$

The characteristic variety is $\Xi_Q = \left(\xi^2 = 0\right) = \Xi_P$.

If there are many functions $u$ of many variables $x$, this $\sigma_P$ is a matrix-valued quadratic form, and the same steps tell us that $\Xi_Q = \Xi_P$ and $\Xi_Q(\mathbb{C}) = \Xi_P(\mathbb{C})$: replacing a system of differential equations by its equivalent first order system doesn't change the characteristic variety.

The *minimal surface equation*

$$\left(1 + u_x^2\right) u_{yy} - 2u_x u_y u_{xy} + \left(1 + u_y^2\right) u_{xx} = 0,$$

representing surfaces which have locally least area with given boundary, has empty characteristic variety at every point. If we pick any embedded real curve in the plane, and along that curve pick functions $f, p, q$, there is a unique solution near that curve so that $u = f, u_x = p, u_y = q$ along the curve.

**B.6** Find the characteristic variety of the *Euler–Tricomi equation* $u_{tt} = tu_{xx}$.

**B.7** If $P = \partial_{xx} + \partial_{yy}$, then the associated equation $Pu = 0$ is the equation of an electrostatic potential energy $u$: the *Laplace equation*, and $\Delta := P(D)$ is the *Laplace operator*. Show that the characteristic variety is cut out by the equation $0 = \xi^2 + \eta^2$, which has $(\xi, \eta) = (0, 0)$ as its only solution.

**B.8** Calculate the linearization, symbol and complex characteristic variety of the minimal surface equation around a solution $u(x, y)$.

**B.9** Suppose that $0 = F(x, u, u_x)$ is an analytic system of partial differential equations, and the symbol of the linearization about any point $(x, u, u_x)$ has a covector $\xi$ for which $\sigma(\xi)$ is a surjective linear map, with kernel of rank $k$ say. Prove that we can locally add $k$ additional equations to make the system determined, and conclude that there are local solutions. What is the correct notion of noncharacteristic initial data?

## Counterexamples

Linear analytic scalar equations which are not determined fail to have either existence or uniqueness of solutions, with analytic initial data [23]. For functions which are not analytic, uniqueness can fail even for determined equations:

*a.* Some smooth linear determined equations have smooth solutions which agree on one side of a noncharacteristic hypersurface and disagree on the other [2].

*b.* Some analytic nonlinear determined equations have smooth solutions which agree on one side of a noncharacteristic hypersurface and disagree on the other [28].

There are some results ensuring analytic local solutions of analytic differential equations with initial data characteristic at some points [11].

Appendix C

# Analytic convergence

## Convergence

A sequence $f_1, f_2, \dots$ of analytic functions of real variables, defined on some open set, *converges* to an analytic function $f$ just when that set lies in an open set in complex variables to which $f$ and all but finitely many $f_i$ extend and where $f_i$ converge uniformly to $f$ on compact sets. By the Cauchy integral theorem [1] p. 120, convergence is uniform with all derivatives on compact sets.

**C.1** Prove: a sequence $f_1, f_2, \dots$ of analytic functions converges to an analytic function $f$ just when we can cover their domain in compact sets, and cover each compact set by an open set in complex variables, to which $f$ and all but finitely many $f_i$ extend and are bounded, and on which $f_i$ converge uniformly to $f$.

## Commutative algebra

A *graded ring* is a ring $R$ which is the direct sum of abelian groups $R_0, R_1, \dots$ so that $R_i R_j \subseteq R_{i+j}$, *elementary* if generated by $R_0$ and $R_1$. An $R$-module $M$ is *graded* if a direct sum of abelian groups $M_i$ with $R_i M_j \subseteq M_{i+j}$; *asymptotically simple* if $R_1 M_k = M_{k+1}$ for all sufficiently large $k$; *elementary* if each $M_k$ is finitely generated as an $R_0$-module.

**Lemma C.1** (Phat Nguyen [30])**.** *An elementary graded module over an elementary graded ring is finitely generated if and only if asymptotically simple.*

*Proof.* Suppose $M$ is finitely generated over $R$, hence by $M_0 + M_1 + \dots + M_k$. For any $j \ge k$,

$$M_{j+1} = R_{j+1} M_0 + R_j M_1 + \dots + R_1 M_j.$$

Since $R$ is elementary

$$M_{j+1} = R_1^{j+1} M_0 + R_1^j M_1 + \dots + R_1 M_j \subseteq R_1 M_j.$$

But $R_1 M_j \subseteq M_{j+1}$, so $M_{j+1} = R_1 M_j$ for all $j \ge k$. Therefore $M$ is asymptotically simple.

Conversely, assume $M$ asymptotically simple: $R_1 M_j = M_{j+1}$ for any $j \ge k$. As an $R$-module, $M$ is generated by $M_0 + M_1 + \dots + M_k$. Because $M$ is elementary, we can select a finite set of generators from each of $M_0, \dots, M_k$ as an $R_0$-module, hence as an $R$-module. □

A module $M$ over a ring $R$ is *Noetherian* if every $R$-submodule is finitely generated; $R$ is *Noetherian* if every ideal is finitely generated.

**Theorem C.2** (Krull Intersection Theorem [30]). *Suppose that $R$ is a Noetherian commutative ring with identity, $I \subseteq R$ an ideal, $M$ a finitely generated $R$-module, and*

$$N := \bigcap_{j=0}^{\infty} I^j M.$$

*Then there is an element $i \in I$ so that $(1+i)N = 0$. If all elements of $1+I$ are units of $R$, then $N = 0$.*

*Proof.* Write $R, I, M, N$ as $R_0, I_0, M_0, N_0$. Define elementary graded ring and modules:

$$\begin{aligned} R &:= R_0 \oplus I_0 \oplus I_0^2 \oplus \dots, \\ M &:= M_0 \oplus I_0 M_0 \oplus I_0^2 M_0 \oplus \dots, \\ N &:= N_0 \oplus N_0 \oplus N_0 \oplus \dots. \end{aligned}$$

Since $M$ is asymptotically simple, $M$ is finitely generated over $R$. Since $R_0$ is Noetherian, $R$ is Noetherian, so $M$ is Noetherian as an $R$-module. Because $N_0 \subset I_0^j M_0$ for $j = 1, 2, \dots$, $N$ is a graded $R$-submodule of $M$, so finitely generated over $R$. By lemma C.1 on the previous page, $N$ is asymptotically simple, so $I_0 N_0 = N_0$.

Take generators $n_a \in N_0$. Then $n_a = \sum_b j_{ab} n_b$ for some $j_{ab} \in I_0$. Let $J = (j_{ab})$, a matrix, and $n = (n_a)$ a vector. So $(I - J)n = 0$. Multiply by the adjugate matrix (the transpose of the matrix of cofactors) of $I - J$ to get $\det(I - J)n = 0$; set $1 + i = \det(I - J)$. □

## Germs

The *germ* of an analytic function at a point is its equivalence class, two functions being equivalent if they agree near that point. A sequence of germs *converges* just when it arises from a convergent sequence of analytic functions. By majorizing, this is just when the Taylor coefficients converge to those of an analytic function. Germs form a ring, under addition and multiplication of functions. Pick a finite dimensional real vector space. Germs of analytic maps to that vector space form a module over that ring.

**Theorem C.3** (Hilbert basis theorem [8] p. 81 theorem 2.7, [18] p. 161 lemma 6.3.2, theorem 6.3.3, [34] p. 44 theorem 3.4.1). *Every submodule of that module is finitely generated.*

**Theorem C.4** (Henri Cartan [5] p. 194 corollaire premiere, [12] p. 46, [15] p. 85 theorem 3, [34] p. 274 theorem 11.2.6). *Every submodule of that module is closed under convergence.*

*Proof.* Pick a submodule $S$. Take a convergent sequence $f_i \to f$ of map germs, with $f_i \in S$. We need to prove that $f \in S$. Let $R$ be the the ring of analytic function germs at the origin, $I \subset R$ the germs vanishing at the origin, $T$ the $R$-module of analytic map germs, which we identify with their Taylor series. The $k$-th order Taylor series are the elements of $T/I^{k+1}T$, and form a vector space of finite dimension, containing a linear subspace $S/I^{k+1}S$, which is therefore a closed subset, so contains the $k$-th order Taylor coefficient of $f$. Let $M := T/S$, so $f \in I^k M$ for all $k = 1, 2, \dots$. The set $1 + I$ is the set of germs of functions equal to 1 at the origin, so units in $R$. By Krull's intersection theorem, $f = 0$ in $M$, i.e. $f \in S$. □

Appendix D

# The moving frame

In most chapters of these notes, we study geometry that is invariant under diffeomorphisms. In this chapter, we study geometry that is invariant under rigid motions of 3-dimensional Euclidean space. We want examples of exterior differential systems. For simplicity, we look to surface theory: differential equations governing surfaces in 3-dimensional Euclidean space. We look for differential equations invariant under rigid motion. But there are no differential forms on Euclidean space invariant under rigid motion. We will instead find invariant differential forms on the *frame bundle* of Euclidean space.

## Notation

*Warning*: We use the notation $\mathbb{E}^3$ instead of $\mathbb{R}^3$ to mean Euclidean space, i.e. $\mathbb{R}^3$ with its usual metric.

*Warning*: In this chapter, $e_1, e_2, e_3$ denotes an *arbitrary* orthonormal basis of $\mathbb{E}^3$, *not* necessarily the standard basis.

## Orthonormal frames

An *orthonormal frame* is a pair $(x, e)$ where $x \in \mathbb{E}^3$ and

$$e = (e_1\ e_2\ e_3)$$

is an orthonormal basis of $\mathbb{E}^3$ arranged into the columns of a matrix.

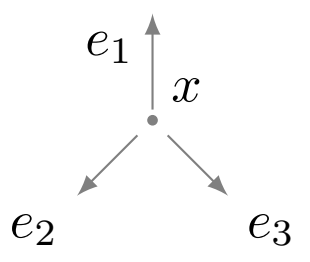

**D.1** For any two frames $(x, e), (x', e')$, there is a unique rigid motion $\varphi$, i.e. a distance preserving transformation, taking one to the other: $x' = \varphi(x)$, $e'_i = \varphi_* e_i$.

Each rigid motion $\varphi$ has associated frame $(x, e)$ which is the frame that arises when we apply $\varphi$ to the standard basis frame at the origin. Let O(3) be the set of $3 \times 3$ orthogonal matrices. To each orthonormal frame $(x, e) \in \ulcorner\mathbb{E}^3$, associate the matrix

$$h = \begin{pmatrix} 1 & 0 \\ x & e \end{pmatrix}$$

giving a map

$$h \colon (x, e) \in \ulcorner\mathbb{E}^3 \mapsto h(x, e) \in \mathbb{E}^3 \ltimes \mathrm{O}(3)\,.$$

**D.2** Prove that: associating to $\varphi$ this matrix $h(x, e)$ above is an isomorphism of groups from the group of rigid motions of Euclidean space to the group $\mathrm{O}(3) \ltimes \mathbb{E}^3$ of such matrices. So $h(x, e)$ is a matrix description of the rigid motion that takes the origin to $x$ and the standard basis to $e$. Prove that

$$h(x, eg) = h(x, e) \begin{pmatrix} 1 & 0 \\ 0 & g \end{pmatrix}$$

for any orthogonal $3 \times 3$ matrix $g$.

## Frames on curves

Take a curve in $\mathbb{E}^3$

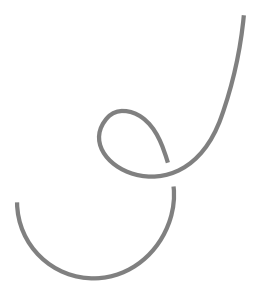

and consider all of the frames $(x, e)$ with $x$ a point of the curve.

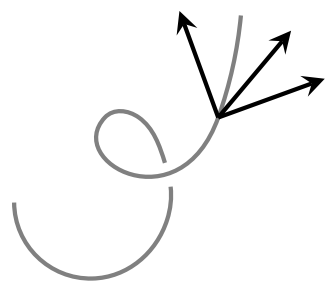

At each point $x$ of the curve, any two such frames agree up to orthogonal transformation; the set of all such frames is a 4-dimensional manifold. If a rigid motion takes one curve to another, it takes these frames along the one curve to those along the other curve. If the curve is immersed, a frame $(x, e)$ is *adapted* if $e_1$ is tangent to the curve.

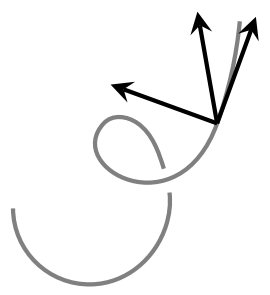

At each point $x$ of the curve, any two adapted frames agree up to changing the direction of $e_1$, and rotating or reflecting $e_2, e_3$ in the plane normal to the tangent line; the set of all such frames is a surface inside that 4-dimensional manifold. This surface of adapted frames is built by attaching to each point of $C$ four circles, representing rotating $e_2, e_3$ and reflecting $e_1$ or the plane or $e_2, e_3$. We will carry out computations directly inside that surface. On a connected curve, that surface has 4 components. Every rigid motion is determined by where it moves any chosen frame. Any rigid motion preserving an immersed curve takes any frame adapted to the curve to another frame adapted to the curve. So we can guess that the group of rigid motions preserving an immersed curve is at most a 2-dimensional submanifold of the 6-dimensional group of rigid motions.

## Frames on surfaces

Similarly, a frame is adapted to a surface if $x$ is a point of the surface and $e_3$ is normal to the surface.

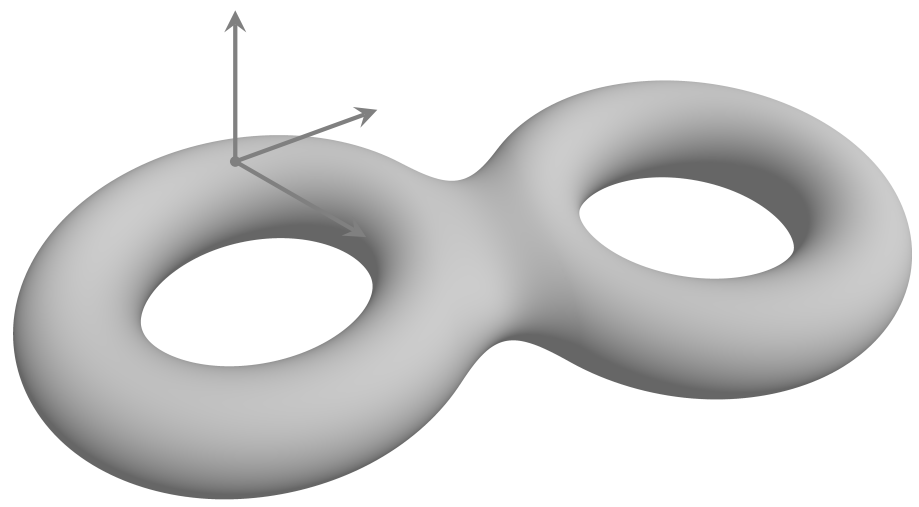

The set of adapted frames $(x, e)$ of a surface is a 3-dimensional manifold, as we can move the point $x$ along the surface, and rotate $e_1, e_2$ around $e_3$. Hence the group of rigid motions preserving an immersed surface is at most a 3-dimensional submanifold of the 6-dimensional group of rigid motions.

## The frame bundle

The *frame bundle* ${}^{\ulcorner}\mathbb{E}^3$ of Euclidean space $\mathbb{E}^3$ is the set of all orthonormal frames. Any rigid motion $\varphi$ of $\mathbb{E}^3$ takes $e_1, e_2, e_3$ to $\varphi_* e_1, \varphi_* e_2, \varphi_* e_3$. If you pick an orthonormal frame $(x, e)$ and I pick another one, say $(x, e')$ at the same point $x$, then yours and mine must agree up to a unique orthogonal matrix, because they are orthonormal bases of the same vector space $\mathbb{E}^3$. So $e'_j = \sum_i g_{ij} e_i$ for a uniquely determined orthogonal $3 \times 3$ matrix $g = (g_{ij})$; denote this $e'$ as $eg$. An orthonormal frame $(x, e)$ with $e = (e_1 \; e_2 \; e_3)$ has $e \in \mathrm{O}(3)$ an orthogonal matrix, so ${}^{\ulcorner}\mathbb{E}^3 = \mathbb{E}^3 \times \mathrm{O}(3)$ and the expression $eg$ is matrix multiplication.

## Adapted frames for a curve

Take a curve $C$ in $\mathbb{E}^3$. At each point $c \in C$, take $e_1$ to be a unit vector tangent to $C$. We can pick an adapted frame $(x, e)$, uniquely up to replacing $(x, e)$ by $(x, eg)$ for any orthogonal matrix

$$g = \begin{pmatrix} \pm 1 & 0 \\ 0 & h \end{pmatrix}, h \in \mathrm{O}(2)\,.$$

The *curvature vector* of $C$ at a point is the acceleration as we move along $C$ at unit speed. If we change the choice of unit speed parameterization (also known as *arclength parameterization*, adding a constant or changing the direction of motion, i.e. the sign, this changes the sign of velocity, but we hit two sign changes when we take two derivatives, so the curvature vector is independent of the choice of unit speed parameterization. In particular, the curvature vector is defined without selecting any orientation of the curve $C$.

Write the (perhaps local) inverse of an arclength parameterization as a function $s$ on an open subset of $C$, uniquely determined up to sign and adding a constant. The identity function $x \mapsto x$ on $\mathbb{E}^3$ restricts to a function $x \colon C \to \mathbb{E}^3$, and clearly

$$e_1 = \frac{dx}{ds},$$

when we take $e_1$ to point in the direction of increasing $s$, i.e. use $s$ to orient $C$. Differentiating again, the curvature vector of $C$ is

$$\kappa := \frac{de_1}{ds}.$$

If $\kappa = 0$ at every point of the curve then the curve is a straight line. If $\kappa \neq 0$ at some point, then we can uniquely determine a unit normal vector $e_2$ to $C$ by

$$e_2 := \frac{\kappa}{|\kappa|},$$

and then define the *curvature* to be $k := |\kappa|$. Having chosen $e_1, e_2$ at each point of our curve, $e_3$ is determined up to sign, and can be chosen smoothly.

**D.3** Calculate that

$$\frac{de_2}{ds} = -ke_1 + te_3, \quad \frac{de_3}{dt} = -te_2,$$

for some function $t$, the *torsion*.

The sign of $t$ depends on the choice of $e_3$, which is unique up to sign, so $te_3$ is defined independent of choice of $e_3$. But $te_3$ still depends on the choice of sign of $e_1$, i.e. the direction of motion. So we could say that the complicated expression $te_3\, ds$, a 1-form valued in the normal line to the curve, is an invariant independent of any choices.

The frame $e_1, e_2, e_3$ of an oriented curve of nonzero curvature is the *Serret–Frenet frame*. Note that we can change sign of $e_1$, i.e. change orientation, and independently change sign of $e_3$. So the frame bundle of $C$ contains four Serret–Frenet frames above each point of $C$, i.e. four curves in the frame bundle. The torsion changes sign if we change the sign of $e_3$, so it is not really a function on the curve $C$. We can pick out one choice of Serret–Frenet curve by orienting both $C$ and $\mathbb{E}^3$. We can guess that the isometries of $\mathbb{E}^3$ which preserve a nonlinear curve are at most 1-dimensional.

**D.4** For a curve with nonzero curvature vector, prove that its torsion vanishes if and only if the curve lies in a plane.

Intuitively, torsion controls the rate at which the curve twists out of its tangent plane, while the curvature controls the rate at which the curve twists in its tangent plane.

With only rigid motions, we can turn a circle around, slide a line along itself, and twist a helix along itself, and allow reflections in direction. The adapted frames of these curves are taken one to another by these rigid motions. So the curvature and torsion functions are constant, except that the torsion changes sign if we change the choice of $e_3$, or of direction of motion along the curve.

**D.5** Compute the curvature and torsion of a helix. Prove that any constant values of curvature and torsion can arise for a suitable choice of helix.

**Theorem D.1.** *Take an interval $I \subset \mathbb{R}$. Given a continuous positive function $k\colon I \to (0, \infty)$ and a continuous function $t\colon I \to \mathbb{R}$, there is a curve $C$ in $\mathbb{E}^3$ with twice continuously differentiable unit speed parameterization $x\colon I \to \mathbb{E}^3$ and with curvature $k$ and torsion $t$. Any two such curves are identified by a unique rigid motion of $\mathbb{E}^3$.*

*Proof.* To each adapted frame $(x\ e)$, associate the matrix

$$g = \begin{pmatrix} 1 & 0 \\ x & e \end{pmatrix}.$$

This $g$ maps each adapted frame $(x\ e)$ of the curve to the rigid motion $g(x\ e)$ which takes the origin to $x$ and the standard basis to $e$. For the adapted frames of a curve, the above gives

$$\frac{dg}{ds} = g \begin{pmatrix} 0 & 0 & 0 & 0 \\ 1 & 0 & -k & 0 \\ 0 & k & 0 & -t \\ 0 & 0 & t & 0 \end{pmatrix}.$$

These equations are linear, so there are global continuously differentiable solutions, by the existence and uniqueness theorem for linear ordinary differential equations, with any initial data. Expand out to

$$\frac{dx}{ds} = e_1,$$

$$\frac{d}{ds}\begin{pmatrix} e_1 & e_2 & e_3 \end{pmatrix} = \begin{pmatrix} e_1 & e_2 & e_3 \end{pmatrix}\begin{pmatrix} 0 & -k & 0 \\ k & 0 & -t \\ 0 & t & 0 \end{pmatrix}$$

If $e_1, e_2, e_3$ are orthonormal at one point $s_0 \in I$, then

$$\frac{d}{ds} e_1 \cdot e_1 = 2e_1 \cdot ke_2 = 0,$$

etc., so remains an orthonormal frame. The choice of initial conditions at a single point determines the curve and the frame. If we translate $x$ or rotate $x, e_1, e_2, e_3$ by a constant rotation, check that this takes us to another solution. Since $e_1(s)$ is continuously differentiable and $dx/ds = e_1$, $x$ is twice continuously differentiable. □

Any curve of constant positive curvature and constant nonzero torsion is a helix, by existence (from the problem above) and uniqueness. Similarly, any curve of constant positive curvature and zero torsion is a circle, and any curve of zero curvature is a line.

**D.6** For a given curve in 3-dimensional Euclidean space, with nowhere vanishing curvature vector, find every other curve so that, parameterized by arc length, points with the same arclength parameter have the same direction of curvature vector.

## Adapted frames for a surface

An *adapted orthonormal frame* (or just *frame* for short) on a surface $S$ in $\mathbb{E}^3$ is an orthonormal frame $(x, e) \in {}^{\ulcorner}\mathbb{E}^3$ so that $x \in S$ and $e_1, e_2$ are tangent to $S$ at $x$, and therefore $e_3$ is normal to $S$.

If $S$ is the unit sphere in $\mathbb{E}^3$, we can write the points of $S$ as

$$x = \begin{pmatrix} x_1 \\ x_2 \\ x_3 \end{pmatrix}$$

with $1 = x_1^2 + x_2^2 + x_3^2$. Away from the north or south pole, we can take

$$e_1 = \frac{1}{\sqrt{x_1^2 + x_3^2}} \begin{pmatrix} x_3 \\ 0 \\ -x_1 \end{pmatrix}, e_2 = \frac{1}{\sqrt{x_1^2 + x_3^2}} \begin{pmatrix} x_1 x_2 \\ -\left(x_1^2 + x_3^2\right) \\ x_2 x_3 \end{pmatrix}$$

tangent and

$$e_3 = \begin{pmatrix} x_1 \\ x_2 \\ x_3 \end{pmatrix}$$

normal. We can choose $e_3$ to point outward from $S$ in either direction, and choose $e_1, e_2$ to have either orientation.

The *frame bundle* $\ulcorner S$ is the set of all adapted orthonormal frames for $S$. The manifold $\ulcorner S$ is a 3-dimensional submanifold of $\ulcorner\mathbb{E}^3 = \mathbb{E}^3 \times \mathrm{O}(3)$.

**D.7** Explain why, if $S = \mathbb{S}^2 \subset \mathbb{E}^3$ is the unit sphere around the origin, then $\ulcorner S$ is diffeomorphic to two disjoint copies of $\mathrm{O}(3)$, given by the equation $x = \pm e_3$.

Each orthogonal $3 \times 3$ matrix $g$ preserving the horizontal plane

$$g = \begin{pmatrix} h & 0 \\ 0 & n \end{pmatrix}$$

with $h \in \mathrm{O}(2)$ and $n \in \mathrm{O}(1) = \{\pm 1\}$, acts on $\ulcorner S$ as $(x, e)g = (x, eg)$.

## The soldering forms

If $v$ is a tangent vector on the manifold $\ulcorner\mathbb{E}^3$, we can write $v$ as $(\dot{x}, \dot{e})$, an infinitesimal motion $\dot{x}$ of the point $x$, and an infinitesimal rotation $\dot{e}$ of the frame $e$. We can write out any tangent vector $\dot{x}$ in terms of the basis $e_1, e_2, e_3$, say as $\dot{x} = a_1 e_1 + a_2 e_2 + a_3 e_3$. The *soldering forms* are the 1-forms $\omega_1, \omega_2, \omega_3$ on $\ulcorner\mathbb{E}^3$ given by $(\dot{x}, \dot{e}) \lrcorner \omega_i = a_i$. So the soldering forms measure, as we move a frame, how the base point of the frame moves, as measured in the frame itself as a basis.

The *identity function* on $\mathbb{E}^3$, which we write as $x$, is defined as $x(p) = p$ for any point $p \in \mathbb{E}^3$. Of course $dx(v) = v$ for any vector $v$, i.e. $dx = I$ is the identity matrix. We define $x$ also on $\ulcorner\mathbb{E}^3$ by $x(p, e) = p$. We can write the 1-forms $\omega_1, \omega_2, \omega_3$ on $\ulcorner\mathbb{E}^3$ as $\omega_1 = e_1 \cdot dx, \omega_2 = e_2 \cdot dx, \omega_3 = e_3 \cdot dx$. On our vector $v = (\dot{x}, \dot{e})$, we have $v \lrcorner \omega_i = e_i \cdot \dot{x} = a_i$.

## The connection forms

When we move a frame $(x, e)$, the soldering forms measure the motion of the underlying point $x$. We want to measure the rotation of the vectors $e_1, e_2, e_3$.

Infinitesimal rotations are complicated. Write the inner product on $\mathbb{E}^3$ as $x \cdot y = x_1 y_1 + x_2 y_2 + x_3 y_3$. So

$$e_i \cdot e_j = \begin{cases} 1 & \text{if } i = j, \\ 0 & \text{if } i \neq j. \end{cases}$$

If we rotate an orthonormal basis $e_1, e_2, e_3$ through a family of orthonormal bases $e_1(t), e_2(t), e_3(t)$, along some curve $x(t)$, these still have the same constant values of $e_i(t) \cdot e_j(t)$ at every time $t$. Differentiate: $0 = \dot{e}_i(t) \cdot e_j(t) + e_i(t) \cdot \dot{e}_j(t)$. Therefore we can write any infinitesimal rotation of frame as $\dot{e}_j = \sum_i a_{ij} e_i$ for an antisymmetric $3 \times 3$ matrix $A = (a_{ij})$. The quantity $a_{ij}$ measures how quickly $e_j$ is moving toward $e_i$.

The *Levi-Civita connection forms* are the 1-forms $\gamma_{ij} = e_i \cdot de_j$, i.e. $v \lrcorner \gamma_{ij} = a_{ij}$: so $\gamma_{ij}$ measures the tendency of $e_j$ to move toward $e_i$ as the frame moves. In particular, $0 = \gamma_{ij} + \gamma_{ji}$. These $\omega_i$ and $\gamma_{ij}$ are defined on ${}^{\ulcorner}\mathbb{E}^3$, *not* on $\mathbb{E}^3$, because they depend on $x$ and $e$. If we move a frame, we said it moves by a velocity vector $v = (\dot{x}, \dot{e})$ with $\dot{e}_i = \sum_j a_{ji} e_j$ for an antisymmetric matrix $a_{ij}$, so

$$v \lrcorner \gamma_{ij} = v \lrcorner e_i \cdot de_j = e_i \cdot \dot{e}_j = a_{ij}.$$

So the 1-forms $\omega_i$ restrict to any "moving frame" $(x(t), e(t))$ to describe how the velocity of the moving point $x(t)$ is expressed in the moving frame $e(t)$, while the 1-forms $\gamma_{ij}$ describe how the infinitesimal rotation of the moving frame $e(t)$ is expressed at each moment in the moving frame $e(t)$.

Write our soldering forms as $\omega$, thought of as a column of 1-forms

$$\omega = \begin{pmatrix} \omega_1 \\ \omega_2 \\ \omega_3 \end{pmatrix}$$

and our connection 1-forms as

$$\gamma = \begin{pmatrix} \gamma_{11} & \gamma_{12} & \gamma_{13} \\ \gamma_{21} & \gamma_{22} & \gamma_{23} \\ \gamma_{31} & \gamma_{32} & \gamma_{33} \end{pmatrix} = \begin{pmatrix} 0 & \gamma_{12} & \gamma_{13} \\ -\gamma_{12} & 0 & \gamma_{23} \\ -\gamma_{13} & -\gamma_{23} & 0 \end{pmatrix}$$

an antisymmetric matrix of 1-forms, the *Levi-Civita connection*.

**D.8** Prove that the soldering and Levi-Civita connection forms satisfying the *structure equations of Euclidean space*:

$$d\omega_i = -\sum_j \gamma_{ij} \wedge \omega_j,$$
$$d\gamma_{ij} = -\sum_k \gamma_{ik} \wedge \gamma_{kj}.$$

**D.9** For $g$ any $3 \times 3$ orthogonal matrix, if we write $r_g(x, e)$ to mean $(x, eg)$, and $g^\top$ for the transpose of $g$, prove that $r_g^* \omega = g^\top \omega$ and $r_g^* \gamma = g^\top \gamma g$. Expanding out, this means $r_g^* \omega_i = \sum_j g_{ji} \omega_j$ and $r_g^* \gamma_{ij} = \sum_{k\ell} g_{ki} \gamma_{k\ell} g_{\ell j}$.

## Curves and forms

In this section, we adopt the notation that $i, j, k, \ell = 2, 3$. Take a curve $C$ in $\mathbb{E}^3$. At each point of $\ulcorner C$, $e_2, e_3$ are perpendicular to $T_x C$, so $0 = \omega_2 = \omega_3$ on $\ulcorner C$. Note that $\omega_1 = e_1 \cdot dx = ds$, if we pick $e_1$ to agree with the direction in which the arclength function $s$ increases. So $d\omega_1 = 0$. From the structure equations of Euclidean space,

$$
\begin{aligned}
0 &= d\omega_2 = -\gamma_{21} \wedge \omega_1. \\
0 &= d\omega_3 = -\gamma_{31} \wedge \omega_1.
\end{aligned}
$$

Therefore $\gamma_{21} = k_2 ds, \gamma_{31} = k_3 ds$ for unique smooth functions $k_2, k_3$ on $\ulcorner C$, but these change sign if we change the direction in which we measure arclength $s$. We can move any adapted frame $(x, e)$ of $C$ with $x$ moving along $C$, or with the $e_i$ rotating among frames of the tangent line. The 1-forms measuring those motions $\omega_1, \gamma_{23}$ are all linearly independent, together forming a *coframing* on $\ulcorner C$, i.e. a collection of linearly independent 1-forms spanning the 1-forms at each point. We obtain the *structure equations* of a curve $C$ in $\mathbb{E}^3$:

$$
\begin{aligned}
d\omega_1 &= 0, \\
d\gamma_{23} &= 0.
\end{aligned}
$$

Each of the 4 curves in the frame bundle, from the 4 Serret–Frenet framings, have their own arclength function $s$, defined to agree in direction with $e_1$, but only defined up to adding a constant. In other words, since $d\omega_1 = 0$, locally $\omega_1 = ds$ for a function $s$ defined up to a constant. On the frame bundle of $C$, we can turn $e_2, e_3$, so have a nonzero 1-form $\gamma_{23}$.

**D.10** Prove that the vector $k_2 e_2 + k_3 e_3$ is actually defined down on $C$, and equals the curvature vector.

If the curvature vector is nonzero, we can define the Serret–Frenet frames to be those on which $k_2 > 0$ and $k_3 = 0$. The set of Serret–Frenet frames consists of 4 curves in the frame bundle. It has $k_3 = 0$, so $\gamma_{31} = 0$, and $\gamma_{21} = k\omega_1 \neq 0$, so differentiate $0 = \gamma_{31}$ to get

$$
\begin{aligned}
0 &= -\gamma_{32} \wedge \gamma_{21}, \\
&= -k\gamma_{32} \wedge \omega_1.
\end{aligned}
$$

so $\gamma_{32} = t\omega_1$ for a unique function $t$ on those 4 curves. We can see a general pattern emerging: write out the structure equations, get the geometry to impose some relation on the 1-forms in the structure equations, and differentiate repeatedly, applying Cartan's lemma when possible, until the invariants pop out, and the structure equations reduce down to having only the linearly independent 1-forms in them.

## Surfaces and forms

The expression $dx \cdot dx$ is the Euclidean inner product in Euclidean space, so is well defined on tangent vectors to any curve or surface in Euclidean space. Take a surface $S$ in $\mathbb{E}^3$. The unit normal $e_3$ is defined up to sign at each point of $S$. Hence $de_3$ is also defined up to sign, as is $de_3 \cdot dx$, and so $e_3(de_3 \cdot dx)$ is defined on $S$, called the *shape operator*, or the *second fundamental form*, denoted $I\!I$. Let us find a different way to write this.

As we move along a curve in $S$, take $e_1$ to be tangent to that curve, and travel at unit speed; pick $e_2, e_3$ to give an adapted frame along the curve. The vector $e_2$ moves in the direction of $e_1$, i.e. the curve turns in the tangent plane to the surface, as the $\omega_1$ component of $e_1 \cdot de_2 = \gamma_{12}$. So this component is the curvature of the curve in the surface. Recall that $e_3$ is normal to the surface; stand on the surface, with $e_3$ pointing straight upward. The vector $e_3$ moves down toward $e_1$, i.e. the surface bends down as we move along the curve, as the $\omega_1$ component of $e_1 \cdot de_3 = \gamma_{13}$. The vector $e_3$ moves toward $e_2$, i.e. the surface twists like a screw, as we move along the curve, as the $\omega_1$ component of $e_2 \cdot de_3 = \gamma_{23}$.

From here on, adopt notation $i, j, k, \ell = 1, 2$. At each point of $\ulcorner S$, $e_3$ is perpendicular to $T_x S$, so $\omega_3 = 0$ on $\ulcorner S$. From the structure equations of Euclidean space,

$$0 = d\omega_3 = -\sum_i \gamma_{3i} \wedge \omega_i.$$

Therefore $\gamma_{3i} = -\sum_j a_{ij}\omega_j$ for unique smooth functions $a_{ij} = a_{ji}$ on $\ulcorner S$. So $a_{ij}$ measures the tendency of $e_i$ to move toward $e_3$ as we move in direction of $e_j$, i.e. we measure how the surface twists like a screw. This measurement is (not obviously) symmetric in $i, j$: if the surface twists as we move in a certain direction, it twists just as much if we move in the orthogonal direction. These $a_{ij}$ are the components of the shape operator:

$$\begin{aligned} I\!I &= e_3(de_3 \cdot dx), \\ &= e_3((e_i \cdot de_3)(e_i \cdot dx)), \\ &= e_3(\gamma_{i3}\omega_i), \\ &= e_3(a_{ij}\omega_j\omega_i). \end{aligned}$$

We can move any adapted frame $(x, e)$ of $S$ with $x$ moving along $S$ and the $e_i$ rotating among frames of the tangent space. The 1-forms measuring those motions $\omega_1, \omega_2, \gamma_{12}$ are all linearly independent, together forming a *coframing* on $\ulcorner S$, i.e. a collection of linearly independent 1-forms spanning the 1-forms at each point. On our surface, $0 = \omega_3$, so often forget $\omega_3$ and let $\omega \coloneqq \omega_1 + i\omega_2$. To simplify notation, let $\alpha \coloneqq \gamma_{12}$ (because the letter $\alpha$ looks like $\gamma$ rotated in the plane). We obtain the *structure equations* of a surface $S$ in $\mathbb{E}^3$:

$$
\begin{aligned}
d\omega &= i\alpha \wedge \omega, \\
d\alpha &= \frac{i}{2} K \omega \wedge \bar{\omega}
\end{aligned}
$$

where $K$ is the *Gauss curvature*. From the structure equations of Euclidean space, $d\alpha = d\gamma_{12} = -\gamma_{13} \wedge \gamma_{32}$, so $K = a_{11}a_{22} - a_{12}^2$.

## Rotate the frame

This section can be omitted without loss of continuity.

Recall that we write each orthogonal $3 \times 3$ matrix $g$ preserving the horizontal plane as

$$
g = \begin{pmatrix} h & 0 \\ 0 & n \end{pmatrix}
$$

with $h \in \mathrm{O}(2)$ and $n = \pm 1$. Write $h$ as a product of a rotation matrix $e^{i\theta}$ and zero or one complex conjugation matrices, which we write as $C$. Expanding out the equation $r_g^* \omega = g^\top \omega$ with $\omega_3 = 0$ gives $r_g^* \omega = h^\top \omega$. In complex notation

$$
\begin{aligned}
r_{e^{i\theta}}^* \omega &= e^{-i\theta} \omega, \\
r_C^* \omega &= \bar{\omega}, \\
r_n^* \omega &= \omega.
\end{aligned}
$$

Expanding out the equation $r_g^* \gamma = g^\top \gamma g$ into $3 \times 3$ matrices gives

$$
\begin{aligned}
r_g^* \alpha &= (\det h) \alpha, \\
r_g^* \begin{pmatrix} \gamma_{13} \\ \gamma_{23} \end{pmatrix} &= n h^\top \begin{pmatrix} \gamma_{13} \\ \gamma_{23} \end{pmatrix}.
\end{aligned}
$$

If

$$
a := \begin{pmatrix} a_{11} & a_{12} \\ a_{21} & a_{22} \end{pmatrix},
$$

then $r_g^* a = n h^\top a h$. In our complex notation, with $\alpha := \gamma_{12}$,

$$
\begin{aligned}
r_{e^{i\theta}}^* \alpha &= \alpha, \\
r_C^* \alpha &= -\alpha, \\
r_n^* \alpha &= \alpha.
\end{aligned}
$$

## Mean curvature

The *mean curvature* is $H = \frac{1}{2}(a_{11} + a_{22})$. The *mean curvature vector* is $He_3$.

**D.11** Prove that the mean curvature vector is defined on the surface, a vector field pointing normal to the surface. Why might the mean curvature not be defined as a function on the surface?

The sphere $\mathbb{E}^3$ of radius $r_0$ has frame bundle consisting of the $(x, e) \in \ulcorner\mathbb{E}^3$ so that $x = \pm r_0 e_3$. Therefore $\gamma_{i3} = e_i \cdot de_3 = \pm e_i \cdot dx/r_0 = \pm\omega_i/r_0$, so $a = \pm I/r_0$ depending on orientation. Gauss curvature is $1/r_0^2$. Mean curvature is $H = \pm 1/r_0$, depending on orientation. The mean curvature vector is $2x/r_0$, independent of orientation.

For a flat plane, $e_3$ lies perpendicular to the plane, while $x, e_1, e_2$ are tangent to it, so $0 = e_3 \cdot dx = e_3 \cdot de_1 = e_3 \cdot de_2$, i.e. $0 = \omega_3 = \gamma_{31} = \gamma_{32}$, so $a = 0$.

**D.12** Compute the shape operator and mean and Gauss curvature of a cylinder of radius $r_0$.

**D.13** We know that $I\!I$ is defined on tangent vectors to $S$, independent of frame; check this again using the rules above for how $a_{ij}$ change when we change frame.

## Principal curvatures

The shape operator of a surface is a symmetric bilinear form on the tangent place $T_x S$ valued in the normal line to the tangent plane, so orthogonally diagonalizable. Its eigenvalues are the *principal curvatures*. If it has only one principal curvature, i.e. $a$ is a multiple of the identity matrix, the point $x$ is an *umbilic* of the surface. On the other hand, if it has two distinct principal curvatures then its two perpendicular eigenlines are the *principal directions*.

## Higher fundamental forms

**D.14** Let

$$Da := da + \alpha \begin{pmatrix} 2a_{12} & a_{22} - a_{11} \\ a_{22} - a_{11} & -2a_{12} \end{pmatrix}.$$

Differentiate the equations $\gamma_{i3} = a_{ij}\omega_j$ to reveal that $Da_{ij} = \sum_k a_{ijk}\omega_k$ where $a_{ijk}$ is symmetric in all lower indices.

The *third fundamental form* is $I\!I\!I := e_3 a_{ijk}\omega_i\omega_j\omega_k$.

## Local picture of a surface

Take a point $x_0$ of a surface $S$ and let $n$ be a unit normal vector to $S$ at $x_0$, say $n = -e_3$ at $x_0$. Fixing that constant value of $n$, take the linear function $f(x) := n \cdot x$. The differential of $f$ is

$$e_i \lrcorner df = e_i \lrcorner d\,(n \cdot x) = n \cdot e_i,$$

for $i = 1, 2$, vanishing at $x_0$. The second derivative matrix of $f$ at $x_0$ is

$$\begin{aligned} f''(e_i, e_j) &= e_i \lrcorner d\,(n \cdot e_j)\,, \\ &= e_i \lrcorner n \cdot de_j, \\ &= -e_i \lrcorner e_3 \cdot de_j, \\ &= -e_i \lrcorner \gamma_{3j}, \\ &= a_{ji}. \end{aligned}$$

Similarly $f'''(e_i, e_j, e_k) = a_{ijk}$. Translate the surface to get $x_0$ the origin, and $T_{x_0}S$ the horizontal plane, and $n$ the unit vertical vector. Then $S$ is locally the graph of $x_3 = f(x_1, x_2)$, and $f$ has a critical zero at the origin, so

$$x_3 = \frac{1}{2} a_{ij} x_i x_j + \frac{1}{6} a_{ijk} x_i x_j x_k + O(x)^4.$$

We can pick the frame $e_1, e_2$ at our point to diagonalize the shape operator, say with eigenvalues (i.e. principal curvatures) $k_1, k_2$:

$$x_3 = \frac{k_1}{2} x_1^2 + \frac{k_2}{2} x_2^2 + \frac{1}{6} a_{ijk} x_i x_j x_k + O(x)^4.$$

**D.15** Prove that, if the Gauss curvature is positive at $x_0$ then $S$ lies on one side of its tangent plane near $x_0$, and if the Gauss curvature is negative then it doesn't.

**D.16** Prove that every compact surface in $\mathbb{E}^3$ has a point of positive Gauss and mean curvature.

Appendix E

# Curvature of surfaces

## Example: surfaces of revolution

Take a curve $C$ in the plane, and rotate it around a line in that plane (the *axis of revolution*) to create a *surface of revolution* $S$. Suppose for simplicity that the curve doesn't touch the axis of revolution. Place the axis of revolution vertically. Write the points of $\mathbb{E}^3$ in cylindrical coordinates $(r\cos\theta, r\sin\theta, z)$. Our plane is $\theta = 0$, and our axis of revolution is $r = 0$. Our plane curve $C$ is arclength parameterized as $x'(s) = (r(s), 0, z(s))$. Take adapted frames $(x', e')$ for the curve in the plane, but swap $e'_2$ and $e'_3$ so that $e'_2$ is normal to the plane and $e'_3$ tangent to the plane but normal to the curve:

$$\begin{aligned} e'_1(s) &= \frac{dx'}{ds} = (\cos\varphi(s), 0, \sin\varphi(s)), \\ e'_2(s) &= (0, 1, 0), \\ e'_3(s) &= (-\sin\varphi(s), 0, \cos\varphi(s)). \end{aligned}$$

Let $\rho(\theta)$ be the rotation matrix rotating $\mathbb{E}^3$ about the axis of revolution by an angle $\theta$. Make an adapted frame at each point of the surface as

$$\begin{aligned} x(s,\theta) &= \rho(\theta)x'(s) = (r\cos\theta, r\sin\theta, z), \\ e(s,\theta) &= \rho(\theta)e'(s), \end{aligned}$$

**E.1** Find all $\omega_i$, $\gamma_{ij}$.

The shape operator:

$$\begin{pmatrix} \gamma_{13} \\ \gamma_{23} \end{pmatrix} = \begin{pmatrix} -\dot\varphi & 0 \\ 0 & -\frac{\sin\varphi}{r} \end{pmatrix} \begin{pmatrix} \omega_1 \\ \omega_2 \end{pmatrix}.$$

So the Gauss curvature of a surface of revolution is

$$K = \frac{\dot\varphi \sin\varphi}{r} = -\frac{\ddot r}{r},$$

while the mean curvature is

$$H = -\frac{1}{2}\left(\dot\varphi + \frac{\sin\phi}{r}\right).$$

Take a function $H(s)$ and try to construct a curve $x'(s)$ in the plane, so that the associated surface of revolution has mean curvature $H(s)$. To do this, we have to solve the equations

$$\dot{r}(s) = \cos\varphi(s),$$
$$\dot{\varphi}(s) = -2H(s) - \frac{\sin\varphi(s)}{r(s)}$$

with $r > 0$. Each initial value $r = r_0 > 0$ and $\varphi = \varphi_0$ at some $s = s_0$ determines a unique local solution, by the existence and uniqueness theorem for ordinary differential equations. Solve for $z$ by $z(s) = \int \sin\varphi(s)\,ds$, again unique up to adding a constant.

**E.2** If we want $H$ to be constant, say equal to a constant $H_0$, can you integrate this coupled system of ordinary differential equations? Hint: is $r\sqrt{1-\dot{r}^2} + H_0 r^2$ constant along solutions?

Similarly, we can prescribe Gauss curvature $K(s)$ arbitrarily by solving the equation $\ddot{r} + Kr = 0$. For example, the surfaces of revolution with zero Gauss curvature are the circular cones given by rotating a line segment. More generally, for constant Gauss curvature $K = K_0$, the solutions are (with $k_0 := \sqrt{|K_0|}$)

$$r(s) = \begin{cases} a\cos k_0 s + b\sin k_0 s, & \text{if } K_0 > 0,\\ a + bs, & \text{if } K_0 = 0,\\ a\cosh k_0 s + b\sinh k_0 s, & \text{if } K_0 < 0.\end{cases}$$

For $K_0 > 0$, translate $s$ to arrange $b = 0$, so write

$$r(s) = r_0 \cos k_0 s,$$

so that

$$\dot{z}^2 = 1 - \dot{r}^2 = 1 - r_0^2 K_0 \sin^2 k_0 s,$$

giving

$$z = \int \sqrt{1 - r_0^2 K_0 \sin^2 k_0 s}\,ds,$$

an elliptic integral, usually denoted $E(k_0 s \,|\, r_0^2 K_0)/k_0$. Similar remarks apply to $K_0 < 0$.

**E.3** Take a curve $C_0$ in the plane $\mathbb{R}^2$. Prove that there is a curve $C_1$ in the plane $\mathbb{E}^2$ whose associated surface of revolution satisfies $(H, K) \in C_0$ at every point.

## Orientations

We could have used only positively oriented frames $e_1, e_2, e_3$ of $\mathbb{E}^3$. We could then look at an oriented surface, and define its *adapted frames* to have $e_1, e_2$ a positively oriented basis of the tangent plane, and $e_1, e_2, e_3$ a positively oriented basis of $\mathbb{E}^3$. The expression

$$g = \begin{pmatrix} h & 0 \\ 0 & n \end{pmatrix}$$

simplifies to having $h$ a rotation matrix, and $n = 1$, so

$$g = \begin{pmatrix} \cos\theta & -\sin\theta & 0 \\ \sin\theta & \cos\theta & 0 \\ 0 & 0 & 1 \end{pmatrix},$$

which we just write as $e^{i\theta}$.

**E.4** Over an oriented surface, prove that the 2-form $dA := \omega_1 \wedge \omega_2$ is the pullback of a unique 2-form $dA$ on the surface, the *area form.*

**E.5** Over an oriented surface, prove that the mean curvature is the pullback of a unique function on the surface.

**E.6** Prove that, on any surface, oriented or not, the shape operator admits a unique expression as

$$\frac{q}{2}\omega^2 + \frac{H}{2}(\omega\bar{\omega} + \bar{\omega}\omega) + \frac{\bar{q}}{2}\bar{\omega}^2.$$

for functions $q \colon {}^{\ulcorner}S \to \mathbb{C}$, $H \colon {}^{\ulcorner}S \to \mathbb{R}$. How are the $a_{ij}$ related to $q, H$?

The expression $q\omega^2$ is the *Hopf differential.*

**E.7** If the surface is oriented, prove that the Hopf differential descends to the surface, as a complex valued quadratic form on tangent vectors.

The *Gauss map* of a surface is $e_3$, mapping the frame bundle to the unit sphere.

**E.8** Prove that the Gauss map of an oriented surface is the composition of the map $(x, e) \in {}^{\ulcorner}S \mapsto x$ with a unique map $x \in S \mapsto e_3(x)$, also called the *Gauss map.*

## Recovering a surface from its inner products and shape operator

This section can be omitted without loss of continuity.

**Theorem E.1** (Bonnet). *If two connected oriented surfaces in 3-dimensional Euclidean space are identified by a diffeomorphism preserving orientation and the inner product and shape operator on every tangent space, then the diffeomorphism is induced by an orientation preserving rigid motion of 3-dimensional Euclidean space.*

*Proof.* Use the orientations of the surfaces to define "identifying" the shape operators: we have a well defined normal vector on each surface, so we treat the shape operator as a symmetric bilinear form on tangent spaces. Such a diffeomorphism $S \to S'$, clearly identifies positively oriented orthonormal frames $e_1, e_2$ of tangent spaces of $S$ with those of $S'$. Taking the usual orientation of $\mathbb{E}^3$, these extend uniquely into positive oriented orthonormal frames $e_1, e_2, e_3 = e_1 \times e_2$, a diffeomorphism of the adapted frames, clearly identifying the soldering 1-forms $\omega_1, \omega_2$ and $\omega_3 = 0$. Differentiate to find that the Levi–Civita connection 1-form $\gamma_{12}$ matches. Having the same shape operator ensures precisely that $\gamma_{31}, \gamma_{32}$ match. It then follows by antisymmetry of $\gamma_{ij}$ that all connection 1-forms agree. Let

$$h = h(x, e) = \begin{pmatrix} 1 & 0 \\ x & e \end{pmatrix}$$

as above, and $h'$ in the same way for the corresponding surface:

$$h^{-1} dh = h'^{-1} dh' = \begin{pmatrix} 0 & 0 \\ \omega & \alpha \end{pmatrix},$$

under the diffeomorphism. So $h$ and $h'$ satisfy the same first order linear ordinary differential equation

$$dh = \begin{pmatrix} 0 & 0 \\ \omega & \alpha \end{pmatrix} h$$

along any curve in the frame bundle $\ulcorner S$. Apply a rigid motion to one of the two surfaces, so that their frame bundles contain at least one frame in common. By existence and uniqueness for linear ordinary differential equations, $h = h'$ along any curve in $\ulcorner S$ starting at that frame. The adapted frames which are positively oriented for $\mathbb{E}^3$ and with $e_1, e_2$ positively oriented on the surface form a connected manifold. Hence $h = h'$. □

## Parallel surfaces

This section can be omitted without loss of continuity.

Start with a surface $S$ with a chosen normal direction, i.e. a unit normal vector $e_3$ at each point, smoothly varying. Let $S'$ be the translation of $S$ in that direction by a distance $t$. So then $x' = x + te_3$, $dx' = dx + t\, de_3$ and $e_3 \cdot dx' = e_3 \cdot dx + te_3 \cdot de_3 = 0$, so $e_3$ is still normal to $S'$ at $x'$. Hence $T_{x'} S' = T_x S$ is constant in $t$, and the Gauss map is unchanged. Identify $\ulcorner S' \cong \ulcorner S$ by

$$\Phi \colon (x, e) \mapsto (x', e') = (x + te_3, e).$$

The soldering 1-forms are

$$\begin{aligned} \omega_i' &= e_i \cdot dx', \\ &= e_i \cdot dx + te_i \cdot de_3, \\ &= \omega_i + t\gamma_{i3}, \\ &= \omega_i + ta_{ij}\omega_j, \end{aligned}$$

which we can write as $\omega' = (I + ta)\omega$. The Levi-Civita connection 1-form is

$$\gamma'_{12} = e_1 \cdot de_2 = \gamma_{12}.$$

The shape operator is

$$\begin{aligned} \gamma'_{i3} &= e_i \cdot de_3, \\ &= \gamma_{i3}, \\ &= a_{ij}\omega_j \end{aligned}$$

If we write this as $\gamma'_{*3} = a\omega$, then

$$\gamma'_{*3} = a\omega = a(I + ta)^{-1}\omega',$$

so

$$a' = a(I + ta)^{-1}.$$

In any frame in which $a$ is diagonalized, so is $a'$:

$$a = \begin{pmatrix} k_1 & 0 \\ 0 & k_2 \end{pmatrix}$$

in principal curvatures, so

$$a' = \begin{pmatrix} k'_1 & 0 \\ 0 & k'_2 \end{pmatrix} = \begin{pmatrix} \frac{k_1}{1+tk_1} & 0 \\ 0 & \frac{k_2}{1+tk_2} \end{pmatrix}.$$

Hence Gauss curvature

$$K' = \frac{\det a}{\det(I + ta)} = \frac{K}{1 + 2Ht + Kt^2}$$

and mean curvature

$$H' = \frac{H + Kt}{1 + 2Ht + Kt^2}.$$

In particular, the parallel surface will not be smooth enough to compute such quantities at any point $x' = x + te_3$ where $1 + tk_1$ or $1 + tk_2$ vanish. But if these don't vanish, then the differential of the map is

$$dx' = dx + t\,de_3 = \omega_i\,e_i + ta_{ij}\omega_j e_i$$

an immersion. Surfaces parallel to a Weingarten surface are also Weingarten, but for a different relation between mean and Gauss curvature.

## Constant mean curvature surfaces

This section can be omitted without loss of continuity.

On an oriented surface, a *complex differential* is an expression of the form $f\omega^k$ on the frame bundle, but which is in fact defined on the surface. In other words, it doesn't transform as we right translate by $e^{i\theta}$. Since $\omega$ transforms by $e^{-i\theta}$, we need $f$ to transform by $e^{ik\theta}$:

$$r^*_{e^{i\theta}} f = e^{ik\theta} f.$$

**E.9** Differentiating this, show that

$$df + ikf\alpha = Df\omega + \bar{D}f\bar{\omega},$$

for unique complex functions $Df, \bar{D}f$.

A *holomorphic differential* is one for which $0 = \bar{D}f$, i.e. $df - ikf\alpha$ is complex linear.

A *holomorphic chart* on an oriented surface in $\mathbb{E}^3$ is a chart on an open subset of the surface with differential $dz$ a holomorphic differential. Any two such charts locally agree up to complex analytic transition maps. The Korn–Lichtenstein theorem says that any oriented surface is covered by domains of holomorphic charts. Hence any oriented surface in $\mathbb{E}^3$ is a Riemann surface.

**E.10** Compute that the Hopf differential satisfies

$$dq + 2iq\alpha = Dq\omega + \bar{D}q\bar{\omega},$$

where $\bar{D}q = 2(H_1 - H_2 i)$ where $H$ is the mean curvature, and $dH = H_1\omega_1 + H_2\omega_2$. Thereby, prove:

**Theorem E.2.** *An oriented surface has constant mean curvature just when its Hopf differential is holomorphic.*

The Riemann–Roch theorem [13] p. 99 shows that a compact Riemann surface of genus $g$ has a space of quadratic differentials of complex dimension:

$$\begin{cases} 0, & g = 0, \\ 1, & g = 1, \\ 3g - 3, & g \ge 2. \end{cases}$$

It is well known that there is precisely one genus zero compact Riemann surface, up to holomorphic diffeomorphism: the Riemann sphere [13] p. 125. Applying Bonnet's theorem:

**Corollary E.3.** *For any constant $H_0$, the sphere has at most one immersion in $\mathbb{E}^3$ with constant mean curvature equal to $H_0$, up to rigid motion: as a round sphere of radius $1/H_0$.*

## Away from umbilics

Recall that principal curvatures of a surface are the eigenvalues of its shape operator. Suppose that they are distinct $k_1 \ne k_2$, i.e. there are no umbilic points. A *Darboux frame* is one in which $e_1, e_2$ are eigenvectors of the shape operator with eigenvalues $k_1 < k_2$. The set of Darboux frames of the surface is a surface, as each Darboux frame is uniquely determined up to signs. On that surface $\gamma_{i3} = k_i\omega_i$, for $i = 1, 2$.

**E.11** Differentiate the equations of a Darboux frame to see that

$$d\begin{pmatrix} k_1 \\ k_2 \end{pmatrix} = \begin{pmatrix} k_{11} & k_{12} \\ k_{21} & k_{22} \end{pmatrix} \begin{pmatrix} \omega_1 \\ \omega_2 \end{pmatrix},$$

for some functions $k_{ij}$, perhaps *not* symmetric in lower indices, while

$$\gamma_{12} = \frac{k_{12}\omega_1 + k_{21}\omega_2}{k_2 - k_1}.$$

Expand $0 = d^2 k_i$ to find that, if we let

$$Dk_{i1} = dk_{i1} + k_{i2}\gamma_{12},$$
$$Dk_{i2} = dk_{i2} - k_{i1}\gamma_{12},$$

then $Dk_{ij} = k_{ijk}\omega_k$ with $k_{ijk} = k_{ikj}$. Recall that Gauss curvature is $K = k_1 k_2$. Plug into the equation

$$d\gamma_{12} = K\omega_1 \wedge \omega_2.$$

to get

$$d\left(\frac{k_{12}\omega_1 + k_{21}\omega_2}{k_2 - k_1}\right) = k_1 k_2 \omega_1 \wedge \omega_2.$$

Expanding this out gives

$$0 = \frac{k_{122} - k_{211}}{k_2 - k_1} + \frac{2k_{12}^2 + 2k_{21}^2}{(k_2 - k_1)^2} + k_1 k_2.$$

**Lemma E.4.** *If a surface bears a point at which the larger principal curvature is maximal and the smaller minimal, then either the point is umbilic or the Gauss curvature is not positive at that point.*

*Proof.* We can assume that $k_2 > k_1$ near that point, and then adapt frames as above. At a critical point of $k_2$ and $k_1$, $0 = dk_1 = dk_2$, so $k_{ij} = 0$ and so $\gamma_{12} = 0$, so $Dk_{ij} = dk_{ij} = k_{ijk}\omega_k$. Since $k_2$ is maximal, along the flow of $e_1$, $k_{211} \le 0$. Since $k_1$ is minimal, along the flow of $e_2$, $k_{122} \ge 0$. Plug in above to get $K = k_1 k_2 \le 0$. □

**Theorem E.5** (Liebmann). *Suppose that $S$ is a connected compact surface in $\mathbb{E}^3$ of (i) constant Gauss curvature or (ii) constant mean curvature with Gauss curvature not having any critical zeroes. Then $S$ is a round sphere.*

*Proof.* By problem D.16 on page 113, there is a point of positive Gauss and mean curvature, so the constant is a positive constant. At any point at which the principal curvatures are as far apart from one another as possible, one is maximal and the other minimal, as their product is constant. Apply lemma E.4 to find a point of nonpositive Gauss curvature where the Gauss curvature is critical. □

**E.12** Suppose that $S$ is a connected surface in $\mathbb{E}^3$ and that one principal curvature is a nonzero constant, and always distinct from the other. Prove that $S$ is a union of circular arcs of constant radius, with centres on a curve.

We define completeness in appendix G; for the moment, we only need know that on a complete surface, every bounded vector field is complete.

**Theorem E.6** (Hilbert). *No complete surface in $\mathbb{E}^3$ has constant negative Gauss curvature.*

*Proof.* Take a surface of constant negative Gauss curvature, say equal to $-1$. We can assume it is connected, and replace it by its universal covering space if needed, so assume it is simply connected. Let $u$ be half the angle between the asymptotic curves, so $0 < u < \pi/2$. Check that the principal curvatures are $-\cot u, \tan u$. A surface of negative Gauss curvature has no umbilics. Differentiate the equations

$$\begin{aligned}\gamma_{13} &= -\cot u\, \omega_1,\\ \gamma_{23} &= \tan u\, \omega_2.\end{aligned}$$

Check that $\omega_1/\sin u$ and $\omega_2/\cos u$ are closed 1-forms; noting that $0 < u < \pi/2$, so the denominators don't vanish. Because the surface is connected and simply connected, there are functions $s_1, s_2$ so that

$$\begin{aligned}ds_1 &= \omega_1/\sin u,\\ ds_2 &= \omega_2/\cos u,\end{aligned}$$

unique up to constant and with $ds_1, ds_2$ linearly independent.

The surface is complete, so $e_1, e_2$ have complete flows. The indefinite integral $\int \omega_1$ along curves $\omega_2 = 0$ is an unbounded function, defined up to constant. But $s_1$ grows at least as fast, so is also unbounded along curves on which $ds_2 = 0$, and vice versa. So the map $(s_1, s_2)$ from our surface to the plane is onto, and a local diffeomorphism.

The vector fields $\partial_{s_1} = \sin u\, e_1$, $\partial_{s_2} = \cos u\, e_2$ are complete and commuting. Flows of these for times $s_1, s_2$ invert our previous map, hence a diffeomorphism.

Let $x := (s_1 - s_2)/2$, $y := (s_1 + s_2)/2$. The level sets of $x$ and $y$ are precisely the asymptotic curves. Check that $\gamma_{12} = -u_x dx + u_y dy$, and that the relation $d\gamma_{12} = K\omega_1 \wedge \omega_2$ becomes the *sine–Gordon equation* for the angle $\theta = 2u$:

$$\theta_{xy} = \sin\theta.$$

Note also that $0 < \theta < \pi$, so $\sin\theta > 0$. The area form is $\omega_1 \wedge \omega_2 = 2\sin\theta dx \wedge dy$.

Apply Stokes's theorem to a rectangle $R$ in the $xy$-plane

$$\int_{\partial R} \gamma_{12} = \int_R d\gamma_{12},$$

to find

$$\theta(x, y) - \theta(0, y) = \theta(x, 0) - \theta(0, 0) + \int \sin\theta\, dx\, dy.$$

The difference in values of $\theta$ along the top of any rectangle is at least as large as the difference along the bottom.

Since $0 < \theta < \pi$, the area of the rectangle, in the geometry of our surface, is never more than $4\pi$. Our coordinates are global, so any measureable subset of the surface has area no more than $4\pi$. Each set of infinite Euclidean area thus contains sets of large measure where $\sin\theta$ gets small, so $\theta$ gets close to 0 or $\pi$.

Suppose that $\theta$ is not constant on some horizontal line. Change the sign of $x$ and $y$ if needed, to arrange that $\theta$ increases along a horizontal line segment. Split that line segment into three line segments. On the leftmost, suppose that $\theta$ increases by some amount $\varepsilon_L$, and on the rightmost by some amount $\varepsilon_R$.

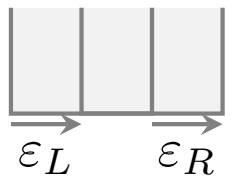

Draw the half strips above the line segments. So $\theta$ increases by at least $\varepsilon_L$ across each parallel line segment across the left half strip. But $\theta > 0$ everywhere, so in the middle strip, $\theta > \varepsilon_L$. Similarly, in the middle strip, $\theta < \pi - \varepsilon_R$. So $\sin\theta$ is bounded away from zero in the middle strip, a contradiction.

Hence $\theta$ is constant on every horizontal line, i.e. $0 = \theta_x$, so $0 = \theta_{xy} = \sin\theta$, a contradiction. $\square$

Appendix F

# The Gauss–Bonnet theorem

We prove that the integral of Gauss curvature over a compact oriented surface in 3-dimensional Euclidean space is a diffeomorphism invariant.

## Length and curvature of curves

**F.1** Suppose that $x(t)$ is a smooth curve in $\mathbb{E}^3$ so that $x'(t) \neq 0$ for every $t$ with $t_0 \le t \le t_1$. Prove that there is a reparameterization $y(t) = x(\tau(t))$ for some function $\tau(t)$ so that $y'(t)$ is a unit length vector, i.e. $y(t)$ has unit speed. Prove that $\tau(t)$ is uniquely determined up to adding a constant. Prove that the length of $y(t)$ for $a \le t \le b$ is equal to the length of $x(t)$ for $\tau(a) \le t \le \tau(b)$.

Take an oriented curve $C$, in other words an oriented 1-dimensional manifold, and a smooth immersion $f : C \to S$ to a surface $S$. An *adapted frame* for $C$ at a point $c \in C$ is a frame $(x, e)$ for $S$, where $x = f(c)$, so that $e_1$ is a positively oriented basis of the tangent line to $C$ at $m$. Let ${}^{\ulcorner}S_C$ be the set of all pairs $(c, e)$ so that $c \in C$ and $e$ is an adapted frame for $C$ at $c$.

**F.2** Prove that the set ${}^{\ulcorner}S_C$ is a curve in $C \times {}^{\ulcorner}S$.

Map $(c, e) \in {}^{\ulcorner}S_C \mapsto (f(c), e) \in {}^{\ulcorner}S$. Using this map to pullback $\omega_1, \omega_2$ and $\alpha = \gamma_{12}$. If we parameterize $C$, say as $x(t)$, then $e_1(t)$ is the unit vector in the direction of $x'(t)$, so that $\omega_1 = e_1(t) \cdot dx = e_1(t) \cdot \frac{dx}{dt}\, dt = |\dot{x}(s)|\, dt$. So $\omega_1$ is a well defined 1-form on the curve $C$ and $\int_C \omega_1$ is the length of the curve $C$. Meanwhile $e_2$ is perpendicular to the tangent line to $C$, so on ${}^{\ulcorner}S_C$, $\omega_2 = e_2 \cdot dx = e_2 \cdot \frac{dx}{dt}\, dt = 0$. Therefore $0 = d\omega_2 = \alpha \wedge \omega_1$.

Applying Cartan's lemma, $\alpha = -\kappa_2\, \omega_1$ for some function $\kappa_2$. The 1-form $\omega_1$ is a basis for the 1-forms on ${}^{\ulcorner}S_C$. The equations $\omega_2 = 0$ and $\alpha = -\kappa_2 \omega_1$ are precisely the linear equations that cut out the tangent spaces of ${}^{\ulcorner}S_C$ inside ${}^{\ulcorner}S$.

**F.3** Prove that $\kappa_2 e_2$, the *geodesic curvature* of the curve $C$ at the point $c$, is the projection of the curvature vector of $C$ (as a space curve in $\mathbb{E}^3$) to the tangent plane of $S$.

**F.4** Prove that the geodesic curvature of an oriented curve is unchanged if we reorient the curve.

A flat plane $P \subset \mathbb{E}^3$ is perpendicular to a constant unit vector $e_3$, so $0 = de_3 = e_3 \cdot dx = e_3 \cdot de_1 = e_3 \cdot de_2$. Therefore on $^{\ulcorner}P$, the shape operator vanishes. Given a curve in the plane, the adapted frames have $\omega_2 = 0$ and $\alpha = -\kappa_2\omega_1$. If the curve is given by $s \mapsto (x_1(s), x_2(s))$ in an arclength parameterization then we write the velocity vector as a

$$e_1 = (\dot{x}_1, \dot{x}_2) = (\cos\phi, \sin\phi)\,,$$

and then, rotating by a right angle,

$$e_2 = (-\dot{x}_2, \dot{x}_1) = (-\sin\phi, \cos\phi)\,,$$

so that

$$-\alpha = e_2 \cdot de_1 = \dot{\phi}(s)\,ds,$$

i.e. $\kappa = \dot{\phi}(s)$, the curvature is the rate of rotation of the unit tangent vector. We can recover a curve from its curvature $\kappa$, by integrating $\kappa$ to find $\phi$, and then integrating

$$\frac{dx}{ds} = (\cos\phi, \sin\phi)$$

to recover the curve. Constants of integration represent rotation and translation of the curve.

To get a curve of constant curvature $\kappa = \kappa_0$ in the plane, we integrate to find $\phi = \kappa s$, up to a constant, and then up to a constant,

$$x(s) = \frac{1}{\kappa_0}\left(-\sin\kappa_0 s, \cos\kappa_0 s\right),$$

a circle.

**Lemma F.1.** *Take a curve $C$ on a surface $S$ in $\mathbb{E}^3$. Suppose that $C$ has curvature vector as a space curve nowhere perpendicular to $S$. Take $\theta$ to be the angle between that curvature vector and its projection to the surface $S$. Its surface invariants (geodesic curvature $\kappa$ and shape operator $a$) are related to its invariants as a space curve (the curvature $k$ and torsion $t$ of $C$) by*

$$\begin{aligned}
\kappa &= k\cos\theta,\\
a(\dot{x}, \dot{x}) &= k\sin\theta,\\
a(\dot{x}, \dot{x}^{\perp}) &= \dot{\theta} + t
\end{aligned}$$

*where $\dot{x}$ is its unit speed velocity in either orientation, and $\dot{x}^{\perp}$ is the unit tangent vector to $S$ normal to $C$ which is closest to the curvature vector.*

*Proof.* Write the Serret–Frenet frame to $C$ as $e_1, e_2', e_3'$, and all adapted frames to the surface $S$ for which $e_1$ is tangent to the curve $C$ as $e_1, e_2, e_3$. Since $e_1$ is the same for both frames, the other two pairs of legs are orthogonal transformations of one another:

$$\begin{aligned} e_2' &= \cos\theta e_2 \mp \sin\theta e_3, \\ e_3' &= \sin\theta e_2 \pm \cos\theta e_3. \end{aligned}$$

The angle $\theta$ is at most $\pi/2$, and equal just when the curvature vector as a space curve is perpendicular to the normal line on $S$, which we suppose never happens. Pick $e_3$ so that $\pm = +$. The Serret–Frenet equations are

$$\begin{aligned} \dot{x} &= e_1, \\ \dot{e}_1 &= k e_2', \\ \dot{e'}_2 &= -k e_1 + t e_3', \\ \dot{e'}_3 &= -t e_2'. \end{aligned}$$

Geodesic curvature is $\kappa ds = -\alpha$. The connection forms are

$$\begin{aligned} \begin{pmatrix} -\alpha \\ -a_{11} ds \\ -a_{12} ds \end{pmatrix} &= \begin{pmatrix} \gamma_{21} \\ \gamma_{31} \\ \gamma_{32} \end{pmatrix}, \\ &= \begin{pmatrix} e_2 \dot{e}_1 \\ e_3 \dot{e}_1 \\ e_3 \dot{e}_2 \end{pmatrix} ds, \\ &= \begin{pmatrix} k \cos\theta \\ \pm k \sin\theta \\ e_3 \frac{d}{ds}\left(\cos\theta e_2' + \sin\theta e_3'\right) \end{pmatrix} ds. \end{aligned}$$

**F.5** Expand out the last line to finish the proof.

□

**F.6** Suppose that a surface has a straight line on it, through every point, in every tangent direction through that point. Prove that it is a plane.

## The Gauss map

Consider first a surface $S$ in $\mathbb{E}^3$. On its frame bundle we have $\omega_1, \omega_2$ and $\gamma_{12}$, so that

$$\begin{aligned} d\omega_1 &= -\gamma_{12} \wedge \omega_2, \\ d\omega_2 &= \gamma_{12} \wedge \omega_1, \\ d\gamma_{12} &= K \omega_1 \wedge \omega_2. \end{aligned}$$

where $K$ is the Gauss curvature of the surface. We can find a geometric description of Gauss curvature. Take the unit vector $e_3$, which is uniquely determined at each point of $S$ up to sign, being normal to $S$. Differentiate this vector, and write its derivative in the frame itself:

$$\begin{aligned} de_3 &= (e_1 \cdot de_3)\, e_1 + (e_2 \cdot de_3)\, e_2 + (e_3 \cdot de_3)\, e_3, \\ &= \gamma_{13} e_1 + \gamma_{23} e_2 + \gamma_{33} e_3, \end{aligned}$$

but $\gamma_{33} = 0$ because $\gamma$ is antisymmetric,

$$= (a_{11}\omega_1 + a_{12}\omega_2)\, e_1 + (a_{12}\omega_1 + a_{22}\omega_2)\, e_2.$$

So if we move in the $e_1$ direction, the unit vector $e_3$ moves by $a_{11}e_1 + a_{12}e_2$, while if we move in $e_2$ direction, the unit vector $e_3$ moves by $a_{12}e_1 + a_{22}e_2$. The unit square $e_1 \wedge e_2$ in the tangent plane $T_xS$ maps by $e_3' = de_3$ to the parallelogram

$$e_{3*}(e_1 \wedge e_2) = \left(a_{11}a_{22} - a_{12}^2\right) e_1 \wedge e_2.$$

But the Gauss curvature is $K = a_{11}a_{22} - a_{12}^2$. Therefore the Gauss curvature is the area stretch factor of the Gauss map $e_3$. While $e_3$ is only defined on the frame bundle of $S$, nonetheless it is defined locally on $S$ up to sign; the sign doesn't alter the value of the area stretch factor. Gauss curvature also doesn't depend on orientation.

## The Gauss–Bonnet theorem

Take a vector field $X$ on an oriented surface $S$. Away from the zeroes of $X$, the unit vector field

$$e_1 := \frac{X}{|X|}$$

is defined. Extend it to an oriented orthonormal basis $e_1, e_2$. This choice of orthonormal frame gives a section of ${}^{\ulcorner}S \to S$ over the points where $X \neq 0$.

Take a point $x_0 \in S$ at which $X$ has an isolated zero. In any positively oriented chart near $x_0$, draw a circle $C_0$ around $x_0$, so that $x_0$ is the only zero of $X$ inside the circle. The *index* of $X$ at $x_0$ is the number of times that the frame $e_1, e_2$ turns around the origin, as measured in some choice of orthonormal frame $e_1', e_2'$ defined near $x_0$. Write

$$\begin{pmatrix} e_1 \\ e_2 \end{pmatrix} = \begin{pmatrix} \cos\theta & -\sin\theta \\ \sin\theta & \cos\theta \end{pmatrix} \begin{pmatrix} e_1' \\ e_2' \end{pmatrix}$$

near $C_0$. The angle $\theta$ varies by a multiple of $2\pi$, say $2\pi n_0$, around $C_0$, and the index of $X$ at $x_0$ is $n_0$. Clearly the index doesn't change if we continuously vary $X$ while the surface $S$ moves continuously in space $\mathbb{E}^3$. The sum of indices is the *Euler characteristic* $\chi_X$.

**F.7** Prove that the vector field

$$X = ax\frac{\partial}{\partial x} + by\frac{\partial}{\partial y}$$

on the plane has index 1 if $a, b > 0$ or $a, b < 0$ and index $-1$ if $a$ and $b$ have opposite signs.

**F.8** Calculate out that

$$\alpha = e_2 \cdot de_1 = e_2' \cdot de_1' - d\theta = \alpha' - d\theta.$$

The *Euler characteristic* of a compact oriented surface $S$ is

$$\chi_S := \frac{1}{2\pi}\int_S K\,dA.$$

**Theorem F.2.** *Suppose that $X$ is a vector field on a compact oriented surface $S$ with only isolated zeroes. Then $\chi_X = \chi_S$.*

*Proof.* Take a set $D_k$ diffeomorphic to a closed disk in some chart around each zero $x_k$ of $X$ and let $C_k = \partial D_k$. The framing $e_1, e_2$ lifts

$$S' := S - \bigcup_k D_k$$

to a compact surface with boundary in $\ulcorner S$. Integrate:

$$\begin{aligned}
\int_{S'} K\,dA &= \int_{S'} d\alpha, \\
&= \int_{\partial S'} \alpha, \\
&= -\sum_k \int_{C_k} \alpha, \\
&= -\sum_k \int_{C_k} e_1 \cdot de_2, \\
&= -\sum_k \int_{C_k} e_1' \cdot de_2' + \sum_k \int_{C_k} d\theta, \\
&= -\sum_k \int_{C_k} K\,dA + \sum_k 2\pi \operatorname{index}_{x_k} X,
\end{aligned}$$

so that

$$\int_S K\,dA = 2\pi \sum_k \operatorname{index}_{x_k} X.$$

□

## Degree

Suppose that $f\colon P \to Q$ is a continuously differentiable map between two oriented manifolds of equal dimension, perhaps with corners. The *sign* of $f$ at a regular point $p \in P$ is 1 if $f$ is orientation preserving near $p$, and $-1$ is $f$ is orientation reversing near $p$. Suppose that $q_0 \in Q$ is a regular value and that there are only finitely points $p \in P$ mapped to $q_0$ by $f$. The *degree* of $f$ above $q_0$ is the sum of signs at these points $p$:

$$\deg_{q_0} f = \sum_{f(p)=q_0} (\text{sign of } f \text{ at } p)\,.$$

The following two theorems are proven in [27].

**Theorem F.3.** *Consider a proper continuously differentiable map $f\colon P \to Q$ between oriented manifolds (perhaps with corners). Suppose that $Q$ is connected. The degree of $f$ is the same at all regular values.*

**Theorem F.4.** *Suppose that $f\colon P \to Q$ is a proper continuously differentiable map between two oriented manifolds (perhaps with corners) of the same dimension, say $n$, and that $Q$ is connected. For any $n$-form $\omega$ on $Q$, which vanishes outside some compact set,*

$$\int_P f^*\omega = \deg f \int_Q \omega.$$

## The Gauss map

For any oriented surface $S$ in $\mathbb{E}^3$, the *Gauss map* $e_3$ taking each point $x \in S$ its unit normal vector, a point of the unit sphere, has degree

$$\deg e_3 = \frac{\int K\,dA}{4\pi} = \frac{\chi_S}{2}.$$

**Theorem F.5.** *On any compact oriented surface $S$ in $\mathbb{E}^3$ with Gauss map $e_3\colon S \to \mathbb{S}^2$, $2\deg e_3 = \chi_X$, for any vector field $X$ on $S$ with only isolated zeroes.*

**Theorem F.6.** *Every compact oriented surface $S$, twice continuously differentiable, has a continuously differentiable vector field $X$ with only finitely many zeroes, each of index 1 or $-1$.*

*Proof.* The Gauss map is onto: take any unit vector $n$, and look at the points $x \in S$ on which the linear function $f(x) := \langle n, x\rangle$ has its maximum value. By Sard's theorem, almost every point $n$ in the unit sphere is a regular value of the Gauss map. For such a value of $n$, the differential of the Gauss map has full rank. Its differential has determinant given by the Gauss curvature, so the Gauss curvature is nonzero at every point $x \in S$ for which $\nu(x) = n$.

The function $f$ restricted to $S$ has critical points at the points $x$ where $n$ is perpendicular to $T_x S$. Again by Sard's theorem, we can arrange that these points are regular values of the Gauss map, i.e. have nonzero Gauss curvature. In any orthonormal framing $e_1, e_2, e_3$, the differential of $f$ is

$$e_i \lrcorner df = e_i \lrcorner (n \cdot x) = n \cdot e_i,$$

for $i = 1, 2$. The second derivative matrix of $f$ at a critical point is

$$\begin{aligned} f''(e_i, e_j) &= e_i \lrcorner d\left(n \cdot e_j\right), \\ &= e_i \lrcorner n \cdot de_j, \\ &= e_i \lrcorner (\pm e_3) \cdot de_j, \\ &= \pm e_i \lrcorner \gamma_{3j}, \\ &= \mp a_{ji}. \end{aligned}$$

So the function $f$ has $f''$ a rank 2 quadratic form on each tangent space where $f$ is critical. In particular, if we write $df = f_1 \omega_1 + f_2 \omega_2$, then the vector field

$$X = \nabla f = f_1 e_1 + f_2 e_2$$

has only finitely many zeroes. Since we can vary our surface continuously without changing indices, we can see that the index doesn't change as we vary $S$ and $f$ so that $f$ remains critical at some point $x_0$ with $f''$ nondegenerate at $x_0$. In particular, we can arrange that $S$ is flat near $x_0$ with $f$ a quadratic function instead, without changing the index. But we saw above that the index of $X = \nabla f$ at $x_0$ is then 1 if $\det f''(x_0) > 0$ and $-1$ if $\det f''(x_0) < 0$. □

**Corollary F.7.** *Diffeomorphic surfaces have equal Euler characteristic.*

Appendix G

# Geodesics on surfaces

Geodesics are locally shortest paths on a surface. We prove that they exist and that they are smooth.

## Geodesic flow

Each point $(x^o, e^o) \in \ulcorner\mathbb{E}^3$ moves along a path $x(t) = x^o + te^o_1, e(t) = e^o$. In other words, don't rotate the frame, but slide the base point $x$ along the direction of the first leg of the frame.

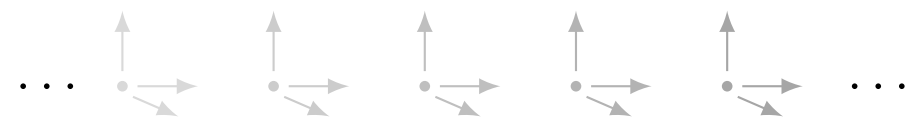

Define a vector field $w(x, e) := (e_1, 0)$ on $\ulcorner\mathbb{E}^3$; these paths are its flow lines.

Imitate this vector field on the frame bundle of a surface $S$ in $\mathbb{E}^3$. Write points of $\ulcorner S \times \mathbb{E}^2$ as $(x, e, y)$ for $y \in \mathbb{E}^2$. The *geodesic spray* of $S$ is the vector field $w$ on $\ulcorner S \times \mathbb{E}^2$ given by

$$w \lrcorner (\omega, \alpha, dy) = (y, 0, 0)$$

This expression does *not* involve $\gamma_{13}$ or $\gamma_{23}$, since these are not linearly independent 1-forms on $\ulcorner S$. The flow of $w$ is the *geodesic flow*. The vector field $w$ moves every frame $(x, e)$, moving $x$ in the direction of $\sum y_i e_i$, and "not rotating" the frame $e$, since $w \lrcorner \gamma_{12} = 0$. Intuitively, $w$ twists frames only in the normal directions (in the $\gamma_{13}, \gamma_{23}$ directions), not tangentially (the $\alpha = \gamma_{12}$ direction), while keeping $e_1, e_2$ tangent to the surface $S$. Map $(x, e) \in \ulcorner S \to (x, e, 1) \in \ulcorner S \times \mathbb{E}^2$, so that $\ulcorner S$ is identified with a submanifold of $\ulcorner S \times \mathbb{E}^2$, making $\ulcorner S$ invariant under the geodesic flow. Write $a^\top$ for the transpose of a matrix $a$. Let $g \in \mathrm{O}(2) \times \pm 1$ act on $\ulcorner S \times \mathbb{R}^p$ by

$$(x, e, y)g = (x, eg, a^\top y).$$

if

$$g = \begin{pmatrix} a & 0 \\ 0 & b \end{pmatrix}.$$

**Lemma G.1.** *On $\ulcorner S \times \mathbb{E}^2$, $r_g^* (\omega, \gamma_{12}, dy) = \left(g^\top \omega, \gamma_{12}, g^\top\, dy\right)$ while $r_{g*} w = w$.*

*Proof.* The reader can check how the 1-forms pull back. The equations defining $w$ imply invariance under $\mathrm{O}(2) \times \pm 1$. $\square$

Let $p\colon (x,e,y) \in {}^{\ulcorner}S \times \mathbb{E}^2 \mapsto (x,v) \in TS$ by $v = \sum y_i e_i$.

**Proposition G.2.** *There is a uniquely determined vector field $w$ on $TS$, also called the* geodesic spray*, so that $p_* w = w$.*

*Proof.* Calculate that at a point $z = (x,e,y)$, for any $g \in \mathrm{O}(2) \times \pm 1$,

$$
\begin{aligned}
p'(zg)w(zg) &= p'(zg)r_g'(z)w(z), \\
&= \left(p \circ r_g\right)'(z)w(z), \\
&= p'(z)w(z).
\end{aligned}
$$

The reader can check that $\mathrm{O}(2) \times \pm 1$ acts transitively on the points $z = (x,e,y)$ lying above a given point $(x,v) = p(z)$. Therefore the vector $p'(z)w(z)$ is independent of the point $z$, depending only on $(x,v) = p(z)$. □

A *geodesic* is a curve in $S$ which is a projection to $S$ of a flow line of $w$ in $TS$.

**G.1** Take a surface $S$ in $\mathbb{E}^3$. Prove that any geodesic $x(t)$ of $S$, as a curve in $\mathbb{E}^3$, has geodesic curvature $a(\dot{x},\dot{x})$. Prove that $a = 0$ just when all geodesics of $S$ are geodesics of $\mathbb{E}^3$, i.e. straight lines, so that $S$ is locally an open subset of a plane.

**Lemma G.3.** *A smooth curve $C$ on a surface $S$ in $\mathbb{E}^3$ has vanishing geodesic curvature just when it is a geodesic.*

*Proof.* The geodesic spray vector field $w$ points along the tangent line cut out by the equations $\omega_2 = 0$ and $\alpha = 0$ and has $w \lrcorner \omega_1 = 1$. The tangent spaces of ${}^{\ulcorner}S_C$ are cut out by the equations $\omega_2 = 0$ and $\alpha = -\kappa_2 \omega_1$. The geodesic spray $w$ lies tangent to ${}^{\ulcorner}S_C$ just when the geodesic curvature of $C$ vanishes, which therefore occurs just when the projection of the flow lines of $w$ lie along $C$, i.e. $C$ is a geodesic. □

## First order deformation

Take a plane curve $(x(s), y(s)) \in \mathbb{R}^2$. Draw a vector field along this curve: a vector $(u(s), v(s))$ at each point $(x(s), y(s))$. The smooth homotopy $s,t \mapsto (x(s) + tu(s), y(s) + tv(s))$ pushes out each point $(x(s), y(s))$ in the direction of our vector field. In the same way, if we have any compactly supported vector field $v(s)$ along a curve $x(s)$ in an open set in $\mathbb{E}^3$, we can let $x_t(s) = x(s) + tv(s)$, and then $x_t(s)$ will be defined for some interval of time $t$. Alternatively, we could create many other smooth homotopies $x_t(s)$ with $x_0(s) = x(s)$ and $\left.\frac{\partial}{\partial t}\right|_{t=0} x_t(s) = v(s)$; we say that $v(s)$ is the *first order deformation* of $x_t(s)$. A *vector field $v$ along* a curve $x(s)$ on a surface $S$ is a smooth choice of vector $v(s) \in T_{x(s)}S$ at each point. As long as the vector field $v(s)$ has compact support in a single chart, we can use the chart to identify an open set of the surface with an open set in the plane, the curve with a plane curve, and each vector field along the curve with a vector field along a plane curve, so create a homotopy with given first order deformation.

**Lemma G.4.** *Take a curve $C$ (perhaps with boundary) on a surface $S$ in $\mathbb{E}^3$, and a compactly supported vector field $v$ along $C$. Then $v$ is a sum $v = v_1 + v_2 + \dots + v_k$ of vector fields along $C$, so that each vector field $v_j$ is the first order deformation of a smooth homotopy of $C$.*

*Proof.* Cover $C$ by charts and take a partition of unity. $\square$

## Smooth homotopies and variation of length

**Lemma G.5.** *On a surface $S$ in $\mathbb{E}^3$, take a vector field $v(s)$ along a unit-speed curve $x(s)$ on $S$. Suppose that $v$ vanishes at $s = s_0$ and at $s = s_1$. Let $\kappa(s)$ be the geodesic curvature of $x(s)$. Suppose that there is a smooth homotopy $x_t(s)$ with first order deformation $v(s)$ along $x_0(s) = x(s)$, so $x_t(s_0) = x(s_0)$ and $x_t(s_1) = x(s_1)$ for all values of $t$. Then*

$$\left.\frac{d}{dt}\right|_{t=0} \operatorname{length} x_t = \int \kappa \cdot v \, ds$$

*Proof.* Our homotopy is a map $x \colon X \to S$ for some domain $X \subset \mathbb{R}^2$ in the $s, t$-plane. Since $x(s)$ has nowhere zero velocity, by continuity the curve $s \mapsto x_t(s)$ has nonzero velocity for $t$ near 0. If we restrict $t$ to a small enough interval, we can ensure that $\frac{\partial x}{\partial s} \neq 0$ for all $(s, t)$. An *adapted frame* for the homotopy $x(s, t)$ is a frame $e$ so that $\frac{\partial x}{\partial s}$ is a positive multiple of $e_1$, and therefore is perpendicular to $e_2, e_3$. So $e_1 = e_1(s, t)$ is a field of vectors, pointing in the direction of each curve in the homotopy. Since $e_1 = e_1(s, t)$ is actually now defined on the $st$-plane, $\omega_1 = e_1 \, dx$ becomes well defined as a form on the $st$-plane as well. At each point $(s, t, e)$ of the adapted frame bundle, we can write the first deformation $v = \frac{\partial x}{\partial t}$ of each of the curves $x_t$ as $\frac{\partial x}{\partial t} = v_1 e_1 + v_2 e_2$. Denote the set of adapted frames of the homotopy by $\ulcorner S_x$; we leave the reader to check that this set is a manifold. On $\ulcorner S_x$, we have $\omega_2 = v_2 \, dt$. Taking exterior derivative of both sides of this equation,

$$\begin{aligned} 0 &= d\omega_2 - dv_2 \wedge dt, \\ &= \alpha \wedge \omega_1 - dv_2 \wedge dt. \end{aligned}$$

By Cartan's lemma (lemma 1.2 on page 6), there are unique smooth functions $\kappa_2, a_2, b_2$ on $\ulcorner S_x$ so that

$$\begin{aligned} -\alpha &= \kappa_2 \omega_1 + a_2 \, dt, \\ dv_2 &= a_2 \, \omega_1 + b_2 \, dt. \end{aligned}$$

Checking dimensions, the tangent spaces of $\ulcorner S_x$ are cut out by the equations $\omega_2 = v_2 \, dt$ and $-\alpha = \kappa_2 \omega_1 + a_2 \, dt$.

Let $\ell(t) = \operatorname{length} x_t$ and let $\pi(s,t) = t$. Write $\pi_*$ for integration over the fibers of $\pi$.

$$\begin{aligned}\ell(t) &= \int |\dot{x}_t|\, ds, \\ &= \int x^* \omega_1, \\ &= \pi_* x^* \omega_1.\end{aligned}$$

Therefore

$$\begin{aligned}d\ell &= \pi_* x^* d\omega_1, \\ &= \pi_* x^* \left(-\alpha \wedge \omega_2\right), \\ &= \sum_i \pi_* x^* \left(\kappa_2\, \omega_1 \wedge v_2\, dt\right), \\ &= \pi_* x^* \left(\kappa \cdot v \omega_1 \wedge dt\right), \\ &= \left(\int \kappa \cdot v \omega_1\right) dt.\end{aligned}$$

Along $t = 0$, since our curve is unit-speed, $\omega_1 = ds$. □

**Lemma G.6.** *If a smooth curve locally minimises length between its points, then it is a geodesic.*

*Proof.* If $C$ is a smooth curve of minimal length between two points, then $0 = \int \kappa \cdot v\, ds$ for any first order deformation $v$. Any vector field $v$ along $C$ is a sum of first order deformations, so again $0 = \int \kappa \cdot v\, ds$. In particular, we can take $v$ a positive multiple of $\kappa$ near some point of $C$, and zero away from there. □

## The exponential map

Write the flow of $w$ on the tangent bundle as $(x, v) \mapsto e^{tw}(x, v)$, defined for all $(x, v)$ with $v$ close enough to zero. The *projection map* is $\pi_S \colon (x, v) \in TS \mapsto x \in S$. For each point $x \in S$, the *exponential map* is $\exp_x \colon v \in \text{open} \subset T_x S \mapsto \pi_S(e^w(x, v)) \in S$. The exponential map might only be defined near $v = 0$.

**G.2** Prove that $\pi_S'(x, v) w(x, v) = v$.

**Lemma G.7.** *The exponential map on a surface $S$ in $\mathbb{E}^3$ has derivative* $\exp'(0) = \operatorname{id} \colon T_0 T_x S = T_x S \to T_x S$.

*Proof.* Going back to the definition of $w$, for any constant $c \in \mathbb{R}$, $cw(x, e, y) = w(x, e, cy)$. Therefore the flow satisfies $e^{tw}(x, e, cy) = e^{ctw}(x, e, y)$ for any con-stant $c \in \mathbb{R}$. In particular,

$$\begin{aligned} w(x,e,y) &= \left.\frac{d}{dt}\right|_{t=0} e^{tw}(x,e,y), \\ &= \left.\frac{d}{dt}\right|_{t=0} e^{w}(x,e,ty). \end{aligned}$$

Taking the projection down to $TS$,

$$v = \left.\frac{d}{dt}\right|_{t=0} \pi_S \circ e^w(x,tv).$$

This gives us

$$\begin{aligned} \exp_x'(0)v &= \left.\frac{d}{dt}\right|_{t=0} \exp_x(tv), \\ &= \left.\frac{d}{dt}\right|_{t=0} \pi_S \circ e^w(x,tv), \\ &= \left.\frac{d}{dt}\right|_{t=0} \pi_S \circ e^{tw}(x,v), \\ &= \pi_S'(x,v)w(x,v), \\ &= v. \end{aligned}$$

□

We need locally uniform control on geodesics. The *injectivity radius* of a surface $S$ in $\mathbb{E}^3$ at a point $x \in S$ is the supremum of numbers $r > 0$ so that $\exp_x$ is a diffeomorphism on the ball of radius $r$ in $T_xS$.

**Lemma G.8.** *For any compact subset $K \subset S$ of a surface $S$ in $\mathbb{E}^3$ there is a lower bound $r = r_K > 0$ on the injectivity radius of all points $x \in K$. If $x \in K$ and $y \in S$ is within distance $r$ of $x$ then there is a unique geodesic of length less than $r$ from $x$ to $y$.*

*Proof.* Consider the map $f\colon TS \to S \times S$ given by $f(x,v) = (x, \exp_x v)$. As above, $f'(x,0)(u,v) = (u, u+v)$, so $f$ is a diffeomorphism of an open neighborhood of $U_{TS} \subset TS$ of $\{ (x,0) \mid x \in S \}$ to a open neighborhood $U_{S\times S}$ of $\{ (x,x) \mid x \in S \} \subset S \times S$. Take any compact set $K \subset S$. Consider a sequence of points $x_j, y_j$ with $x_j \in K$ and $y_j \in S$ with distances between $x_j$ and $y_j$ getting small. By compactness, after replacing with a subsequence, our sequence converges to a point $(x,x)$, which lies in $U_{S\times S}$. But $U_{S\times S}$ is open, so all but finitely many pairs $(x_j, y_j)$ belong to $U_{S\times S}$. Suppose that there are arbitrarily close points $x_j, y_j$ with $x_j \in K$ with no unique shortest geodesic. Then $(x_j, y_j)$ is not in the image of $f$, i.e. not in $U_{S\times S}$, a contradiction. Similarly, we cannot find arbitrarily close points $x_j, y_j$ with $x_j \in K$ where $f^{-1}$ is not defined. But $f^{-1}(x,y) = \exp_x^{-1} y$, so $\exp_x y$ is defined for near enough points $x$ and $y$, and is a diffeomorphism because $f$ is a diffeomorphism. □

## Geodesic polar coordinates

**G.3** Suppose that $\eta_1, \eta_2$ are an orthonormal coframing of an oriented surface $S$ in $\mathbb{E}^3$, i.e. $\eta_1, \eta_2$ yield an orthonormal basis of each cotangent space. Let $u_1, u_2$ be the dual vector fields, also orthonormal, and take $u_3$ to be the positively oriented normal vector field. Define a map $f: x \in S \mapsto (x, u_1(x), u_2(x), u_3(x)) \in {}^{\ulcorner}S$. Prove that $f^*\omega_i = \eta_i$.

For this reason, we prefer to just denote any choice of orthonormal coframing $\eta_1, \eta_2$ as $\omega_1, \omega_2$, and the associated 1-form $f^*\gamma_{12}$ as just $\gamma_{12}$, which we say is the connection 1-form of this coframing.

Take the exponential map $\exp_{x_0}$ and use it to identity $T_{x_0}S$ with an open set around $x_0$ in $S$. Make an isometric linear isomorphism $T_{x_0}S \cong \mathbb{E}^2$. Take polar coordinates on $\mathbb{E}^2$. The direction $\partial_r$ is, by definition of the exponential map, the direction of the geodesics moving at unit speed away from $x_0$, so is a unit vector field. Take $\omega_1, \omega_2$ an orthonormal coframing so that $\omega_1 = 1$ and $\omega_2 = 0$ on $\partial_r$. So

$$\begin{aligned} \omega_1 &= dr + a(r,\theta)\, d\theta, \\ \omega_2 &= h(r,\theta)\, d\theta. \end{aligned}$$

Differentiate, denoting $\frac{\partial h}{\partial r}$ as $h_r$, and so on, to find that

$$\alpha = -\frac{a_r \omega_1 + h_r \omega_2}{h}.$$

The curves $\omega_2 = 0$ have geodesic curvature $a_r/h$, which vanishes because they are geodesics, so $a = a(\theta)$.

**G.4*** Prove that $a(\theta) = 0$ and $h(r,\theta) \to 0$ and $h_r(0,\theta) \to 1$ as $r \to 0$. Differentiate $\alpha$ to find the Gauss curvature $K$.

Summing up:

**Lemma G.9.** *Any surface, in geodesic normal coordinates, has orthonormal coframing, connection form and curvature:*

$$\begin{aligned} \omega &= dr + i h(r,\theta) d\theta, \\ \alpha &= -h_r d\theta, \\ K &= -\frac{h_{rr}}{h}, \\ h(r,\theta) &\to 0 \text{ as } r \to 0, \\ h_r(0,\theta) &\to 1 \text{ as } r \to 0, \\ h(r,\theta + 2\pi) &= h(r,\theta). \end{aligned}$$

The equation $h_{rr} + Kh = 0$ is linear in $h$, and uniquely determines $h$ subject only to $h = 0$ and $h_r = 1$ at $r = 0$. By uniqueness $h(r, \theta + 2\pi) = h(r,\theta)$.

For example, if $K = 1$ is constant along a geodesic ray, then $h = \sin r$ along that ray. (Note that $h$ cannot be zero as $\omega_2 = h\,dr$, so the exponential map is not a diffeomorphism and geodesic polar coordinates are not defined for $r \ge \pi$. Similarly, if $K = 0$ along a geodesic ray, then $h = r$ along that ray, just as on a flat plane. Similarly, if $K = -1$ along a geodesic ray, then $h = \sinh r$ along that ray.

A *length preserving diffeomorphism* of surfaces is a diffeomorphism preserving the lengths of all curves.

**Theorem G.10.** *Any surface of constant Gauss curvature $K = K_0$ is locally length preserving diffeomorphic to any other.*

*Proof.* In our geodesic normal coordinates,

$$\omega_1 = dr, \omega_2 = h(r)\,d\theta$$

where $k_0 := \sqrt{|K_0|}$ and

$$h(r) = \begin{cases} \frac{\sin k_0 r}{k_0}, & \text{if } K_0 > 0, \\ r, & \text{if } K_0 = 0, \\ \frac{\sinh k_0 r}{k_0}, & \text{if } K_0 < 0. \end{cases}$$

Matching up the values of $r, \theta$ on our two surfaces then matches up $\omega_1, \omega_2$ and so matches up lengths of any curves, as the length of a curve $(r(t), \theta(t))$ is $\int \sqrt{\omega_1^2 + \omega_2^2}$. □

More generally, the function $h(r, \theta)$ is constant in $\theta$ just when the surface admits "rotations", i.e. rotation of $\theta$, as symmetries, fixing the point about which the geodesic coordinates are computed. One easily checks that the surface is then length preserving diffeomorphic near that point to a surface of revolution, with that point fixed under the revolution: an intrinsic description of the shape of a surface of revolution. The length of the circle of radius $r$ around the fixed point is $\ell(r) = \int \omega_2 = 2\pi h(r)$, and the curvature is $K = -\ell''(r)/\ell(r)$: two surfaces of revolution are length preserving diffeomorphic near fixed points just when their curvatures are equal as functions of distance from fixed point.

## Geodesics are locally minimal

**Theorem G.11.** *Geodesics are locally minimal, and every locally minimal curve is a geodesic.*

*Proof.* In geodesic polar coordinates, $r$ increases at unit rate as we move along the radial curves given by constant $\theta$, and at less than unit rate as we move at unit speed in any other direction. Note that these are geodesics, by definition of the coordinates. Hence the length of a piece of curve increases by at least the increase in $r$ across its end points, and exactly by this amount precisely for the radial curves. So the shortest curves of given increase in $r$, i.e. shortest paths between circles of constant $r$, are these radial curves. For any path, any two points on it sufficiently close together are each in the injectivity radius of the other, so the geodesic polar coordinates are defined, and the unique shortest path between those points is the unique geodesic between them in those coordinates. □

The *length metric distance* between two points $p, q \in S$ of a surface $S$ in $\mathbb{E}^3$ is the infimum of lengths of paths from $p$ to $q$. As above, the geodesics between sufficiently close points have minimal length, so their length achieves the length metric distance.

**G.5** Prove that if the length metric distance between points goes to zero, then the distance in $\mathbb{E}^3$ also goes to zero.

**Corollary G.12.** *Every invertible map $f: P \to Q$ of surfaces $P, Q$ in $\mathbb{E}^3$ which preserves lengths of curves is a diffeomorphism.*

*Proof.* The map takes geodesics to geodesics, since they are precisely the locally shortest paths. The map identifies each point $p_0 \in P$ with a point $q_0 \in Q$, and the circle around $p_0$ of a given radius with the corresponding circle around $q_0$ of the same radius. The lengths of geodesic segments are preserved, so if a point moves at unit rate along a geodesic on $P$ then its image in $Q$ does the same. Indeed, the same is true for a point moving along any curve (for example, along a circle), because the path is approximated by piecewise geodesic paths. Therefore at every point near $p_0$, except at $p_0$ itself, $f_* e_1 = e_1$ and $f_* e_2 = e_2$. Therefore $f$ is continuously differentiable at all those points. By lemma G.8 on page 137, every point $p \in P$ is near enough but not equal to some such point $p_0$ to use this argument to prove that $f$ is continuously differentiable at $p$. □

**Corollary G.13.** *On any connected surface $S$ in $\mathbb{E}^3$ the following are equivalent:*

a. *The surface $S$ is complete in the length metric, i.e. length metric Cauchy sequences converge.*

b. *The surface $S$ is proper in the length metric, i.e. all length metric balls are compact.*

c. *The exponential map on $S$ has flow defined for all time through all points.*

d. *Every geodesic on $S$ defined on an open interval extends uniquely to a geodesic defined on a closed interval.*

e. *There is a point $x_0 \in S$ so that every shortest path $x: [a, b) \to S$ with $x(a) = x_0$ extends uniquely to a continuous path $x: [a, b] \to S$.*

*In particular, if $S$ is a closed subset of $\mathbb{E}^3$, then all of these hold.*

*Proof.* Every geodesic on $S$, parameterised by arc length, is defined on a maximal open interval, since the geodesics are the solutions of an ordinary differential equation.

We recall some facts about metric spaces [14]. A *length space* is a metric space in which the distance between points is the infimum of lengths of paths between them. A *minimal geodesic* in a length space $X$ is a path $x(t)$ on which the distance between points $x(t_0), x(t_1)$ is $|t_1 - t_0|$. A *geodesic* is a path which is locally a minimal geodesic. By the Hopf–Rinow theorem for metric spaces [14, 27], any locally compact length space is complete just when any geodesic defined on an open interval extends to be defined on a closed interval.

In any complete surface, geodesics then extend further to be defined on larger open intervals, since they continue to solve their differential equation. So if $S$ is complete, geodesics are defined for all time. But if geodesics are defined for all time, they extend, so by Hopf–Rinow $S$ is complete. □

**G.6** Prove that, on any complete surface, vector fields of bounded norm are complete.

**Corollary G.14.** *Any two local isometries $f, g\colon P \to Q$ of the length metric between two connected surfaces are identical just when there is a point $p_0 \in P$ at which both $f(p_0) = g(p_0)$ and $f'(p_0) = g'(p_0)$.*

*Proof.* Each vector $v \in T_{p_0}P$ is taken to the same vector $f'(p_0)v = g'(p_0)v$, so exponentiates to the same point of $Q$, if $v$ is small enough. Because $f$ and $g$ are local isometries, they locally identify the geodesics of $P$ and $Q$, and the exponential maps, so

$$f \circ \exp_{p_0} = \exp_{q_0} \circ f'(p_0) = \exp_{q_0} \circ g'(p_0) = g \circ \exp_{p_0},$$

and since the exponential map is a local diffeomorphism, $f = g$ near $p_0$. So the set of points where $f = g$ is both open and closed. □

## Theorema Egregium

A *circle* on a surface is the set of points at given length metric distance from a given point of the surface. We would like to see that the Gauss curvature controls the speed at which circles grow with radius. Along each geodesic coming out of a given point $o$, we erect a frame with $e_1$ tangent to the geodesic and $e_2$ normal, a point of the frame bundle associated to each point near (but not equal to) $o$. Motion along $e_1$ takes us away from $o$ at unit speed, while motion along $e_2$ spins us around the circles around $o$ at unit speed.

$$h(r, \theta) = r - \frac{K(o)}{2} r^2 + o(r)^2 .$$

The length of the circular arc of radius $r$ about $o$ between two angles $\theta_0, \theta_1$ is

$$\int \omega_2 = \int_{\theta_0}^{\theta_1} h(r, \theta)\, d\theta = |\theta_1 - \theta_0| r \left(1 - \frac{K(o)\, r}{2}\right) + o(r)^2 .$$

**Theorem G.15** (Theorema Egregium)**.** *Any isometry of length metric between analytic surfaces is analytic and preserves Gauss curvature.*

*Proof.* Distances are preserved, so circles (points of given distance) are preserved. Lengths are preserved, so lengths of circles are preserved, so Gauss curvature is preserved. Geodesics coming out of a given point are mapped to geodesics out of the corresponding point. The geodesics coming out of $o$ have constant values of $\theta$, and these are preserved. So $\theta$ changes under isometry only as a function of $\theta$, independent of $r$. Since $r$ is preserved, lengths of circular arcs are preserved, so the absolute differences of values of $\theta$ are preserved. So all isometries have the form $r, \theta \mapsto r, \pm\theta + c$, for some constant $c$. Arrange by choice of coordinates that the rays $\theta = 0$ are matched, and the direction of increasing $\theta$ is matched: $r, \theta \mapsto r, \theta$. $\square$

**G.7** Suppose that $S$ is a surface in $\mathbb{E}^3$ diffeomorphic to a triangle in the plane, say with interior angles $\alpha, \beta, \gamma$, and with geodesic sides (a "geodesic triangle"). Prove that

$$\alpha + \beta + \gamma = \pi + \int_S K\, dA.$$

In particular, negatively (positively) curved geodesic triangles have larger (smaller) angle sums than Euclidean triangles.

## The Hessian

We will prove:

**Theorem G.16.** *Take an open set on a surface in $\mathbb{E}^3$, with smooth boundary of nonnegative geodesic curvature. Whenever a geodesic enters the interior of that open set, it is not tangent to the boundary.*

Take a function $f \colon S \to \mathbb{R}$. Pulling $f$ back to the frame bundle of $S$, because $df$ is semibasic, $df = f_1\omega_1 + f_2\omega_2$ for some functions $f_1, f_2$. Differentiate to find

$$d\begin{pmatrix} f_1 \\ f_2 \end{pmatrix} + \begin{pmatrix} 0 & \alpha \\ -\alpha & 0 \end{pmatrix}\begin{pmatrix} f_1 \\ f_2 \end{pmatrix} = \begin{pmatrix} f_{11} & f_{12} \\ f_{21} & f_{22} \end{pmatrix}\begin{pmatrix} \omega_1 \\ \omega_2 \end{pmatrix}$$

for unique functions $f_{ij} = f_{ji}$. Check that the operator

$$D^2 f := f_{ij}\omega_i\omega_j$$

is a symmetric bilinear form on tangent spaces of $S$, the *Hessian* of $f$. The differential $df$ restricts to any curve to become

$$df = \frac{df}{ds} ds = f_1\omega_1$$

measuring the rate at which $f$ increases along the curve, and the Hessian becomes

$$\frac{d^2 f}{ds^2} = \frac{df_1}{ds} = f_{11},$$

the acceleration of $f$ along the curve. Any curve on which the Hessian is positive has $f$ accelerating, so at any point of that curve, if $df$ vanishes at that point of the curve, i.e. the curve is tangent to a level set of $f$, then $f$ increases on either side of that point, i.e. the curve stays outside of the sublevel set.

If $df \neq 0$, we can pick out a particular choice of frame: make $e_1$ tangent to level sets of $f$. So $df = f_2\omega_2$, i.e. $f_1 = 0$, on this subbundle of the frame bundle, adapted to $f$. We can further orient the surface by requiring $f_2 > 0$. Our equations simplify to

$$\begin{aligned} f_2\alpha &= f_{11}\omega_1 + f_{12}\omega_2, \\ df_2 &= f_{21}\omega_1 + f_{22}\omega_2. \end{aligned}$$

On any one level set $\omega_2 = 0$ so

$$\alpha = -\kappa_2\omega_1 = -\frac{f_{11}}{f_2},$$

so the geodesic curvature of the level set is

$$\kappa_2 = \frac{f_{11}}{f_2},$$

so the sign of the geodesic curvature of each level set is the sign of the acceleration of $f$ along that level set. In particular, the level sets are geodesic just when the Hessian vanishes along all of them.

Given any curve $C$, we can locally make a function $f$ with that curve as level set. Take two curves $C_0$ and $C_1$ tangent at a point. Suppose they are both oriented, and lie on an oriented surface. Write each as a level set of functions $f_0$ and $f_1$, with differentials pointing on the positive side of the oriented tangent line. At the point of tangency, both functions have vanishing differential on both curves. Clearly $f_0$ has larger Hessian on that tangent line just when $C_0$ has larger curvature. So then $f_0 - f_1$ has positive Hessian on that tangent line, and so $f_0 - f_1$ grows nearby on $C_1$. But $f_1$ is constant on $C_1$. So $f_0$ grows on $C_1$ at all points nearby. So $C_1$ has just one point of tangency with $C_0$ near that point, and everywhere else nearby leaves the sublevel set of $f_0$. In particular, if $C_1$ is a geodesic, $f_1$ has vanishing Hessian along $C_1$, proving theorem G.16 on the facing page.

**G.8** Suppose that $M$ is a complete surface in 3-dimensional Euclidean space, and $f$ is a smooth proper function on $M$ whose only critical point is a single local minimum. Suppose that all level sets of $f$ are positively curved. Prove that every geodesic on $M$ enters each compact set only finitely many times.

## Abstract surfaces

A *surface* is a manifold of dimension 2. A *Riemannian metric* on a surface is a choice of inner product in each tangent space, so that inner products of locally defined analytic vector fields are analytic functions. For simplicity of notation, write inner products as dot products. An *orthonormal frame* $(x, e)$ is a choice of point $x$ of $S$ and a choice of orthonormal basis $e_1, e_2$ of the tangent plane $T_x S$. The *tangent frame bundle* of $S$ is the collection of orthonormal frames of tangent spaces of $S$. Note that for a surface $S$ in $\mathbb{E}^3$, this is *not* the same as the frame bundle we constructed previously in appendix D on page 101, which we could refer to as the *bundle of ambient adapted frames*. The tangent frame bundle of a surface in $\mathbb{E}^3$ is the quotient of the ambient adapted frame bundle by ignoring $e_3$, giving a $2-1$ covering map. Nonetheless, we will use the same notation $\ulcorner S$ for the tangent frame bundle.

Take any two linearly independent vector fields on an open subset $U \subseteq S$ of our surface. Apply the Gram–Schmidt orthogonalization procedure to them to construct orthonormal vector fields, say $e'_1, e'_2$, defined on $U$. Over $U$, any orthonormal frame $e_1, e_2$ is uniquely expressed as

$$
\begin{aligned}
e_1 &= \cos\theta\, e'_1 \pm \sin\theta\, e'_2, \\
e_2 &= -\sin\theta\, e'_1 \pm \cos\theta\, e'_2.
\end{aligned}
$$

In this way, we identify $\ulcorner U \cong U \times (S^1 \cup S^1)$, so the tangent frame bundle becomes a manifold. Let $p \colon (x, e) \in \ulcorner S \mapsto x \in S$. The *soldering forms* $\omega_1, \omega_2$ on $\ulcorner S$ are defined by

$$
p'(m, e_1, e_2)v = \omega_1(v)e_1 + \omega_2(v)e_2.
$$

If the dual 1-forms to $e'_1, e'_2$ on $U$ are $\omega'_1, \omega'_2$ then

$$
\begin{pmatrix} \omega_1 \\ \omega_2 \end{pmatrix} = \begin{pmatrix} \cos\theta & \pm\sin\theta \\ -\sin\theta & \pm\cos\theta \end{pmatrix} \begin{pmatrix} \omega'_1 \\ \omega'_2 \end{pmatrix}.
$$

Suppose henceforth that the surface $S$ is oriented. We pick out only the positively oriented frames, to get rid of $\pm$ so (using the same notation to also denote the bundle of positively oriented orthonormal frames) $\ulcorner U \cong U \times S^1$. Let $\omega := \omega_1 + i\omega_2$, and write this as

$$
\omega = e^{i\theta}\omega'.
$$

The area form $dA := \omega'_1 \wedge \omega'_2$ pulls back to $dA = \omega_1 \wedge \omega_2$, independent of the choice of $e'_1, e'_2$, so globally defined on $S$. Since $dA$ is a nonzero area form, $d\omega'_1, d\omega'_2$ are multiples of it, say

$$
\begin{aligned}
d\omega'_1 &= -a_1 dA, \\
d\omega'_2 &= -a_2 dA.
\end{aligned}
$$

Let $\alpha' := a_1\omega'_1 + a_2\omega'_2$. Take exterior derivative of both sides to see that $d\omega' = i\alpha' \wedge \omega'$ on $U$, and that

$$
d\omega = i\,(\alpha' + d\theta) \wedge \omega.
$$

The real-valued 1-form $\alpha'$ on $U$ is uniquely determined by $0 = d\omega' - i\alpha' \wedge \omega'$. So the real-valued 1-form $\alpha := \alpha' + d\theta$ is as well, giving $0 = d\omega - i\alpha \wedge \omega$ for a unique real-valued *Levi-Civita connection* 1-*form* $\alpha$. Take exterior derivative of that equation to see that $d\alpha = K\omega_1 \wedge \omega_2$, for a unique function $K$ on the tangent frame bundle. But on $\ulcorner U$, $d\alpha = d(\alpha' + d\theta) = d\alpha'$ is the same for any choice of framing, i.e. $K$ is a function well defined on the surface $S$.

**G.9** Suppose that $S$ is not orientable. Prove that $K$ is nonetheless well defined on $S$, while the area form $dA = \omega_1 \wedge \omega_2$ is not. Prove that the curvature form $K\,dA$ is well defined on $S$ if and only if $K = 0$.

**Theorem G.17** (Theorema egregium)**.** *The Levi-Civita connection* 1*-form and the Gauss curvature of a Riemannian metric on a surface depend only on the Riemannian metric, not on any choice of isometric immersion into Euclidean space.*

Geodesic normal coordinates work out exactly the same. But on abstract surfaces, we can reverse their construction. Take any function $K(r, \theta)$; there is a unique solution $h$ to $h_{rr} + Kh = 0$ with $h = 0$ and $h_r = 1$ at $r = 0$: every function $K$ defined near the origin of the plane is the Gauss curvature of some abstract Riemannian metric $dr^2 + h^2 d\theta^2$ near the origin.

**G.10** Prove that this metric is smooth at the origin.

If $K < 0$ then $h_{rr} = -Kh > 0$ for $h > 0$, so $h > r$; if $K$ is defined in the entire plane, $r, \theta$ are geodesic normal coordinates, so this metric is complete.

**Theorem G.18.** *Any connected surface bearing a complete metric of negative Gauss curvature has exponential map at any point a smooth universal covering map from the plane.*

*Proof.* The exponential map is defined in the plane by completeness, but $h$ has no zeroes so the exponential map is a local diffeomorphism. The pullback metric is complete:

$$\omega_1^2 + \omega_2^2 = dr^2 + h^2 d\theta^2 \geq dr^2 + r^2 d\theta^2.$$

Take a closed disk on the surface, say of radius $r$. Its preimage is a union of disks in the pullback metric, closed by completeness. Each sits in a Euclidean disk of smaller radius. If, no matter how small we make $r$, some of these disks overlap, then there are points of arbitrarily small Euclidean distance mapping to the same point of the surface. Take a convergent sequence of these, a contradiction. So the exponential map covers each small enough embedded closed disk on the surface with disjoint embedded closed disks, and so is a covering map. □

Smoothness of isometries and the Gauss–Bonnet theorem follow, with the same proof, for any oriented surface with Riemannian metric.

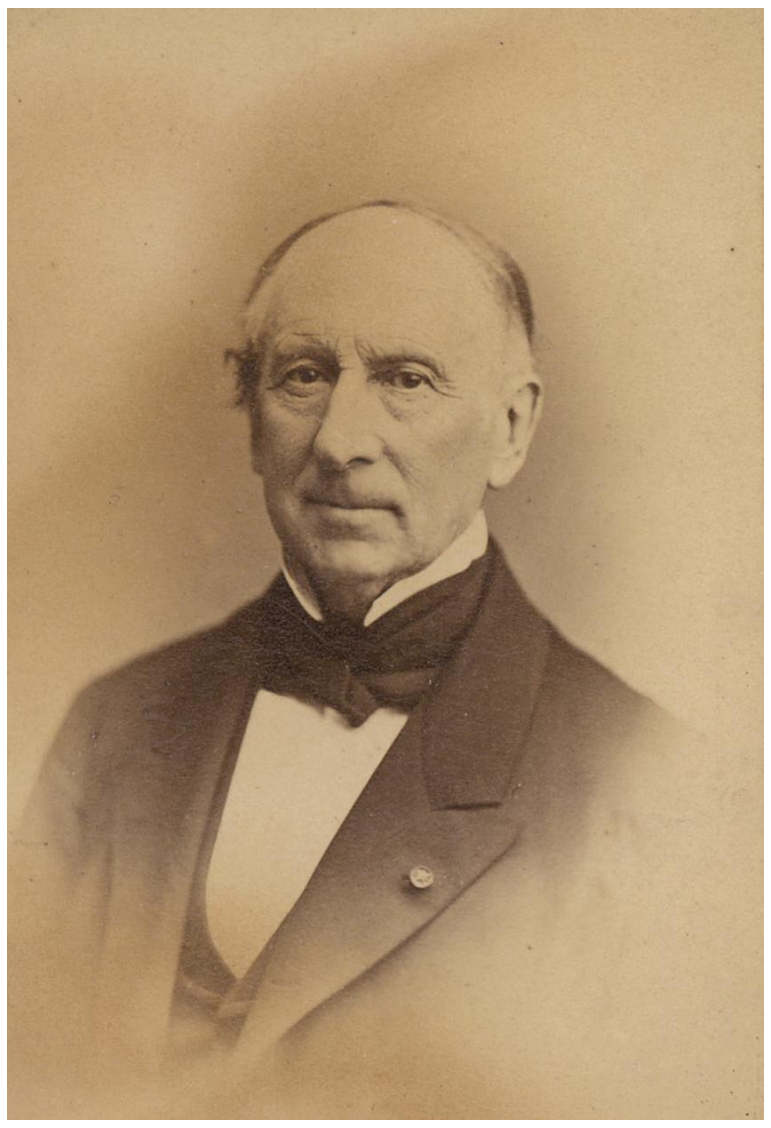

Augustin Louis Cauchy

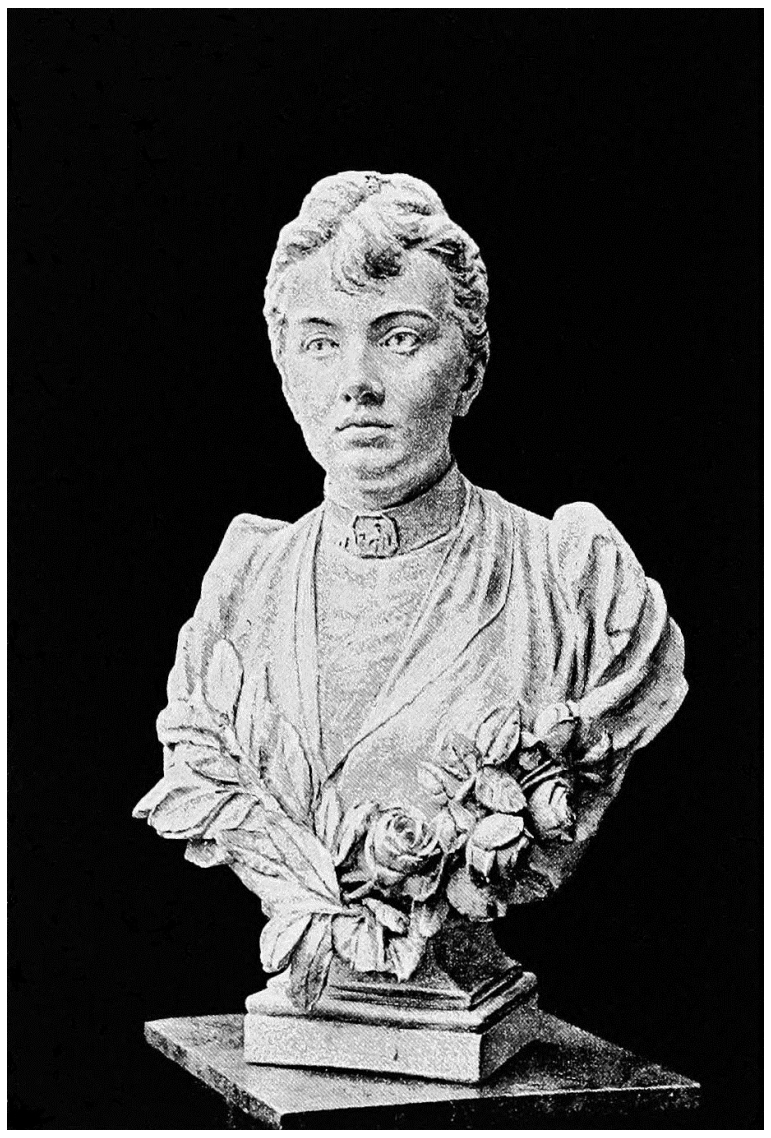

Sofia Vasilyevna Kovalevskaya

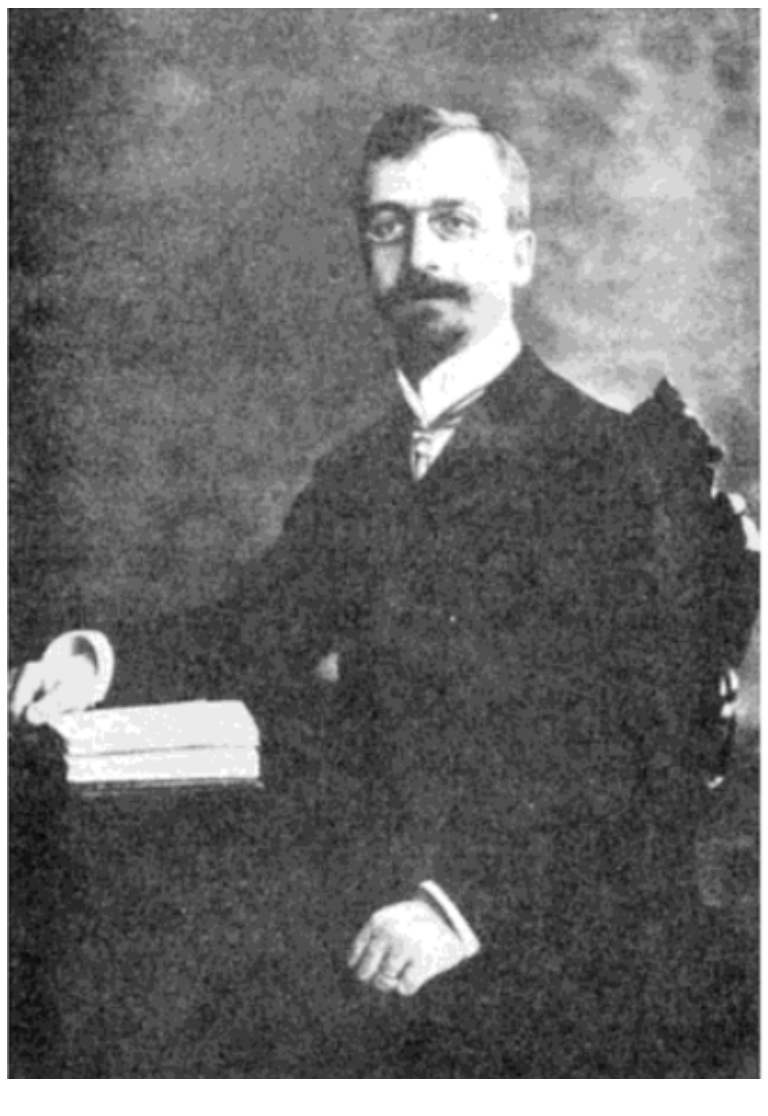

Élie Cartan

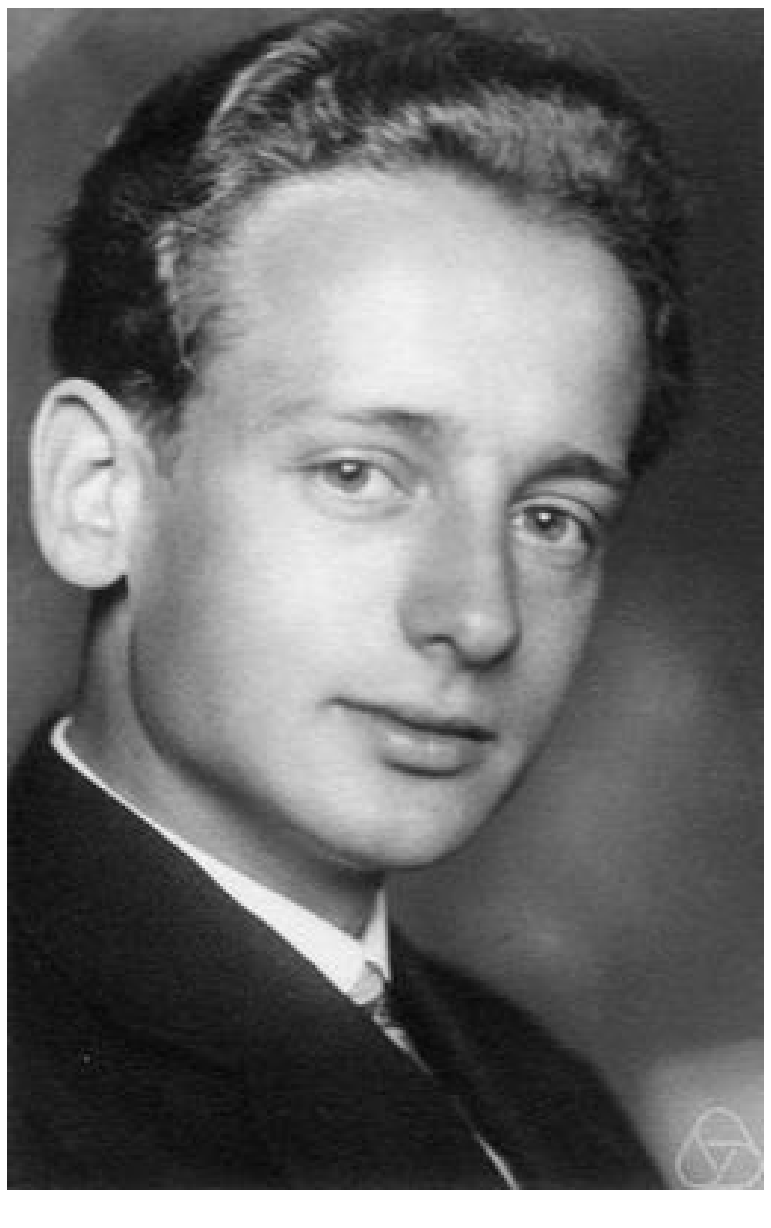

Erich Kähler

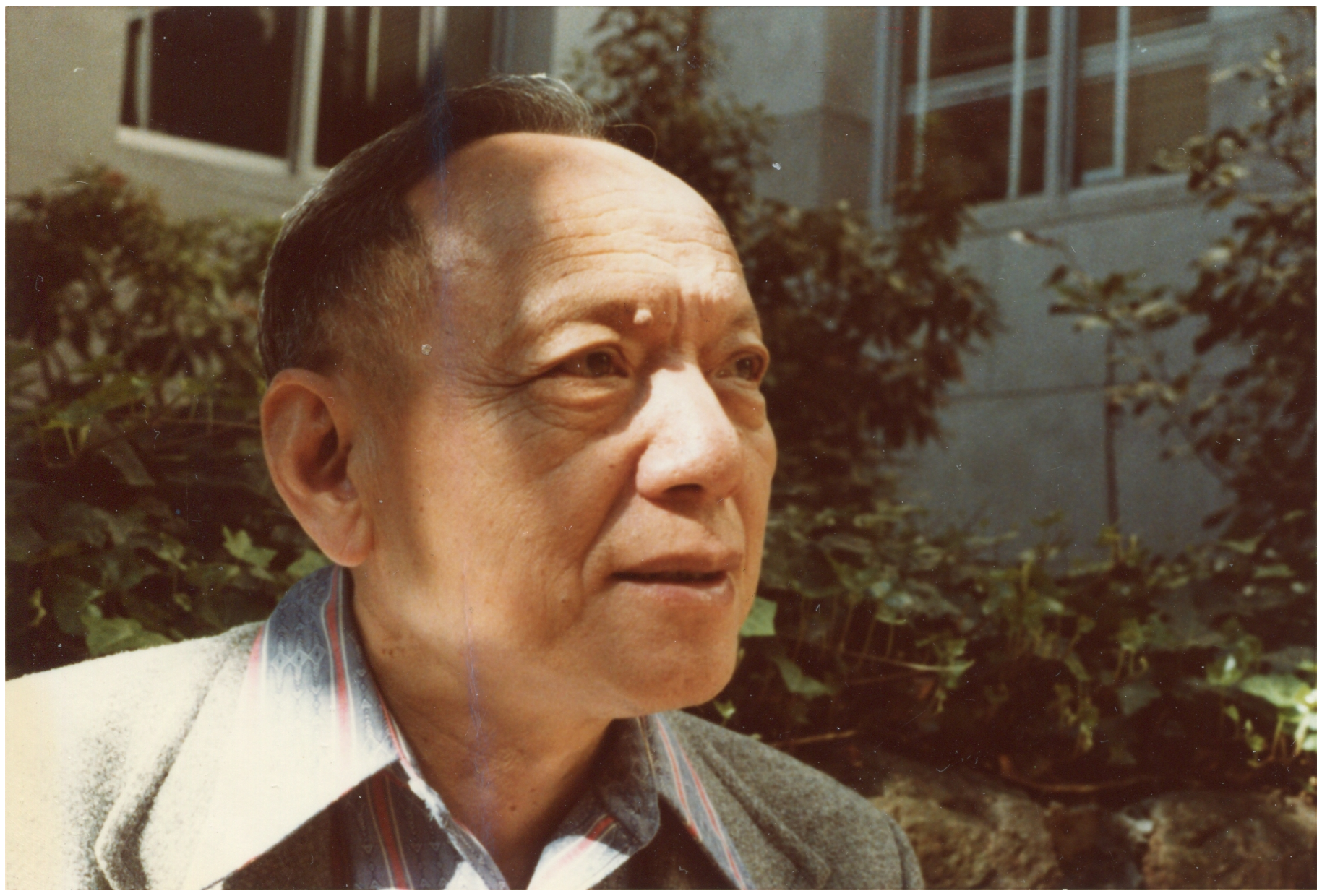

Shiing-Shen Chern[b]

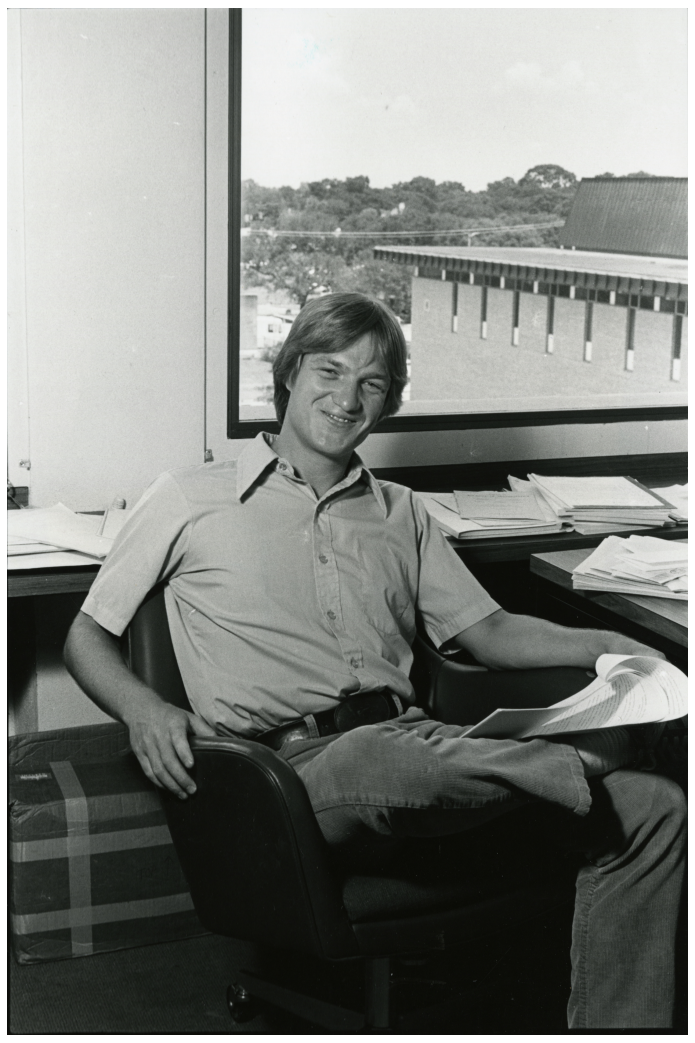

Robert Bryant[a]

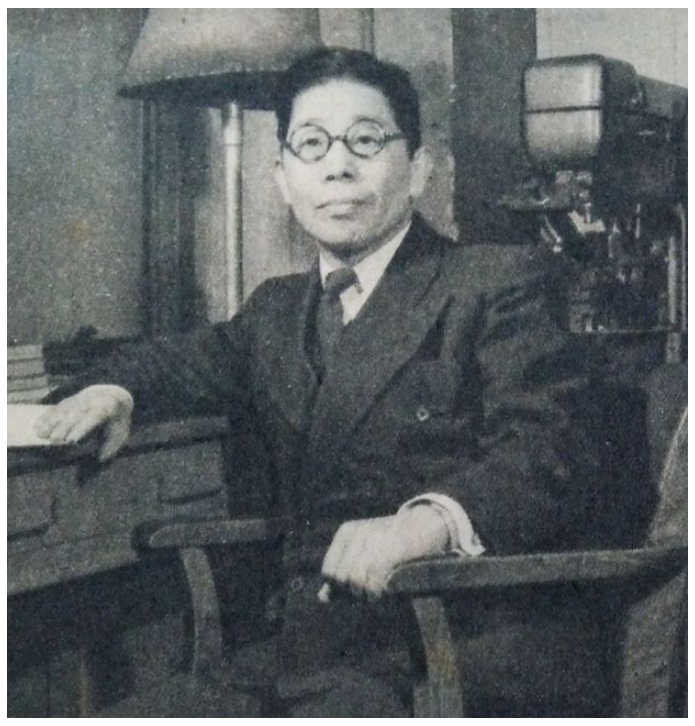

Masatsugu Kuranishi

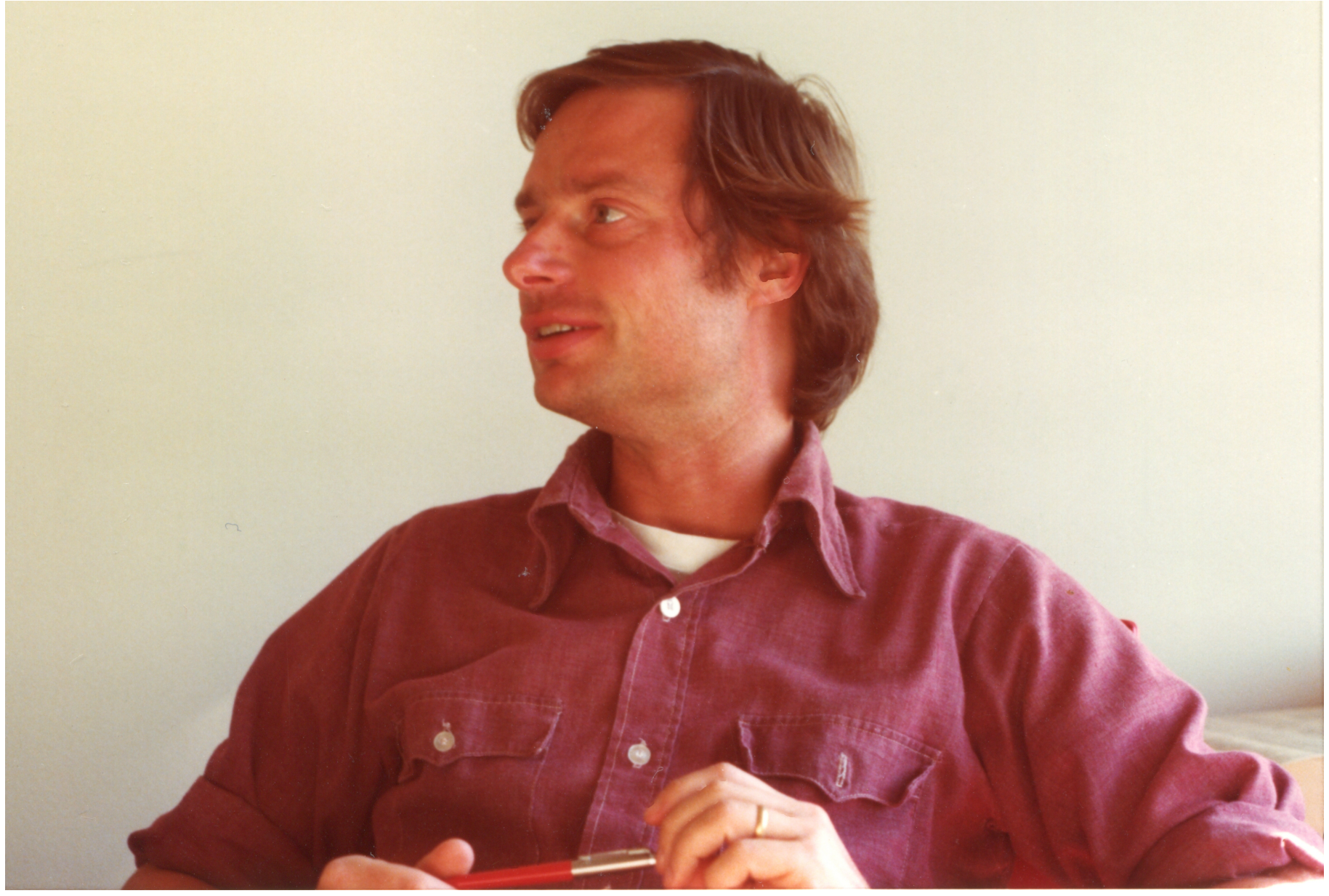

PHILLIP GRIFFITHS[c]

*a.* Aronosky, James (photographer). “Robert Bryant, Mathematics, Rice University.” (1980) Rice University: https://hdl.handle.net/1911/90415. This digital version is licensed under a Creative Commons Attribution 3.0 Unported

*b.* This file is licensed under the Creative Commons Attribution-Share Alike 4.0 International license. Attribution: George M. Bergman

*c.* This file is licensed under the Creative Commons Attribution-Share Alike 4.0 International license. Attribution: George M. Bergman

# Hints

*It is not enough for me to read that the sand on the seashore is soft. My bare feet must feel it.*

— André Gide
Les nourritures terrestres

**1.1.**

*a.* All submanifolds.

*b.* Any discrete set of points.

*c.* Take a 0-dimensional submanifold of the image. Its preimage is a submanifold. Take any submanifold of that submanifold.

*d.* If $dx^1 \wedge dx^2 \neq 0$, these are locally $y_i = \frac{\partial S}{\partial x_i}$; in general the integral surfaces are the Lagrangian surfaces.

*e.* If $dx \neq 0$ these are $y = y(x)$, $z = dy/dx$; in general they are Legendre curves.

**1.2.**

*a.* All linear subspaces of all tangent spaces.

*b.* The zero subspaces of all tangent spaces.

*c.* The linear subspaces tangent to the fibers.

*d.* If $dx^1 \wedge dx^2 \neq 0$, these are $dy_i = a_{ij} dx_j$ with $a_{ij} = a_{ji}$; in general they are zero subspaces, lines, or Lagrangian planes.

*e.* If $dx \neq 0$ these are $dy = z\,dx$, $dz = a\,dx$ for any real number$a$.

**1.4.** $M = \mathbb{R}$, $\mathcal{I} = (x\,dx)$, $X = \{ 0 \}$

**1.8.** Take any local coordinates on $M$, say $x^i, u^a$, near a point $m_0$. Every $p$-plane $E \subset T_m M$, with $m$ near $m_0$, and with $dx^1 \wedge \cdots \wedge dx^p \neq 0$ on $E$, is the set of tangent vectors on which $du^a = q^a_i\,dx^i$, for uniquely determined numbers $q^a_i$. As a short hand, write this as $du = q\,dx$. The functions $x^i, u^a, q^a_i$ on $\mathrm{Gr}_p M$ are local coordinates. The map $(m, E) \in \mathrm{Gr}_p M \mapsto m \in M$ is $(x, u, q) \mapsto (x, u)$, a submersion. We leave the change of coordinates to the reader.

**1.9.** If $\theta = du - q^2\,dx$ then $d\theta = -2q\,dq \wedge dx$.

**1.12.** Suppose that, to each point $x \in M$ of a manifold $M$, we have an associated $p \times q$ matrix $\varphi(x)$, with entries analytic functions on $M$. Let $K$ be the set of pairs $(x, v)$ so that $x \in M$ and $v$ is a vector in the kernel of $\varphi(x)$. If $\varphi(x)$ has rank independent of $x$, let us prove that $K$ is an embedded submanifold of $M \times \mathbb{R}^q$. Let's also prove that, for any map $w\colon M \to \mathbb{R}^p$, the following are equivalent:

*a.* At every point $x \in M$, there is a vector $v$ so that $\varphi(v) = w(x)$.

*b.* Near every point $x \in M$, there is an analytic vector valued function $v$ on an open subset of $M$, so that $\varphi(v(x)) = w(x)$; if $\varphi$ is 1-1 at every point $x \in M$, then $v$ is unique.

We can change $\varphi$ by multiplying by arbitrary invertible matrices of analytic functions. Since $\varphi$ has constant rank, for any chosen point of $M$, we can permute rows and columns to get the upper left corner to be invertible at that point, and hence near that point, of the same rank as $\varphi$:

$$\varphi = \begin{pmatrix} A & B \\ C & D \end{pmatrix},$$

with $A$ invertible. Multiply by

$$\begin{pmatrix} A^{-1} & 0 \\ 0 & I \end{pmatrix},$$

to get $A = I$. Having rank exactly that of the upper left corner is precisely $D = CB$. Multiply:

$$\begin{pmatrix} I & 0 \\ -C & I \end{pmatrix} \varphi \begin{pmatrix} I & -B \\ 0 & I \end{pmatrix} = \begin{pmatrix} I & 0 \\ 0 & I \end{pmatrix}.$$

The reader familiar with vector bundles [6] may generalize.

**1.13.** We prove it in each tangent space separately. Real analyticity then follows from the uniqueness of solutions of the linear equations by problem 1.12 on page 6. So we prove it for constant coefficient 1-forms in $\mathbb{R}^p$. Choose basis in $\mathbb{R}^p$ so that in coordinates

$$x_1, x_2, \dots, x_k, y_1, y_2, \dots, y_{p-k}$$

our 1-forms are $\xi_i = dx_i$. Write $\alpha_i = \sum a_{ij} dx_j + \sum b_{i\ell} dy_\ell$, say. The equation $0 = \alpha_1 \wedge \xi_1 + \alpha_2 \wedge \xi_2 + \dots + \alpha_k \wedge \xi_k$ is precisely

$$\begin{aligned} 0 &= \sum_j a_{ij} dx_j \wedge dx_i + \sum_\ell b_{i\ell} dy_\ell \wedge dx_i, \\ &= \frac{1}{2} \sum_j (a_{ij} - a_{ji})\, dx_j \wedge dx_i + \sum_\ell b_{i\ell} dy_\ell \wedge dx_i, \end{aligned}$$

Plugging in the unit vectors pointing along various coordinate axes, we find $a_{ij} = a_{ji}$ and $b_{i\ell} = 0$.

**1.15.** This is already clear geometrically: the shape operator is the curvature of geodesics. Let $M := {}^{\ulcorner}\mathbb{E}^3$ be the frame bundle of $\mathbb{E}^3$. The frame bundle ${}^{\ulcorner}S$ of any surface $S$ in $\mathbb{E}^3$ satisfies $\omega_3 = 0$ and

$$\begin{pmatrix} \gamma_{13} \\ \gamma_{23} \end{pmatrix} = \begin{pmatrix} a_{11} & a_{12} \\ a_{21} & a_{22} \end{pmatrix} \begin{pmatrix} \omega_1 \\ \omega_2 \end{pmatrix}$$

for some smooth functions $a_{11}, a_{12} = a_{21}, a_{22}$. To have vanishing shape operator, $0 = a_{11} = a_{12} = a_{22}$, so $0 = \omega_3 = \gamma_{13} = \gamma_{23}$.

Consider on $M$ the exterior differential system with equations $\vartheta_0 = \omega_3, \vartheta_1 = \gamma_{13}, \vartheta_2 = \gamma_{23}$. The frame bundle ${}^{\ulcorner}S$ of any surface with vanishing shape operator is a 3-dimensional integral manifold. To be precise, each component of ${}^{\ulcorner}S$ is an integral manifold, since $S$ might have more than one orientation, so ${}^{\ulcorner}S$ might have more than one component. Check that the exterior differential system is Frobenius, so there is a unique integral manifold through each point of $M$. But we already have an example of such an integral manifold, so that must be the only one.

**1.16.** Again take $M = {}^{\ulcorner}\mathbb{E}^3$ with the exterior differential system

$$\omega_3, \gamma_{13} - c_0\omega_1, \gamma_{23} - c_0\omega_2,$$

Again the system is Frobenius, so there is only one 3-dimensional integral manifold through each point. Rotate and translate a sphere of suitable radius to see one such integral manifold through each point of $M$, the frame bundle of that sphere.

**1.17.** Problem 1.16 on page 7 handles the case of one principal curvature, so assume that there are two principal curvatures. After picking a frame on a surface $S$ so that $e_1$ and $e_2$ are in the principal directions, we find $0 = \omega_3 = \gamma_{13} - k_1\omega_1 = \gamma_{23} - k_2\omega_2$, for constants $k_1, k_2$. Differentiate all three equations to find that $0 = (k_1 - k_2)\,\gamma_{12} \wedge \omega_2 = (k_1 - k_2)\,\gamma_{12} \wedge \omega_1$, forcing $\gamma_{12} = 0$ by Cartan's lemma. Differentiating the equation $\gamma_{12} = 0$, we find $0 = k_1k_2\omega_1 \wedge \omega_2$ so that either $k_1 = 0$ or $k_2 = 0$. We can assume that $k_1 = 0$, i.e. $\gamma_{13} = 0$. But then

$$\begin{aligned} de_1 &= (e_1 \cdot de_1)\, e_1 + (e_2 \cdot de_1)\, e_2 + (e_3 \cdot de_1)\, e_3, \\ &= \gamma_{11}e_1 + \gamma_{21}e_2 + \gamma_{31}e_3, \\ &= 0. \end{aligned}$$

Therefore $e_1$ is constant as we travel along the surface. So the surface is a collection of straight lines, in this $e_1$ direction, all placed perpendicular to a curve in the $e_2, e_3$-plane. On that curve, $\omega_1 = 0$, and $d\omega_2 = 0$ so we can write locally $\omega_2 = ds$ for some function $s$. Then we find

$$\begin{aligned} de_2 &= (e_1 \cdot de_2)\, e_1 + (e_2 \cdot de_2)\, e_2 + (e_3 \cdot de_2)\, e_3, \\ &= \gamma_{12}e_1 + \gamma_{22}e_2 + \gamma_{32}e_3, \\ &= -k_2 ds e_3. \end{aligned}$$

and

$$\begin{aligned} de_3 &= (e_1 \cdot de_3)\, e_1 + (e_2 \cdot de_3)\, e_2 + (e_3 \cdot de_3)\, e_3, \\ &= \gamma_{13}e_1 + \gamma_{23}e_2 + \gamma_{33}e_3, \\ &= k_2 ds e_2. \end{aligned}$$

Check that the vectors $E_2 = k_2\cos(s)e_2 + k_2\sin(s)e_3$, $E_3 = -k_2\sin(s)e_2 + k_2\cos(s)e_3$ are constant. Rotate so that $E_1 = e_1, E_2, E_3$ are the standard basis vectors of $\mathbb{E}^3$, to see that the surface $S$ is a right circular cylinder of radius $1/\left|k_2\right|$.

**1.21.** Start by looking at the geometry of the given foliation. Suppose that the foliation is of an open subset $U$ of 3-dimensional Euclidean space. Consider the 4-dimensional manifold $M$ of all orthonormal frames $e_1, e_2, e_3$ at points $x \in U$ for which $e_3$ is perpendicular to the leaf through $x$ of the given foliation. By the Frobenius theorem,

$$\begin{aligned} 0 &= \omega_3 \wedge d\omega_3, \\ &= \omega_3 \wedge (\gamma_{13} \wedge \omega_1 + \gamma_{23} \wedge \omega_2), \\ &= \gamma_{13} \wedge \omega_1 \wedge \omega_3 + \gamma_{23} \wedge \omega_2 \wedge \omega_3. \end{aligned}$$

So on $M$,

$$\begin{pmatrix}\gamma_{13}\\ \gamma_{23}\end{pmatrix} = \begin{pmatrix}a_{11} & a_{12} & a_{13}\\ a_{21} & a_{22} & a_{23}\end{pmatrix}\begin{pmatrix}\omega_1\\ \omega_2\\ \omega_3\end{pmatrix}$$

for some functions $a_{ij}$ on $M$. Plugging these in above, we find that $a_{12} = a_{21}$. Clearly $a_{11}\omega_1^2 + a_{12}\omega_1\omega_2 + a_{21}\omega_2\omega_1 + a_{22}\omega_2^2$ is the shape operator of each leaf of the foliation.

Now consider the problem of constructing a triply orthogonal web incorporating this foliation as one of its three. We need to find a choice of $e_1, e_2$ at each point to construct a frame, so that $e_1$ will be perpendicular to the leaves of the first foliation, and $e_2$ perpendicular to the leaves of the second foliation, and the third foliation will be the one we started with, already perpendicular to $e_3$. So we will need to solve the exterior differential system

$$0 = \omega_1 \wedge d\omega_1 = \omega_2 \wedge d\omega_2$$

on $M$, i.e. with the equations

$$\begin{pmatrix}\gamma_{13}\\ \gamma_{23}\end{pmatrix} = \begin{pmatrix}a_{11} & a_{12} & a_{13}\\ a_{21} & a_{22} & a_{23}\end{pmatrix}\begin{pmatrix}\omega_1\\ \omega_2\\ \omega_3\end{pmatrix}$$

already in force.

Note that on $M$, $\omega_1, \omega_2, \omega_3, \gamma_{12}$ are linearly independent 1-forms. We are looking for an integral 3-manifold $X$ of that exterior differential system, on which we want $\omega_1, \omega_2, \omega_3$ to be linearly independent, i.e. $X$ projects by local diffeomorphism to 3-dimensional Euclidean space.

It might be simpler to write our foliation shape operator as

$$\begin{pmatrix}\gamma_{13}\\ \gamma_{23}\end{pmatrix} = \begin{pmatrix}a_{131} & a_{132} & a_{133}\\ a_{231} & a_{232} & a_{233}\end{pmatrix}\begin{pmatrix}\omega_1\\ \omega_2\\ \omega_3\end{pmatrix}$$

and then we are imposing only the relations $a_{132} = a_{231}$, i.e. symmetry in these two outer indices.

Differentiate the equations of our exterior differential system to find that on any integral 3-manifold $X$:

$$\begin{aligned} 0 &= \omega_i \wedge d\omega_i, \\ &= -\omega_i \wedge \sum_j \gamma_{ij} \wedge \omega_j, \end{aligned}$$

which forces $\gamma_{ij}$ to be a linear combination

$$\gamma_{ij} = \sum_k a_{ijk}\omega_k,$$

with $a_{ijk} = -a_{jik}$ since $\gamma_{ij} = -\gamma_{ji}$. But plug in to get

$$0 = \omega_i \wedge \sum_{jk} (a_{ijk} - a_{ikj})\omega_k \wedge \omega_j$$

so that $a_{ijk} = a_{ikj}$ if $i, j, k$ are all distinct. So for distinct indices, $a_{ijk}$ is symmetric in $jk$, but antisymmetric in $ij$. These two involutions generate the permutation group, and so the sign of $a_{ijk}$ is a representation of the permutations on 3 letters. But any two involutions are conjugate in the permutation group, so they must force the same sign change. Hence $a_{ijk} = 0$ for $i, j, k$ distinct. This is precisely the demand that the shape operators of the leaves of all three foliations are thus diagonal in the frame $e_1, e_2, e_3$. In other words, each leaf of each foliation lies normal to each leaf of each other foliation, and intersects tangent to a principal direction, i.e. along a principal curve.

We can assume that $U$ is connected. Either

*a.* the leaves of the given foliation are everywhere umbilic, hence each leaf is an open subset of a sphere or plane, and so $a_{132} = 0$ and $a_{131} = a_{232}$, or else

*b.* the vectors $e_1, e_2$ have to be chosen in principal directions and this determines them up to 4 choices, at least on a dense open subset of $U$.

Suppose that the leaves are spheres or planes, so $a_{132} = 0$ and $a_{131} = a_{232}$, say. The exterior differential system is generated in dimension 3 by $\gamma_{12} \wedge \omega_1 \wedge \omega_2 = 0$, so every line or plane in any tangent space of $M$ is an integral element. Planes, i.e. integral planes, generically have $\omega_1, \omega_2$ linearly independent on them, so are generically of the form

$$\gamma_{12} = p_1\omega_1 + p_2\omega_2,$$
$$\omega_3 = q_1\omega_1 + q_2\omega_2,$$

so lie in a unique 3-dimensional integral element $\gamma_{12} = p_1\omega_1 + p_2\omega_2$. The polar equations of any integral point or line are trivial, but the generic integral plane has polar equation $\gamma_{12} - p_1\omega_1 - p_2\omega_2$, so $s_1 = 0, s_2 = 1, s_3 = 0$, solutions depend on 1 function of 2 variables.

We can explicitly construct the web: take any leaf of our given foliation, the *initial leaf*, and draw on it any foliation by curves. On that same surface, draw the orthogonal foliation by curves. Drag each leaf of each of those foliations along the flow of $e_3$, through space, to trace out a surface.

We want to see that this construction always creates a triply orthogonal web. Any triply orthogonal web containing the given foliation has to arise from this construction: the leaves of the other two foliations intersect the initial leaf in curves, and the vector field $e_3$ is tangent to every leaf of the other two foliations.

Conversely, take a foliation of the initial leaf by curves. Locally pick any orthonormal vector fields $e_1, e_2$ tangent to that leaf, with $e_1$ tangent and $e_2$ perpendicular to that curve foliation. Drag, as above. We produce two more foliations of Euclidean space, defined near that leaf. But it is not clear whether they remain perpendicular. The leaves of the two new foliations both contain $e_3$, and they start off perpendicular along the initial leaf.

Compute the change in the Euclidean metric along the flow of $e_3$:

$$\mathcal{L}_{e_3}\omega_i\omega_i = 2a_{i3j}\omega_i\omega_j.$$

So umbilicity of the given leaves is precisely the condition that the Euclidean metric varies only by scaling as we flow along $e_3$, on the perpendicular vectors to $e_3$, and so any pair of planes in a tangent space which start perpendicular will remain so. Hence the construction always succeeds.

Suppose now that the leaves are nowhere umbilic. We will see that each foliation by nonumbilic surfaces lies in at most one, and typically no, triply orthogonal web.

We need $e_1, e_2$ to diagonalize the shape operator, i.e. to lie in principal directions. Imagine drawing the principal curves, i.e. the curves in those directions, on each leaf. Our given foliation by surfaces has now, on each surface, two foliations by perpendicular curves. We want to see whether, when we flow along the unit normal vector field $e_3$ of the foliation, these principal curves flow into one another. For a generic foliation by surfaces, the flow of $e_3$ will spin the principal curves of one leaf around into various curves on other leaves, not necessarily principal.

We need $a_{132} = 0$, so the shape operator is diagonalized, i.e. $e_1, e_2$ point in principal directions, and we suppose the leaves are not umbilic, so $a_{131} \neq a_{232}$. Our frames form a 3-manifold $X$ in the frame bundle. We have to decide whether $X$ is an integral manifold. On $X$, $\gamma_{ij} = \sum_k a_{ijk}\omega_k$, for some functions $a_{ijk} = -a_{jik}$. To have a triply orthogonal web, we need precisely that $a_{ijk} = 0$ if $i, j, k$ are distinct, i.e. the shape operators are diagonalized. This is precisely the condition that the prinicipal curves flow into one another.

**2.1.** We can ask that $e_1$ be perpendicular to the leaves of our foliation. (If you use $e_3$ instead of $e_1$ here, you will find it more complicated to write out the tableau.) So then $\omega_1 \wedge d\omega_1 = 0$ on the foliation. Expand out

$$d\omega_1 = -\gamma_{12} \wedge \omega_2 - \gamma_{13} \wedge \omega_3$$

to arrive at the tableau

$$\begin{pmatrix} \gamma_{12} & \gamma_{13} \end{pmatrix} \wedge \begin{pmatrix} \omega_1 \wedge \omega_2 \\ \omega_1 \wedge \omega_3 \end{pmatrix}.$$

so $s_1 = 0, s_2 = 1, s_3 = 1$, foliations of open sets of 3-dimensional Euclidean space depend on 1 function of 3 variables. We can see them as the level sets of 1 function of 3 variables, but the function is defined only up to composition with a strictly increasing or strictly decreasing function.

**2.2.** The polar equations of $E_j$ are given by setting $\omega^i = 0$ for $i > j$, and plugging in, with $\pi^\alpha = 0$ on $E_j$, so all terms with 2 or more $\pi$ vanish, i.e. the polars in grades $0, 1, 2, \dots, j$. Hence the characters of $E_p$ are the numbers of polars in each grade.

**2.5.** Take $\mathbb{R}^{15}$ with coordinates $x^i, u^i, u^i_j$ for $i, j = 1, 2, 3$. Note that our differential equations, spelled out as algebraic equations

$$\begin{aligned} u^3_2 - u^2_3 &= f^1 - u^1, \\ u^1_3 - u^3_1 &= f^2 - u^2, \\ u^2_1 - u^1_2 &= f^3 - u^3, \end{aligned}$$

cut out a submanifold $M \subset \mathbb{R}^{15}$ of dimension 12. Take $\omega^i := dx^i$, $\theta^i := du^i - u^i_j\,dx^j$, and $\pi^i_j := du^i_j$. It will help to denote $\frac{\partial f^i}{\partial x^j}$ by $f^i_j$. The equations of $M$ force 3 linear relations among the $\pi^i_j$:

$$du^3_2 - du^2_3 = \left(f^1_i - u^1_i\right) dx^i,$$

modulo $\theta^1, \theta^2, \theta^3$, and so on, i.e.

$$\begin{aligned} \pi^2_3 &= \pi^3_2 - \left(f^1_i - u^1_i\right)\omega^i, \\ \pi^1_3 &= \pi^3_1 + \left(f^2_i - u^2_i\right)\omega^i, \\ \pi^1_2 &= \pi^2_1 - \left(f^3_i - u^3_i\right)\omega^i. \end{aligned}$$

Our tableau: modulo $\theta^1, \theta^2, \theta^3$,

$$d\begin{pmatrix}\theta^1\\ \theta^2\\ \theta^3\end{pmatrix} = -\begin{pmatrix}\pi^1_1 & \pi^1_2 & \pi^1_3\\ \pi^2_1 & \pi^2_2 & \pi^2_3\\ \pi^3_1 & \pi^3_2 & \pi^3_3\end{pmatrix}\wedge\begin{pmatrix}\omega^1\\ \omega^2\\ \omega^3\end{pmatrix}$$

$$= -\begin{pmatrix}\pi^1_1 & \pi^2_1 & \pi^3_1\\ \pi^2_1 & \pi^2_2 & \pi^3_2\\ \pi^3_1 & \pi^3_2 & \pi^3_3\end{pmatrix}\wedge\begin{pmatrix}\omega^1\\ \omega^2\\ \omega^3\end{pmatrix} + \begin{pmatrix}\tau^1\\ \tau^2\\ 0\end{pmatrix}$$

where the torsion is

$$\begin{pmatrix}\tau^1\\ \tau^2\end{pmatrix} = \begin{pmatrix}f^3_1 - u^3_1\\ 0\end{pmatrix}\omega^{12} + \begin{pmatrix}u^2_1 - f^2_1\\ f^1_1 - u^1_1\end{pmatrix}\omega^{13} + \begin{pmatrix}u^2_2 - f^2_2 + u^3_3 - f^3_3\\ f^1_2 - u^1_2\end{pmatrix}\omega^{23}$$

We can try to absorb torsion, for example by using

$$\begin{pmatrix}\pi'^1_1\\ \pi'^2_1\\ \pi'^2_2\end{pmatrix} := \begin{pmatrix}\pi^1_1\\ \pi^2_1\\ \pi^2_2\end{pmatrix} + \begin{pmatrix}0 & 0 & u^2_1 - f^2_1\\ u^3_1 - f^3_1 & 0 & 0\\ 0 & 0 & f^1_2 - u^1_2\end{pmatrix}\begin{pmatrix}\omega^1\\ \omega^2\\ \omega^3\end{pmatrix},$$

which we denote as $\pi$ instead of $\pi'$ to simplify notation. Our tableau: modulo $\theta^1, \theta^2, \theta^3$,

$$d\begin{pmatrix}\theta^1\\ \theta^2\\ \theta^3\end{pmatrix} = -\begin{pmatrix}\pi^1_1 & \pi^2_1 & \pi^3_1\\ \pi^2_1 & \pi^2_2 & \pi^3_2\\ \pi^3_1 & \pi^3_2 & \pi^3_3\end{pmatrix}\wedge\begin{pmatrix}\omega^1\\ \omega^2\\ \omega^3\end{pmatrix} + \begin{pmatrix}u^i_i - f^i_i\\ 0\\ 0\end{pmatrix}\omega^{23}.$$

The torsion is $u^i_i - f^i_i$ (Einstein notation: implicitly summed over $i$). Take the submanifold $M' \subset M$ cut out by the equation $u^i_i = f^i_i$. For simplicity, denote this submanifold as $M$ henceforth. On $M$,

$$0 = d(u^i_i - f^i_i) = \pi^i_i - f^i_{ij}\omega^j.$$

Our tableau: modulo $\theta^1, \theta^2, \theta^3$,

$$d\begin{pmatrix}\theta^1\\ \theta^2\\ \theta^3\end{pmatrix} = -\begin{pmatrix}\pi^1_1 & \pi^2_1 & \pi^3_1\\ \pi^2_1 & \pi^2_2 & \pi^3_2\\ \pi^3_1 & \pi^3_2 & -\pi^1_1 - \pi^2_2 + f^i_{ij}\omega^j\end{pmatrix}\wedge\begin{pmatrix}\omega^1\\ \omega^2\\ \omega^3\end{pmatrix}.$$

Let

$$\begin{pmatrix}\pi'^3_1\\ \pi'^3_2\end{pmatrix} = \begin{pmatrix}\pi^3_1\\ \pi^3_2\end{pmatrix} - \begin{pmatrix}f^i_{i1}\\ f^i_{i2}\end{pmatrix}\omega^3,$$

and once again just write these as $\pi$ instead of $\pi'$. Our tableau: modulo $\theta^1, \theta^2, \theta^3$,

$$d\begin{pmatrix}\theta^1\\ \theta^2\\ \theta^3\end{pmatrix} = -\begin{pmatrix}\pi^1_1 & \pi^2_1 & \pi^3_1\\ \pi^2_1 & \pi^2_2 & \pi^3_2\\ \pi^3_1 & \pi^3_2 & -\pi^1_1 - \pi^2_2\end{pmatrix}\wedge\begin{pmatrix}\omega^1\\ \omega^2\\ \omega^3\end{pmatrix}.$$

| $s_1$ | $s_2$ | $s_3$ |
|---|---|---|
| 3 | 2 | 0 |

Integral elements $\pi^i_j = p^i_{jk}\omega^k$. We highlight certain coefficients, to be discussed in chapter 5:

$$
\begin{aligned}
p^1_{12} &= p^2_{11} && (1)\\
p^1_{13} &= p^3_{11} && (2)\\
p^2_{13} &= p^3_{12} && (3)\\
p^2_{12} &= p^2_{21} && (4)\\
p^2_{13} &= p^3_{21} && (5)\\
p^3_{12} &= p^3_{21} && (6)\\
p^3_{13} &= -p^1_{11} - p^2_{21} && (7)\\
p^3_{23} &= -p^1_{12} - p^2_{22}. && (8)
\end{aligned}
$$

These coefficients are solved for in terms of others, except for the 3rd and 8th equations. But we use the 6th equation to fix up the 3rd, and the 1st equation to fix up the 8th, to solve for highlighted coefficients in terms of others. Hence the space of integral elements at each point of the 11-dimensional manifold $M'$ has dimension given by counting the other coefficients: 7 dimensions of integral element at each point. Involution, with the general solution depending on 2 functions of 2 variables.

**2.8.** $q^i_{jk} = c^i_{m\ell}\left(g^m_k g^\ell_j - g^m_j g^\ell_k\right)$

**2.11.** Choose the third fundamental form so that the Gauss curvature and the squared mean curvature have linearly independent differentials. Replace the surface by an open subset on which the Gauss curvature and the squared mean curvature are global coordinates invariant under rigid motion of the surface. Pick the eigenvalues so that the mean curvature is not zero; the surface is not symmetric under reflection in the tangent plane. So rigid motions fix every point of the surface, and also fix a normal direction, so are trivial.

**3.5.** If $u_\mu, v_\mu$ are the vector fields dual to the real and imaginary parts of the $\omega^\mu$, the *torsion tensor*, or *Nijenhuis tensor*

$$
T := (u_\mu + i v_\mu)\tau\omega^\mu
$$

depends only on the almost complex structure, a section of $\Omega^{0,2}_M \otimes TM \otimes \mathbb{C}$.

**3.7.** Consider left action of $\mathrm{SU}_3$ on itself: $L_h z = hz$. The identity function $g(z) = z$ behaves like $(L_h^* g)(z) = g(L_h z) = g(hz) = hz = hg(z)$, so $L_h^* g = hg$. Thus for any constant matrix $h \in \mathrm{SU}_3$,

$$
\begin{aligned}
L_h^*\omega &= L_h^*(g^{-1}dg),\\
&= (L_h^* g)^{-1} d L_h^* g,\\
&= (hg)^{-1} d(hg),\\
&= g^{-1}h^{-1}h\,dg,\\
&= g^{-1}\,dg,\\
&= \omega.
\end{aligned}
$$

**3.8.** Differentiating $\omega = g^{-1}\,dg$, i.e. $dg = g\omega$,

$$
\begin{aligned}
0 &= dg \wedge \omega + g\,d\omega, \\
&= g\omega \wedge \omega + g\,d\omega,
\end{aligned}
$$

we find

$$d\omega = -\omega \wedge \omega.$$

**3.11.** [36] pp. 36–37

**5.3.** Take a $\vartheta$-integral element

$$E_0 = \langle e_1, e_2, \dots, e_p \rangle .$$

Move it:

$$E_t = \left\langle e_1 + tw_1 + O(t)^2, \dots, e_p + tw_p + O(t)^2 \right\rangle .$$

Every motion through $p$-dimensional subspaces has this form locally. Expand:

$$\vartheta|_{E_t} = t \sum_i \vartheta(e_1, \dots, e_{i-1}, w_i, e_{i+1}, \dots, e_p) + O(t)^2.$$

The differentials, i.e. linear terms, are sums of polar equations. Set all but one $w_i$ to zero: every polar equation is a differential. In coordinates, $E_0 = (du = 0)$, $E_t = (du = q(t)\,dx)$, let $w_i(t) := q_i^a(t)\partial_{u^a}$.

**5.4.** If an integral is not ordinary, we can find integral elements nearby with characters "borrowed" downward, so larger $s_0$, or the same $s_0$ but larger $s_1$, or some such. As $p + s_0 + \dots + s_p = \dim M$, borrowing raises one character and lowers some later character, decreasing the dimension

$$\dim M + s_1 + 2s_2 + \dots + ps_p,$$

of the submanifold containing all nearby integral elements.

**5.5.** Generic linear subspaces are regular, so ordinary, so involutive. The polar equations along a generic integral element are the same.

**5.6.** Pick a tableau at $m_0$ which is generic in the usual sense: having maximal $s_1$ among all tableaux for the same system, and maximal $s_2$ subject to that $s_1$, and so on. Let $\mathcal{I}' \subseteq \mathcal{I}$ be the exterior differential system generated by the $\theta^a$ and the rows of the tableau. It is possible that various forms in $\mathcal{I}$ vanish at $m_0$, and near $m_0$ lie outside $\mathcal{I}'$. Note that the characters of $\mathcal{I}'$ at any $\mathcal{I}$-integral element are no larger than those of $\mathcal{I}$ since $\mathcal{I}' \subseteq \mathcal{I}$.

Suppose now that $\mathcal{I}$ is involutive. The characters of $\mathcal{I}$ are locally constant by problem 5.5 on page 41. The characters of $\mathcal{I}'$ are at least as large as those of the tableau, which match those of $\mathcal{I}$ at $m_0$. So $\mathcal{I}'$ and $\mathcal{I}$ have the same characters, and both are involutive on all of their integral elements, with the same dimensions of spaces of integral elements. Hence they have the same integral elements where the tableau maintains those same characters.

**5.7.** "Generic" is just a condition on linear independence of $\pi$ 1-forms, so holds locally. Torsion is absorbable at a point just when there is a torsion free integral element coframed by $\omega^1, \dots, \omega^p$. Absorbing torsion is solving a system of linear equations of constant rank; apply problem 1.12 on page 6

**6.2.** Take local coordinates $x^1, x^2, \dots, x^p, y^1, y^2, \dots, y^q$, where $E$ is the graph of $dy = 0$. If $\vartheta = c_{IA} dx^I \wedge dy^A$, we will see that

$$\left.\mathcal{L}_v \vartheta\right|_E = \frac{\partial v^a}{\partial x^j} c_{Ia} dx^{Ij} + v^a \frac{\partial c_I}{\partial y^a} dx^I.$$

Let $(-1)^I$ mean $(-1)^m$ if $I$ consists of $m$ indices. Write

$$v = v^j \frac{\partial}{\partial x^j} + v^b \frac{\partial}{\partial y^b}.$$

Note that

$$v \lrcorner c_I dx^I = (-1)^J v^i c_{JiK} dx^{JK}.$$

Commuting with exterior derivative,

$$\begin{aligned} \mathcal{L}_v dx^i &= \frac{\partial v^i}{\partial x^j} dx^j + \frac{\partial v^i}{\partial y^b} dy^b, \\ \mathcal{L}_v dy^a &= \frac{\partial v^a}{\partial x^j} dx^j + \frac{\partial v^a}{\partial y^b} dy^b. \end{aligned}$$

By the Leibnitz rule,

$$\begin{aligned} \mathcal{L}_v \vartheta = {} & v^i \frac{\partial c_{IA}}{\partial x^i} dx^I \wedge dy^A + v^a \frac{\partial c_{IA}}{\partial y^a} dx^I \wedge dy^A \\ & + c_{JiKA} dx^J \wedge \left( \frac{\partial v^i}{\partial x^j} dx^j + \frac{\partial v^i}{\partial y^b} dy^b \right) \wedge dx^K \wedge dy^A \\ & + c_{IBaC} dx^I \wedge dy^B \wedge \left( \frac{\partial v^a}{\partial x^j} dx^j + \frac{\partial v^a}{\partial y^b} dy^b \right) \wedge dy^C. \end{aligned}$$

On $E$, $dy = 0$ so

$$\begin{aligned} \left.\mathcal{L}_v \vartheta\right|_E = {} & v^i \frac{\partial c_I}{\partial x^i} dx^I + v^a \frac{\partial c_I}{\partial y^a} dx^I \\ & + c_{JiK} dx^J \wedge \frac{\partial v^i}{\partial x^j} dx^j \wedge dx^K \\ & + c_{Ia} dx^I \wedge \frac{\partial v^a}{\partial x^j} dx^j. \end{aligned}$$

This is not quite the same as

$$\begin{aligned} \left. v \lrcorner d\vartheta \right|_E = {} & v^i \frac{\partial c_I}{\partial x^i} dx^I + v^a \frac{\partial c_I}{\partial y^a} dx^I \\ & + (-1)^{jI} v^\ell \frac{\partial c_{I\ell k}}{\partial x^j} dx^{jIK}. \end{aligned}$$

Write the tangent part of $v$ as

$$v' = v^i \frac{\partial}{\partial x^i}.$$

Let $A$ be the linear map $A \colon E \to E$ given by

$$A^i_j = \frac{\partial v^i}{\partial x^j}$$

and apply this by derivation to forms on $E$ as

$$(A\xi)(v_1,\dots,v_k) = \xi(Av_1,v_2,\dots,v_k) - \xi(v_1,Av_2,v_3,\dots,nv_k) + \dots.$$

Then

$$\mathcal{L}_v\vartheta|_E = v' \lrcorner\, d\vartheta\big|_E + v^a\frac{\partial c_I}{\partial y^a}dx^I + A\vartheta|_E + c_{Ia}dx^I\wedge\frac{\partial v^a}{\partial x^j}dx^j.$$

Since $0 = \vartheta|_E = d\vartheta|_E$, we find

$$\mathcal{L}_v\vartheta|_E = v^a\frac{\partial c_I}{\partial y^a}dx^I + c_{Ia}dx^I\wedge\frac{\partial v^a}{\partial x^j}dx^j.$$

**6.3.** Take local coordinates $x^1,x^2,\dots,x^p,y^1,y^2,\dots,y^q$, where $E$ is the graph of $dy = 0$. Write $\vartheta = c_{IA}dx^I\wedge dy^A$. Write $X$ as the graph $y = y(x)$, so

$$\vartheta|_X = \sum c_{IA}(x,y)dx^I\wedge\frac{\partial y^{a_1}}{\partial x^{j_1}}dx^{j_1}\wedge\dots\wedge\frac{\partial y^{a_\ell}}{\partial x^{j_\ell}}dx^{j_\ell}.$$

so

$$\mathcal{L}_v\vartheta|_E = v^a\frac{\partial c_I}{\partial y^a}dx^I + c_{Ia}dx^I\wedge\frac{\partial v^a}{\partial x^j}dx^j.$$

depends only on knowledge of the point $m$ where we compute coefficients of $\vartheta$ and of $E = T_mX$, so that we drop $dy$ terms.

**6.4.** In coordinates, the linearization of an exterior differential system $\mathcal{I}$ about an integral element $E\subset T_mM$ does not involve any $dy\wedge dy$ terms but depends on the terms with no $dy$ and with one $dy$ explicitly. But $P = \left\langle E^\perp\right\rangle = \langle dy^a\rangle$.

**6.6.**

$$\mathcal{L}_v\vartheta|_E = c_{Ia}dx^I\wedge\frac{\partial v^a}{\partial x^j}dx^j$$

**6.8.** Suppose $E\subset E_+$ is a noncharacteristic hyperplane. By uniqueness of extension $E_+$, the polar equations of $E$ cut out precisely $E_+$ inside $T_mM$, i.e. the dimension of polar equations of $E$ is the dimension of $T_mM/E_+$. The polar equations of any regular hyperplane in $E_+$ are satisfied on $E_+$, so have same rank. So the regular integral elements are also noncharacteristic. If $p := \dim E_+$, the rank of polar equations of $E$ is $s_0+s_1+\dots+s_{p-1}$, while the dimension of $T_mM/E_+$ is $s_0+\dots+s_p$. So $s_p = 0$ on $E_+$ just when every regular hyperplane in $E_+$ is noncharacteristic.

**6.9.** Any hyperplane $E_-\subset E$ is $\mathcal{I}$-characteristic just when it lies in some other $\mathcal{I}$-integral element $E'$ of same dimension as $E$. But then $E'$ is also a $\mathcal{J}$-integral element, so $E_-$ is $\mathcal{J}$-characteristic: $\Xi^{\mathcal{I}}_E\subset\Xi^{\mathcal{J}}_E$.

**6.10.** At each point, there is a unique maximal integral element, so every hyperplane in it lies in that unique maximal integral element: $\Xi_E$ is empty.

**6.12.** The choices of $p$-dimensional integral element arise from the semipositive grade coefficients: $k$ coefficients of each polar in grade $k$, so $ks_k$ in that grade in all. So if $s_p = 0$ then there is a unique $p$-dimensional integral element containing $E$. If $E_+$ is involutive, then the semipositive grade coefficients $p^\alpha_i$ are arbitrary. In an adapted tableau, $E$ is determined by $p^\alpha_i$ for $i<p$, which are all of the coefficients just when $s_p = 0$.

**6.13.** In terms of the tableau we gave previously in solving problem 7.1 on page 54, the complex points in the projective plane satisfying the equations of the minors are

$$[0,1,0],[i,1,0],[-i,1,0],[i,0,1],[-i,0,1],[1,1,\sqrt{2}i],[1,1,-\sqrt{2}i].$$

The real ones for the characteristic variety as defined above, i.e. just $[0,1,0]$ corresponding to the hyperplane $0 = dx^2$.

**6.14.** For the harmonic function exterior differential system $\mathcal{I}$,

$$\mathcal{I}^2_{(x,y,u_x,u_y)} = \langle \Theta, d\vartheta, \vartheta \wedge (dx, dy, du_x, du_y) \rangle$$

has dimension 6. The characters are $s_0, s_1, s_2 = 1, 2, 0$, summing to 3, as we see from the tableau

$$\begin{pmatrix} du_y & -du_x \\ du_x & du_y \end{pmatrix} \wedge \begin{pmatrix} dx \\ dy \end{pmatrix}$$

So the symbol matrix is not square: $s_0 + s_1 + s_2 = 3 < 6 = \dim \mathcal{I}^2_{(x,y,u_x,u_y)}$.

Here is a way to "fix" this problem. Let $\mathcal{I}'$ be the ideal generated by the 2-forms $\vartheta \wedge dx, \Theta, d\vartheta$. The tableau of $\mathcal{I}'$ is

$$\begin{pmatrix} \vartheta & 0 \\ du_y & -du_x \\ du_x & du_y \end{pmatrix} \wedge \begin{pmatrix} dx \\ dy \end{pmatrix}$$

with characters $s'_0, s'_1, s'_2 = 0, 3, 0$. But now

$$\mathcal{I}'^2_{(x,y,u_x,u_y)} = \langle \Theta, d\vartheta, \vartheta \wedge dx \rangle$$

has dimension 3. Take an $\mathcal{I}'$-integral surface $S$ coframed by $dx, dy$ and containing an $\mathcal{I}$-integral curve $C$ on which $dx \neq 0$, say with $y = y(x)$. On $S$, $\vartheta \wedge dx = 0$ so $\vartheta = f(x,y)\, dx$ locally. On $C$, $0 = \vartheta = f(x, y(x))\, dx$. So $f(x, y(x)) = 0$. On $S$, $0 = d\vartheta = df \wedge dx$, so $f = f(x)$ and $f(x,y) = f(x) = f(x, y(x)) = 0$, hence $\vartheta = 0$ on $S$. We conclude that any $\mathcal{I}$-integral curve $C$ on which $dx \neq 0$ lies in a locally unique integral surface $S$. The problem with directly applying the Cauchy–Kovalevskaya theorem is that $\mathcal{I}$-integral curves are not arbitrary curves.

**7.1.** The exterior differential system lies on an 11-dimensional manifold $M$ with tableau modulo $\theta^1, \theta^2, \theta^3$

$$d \begin{pmatrix} \theta^1 \\ \theta^2 \\ \theta^3 \end{pmatrix} = - \begin{pmatrix} \pi^1_1 & \pi^2_1 & \pi^3_1 \\ \pi^2_1 & \pi^2_2 & \pi^3_2 \\ \pi^3_1 & \pi^3_2 & -\pi^1_1 - \pi^2_2 \end{pmatrix} \wedge \begin{pmatrix} \omega^1 \\ \omega^2 \\ \omega^3 \end{pmatrix}.$$

Denote our manifold as $M^{11}$ to indicate that it is 11-dimensional. Take the flag $M^8_0 \subset M^9_1 \subset M^{10}_2 \subset M^{11}$ given by

$$\begin{aligned} M^8_0 &:= (0 = x_1 = x_2 = x_3), \\ M^9_1 &:= (0 = x_2 = x_3), \\ M^{10}_2 &:= (0 = x_3). \end{aligned}$$

We have $\omega^i := dx^i$, $\theta^i = du^i - u^i_j dx^j$, and $\pi^i_j = du^i_j$ modulo $\theta^a, \omega^i$.

Pick a submanifold $R^8_2 \subset M^{10}_2$ of codimension $s_2 = 2$ given by equations

$$\begin{aligned} u^2_2 &= u^2_2(x^1, x^2), \\ u^3_2 &= u^3_2(x^1, x^2), \end{aligned}$$

a submanifold $R_1^4 \subset R_2^8 \cap M_1^9$ of codimension $s_1 = 3$ given by equations

$$\begin{aligned} u_1^1 &= u_1^1(x^1), \\ u_1^2 &= u_1^2(x^1), \\ u_1^3 &= u_1^3(x^1), \end{aligned}$$

and a point $R_0^0 \subset R_1^4 \cap M_0^8$ of codimension $s_0 = 3$ given by equations

$$\begin{aligned} u^1 &= c^1, \\ u^2 &= c^2, \\ u^3 &= c^3, \end{aligned}$$

for functions $u_j^i$ and constants $c^1, c^2, c^3$. Since there are no free derivatives to restrain in the final column of the tableau, there is no further restraining manifold.

So $R_0^0$ is a point. On $R_1^4$, the equations $0 = \theta^1 = \theta^2 = \theta^3$ become ordinary differential equations for functions $u^1, u^2, u^3$ of $x^1$, which have a unique solution through the point $R_0^0$. Check that on $R_1^4$, all of the tableau vanishes, so we don't have to solve any other equations than $0 = \theta^1 = \theta^2 = \theta^3$ to produce an integral curve.

One $R_2^8$, the tableau expands out to give derivatives in $x^2$:

$$\begin{aligned} \frac{\partial u_1^1}{\partial x^2} &= \frac{\partial u_1^2}{\partial x^1} + u_1^3 - f_1^3, \\ \frac{\partial u_1^2}{\partial x^2} &= \frac{\partial u_2^3}{\partial x^1}, \\ \frac{\partial u_1^3}{\partial x^2} &= \frac{\partial u_2^3}{\partial x^1}. \end{aligned}$$

The last two have right hand sides expressed in terms of the restraining manifold data:

$$\begin{aligned} u_1^2(x^1, x^2) &= u_1^2(x^1) + \int \frac{\partial u_2^3}{\partial x^1} dx^2, \\ u_1^3(x^1, x^2) &= u_1^3(x^1) + \int \frac{\partial u_2^3}{\partial x^1} dx^2. \end{aligned}$$

With these solved for, the first equation then solves:

$$u_1^1(x^1, x^2) = \int \left( \frac{\partial u_1^2}{\partial x^1} + u_1^3 - f_1^3 \right) dx^2.$$

It is not clear that the graph of these functions is an integral surface.

Finally, we expand out the tableau to find equations for derivatives in $x^3$:

$$\begin{aligned} \frac{\partial u_1^1}{\partial x^3} &= \frac{\partial u_1^3}{\partial x^1}, \\ \frac{\partial u_1^2}{\partial x^3} &= \frac{\partial u_2^3}{\partial x^1}, \\ \frac{\partial u_2^2}{\partial x^3} &= \frac{\partial u_2^3}{\partial x^2} + u_2^1 - f_2^1, \\ \frac{\partial u_1^2}{\partial x^3} &= -\frac{\partial u_1^1}{\partial x^1} - \frac{\partial u_2^2}{\partial x^1} - \frac{\partial f_i^i}{\partial x^1}, \\ \frac{\partial u_2^3}{\partial x^3} &= -\frac{\partial u_1^1}{\partial x^2} - \frac{\partial u_2^2}{\partial x^2} + \frac{\partial f_i^i}{\partial x^2}. \end{aligned}$$

Any 3-dimensional integral manifold will arise by solving this determined system, with initial data from the integral surface. The proof of the Cartan–Kähler theorem shows that this procedure always constructs a 3-dimensional integral manifold, which is not obvious.

**7.2.** Two proofs:

*a.* Writing out $\mathcal{I}$ in local coordinates as differential equations:

$$\frac{\partial u^i}{\partial x^{>i}} = \text{some function}\left(x, u, \frac{\partial u^j}{\partial x^{\le j}}\right),$$

as on page 42, $\mathcal{I}'$ yields only the subset of those differential equations which have the form

$$\frac{\partial u^i}{\partial x^p} = \text{some function}\left(x, u, \frac{\partial u^j}{\partial x^{\le j}}\right).$$

Apply the Cauchy–Kovalevskaya theorem to these equations.

*b.* The $\mathcal{I}'$-polar equations on $E_{p-1}$ are precisely those of $\mathcal{I}$. The $\mathcal{I}'$-integral elements containing $E_{p-1}$ are therefore those of $\mathcal{I}$. So $E_{p-1}$ lies in a unique integral element of $\mathcal{I}'$, which is $E_p$. The characters of $\mathcal{I}'$ are

| $s'_0$ | $s'_1$ | $\dots$ | $s'_{p-2}$ | $s'_{p-1}$ | $s'_p$ |
|---|---|---|---|---|---|
| 0 | 0 | $\dots$ | 0 | $s_0 + \dots + s_{p-1}$ | 0 |

Apply theorem 6.2 on page 49.

**7.4.** For $k \le p$, if $\mathcal{I}^k = 0$ then $\mathcal{I}^{k-1} \wedge \Omega^1 \subseteq \mathcal{I}^k$, so $0 = \mathcal{I}^{k-1} \wedge \Omega^1$ and so $0 = \mathcal{I}^{k-1}$.

**7.5.** If $\omega$ is identically zero, then $\mathcal{I}^{p-1} \wedge \Omega^1 = 0$ so $\mathcal{I}^{p-1} = 0$. Assume that $\omega$ is not identically zero. Write $X$ locally as $X = (0 = f)$ for some function $f$ with $df \ne 0$. If $\omega = 0$ at every point of $X$, replace $\omega$ by $\omega/f$, and repeat until $\omega \ne 0$ at some point of $X$. Replace $X_+$ by the open subset of $X_+$ containing a point of $X$ at which $\omega \ne 0$.

**7.6.** Since $e_i$ form a basis, $[e_i, e_j] = c^k_{ij} e_k$ for some functions $c^k_{ij}$. Write $e, e_{\hat\imath}, e_{\hat\imath\hat\jmath}$ to denote

$$e_1, \dots, e_p, \quad e_1, \dots, \hat{e}_i, \dots, e_p, \quad e_1, \dots, \hat{e}_i, \dots, \hat{e}_j, \dots, e_p.$$

By lemma 7.5 on page 55,

$$d\phi(e) = (-1)^{i+1} e_i \left(\phi(e_{\hat\imath})\right) + \sum_{i<j} \phi([e_i, e_j], e_{\hat\imath\hat\jmath}).$$

Plug in $d\phi = h\phi \wedge \omega^p$ and $\phi(e_{\hat\imath}) = (-1)^{p-i} \phi \wedge \omega^i(e)$ and similarly that

$$\begin{aligned} &\phi(e_k, e_{\hat\imath\hat\jmath}) \\ &= (-1)^{p+i+j} \left(\delta_{j=k} \phi \wedge \omega^i - \delta_{i=k} \phi \wedge \omega^j\right)(e) \end{aligned}$$

to get

$$e_p f = -\sum_{i<p} e_i(h^i f) + Hf,$$

where

$$H := (-1)^{p+1} h + \sum_{i<j} (c^j_{ij} h^i - c^i_{ij} h^j).$$

**8.1.** Pick a noncompact manifold $M$ of positive dimension, and a discrete infinite set $D \subset M$. Let $\mathcal{J}^0$ be the set of functions $f\colon M \to \mathbb{R}$ so that $f$ vanishes at all but finitely many points of $D$. Let $\mathcal{J}^k = \Omega^k_M$ for $k \ge 1$.

**8.2.** On $M := \mathbb{R}$, let $\mathcal{I}^0$ be the analytic functions vanishing at all but finitely many integers, and $\mathcal{I}^1 := \Omega^1_M$. Recall the infinite product expansion

$$\frac{\sin \pi x}{\pi x} = \prod_{k=1}^{\infty} \left(1 - \frac{x^2}{k^2}\right),$$

convergent in the complex plane [1] p. 197, [37] p. 239 12·14. So

$$f_n(x) := \prod_{k=n+1}^{\infty} \left(1 - \frac{x^2}{k^2}\right) \to 1$$

as analytic functions, but $f_n(x) = 0$ for $x$ integer, except at $x = -n, -n+1, \dots, -1, 0, 1, \dots, n-1, n$.

**8.3.** Where $v \ne 0$, straighten out, i.e. take coordinates in which $v = \partial_{x^1}$. Where $v = 0$, add a small multiple of a nonzero vector field and take a limit. For a more detailed proof: recall that for any function $f$,

$$\frac{\partial}{\partial t} e^{tv*} f = e^{tv*} \mathcal{L}_v f.$$

Apply induction, to get

$$\frac{\partial^k}{\partial t^k} e^{tv*} f = e^{tv*} \mathcal{L}_v^k f.$$

Taylor expand in $t$ in any coordinates, so the result holds for any function $f$. Taking exterior derivative, the result holds for $df$. If the result holds for two differential forms, then it holds for their wedge product: expand. Any differential forms are locally obtained by repeated wedging and exterior differentiating on functions.

**8.4.** If $v$ is a symmetry then

$$\mathcal{L}_v \vartheta = \left.\frac{d}{dt}\right|_{t=0} e^{tv*} \vartheta$$

so $\mathcal{L}_v \mathcal{I} \subseteq \mathcal{I}$.

If $\mathcal{L}_v \mathcal{I} \subseteq \mathcal{I}$, then in problem 8.3 on page 60, each term lies in the ideal. By closure under convergence (theorem 8.1 on page 60), $v$ is a symmetry.

**8.5.** For simplicity, let us just consider the system $\mathcal{I}$ on the plane $M = \mathbb{R}^2_{x,y}$ generated by $dy$. A vector field $v = a(x,y)\partial_x + b(x,y)\partial_y$ is a symmetry vector field just when $\mathcal{L}_v dy = f\, dy$ for some function $f$.

$$\begin{aligned}
\mathcal{L}_v dy &= d\mathcal{L}_v y, \\
&= db, \\
&= b_x\, dx + b_y\, dy,
\end{aligned}$$

we see that $b(x,y)$ depends only on $y$, i.e. $v = a(x,y)\partial_x + b(y)\partial_y$. Geometrically, $v$ flows points with equal $y$-value to points with equal $y$-value, i.e. its $y$-component depends only on $y$.

**8.7.** Take any nonzero compactly supported function $f$ on any manifold $M$. (We could even take $M = \mathbb{R}$.) Take any nonzero complete vector field $v$, with the support of $f$ contained in the interior of the support of $v$, so that the flow of $v$ takes some point in the support of $f$ outside of the support of $f$. (We could even take $v = \partial_x$.) Generate $\mathcal{I}$ with $f, \mathcal{L}_v f, \dots$. Each of these functions is supported in the support of $f$, so under the flow of $v$ is taken out of $\mathcal{I}$.

**8.8.** On $M = \mathbb{R}$, for each open set $U \subseteq M$, let $\mathcal{I}_U$ be generated by all analytic functions vanishing at the points $x = 1, 1/2, 1/3, \dots$ which lie in $U$. So if $0 \in U$, then $\mathcal{I}_U = 0$. A vector field $v$ vanishing at those points, and at the origin, and having compact support, is a symmetry of the smooth exterior differential system. Any analytic symmetry of $\mathcal{I}$ has to vanish at all of those points, so vanishes to all orders at the origin, so vanishes.

**8.10.** Pick local generators $\vartheta_i$. Write $\mathcal{L}_v \vartheta_i$ in those generators: $\mathcal{L}_v \vartheta_i = a_i^j \vartheta_i$. Pick an embedded smooth hypersurface $H$ on which $v \ne 0$. Make a function $g = I$ on $H$ and extend $g$ off of $H$ as a local solution of $\mathcal{L}_v g = -ga$. By straightening out, such a function $g$ exists, at least near each point of $H$. Check that $\mathcal{L}_v(g\vartheta) = 0$, so $g\vartheta$ is a $v$-invariant collection of local generators of $\mathcal{I}$.

For a symmetry $v$,

$$\mathcal{L}_v \vartheta = \left.\frac{d}{dt}\right|_{t=0} e^{-tv*}\vartheta.$$

**8.11.** Take any manifold $M$ of positive dimension and a foliation defined in some open set. Take an open set $W$ so that the foliation is defined on the closure of $W$, and a point $m_0$ on the boundary of $W$ near which the boundary of $W$ is an analytic hypersurface of $M$. For any open set $U \subseteq M$, define $\mathcal{I}_U$ to be

*a.* the forms pulling back to zero on each leaf of the foliation, if $U$ intersects $W$, and

*b.* $\mathcal{I}_U := \Omega_U^*$ otherwise.

Every point of $W$ lies in a unique integral manifold: its leaf.

Every open set $U \subseteq M$ containing $m_0$ intersects $W$, so $\mathcal{I}_U$ consists of the forms vanishing on the leaves. The foliation extends beyond $W$. If the leaf through $m_0$ is tangent to the boundary of $W$ near $m_0$, then it is an integral manifold near $m_0$. Otherwise there is no integral manifold through $m_0$, although there is an integral manifold with boundary. At every point outside the closure of $W$, there are no integral manifolds of dimension $p$. So the union of the integral manifolds might be neither open nor closed.

The space of integral elements is a manifold with boundary near $m_0$, so the system is *not* involutive, as the definition of involution requires the integral elements to form a manifold without boundary.

**8.14.** Let $M := \mathbb{R}^3_{x,u,v}$, $W \subset M$ an open subset not equal to $M$. For any open set $U \subseteq M$, define $\mathcal{I}_U$ to be

*a.* generated by $du \wedge dx, dv \wedge dx$ if $U$ intersects $W$,

*b.* generated by $du \wedge dx, dv \wedge dx, du \wedge dv$ otherwise.

For any point on the boundary of $W$, and open set $U$ containing that point, we don't have $du \wedge dv$ in $\mathcal{I}_U$, so for an open set $U' \subset U$ not intersecting $W$, $du \wedge dv \in \mathcal{I}_{U'}$ is not generated by any generators of $\mathcal{I}_U$. The 1-dimensional integral elements $du = u'\,dx$,$dv = v'\,dx$ are involutive.

**8.17.** $x\,dx$

**8.18.** The problem is local: we can assume that $\mathcal{I}$ is globally generated by pointwise linearly independent differential forms $\vartheta_i$. The flow of $v$ acts on the differential form bundle, as linear transformations of its fibers. We need to prove that the flow of $v$ preserves the vector subbundles. Take a pointwise basis $\vartheta^i, \vartheta^I$ of the differential forms. Write $\mathcal{L}_v\vartheta = a\vartheta$, for a block matrix

$$a = \begin{pmatrix} a^i_j & 0 \\ a^i_J & a^I_J \end{pmatrix}.$$

Then

$$e^{tv*}\vartheta = g\vartheta,$$

for a unique smooth function $g(x,t)$ for $(x,t$ in some open subset of $M \times \mathbb{R}$:

$$\begin{aligned} \frac{\partial}{\partial t}\vartheta &= \frac{\partial}{\partial t}e^{tv*}\vartheta, \\ &= \mathcal{L}_v e^{tv*}\vartheta, \\ &= e^{tv*}\mathcal{L}_v\vartheta, \\ &= e^{tv*}(a\vartheta), \\ &= (e^{tv*}a)g\vartheta, \end{aligned}$$

so that

$$g^{-1}\frac{\partial g}{\partial t} = e^{tv*}a,$$

has derivative lying in a block matrix of the prescribed form. Since $g(0) = I$ is also such a block matrix, $g(t)$ is such a block matrix for all $t$.

**8.19.** $du \wedge dx, dv \wedge dx, u\, du \wedge dv$ on $\mathbb{R}^4_{x,y,u,v}$

**8.21.** Apply the Cartan formula $\mathcal{L}_v\vartheta = d(v \lrcorner \vartheta) + v \lrcorner d\vartheta$ to a Cauchy characteristic vector field $v$, so see that $v$ is a symmetry vector field. If $v$ and $w$ are Cauchy characteristic vector fields, then

$$\begin{aligned} [v,w] \lrcorner \vartheta &= \mathcal{L}_v\mathcal{L}_w\vartheta - \mathcal{L}_w\mathcal{L}_v\vartheta, \\ &= \mathcal{L}_v(d(w \lrcorner \vartheta) + w \lrcorner d\vartheta) - \dots \end{aligned}$$

we expand out $\mathcal{L} = d\lrcorner + \lrcorner d$.

**8.22.** Fix a point $m$. Take a set of forms $\vartheta^a \in \mathcal{I}_U$ on some open set $U$ containing $m$, giving a basis of $\mathcal{I}_m$. In particular, the forms are linearly independent at $m$, so span a vector subbundle of the differential form bundle near $m$. Pick additional forms $\vartheta^\mu$ so that $\vartheta^a, \vartheta^\mu$ is a basis of the exterior algebra at $m$. These forms remain linearly independent nearby, so form a basis of the differential forms near $m$. Every tangent vector $v$ has $v \lrcorner \vartheta^a = \lambda^a_b(v)\vartheta^b + \lambda^a_\mu(v)\vartheta^\mu$, for unique $\lambda^a_b, \lambda^a_\mu \in T^*_m M$, hence linearly independent. The Cauchy characteristic subspace at $m$ is exactly the kernel of the various $\lambda^a_\mu$. For nearby points, the same linearly independent forms have a kernel of the same rank, containing the Cauchy characteristics. By constancy of dimension of Cauchy characteristics, this kernel is still the space of Cauchy characteristics.

**8.23.** Suppose constant rank. As in problem 8.22 on page 62, the equation of Cauchy characteristic vectors becomes a kernel of a constant rank vector bundle map, and so the local sections of that vector bundle are the local sections in the kernel of that map. The vector bundle map is the map quotienting $v \lrcorner \vartheta^a$ by $\vartheta^b$, so the kernel lies inside the set of Cauchy characteristic vector fields. But each Cauchy characteristic vector field lies in the kernel.

**8.24.** $dx \wedge du, dx \wedge dv, y\, dy \wedge du \wedge dv$ has involutive integral plane $\langle \partial_x, \partial_y \rangle$.

**8.29.** Note that $\pi_* v = 0$ just when $v \lrcorner \pi^* \vartheta = 0$ for any $\vartheta \in \Omega^*_{\bar{M}}$, so $v$ is a Cauchy characteristic of $\mathcal{I}$.

**8.30.** Consider a push forward. Take a generic tableau for the push forward, and pull it back. Add forms to it as needed, and make their rows generic as well. So we treat it, locally, as a part of the tableau for the original system. Take an integral manifold for the push forward. Restrict the original system to the preimage of the integral manifold. The rows in the tableau that were pulled back are now zero. But the polars of the other rows are still linearly independent.

**8.31.** By problem 1.12 on page 6, the space of Cauchy characteristics is a vector subbundle of the tangent bundle. Above we saw that Cauchy characteristics are bracket closed. Problem 1.14 on page 7 shows that the retracting space is therefore Frobenius, so $\pi$ exists locally.

**8.33.** Invariance of our generators under the flows of the Cauchy characteristic vector fields extends their definition to the largest set invariant under those flows and containing their domain. Each diffeomorphism $\phi$ preserving $\mathcal{I}$ and $\pi$ interchanging components allows us to extend that domain further, so that it becomes the preimage of an open set in $\bar{M}$.

**8.35.** Let $\mathcal{I} := \pi^* \bar{\mathcal{I}}$ and $\bar{\mathcal{J}} := \pi_* \mathcal{I}$. Easily $\bar{\mathcal{I}} \subseteq \bar{\mathcal{J}}$. Take some $\bar{\vartheta} \in \bar{\mathcal{J}}$, a differential form on $\bar{M}$. On $M$, $\bar{\vartheta} = \phi_i \wedge \bar{\vartheta}^i$ for some forms $\phi_i$ and some forms $\bar{\vartheta}_i$ from $\bar{\mathcal{I}}$.

Suppose that $\pi$ is a surjective submersion. Take coordinates $x^i$ on some open subset of $\bar{M}$. Pullback and extend to coordinates $x^i, y^a$ on some open subset of $M$. So

$$\bar{\vartheta} = f_I(x) dx^I = \phi_i \wedge \bar{\vartheta}^i = g_{JA}(x,y) dx^J \wedge dy^A \wedge f^i_I(x) dx^I.$$

Average over $y$: all coefficients are functions of $x$ only. Drop any terms with $dy$ in them, as they must cancel out.

**8.38.** For example, on $M = \mathbb{R}^4_{x,y,z,w}$, the exterior differential system $\mathcal{I}$ generated by $dy - z^2\, dx$ is generated by pullbacks of forms in the exterior differential system $\bar{\mathcal{I}}$ generated by $dy - z^2\, dx$ on $\bar{M} = \mathbb{R}^3_{x,y,z}$. But $\mathcal{I}^2$ consists of the 2-forms $z\, f(x,y,z)\, dz \wedge dx$, all of which vanish at $z = 0$.

**8.40.** The characteristic variety for the isometric immersion problem emerges from plugging $Da = v^a \xi$, $Db = v^b \xi$ and $Dc = v^c \xi$ into the tableau, but with $a\, Dc + c\, Da = 2b\, Db$, so $a\, v^c + c\, v^a = 2b\, v^b$, giving

$$\begin{aligned}
0 &= \begin{pmatrix} v^a & v^b \\ v^b & v^c \end{pmatrix} \begin{pmatrix} \xi \wedge \omega_1 \\ \xi \wedge \omega_2 \end{pmatrix}, \\
&= \begin{pmatrix} v^a(-\xi_2) + v^b \xi_1 \\ v^b(-\xi_2) + v^c \xi_1 \end{pmatrix} \omega_{12}, \\
&= \begin{pmatrix} -v^a \xi_2 + v^b \xi_1 \\ -v^b \xi_2 + \frac{1}{a}(2b v^b - c v^a) \xi_1 \end{pmatrix} \omega_{12}, \\
&= \begin{pmatrix} -\xi_2 & \xi_1 \\ -\frac{c}{a} \xi_1 & -\xi_2 + \frac{2b}{a} \xi_1 \end{pmatrix} \begin{pmatrix} v^a \\ v^b \end{pmatrix} \omega_{12}.
\end{aligned}$$

After we drop the Cauchy characteristics, the characteristic variety is the determinant of this matrix, i.e.

$$0 = -c\xi_1^2 + 2b\xi_1\xi_2 - a\xi_2^2.$$

Recall that the characteristic variety consists of the lines

$$0 = \sum_i \xi_i \omega_i = 0,$$

satisfying these equations. A vector $v = v_1 e_1 + v_2 e_2$ lies in such a characteristic line just when $0 = \xi_1 v_1 + \xi_2 v_2$, so then, up to scaling

$$(\xi_1, \xi_2) = (v_2, -v_1).$$

Plug this in to see that the characteristics are the curves whose velocities satisfy

$$0 = av_1^2 + 2bv_1v_2 + cv_2^2.$$

**9.3.**

$$\begin{aligned}
\omega_1' &= e_1' \cdot dx', \\
&= e^{-u}\left(\phi'(x)e_1\right) \cdot d\phi(x), \\
&= e^{-u}\left(\phi'(x)e_1\right) \cdot \phi'(x)\, dx, \\
&= e^{-u} e_1 \cdot dx, \\
&= e^{-u}\omega_1,
\end{aligned}$$

and similarly $\omega_2' = e^{-u}\omega_2$. So $\omega' = e^{-u}\omega$. Differentiate to get

$$d\omega' = i\alpha' \wedge \omega' = d(e^{-u}\omega),$$

which we expand out to find

$$(du + i(\alpha' - \alpha) \wedge \omega' = 0.$$

By the complex linear form of Cartan's lemma (which the reader can state and prove), we get the result.

**10.1.** If the eigenvalues of the shape operator at a point are $\lambda_1, \lambda_2$, then

$$H = \frac{1}{2}(\lambda_1 + \lambda_2), K = \lambda_1\lambda_2,$$

so

$$0 \le (\lambda_1 - \lambda_2)^2 = 4\left(H^2 - K\right).$$

**10.3.**

$$\begin{aligned}
\frac{\partial f}{\partial a_{22}} &= f_H \frac{\partial H}{\partial a_{22}} + f_K \frac{\partial K}{\partial a_{22}}, \\
&= \frac{f_H}{2} + f_K a_{11}, \\
&= \frac{f_H}{2} - e_3 \cdot \frac{de_1}{dt} f_K, \\
&\ne 0.
\end{aligned}$$

**A.2.** [31], p. 68.

**A.3.** [31], p. 69.

**B.1.**

$$P(x, y, \partial_x, \partial_y) = \begin{pmatrix} \partial_x & -\partial_y \\ \partial_y & \partial_x \end{pmatrix}$$

so turning derivatives $\partial_x, \partial_y$ into Greek variables $\xi, \eta$:

$$P(x, y, \xi, \eta) = \begin{pmatrix} \xi & -\eta \\ \eta & \xi \end{pmatrix}.$$

**B.3.** If we write each antisymmetric matrix

$$\begin{pmatrix} 0 & z & -y \\ -z & 0 & x \\ y & -x & 0 \end{pmatrix}$$

as $[x, y, z]$, Maxwell's equations become

$$\partial_t E = [\partial_x, \partial_y, \partial_z]H,$$
$$\partial_t H = -[\partial_x, \partial_y, \partial_z]E.$$

As an operator

$$P \begin{bmatrix} E \\ H \end{bmatrix} = \begin{pmatrix} \partial_t I & -[\partial_x, \partial_y, \partial_z] \\ [\partial_x, \partial_y, \partial_z] & \partial_t I \end{pmatrix}$$

so if we write 1-forms as $\xi\, dx + \eta\, dy + \zeta\, dz + \tau\, dt$, the symbol matrix is

$$\sigma_P = \begin{pmatrix} \tau I & -[\xi, \eta, \zeta] \\ [\xi, \eta, \zeta] & \tau I \end{pmatrix}$$

a square matrix with determinant

$$\left(\tau^2 - \xi^2 - \eta^2 - \zeta^2\right)^2 \tau^2.$$

The characteristic variety consists of the hyperplanes tangent to the light cone and the hyperplane $\tau = 0$ tangent to space at constant time.

**D.3.** Since $e$ is orthogonal, its transpose is $e^{\top} = e^{-1}$, i.e. $e^{\top} e = I$, i.e. $e_i \cdot e_j = 1$ if $i = j$, 0 otherwise. Differentiate: $\dot{e}^{\top} e + e^{\top} \dot{e} = 0$, i.e. $e^{\top} \dot{e} = -\dot{e}^{\top} e = -(e^{\top} \dot{e})^{\top}$, i.e. $e^{\top} \dot{e}$ is antisymmetric, with entries $e_i \cdot \dot{e}_j$. Since $\dot{e}_1 = k e_2$,

$$e^{\top} \dot{e} = \begin{pmatrix} 0 & ? & ? \\ k & ? & ? \\ 0 & ? & ? \end{pmatrix}.$$

By antisymmetry,

$$e^{\top} \dot{e} = \begin{pmatrix} 0 & -k & 0 \\ k & 0 & -t \\ 0 & t & 0 \end{pmatrix}$$

for some $t$.

**D.5.** For any constants $\alpha, \beta, \gamma$, let

$$\begin{pmatrix} x_1 \\ x_2 \\ x_3 \end{pmatrix} := \begin{pmatrix} \alpha \cos \beta s \\ \alpha \sin \beta s \\ \gamma s \end{pmatrix}.$$

so

$$\begin{pmatrix} \dot{x}_1 \\ \dot{x}_2 \\ \dot{x}_3 \end{pmatrix} = \begin{pmatrix} -\alpha\beta \sin \beta s \\ \alpha\beta \cos \beta s \\ \gamma \end{pmatrix},$$

$$|\dot{x}|^2 = \alpha^2\beta^2 + \gamma^2.$$

So we need this to equal 1. Clearly $\alpha = 0$ or $\beta = 0$ is a line, and $\gamma = 0$ is a circle. We ignore those cases and so the equation $\alpha^2\beta^2 + \gamma^2 = 1$ can be solved for any of these three constants in terms of the other two. We can reflect in the $x_1x_2$ plane to arrange that $\alpha, \beta > 0$ and in the $x_3$ axis to arrange that $\gamma > 0$. We get

$$e_1 = \begin{pmatrix} -\alpha\beta \sin \beta s \\ \alpha\beta \cos \beta s \\ \gamma \end{pmatrix}$$

so

$$\dot{e}_1 = k e_2 = \begin{pmatrix} -\alpha\beta^2 \cos \beta s \\ -\alpha\beta^2 \sin \beta s \\ 0 \end{pmatrix},$$

so

$$k = \alpha\beta^2$$

and

$$e_2 = -\begin{pmatrix} \cos \beta s \\ \sin \beta s \\ 0 \end{pmatrix}.$$

Differentiate

$$\dot{e}_2 = \begin{pmatrix} \beta \sin \beta s \\ -\beta \cos \beta s \\ 0 \end{pmatrix}$$

so

$$t e_3 = \dot{e}_2 + k e_1 = t e_3 = \begin{pmatrix} \beta\gamma^2 \sin \beta s \\ -\beta\gamma^2 \cos \beta s \\ \alpha\beta^2\gamma \end{pmatrix}.$$

Finally,

$$t = \pm\beta\gamma,$$

and

$$e_3 = \pm \begin{pmatrix} \gamma \sin \beta s \\ -\gamma \cos \beta s \\ \alpha\beta \end{pmatrix}.$$

So if we set

$$\beta := \sqrt{k^2 + t^2},$$

$$\alpha := \frac{k}{\beta^2},$$

$$\gamma := \pm\frac{t}{\beta},$$

we get a helix with prescribed constant values of $k, t$.

**D.6.** If $e_2' = e_2$, then $e_1, e_3$ and $e_1', e_3'$ are orthonormal bases of the same plane $e_2^\perp$. Choose the sign of $e_3$ to get

$$e_1' + ie_3' = e^{i\theta}(e_1 + ie_3)$$

for some angle $\theta$. Differentiate to get

$$de_1' + i\,de_3' = \dot{\theta}e^{i\theta}(e_1 + ie_3) + e^{i\theta}(\dot{e}_1 + i\dot{e}_3).$$

Expand to get

$$(k' - it')e_2 = \dot{\theta}e^{i\theta}(e_1 + ie_3) + e^{i\theta}(k - it)e_2$$

Differentiate $e_2 = e_2'$ to get $-k'e_1' + t'e_3' = -ke_1 + te_3$. Hence

$$k' + it' = e^{-i\theta}(k + it).$$

So $\dot{\theta} = 0$, a constant rotation.

On the other hand, we can now pick any constant $\theta$, and define:

$$\begin{aligned} e_1' + ie_3' &:= e^{i\theta}(e_1 + ie_3), \\ e_2' &:= e_2, \\ k' + it' &:= e^{-i\theta}(k + it), \\ x' &= \int e_1'\,ds \end{aligned}$$

and check that the Serret–Frenet equations are satisfied, so these yield the Serret–Frenet frame of a curve.

**D.8.** The proof requires some unwinding of notation: the expression $\omega_i = e_i \cdot dx$ means that $\omega_i = \sum_j e_{ji} dx_j$, which allows us to unwind the following formal steps:

$$\begin{aligned} d\omega_i &= d\left(e_i \cdot dx\right), \\ &= de_i \wedge dx, \\ &= \sum_j \left(e_j \cdot de_i\right) \wedge \left(e_j \cdot dx\right), \\ &= \sum_j \gamma_{ji} \wedge \omega_j. \end{aligned}$$

Similarly, the expression $\gamma_{ij} = e_i \cdot de_j$ means that $\gamma_{ij} = \sum_k e_{ki} de_{kj}$, so:

$$\begin{aligned} d\gamma_{ij} &= d\left(e_i \cdot de_j\right), \\ &= de_i \wedge de_j, \\ &= \sum_k \left(e_k \cdot de_i\right) \wedge \left(e_k \cdot de_j\right), \\ &= \sum_k \gamma_{ki} \wedge \gamma_{kj}. \end{aligned}$$

**D.9.** Here are two proofs:

*a.* If we think of $x$ and $e$ as functions on $^{\ulcorner}\mathbb{E}^3$, then $r_g^* x = x$, $r_g^* e = eg$. Hence $r_g^* dx = dx$ and $r_g^* de = (de)g$. So $r_g^* \omega = r_g^*(e^\top dx) = (eg)^\top dx = g^\top e^\top dx = g^\top \omega$ and $r_g^* \gamma = r_g^*(e^\top de) = (eg)^\top d(eg) = g^\top e^\top de\, g = g^\top \gamma g$.

*b.* The action is $r_g(x, e) = (x, eg)$, where $(eg)_i = \sum_j g_{ji} e_j$. Hence

$$r'_g(x, e)(\dot{x}, \dot{e}) = (\dot{x}, \sum_j g_{ji} \dot{e}_j).$$

$$\begin{aligned}
(r_g^* \omega)_{(x,e)}(\dot{x}, \dot{e}) &= \omega_{r_g(x,e)} r'_g(x, e)(\dot{x}, \dot{e}), \\
&= \omega_{(x,eg)}(\dot{x}, \sum_j g_{ji} \dot{e}_j), \\
&= (eg) \cdot \dot{x}, \\
&= \sum_j g_{ji} e_j \cdot \dot{x}, \\
&= \sum_j g_{ji} \omega_{(x,e)}(\dot{x}, \dot{e}).
\end{aligned}$$

**D.13.** For any tangent vectors $u, v \in T_x S$,

$$I\!I(u, v) = I\!I(v, u) = \sum_{ij} a_{ij} u_i v_j e_3,$$

We can write $I\!I$ as $I\!I(u, v) = \sum a_{ij} \omega_i(u) \omega_j(v) e_3$, i.e. $I\!I = \omega^\top a \omega e_3$, so

$$r_g^* I\!I = r_g^* \omega^\top a \omega e_3 = (h^\top \omega)^\top (n h^\top a h)(h^\top \omega)(n e_3) = \omega^\top a \omega e_3 = I\!I,$$

i.e. $r_g^* I\!I = I\!I$ so $I\!I$ is invariant.

**D.16.** Take any point $p_0 \in \mathbb{E}^3$. At a point $x_0$ where the surface acheives maximal distance from some point $p_0 \in \mathbb{E}^3$, differentiate distance to see that $T_{x_0} S$ is perpendicular to the ray from $p_0$ to $x_0$. Translate $x_0$ to the origin, and rescale and rotate to get $p_0$ a unit vector above the origin. Apply our local picture of surfaces to see that $S$ being inside the unit sphere around $p_0$ forces $S$ to be locally the graph of a function with critical zero at the origin, and eigenvalues of the second derivative both at least 1.

**E.1.** Note that $\rho$ rotates the point $x'$ of our plane around on a circle perpendicular to $e_1, e_3$ and tangent to $e_2$, with radius $r$, so that

$$\frac{d\rho}{d\theta} x' = r e_2.$$

Therefore

$$\begin{aligned}
\omega_1 &= e_1 \cdot dx, \\
&= e_1 \cdot d\left(\rho x'\right), \\
&= e_1 \cdot \rho\, dx' + e_1 \cdot \frac{d\rho}{d\theta} x'\, d\theta, \\
&= e_1 \cdot \rho\, e'_1 ds + e_1 \cdot r e_2\, d\theta, \\
&= e_1 \cdot e_1 ds \\
&= ds.
\end{aligned}$$

Similarly, $\omega_3 = 0$ and

$$\begin{aligned}\omega_2 &= e_2 \cdot dx, \\ &= e_2 \cdot d\left(\rho x'\right), \\ &= e_2 \cdot \rho\, dx' + e_2 \cdot \frac{d\rho}{d\theta}\, x'\, d\theta, \\ &= e_2 \cdot \rho\, e_1' + e_2 \cdot r e_2\, d\theta, \\ &= r\, d\theta.\end{aligned}$$

Note that

$$\begin{aligned}de_1' &= \dot\varphi\, e_3' ds, \\ de_2' &= 0, \\ de_3' &= -\dot\varphi\, e_1' ds,\end{aligned}$$

and

$$d\rho = \rho \begin{pmatrix} 0 & -1 & 0 \\ 1 & 0 & 0 \\ 0 & 0 & 0 \end{pmatrix} d\theta.$$

Compute

$$\begin{aligned}de_1 &= d(\rho e_1'), \\ &= d\rho\, e_1' + \rho\, de_1', \\ &= \rho \begin{pmatrix} 0 & -1 & 0 \\ 1 & 0 & 0 \\ 0 & 0 & 0 \end{pmatrix} e_1' d\theta + \rho\dot\varphi e_3' ds, \\ &= \rho \begin{pmatrix} 0 & -1 & 0 \\ 1 & 0 & 0 \\ 0 & 0 & 0 \end{pmatrix} \begin{pmatrix} \cos\varphi \\ 0 \\ \sin\varphi \end{pmatrix} d\theta + \dot\varphi e_3 ds, \\ &= \rho \begin{pmatrix} 0 \\ \cos\varphi \\ 0 \end{pmatrix} d\theta + \dot\varphi e_3 ds, \\ &= \cos\varphi e_2 d\theta + \dot\varphi e_3 ds.\end{aligned}$$

and similarly

$$\begin{aligned}de_2 &= d(\rho e_2'), \\ &= d\rho\, e_2', \\ &= \rho \begin{pmatrix} 0 & -1 & 0 \\ 1 & 0 & 0 \\ 0 & 0 & 0 \end{pmatrix} e_2' d\theta \\ &= -\rho \begin{pmatrix} 1 \\ 0 \\ 0 \end{pmatrix} d\theta \\ &= -\begin{pmatrix} \cos\theta \\ \sin\theta \\ 0 \end{pmatrix} d\theta.\end{aligned}$$

Finally,

$$\begin{aligned} de_3 &= d(\rho e_3'), \\ &= d\rho\, e_3' + \rho\, de_3', \\ &= \rho \begin{pmatrix} 0 & -1 & 0 \\ 1 & 0 & 0 \\ 0 & 0 & 0 \end{pmatrix} e_3' d\theta - \rho\dot{\varphi} e_1' ds, \\ &= \rho \begin{pmatrix} 0 & -1 & 0 \\ 1 & 0 & 0 \\ 0 & 0 & 0 \end{pmatrix} \begin{pmatrix} -\sin\varphi \\ 0 \\ \cos\varphi \end{pmatrix} d\theta - \dot{\varphi} e_1 ds, \\ &= \rho \begin{pmatrix} 0 \\ -\sin\varphi \\ 0 \end{pmatrix} d\theta - \dot{\varphi} e_1 ds, \\ &= -\sin\varphi\, e_2 d\theta - \dot{\varphi} e_1 ds. \end{aligned}$$

Therefore

$$\begin{aligned} \gamma_{12} &= -\gamma_{21}, \\ &= -e_2 \cdot de_1, \\ &= -e_2 \cdot (\cos\varphi\, e_2 d\theta + \dot{\varphi} e_3 ds), \\ &= -\cos\varphi\, d\theta. \end{aligned}$$

$$\begin{aligned} \gamma_{13} &= e_1 \cdot de_3, \\ &= -\dot{\varphi}\, ds, \end{aligned}$$

$$\begin{aligned} \gamma_{23} &= e_2 \cdot de_3, \\ &= -\sin\varphi\, d\theta. \end{aligned}$$

**E.2.** If $\sin\varphi = 0$ everywhere, the curve is a horizontal line and the surface of revolution an annulus. Suppose that $\sin\varphi$ is not everywhere zero; pick an interval on which $\sin\varphi \neq 0$, and let $\varepsilon$ be the sign of $\sin\varphi$. Differentiate to find that the function

$$\beta := \varepsilon r \sqrt{1 - \dot{r}^2} + H_0 r^2$$

is constant along solutions, say equal to $\beta_0$. Solve for $\dot{r}$:

$$\dot{r}^2 = 1 - \frac{\beta_0 - H_0 r^2}{r^2}.$$

Substitute $u = r^2$:

$$\dot{u} = 2\sqrt{H_0 u^2 + u - \beta_0}.$$

Integrate:

$$\int \frac{du}{2\sqrt{H_0 u^2 + u - \beta_0}} = s.$$

This integral can be solved in elementary functions giving $s = s(r)$; for example, if $H_0 > 0$,

$$s = \frac{(4\beta_0 H_0 + 1)\log\left(2\sqrt{H_0}\sqrt{H_0 u^2 + u - \beta_0} + 2H_0 u + 1\right)}{16 H_0^{3/2}},$$
$$= \frac{(4\beta_0 H_0 + 1)\log\left(2\sqrt{H_0}\sqrt{H_0 r^4 + r^2 - \beta_0} + 2H_0 r^2 + 1\right)}{16 H_0^{3/2}},$$

The surfaces of revolution of constant positive mean curvature are obtained by solving implicitly for $r = r(s)$.

**E.6.** Taking real and imaginary parts $q = q_1 + iq_2$:

$$\begin{pmatrix} a_{11} & a_{12} \\ a_{12} & a_{22} \end{pmatrix} = \begin{pmatrix} H + q_1 & -q_2 \\ -q_2 & H - q_1 \end{pmatrix}.$$

So $H = (a_{11} + a_{22})/2$ is the mean curvature, and

$$q = \frac{a_{11} - a_{22} - 2i\, a_{12}}{2}.$$

Note that $|q|^2 = q\bar{q} = H^2 - K$, so $q = 0$ precisely at umbilic points.

**E.9.** In our complex notation, we are identifying the matrix

$$\begin{pmatrix} 0 & -1 & 0 \\ 1 & 0 & 0 \\ 0 & 0 & 0 \end{pmatrix}$$

with the complex number $i$. The fibers of ${}^{\ulcorner}S \to S$ have the form

$$(x(t), e(t)) = (x_0, e_0 e^{it}).$$

So if we write out $\gamma = e^{\top} de$ in our complex notation, on each fiber,

$$\gamma = i\, dt,$$

so $\alpha = \gamma_{12} = -dt$. We are given that

$$f(x(t), e(t)) = f(x_0, e_0 e^{it}) = e^{ikt} f(x_0, e_0),$$

so on the fiber

$$\frac{df}{dt} = ikf,$$

i.e.

$$df + ikf\alpha = 0.$$

So in complex notation, $df + ikf\alpha$ vanishes on vertical vectors for ${}^{\ulcorner}S \to S$, i.e. is a multiple of $\omega, \bar{\omega}$.

**E.10.**

$$Dq = a_{111} - 3a_{122} + i(a_{222} - 3a_{112}),$$
$$\bar{D}q = 2(H_1 - iH_2).$$

**E.11.**

$$
\begin{aligned}
0 &= \begin{pmatrix} d(\gamma_{13} - k_1\omega_1) \\ d(\gamma_{23} - k_2\omega_2) \end{pmatrix} \\
&= \begin{pmatrix} -\gamma_{12} \wedge \gamma_{23} - dk_1 \wedge \omega_1 + k_1\gamma_{12} \wedge \omega_2 \\ -\gamma_{21} \wedge \gamma_{13} - dk_2 \wedge \omega_2 + k_2\gamma_{21} \wedge \omega_1 \end{pmatrix} \\
&= \begin{pmatrix} -\gamma_{12} \wedge k_2\omega_2 - dk_1 \wedge \omega_1 + k_1\gamma_{12} \wedge \omega_2 \\ \gamma_{12} \wedge k_1\omega_1 - dk_2 \wedge \omega_2 - k_2\gamma_{12} \wedge \omega_1 \end{pmatrix} \\
&= - \begin{pmatrix} dk_1 \wedge \omega_1 + (k_2 - k_1)\gamma_{12} \wedge \omega_2 \\ (k_2 - k_1)\gamma_{12} \wedge \omega_1 + dk_2 \wedge \omega_2 \end{pmatrix}
\end{aligned}
$$

**E.12.** Arrange $e_1, e_2$ to diagonalize the shape operator, so

$$
\begin{aligned}
\gamma_{13} &= k_1\omega_1, \\
\gamma_{23} &= k_2\omega_2,
\end{aligned}
$$

with $k_1$ constant. Differentiate to find that $\gamma_{12}$ is a multiple of $\omega_2$. So then along the flow of $e_1$, $\gamma_{12} = 0$ and $\gamma_{23} = 0$ i.e. $e_2$ is constant, so $e_1, e_3$ rotate in the plane, with $\gamma_{13} = k_1$, so at a constant rate, i.e. on a circle.

**F.5.**

$$
\begin{aligned}
e_3 \frac{d}{ds}\left(\cos\theta e_2' + \sin\theta e_3'\right) &= e_3\left(-\dot\theta \sin\theta e_2' + \cos\theta \dot{e'}_2 + \dot\theta\cos\theta e_3' + \sin\theta \dot{e'}_3\right), \\
&= e_3\left(-\dot\theta\sin\theta e_2' + \cos\theta(-ke_1 + te_3') + \dot\theta\cos\theta e_3' - t\sin\theta e_2'\right), \\
&= e_3\left(-\dot\theta\sin\theta e_2' + t\cos\theta e_3' + \dot\theta\cos\theta e_3' - t\sin\theta e_2'\right), \\
&= \dot\theta\sin^2\theta t\cos^2\theta\dot\theta\cos^2\theta + t\sin^2\theta, \\
&= \dot\theta + t.
\end{aligned}
$$

**F.6.** We get $\dot{e}_1 = 0$ along such a curve, at that point, so $0 = \gamma_{21} = \gamma_{31}$, so $a_{11} = 0$. Since this occurs in all directions, $a = 0$.

**G.4.** The fact that $\exp'_{x_0} = I$ means that at the origin,

$$
\omega_1^2 + \omega_2^2 = dx^2 + dy^2,
$$

in rectangular coordinates, i.e. near the origin

$$
\omega_1^2 + \omega_2^2 - dr^2 - r^2\, d\theta^2
$$

is smooth. Expand out, and plug in

$$
r\, d\theta = \cos\theta\, dy - \sin\theta\, dx
$$

to see that

$$
\begin{aligned}
\omega_1^2 + \omega_2^2 &= dr^2 + r^2\, d\theta^2 + \left(\left(\frac{h}{r}\right)^2 - 1\right) r^2 d\theta^2, \\
&= dx^2 + dy^2 + \left(\left(\frac{h}{r}\right)^2 - 1\right)(\cos\theta\, dy - \sin\theta\, dx).
\end{aligned}
$$

In particular,

$$\frac{h}{r} \to 1$$

as $r \to 0$. Hence $h \to 0$ as $r \to 0$ and

$$\partial_r h \to \lim_{r \to 0} \frac{h}{r} = 1.$$

**G.7.** Gauss–Bonnet

**G.8.** If our geodesic is $x(t)$, then $f(x(t))$ has local minima as its only critical points. Moreover, they are isolated, as $f(x(t)$ is critical only where $x(t)$ reaches the minimum of $f$ or becomes tangent to a level set of $f$. If $f(x(t)$ stays bounded as $t \to \infty$, then $x(t)$ stays in a compact sublevel set of $f$. But $f(x(t))$ increases for large enough $t$, so $x(t)$ approaches the boundary of that sublevel set, i.e. a level set. Geodesics near each level set enter or leave soon, a contradiction.

**G.10.** The differential equation forces $h = r - Kr^3/3! + O(r)^4$, and inductively forces $\partial_r^k \partial_\theta^\ell h = O(r)^\ell$. Note that

$$\begin{aligned}
\partial_x &= \cos\theta \partial_r - \frac{\sin\theta}{r}\partial_\theta, \\
\partial_y &= \sin\theta \partial_r + \frac{\cos\theta}{r}\partial_\theta.
\end{aligned}$$

So $h^2 - r^2$ is a smooth function of $x, y$ at the origin.

$$\begin{aligned}
dr^2 + h^2 d\theta^2 &= dr^2 + r^2 d\theta^2 + (h^2 - r^2) d\theta^2, \\
&= dr^2 + r^2 d\theta^2 + \frac{(h^2 - r^2)(-x\,dy + y\,dx)}{r^2}, \\
&= dx^2 + dy^2 - \frac{K(x^2 + y^2) + O(r)^3}{3}(-x\,dy + y\,dx)
\end{aligned}$$

is smooth at the origin.

# Index